\documentclass[a4paper,10pt]{article}
\usepackage[english]{babel}
\usepackage{amssymb}
\usepackage{amsmath}
\usepackage{indentfirst}
\usepackage[dvips]{graphicx}
\usepackage{epsfig}
\usepackage{lscape}
\usepackage{hyperref}
\usepackage{color}
\usepackage{amsthm}
\usepackage{verbatim}

\newtheorem{defi}{Definition}[section]
\newtheorem{cor}{Corollary}[section]
\newtheorem{prop}[cor]{Proposition}
\newtheorem{lemma}[cor]{Lemma}
\newtheorem{teo}[cor]{Theorem}

\newtheorem{rem}[cor]{Remark}
\newtheorem{prob}[cor]{Problem}
\newtheorem{ass}[cor]{Assumption}
\newtheorem{ex}[cor]{Example}

\newcommand{\R}{\mathbb{R}}
\newcommand{\N}{\mathbb{N}}

\def\mathref#1{\ifmmode\mathrm{(\ref{#1})}\else{\rm(\ref{#1})}\fi} 
\def\nref#1{\ifmmode\mathrm{\ref{#1}}\else{\rm\ref{#1}}\fi}

\numberwithin{equation}{section}
\allowdisplaybreaks

\title{Minimal horizontal triods: Analysis and computation}

\author{Robert N\"urnberg%
\thanks{Dipartimento di Matematica, Universit\`a di Trento, 
38123 Trento, Italy, 
\url{robert.nurnberg@unitn.it}}
~and 
Paola Pozzi%
\thanks{Fakult\"at f\"ur Mathematik, Universit\"at Duisburg-Essen,
Thea-Leymann-Stra\ss e 9,
45127 Essen, Germany, \url{paola.pozzi@uni-due.de}} 
}
\date{}

\begin{document}

\maketitle

\begin{abstract}
In this article we investigate the question of finding a network configuration of minimal length connecting three given points in the Heisenberg group. After proving existence of (possibly degenerate) minimal horizontal triods, we investigate their characterization. We then formulate a horizontal curve shortening flow that deforms any given suitable initial triod into a critical point for the length functional. Numerical experiments based on a stable fully discrete finite element scheme provide useful insights into the rich landscape of this sub-Riemannian geometry. 
\end{abstract}

\bigskip
 

\section{Introduction}

In this article we investigate the classical Steiner problem in the Heisenberg group. The (first) Heisenberg group is the simplest example of a sub-Riemannian space. Thus the geometry is constrained, in the sense that  motion is only allowed along a given set of directions that change from point to point.

The classical Euclidean problem of finding the path of shortest length connecting three given distinct points $P_{1}, P_{2},P_{3} \in \R^{3}$, see e.g.\ \cite{Tromba,JarnikK34,GP68,HwangRW92}, is now posed  as the question of existence of (possibly degenerate) shortest horizontal triods of minimal length. 
A horizontal triod is the union of three horizontal curves connecting the points $P_{\alpha}$, $\alpha=1,2,3$, to a common junction point $\Sigma$. For a curve to be horizontal, its tangent vector must lie in a specific plane that changes from point to point. The triod is said to be degenerate if $\Sigma $ coincides with one of the $P_{\alpha}$ and the corresponding curve collapses to a point. The length of the triod is given as the sum  of the lengths of its component curves in an appropriate metric.
 Unlike the Euclidean case, where triods are unions of straight line segments, and where the Steiner problem can be treated with elementary geometry (see the beautiful presentation given in \cite{Tromba}), we need here some more refined mathematics.

In this work   we show existence of minimal horizontal (possibly degenerate) triods and give a characterization of the minimal configurations using direct methods of calculus of variations and using some ideas from the field of minimal surfaces. 
We find that: if a minimal triod is composed of three curves then its projected planar curves meet at a 120 degree angle and have curvatures that sum up to zero. In particular, the projected triod is composed of lines and/or arcs of circles. If the minimal triod is degenerate, then the remaining projected curves have constant curvature but no angle condition needs to hold. 
Uniqueness can in general not be expected, since even a curve of minimal length connecting two given distinct points need not to be unique in the Heisenberg group (see for instance \cite[\S~2.3]{CapognaDanielli} or the direct computations provided in Remark~\ref{rem1.10} in the Appendix).

 Additionally we  formulate  a flow, which we call horizontal curve shortening flow,  that lets a given initial horizontal triod evolve into a critical configuration for the length functional. Such evolutions are discretized by  finite elements with an unconditionally stable algorithm. The simulations  we propose allow for  a deeper understanding of this fascinating anisotropic space.
 In particular, the numerical experiments show that the landscape of degenerate triods is much richer than in the Euclidean case. While in the Euclidean setting it is known that as soon as the largest angle of the triangle formed by the three points $P_{\alpha}$, $\alpha=1,2,3$, is larger or equal to 120 degrees, then the minimal configuration is given by two straight segments, see e.g.\ \cite{Tromba,GP68,HwangRW92},
we could not find a similar condition experimentally for the horizontal set-up (see for instance Experiments~2,3 and 4 in Section~\ref{sec:nr}). 
 Moreover, degenerate triods are clearly non-unique:  Experiment~10  in Section~\ref{sec:nr} indicates the existence of a degenerate solution made up of the horizontal lift of two circles  connecting three points lying on the $z$-axis. Any rotation around the $z$-axis yields another critical configuration of the same length, hence providing an infinite family of degenerate solutions.

\bigskip
The remainder of this article is organized as follows: in Section~\ref{sec2} we recall all necessary notions needed to understand the geometry of the Heisenberg group. In Section~\ref{sec:Htriods} we formulate the problem and provide proof of existence of minimal triods.  After deriving
the Euler--Lagrange equations for minimal and critical configurations, the horizontal curve shortening flow for triods is formulated. In Section~\ref{sec:FEM} a finite element scheme is proposed and its stability properties are proved. Several interesting numerical simulations are given in Section~\ref{sec:nr}.

\bigskip

Throughout this paper, we use superscripts to denote the components of a vector in $\R^2$ or $\R^3$, so that e.g.\ $P=(P^{1}, P^{2}, P^{3})^t$ for $P\in\R^3$. For a vector $P\in\R^2$, we define the perpendicular vector by $P^\perp = \binom{P^1}{P^2}^{\perp}=\binom{-P^2}{P^1}$.

\section{Preliminaries}\label{sec2}

We consider the Heisenberg group $\mathcal{H}= (\R^{3}, \circ)$, where $\circ$ is the group operation
given by
\[
\begin{pmatrix} x \\ y \\ z \end{pmatrix}
\circ 
\begin{pmatrix} \hat{x} \\ \hat{y} \\ \hat{z} \end{pmatrix} 
= 
\begin{pmatrix} x+ \hat{x} \\ y+ \hat{y} \\ 
z+ \hat{z} + \frac{1}{2}(x \hat{y}-y \hat{x}) \end{pmatrix}.
\]
For comprehensive information about $\mathcal{H}$ we refer to \cite{Calin}, \cite{CapognaDanielli}, and the references given there.
Here we just  briefly recall that $\mathcal{H}$ is a Lie group with Lie algebra generated by the left invariant vector fields
\begin{align*}
X_{1}&=\frac{\partial}{\partial x}- \frac{y}{2}\frac{\partial}{\partial z},\qquad
X_{2}=\frac{\partial}{\partial y}+ \frac{x}{2}\frac{\partial}{\partial z},\qquad
X_{3} =\frac{\partial}{\partial z}.
\end{align*}
Note that $X_{3}=[X_{1},X_{2}]$  and that $\mathcal{H}$ is a nilpotent Lie group of step 2.

With some abuse of notation we set
\begin{align*}
X_{1}(x,y,z)=\begin{pmatrix} 1\\0 \\ -\frac{y}{2} \end{pmatrix}, \qquad
X_{2}(x,y,z)=\begin{pmatrix} 0\\1 \\ \frac{x}{2} \end{pmatrix}, \qquad
X_{3}(x,y,z)=\begin{pmatrix} 0\\0 \\ 1 \end{pmatrix},
\end{align*}
i.e.\ we identify the vector fields with the associated  differential operators.
At each point $p=(x,y,z)^t$, the three vectors $X_{i}(x,y,z)$ are linear independent.
We denote by $g_{p}(\cdot, \cdot): \R^{3} \times \R^{3} \to [0, \infty)$ the inner product for which the vectors $X_{i}(p)$ build an orthonormal basis and set
\begin{align*}
|v|_{g}=\sqrt{g_{p}(v,v)} \qquad \text{ for } v \in \mbox{span} \{X_{1}(p),X_{2}(p),X_{3}(p)\} .
\end{align*}
The horizontal tangent space at a point $p=(x,y,z)^t$ is linearly spanned by $X_{1}(p)$ and $X_{2}(p)$.

\begin{defi}
A (sufficiently smooth) curve $\gamma:[0,1] \to \mathcal{H}$, $\gamma=\gamma(u)$, is said to be \emph{horizontal} if its tangent vector $\gamma'(u)$ is a linear combination of $X_{1}(\gamma(u))$ and $X_{2}(\gamma(u))$ for each $u \in [0,1]$. 
\end{defi}

\newcommand{\oldG}{\mathfrak{G}}
\begin{rem}\label{rem1.1}
Note that since for $\gamma:[0,1] \to \R^{3}$, $\gamma=\gamma(u)$ (sufficiently smooth), it holds
\begin{align*}
\gamma'(u)&=\gamma_{u}(u) = \begin{pmatrix}
\gamma^{1}_{u}(u)\\\gamma^{2}_{u}(u) \\ \gamma^{3}_{u}(u)
\end{pmatrix}\\
&= \gamma^{1}_{u}(u)X_{1}(\gamma(u)) +\gamma^{2}_{u}(u)X_{2}(\gamma(u))  +\left(\gamma^{3}_{u}(u) +\frac{1}{2}\gamma^{2}(u)\gamma^{1}_{u}(u) -\frac{1}{2}\gamma^{1}(u)\gamma^{2}_{u}(u) \right)X_{3}(\gamma(u)),
\end{align*}
we have that $\gamma$ is horizontal if and only if
\begin{equation} \label{eq:iff}
\gamma^{3}_{u} =-\frac{1}{2}\gamma^{2}\gamma^{1}_{u} +\frac{1}{2}\gamma^{1}\gamma^{2}_{u} \quad \text{in } [0,1].
\end{equation}
In particular, for a horizontal curve its third space dimension component $\gamma^{3}$ is  determined by the other two space components via
\begin{equation} \label{g3}
\gamma^{3}(u) = \gamma^{3}(0) + \int_{0}^{u} \left(-\frac{1}{2}\gamma^{2}(\xi)\gamma^{1}_{u}(\xi) +\frac{1}{2}\gamma^{1}(\xi)\gamma^{2}_{u}(\xi) \right) d\xi,
\end{equation}
or, equivalently,
\begin{equation} \label{g3reverse}
\gamma^{3}(u) = \gamma^{3}(1) - \int_{u}^{1} \left(-\frac{1}{2}\gamma^{2}(\xi)\gamma^{1}_{u}(\xi) +\frac{1}{2}\gamma^{1}(\xi)\gamma^{2}_{u}(\xi) \right) d\xi.
\end{equation}
Finally, note that for a horizontal curve $\gamma$ we have
\[
|\gamma_{u}(u)|_{g}=\sqrt{ (\gamma^{1}_{u}(u))^{2} + (\gamma^{2}_{u}(u))^{2}},
\]
and its length is given by
\[
L(\gamma)=\int_{0}^{1} |\gamma_{u}(u)|_{g} du= \int_{0}^{1} \sqrt{ (\gamma^{1}_{u}(u))^{2} + (\gamma^{2}_{u}(u))^{2}} du.
\]
The length functional is a geometric functional in the sense that it is invariant under reparametrizations of the curve $\gamma$. The same holds for the functional (cf.\ with \eqref{g3})
\begin{equation} \label{eq:G}
\oldG(\gamma):=\int_{0}^{1} \left(-\frac{1}{2}\gamma^{2}(u)\gamma^{1}_{u}(u) +\frac{1}{2}\gamma^{1}(u)\gamma^{2}_{u}(u) \right) du
\end{equation}
when restricting to reparametrizations that keep the  same orientation of the curve.
\end{rem}

\begin{defi}[Projected curve and horizontal lift]
Given a horizontal curve $\gamma:[0,1] \to \mathcal{H}$, we denote by $c:[0,1] \to \R^{2}$, $c=(\gamma^{1}, \gamma^{2})^{t}$, the projected curve of $\gamma$ in the plane $\{ z=0\}.$

Vice-versa, given a planar curve $c:[0,1] \to \R^{2}$ we call \emph{horizontal lift} of $c$ any map
$\gamma:[0,1] \to \mathcal{H}$, $\gamma=(\gamma^{1}, \gamma^{2}, \gamma^{3})^{t}$, with $(\gamma^{1}, \gamma^{2})^{t}=c$ and $\gamma^{3}$ satisfying \eqref{g3}.
\end{defi}

\begin{rem}
Clearly, the lift of $c$ is only unique up to specifying $\gamma^{3}(0)$, for
example. Hence we say that $\gamma$ is the horizontal lift of $c$ starting
at $P \in \R^3$ if in addition $\gamma(0) = P$, recall \eqref{g3}. 
Similarly, we say that $\gamma$
is the horizontal lift of $c$ ending at $Q \in \R^3$ if $\gamma(1) = Q$,
recall \eqref{g3reverse}.
\end{rem}

\begin{rem}[Geometrical interpretation for the $z$-component of a horizontal curve]\label{rem1.2}
Note that for a horizontal curve joining the points $P=\gamma(0)$ and $Q=\gamma(1) \in \R^{3}$, equation~\eqref{g3}  yields that
\[
|Q^{3}-P^{3}|=\left|\int_{0}^{1}\frac{1}{2} \binom{\gamma^{1}(u)}{\gamma^{2}(u)}\cdot \binom{-\gamma^{2}_{u}(u)}{\gamma^{1}_{u}(u)} du \right|
= \left|\int_{0}^{1}\frac{1}{2} c(u)\cdot c_{u}^\perp(u) du \right|
=| \text{oriented area}(S)|,
\]
where $S$ is the portion in the plane swept by the vector $c(u)$ of the projected curve $c$.

For example, a planar semicircle of radius $R>0$, connecting, say, the points $P=0$ and $Q=(2R,0,0)^t$, is \emph{not} a horizontal curve, since $|Q^{3}-P^{3}|=0$, while the swept region has area $\frac\pi2 R^2$. On the other hand, any planar straight line segment that goes through the origin is horizontal.
\end{rem}

\begin{lemma}\label{leftT}
For $a=(a^{1},a^{2},a^{3})^t \in \R^{3}$ let $\tau_{a}(x,y,z):=(a^{1},a^{2},a^{3})^t \circ (x,y,z)^t$ denote the left translation. If $\gamma:[0,1] \to \mathcal{H}$ is horizontal, then so is the curve $u \mapsto \tau_{a}(\gamma(u))$.
Moreover, $L(\gamma)=L(\tau_{a}(\gamma))$.
\end{lemma}
\begin{proof}
The claims follow by direct computation.
\end{proof}
\begin{rem}\label{rem1.4}
Left translation  of a horizontal curve $\gamma:[0,1] \to \mathcal{H}$ by $a=(-\gamma^{1}(0), -\gamma^{2}(0),-\gamma^{3}(0))^t$ ensures that the translated curve  $\tau_{a}(\gamma)$ starts at the origin.
\end{rem}
\begin{rem}
For $a=(0,0,a^{3})^t$ we have that $\tau_{a}(\gamma(u)) =\gamma(u)+(0,0,a^{3})^{t}$.
\end{rem}
We recall the following important connectivity theorem.
\begin{prop}
Any two points in $\mathcal{H}$ can be connected by a horizontal smooth curve.
\end{prop}
\begin{proof}
See, for instance, \cite[\S1.3]{Calin}.
\end{proof}

\subsection{On tangent, normal, and curvature vectors}
\begin{defi}\label{def:gammareg}
We say that a (sufficiently smooth) horizontal curve $\gamma:[0,1] \to \mathcal{H}$ is regular if $|\gamma_{u}(u)|_{g} \neq 0$ for all $u \in [0,1]$.
\end{defi}
To motivate the definition of curvature vector for horizontal curves, let us first recall the following lemma:
\begin{lemma}\label{keuclidea}For a sufficiently smooth regular horizontal curve $\gamma:[0,1] \to \R^{3}$ let
\[
\vec{k}(u):=\frac{1}{|\gamma_{u}(u)|_{g}} \left( \frac{1}{|\gamma_{u}(u)|_{g}}  
\binom{\gamma^{1}_{u}(u)}{\gamma^{2}_{u}(u)}
\right)_{u} \quad \in \R^{2}.
\]
With 
\begin{equation} \label{eq:k}
k(u):= \frac{\gamma^{2}_{uu}(u) \gamma^{1}_{u}(u) - \gamma^{2}_{u}(u) \gamma^{1}_{uu}(u)}{((\gamma^{1}_{u}(u))^{2} + (\gamma^{2}_{u}(u))^{2})^{\frac{3}{2}}}
\end{equation}
there holds
\[
\vec{k}(u)= k(u) \frac{1}{|\gamma_{u}(u)|_{g}} 
\binom{-\gamma^{2}_{u}(u)}{\gamma^{1}_{u}(u)}
.
\]
\end{lemma}
\begin{proof}
The claim follows by direct computation. 
\end{proof}
Note that $\vec{k}$ is the standard Euclidean curvature vector of the projected planar curve $c:I \to \R^{2}$, $c(u):=(\gamma^{1}(u), \gamma^{2}(u))^{t}$, whose arc length  is given by $ds=\sqrt{ (\gamma^{1}_{u}(u))^{2} + (\gamma^{2}_{u}(u))^{2}} du$. Thus $\vec{k}=\partial_{ss} c$.

In analogy to the Euclidean case we provide the following definition of \emph{horizontal curvature vector}.
\begin{lemma}
For a sufficiently smooth regular horizontal curve $\gamma:[0,1] \to \R^{3}$ we define the horizontal unit tangent
\begin{align}\label{eq:T}
T(u) & :=\frac{1}{|\gamma_{u}(u)|_{g}} 
\begin{pmatrix}
\gamma^{1}_{u}(u)\\ \gamma^{2}_{u}(u) \\ \gamma^{3}_{u}(u)
\end{pmatrix} =  \frac{1}{|\gamma_{u}(u)|_{g}}\begin{pmatrix}
\gamma^{1}_{u}(u)\\ \gamma^{2}_{u}(u) \\ -\frac{1}{2}\gamma^{2}(u)\gamma^{1}_{u}(u) +\frac{1}{2}\gamma^{1}(u)\gamma^{2}_{u}(u)
\end{pmatrix} \notag \\
&=\frac{\gamma^{1}_{u}(u)}{|\gamma_{u}(u)|_{g}}X_{1}(\gamma(u)) + \frac{\gamma^{2}_{u}(u)}{|\gamma_{u}(u)|_{g}}X_{2}(\gamma(u)) ,
\end{align}
the horizontal normal vector
\begin{equation} \label{eq:N}
N(u):=-\frac{\gamma^{2}_{u}(u)}{|\gamma_{u}(u)|_{g}}X_{1}(\gamma(u)) + \frac{\gamma^{1}_{u}(u)}{|\gamma_{u}(u)|_{g}}X_{2}(\gamma(u)) 
= \frac{1}{|\gamma_{u}(u)|_{g}}
\begin{pmatrix}
-\gamma^{2}_{u}(u)\\ \gamma^{1}_{u}(u) \\ \frac{1}{2}\gamma^{2}(u)\gamma^{2}_{u}(u) +\frac{1}{2}\gamma^{1}(u)\gamma^{1}_{u}(u)
\end{pmatrix},
\end{equation}
and the horizontal curvature vector
\[
\vec{k}_{g}(u) := \frac{1}{|\gamma_{u}(u)|_{g}}T_{u}(u).
\]
Then we have that
\begin{align}\label{kg}
\vec{k}_{g}(u) &= k(u)N(u) \notag \\
&= \begin{pmatrix}
\vec{k}(u)\\ \frac{k(u)}{|\gamma_{u}(u)|_{g}} (\frac{1}{2}\gamma^{2}(u)\gamma^{2}_{u}(u) +\frac{1}{2}\gamma^{1}(u)\gamma^{1}_{u}(u))
\end{pmatrix}
= \frac{k(u)}{|\gamma_{u}(u)|_{g}}\begin{pmatrix}
 -\gamma^{2}_{u}(u)\\\gamma^{1}_{u}(u) \\  \frac{1}{2}\gamma^{2}(u)\gamma^{2}_{u}(u) +\frac{1}{2}\gamma^{1}(u)\gamma^{1}_{u}(u)
\end{pmatrix},
\end{align}
where $\vec{k}$ and $k$ are as in Lemma~\ref{keuclidea}.
\end{lemma}
\begin{proof}
It follows from Lemma~\ref{keuclidea} that we only need to verify the third component in \eqref{kg}. On recalling \eqref{eq:iff}, we have that
\begin{equation} \label{eq:star}
\gamma^{3}_{u}=\frac{1}{2}( \gamma^{1}\gamma^{2}_{u}-\gamma^{2}\gamma^{1}_{u}), 
\qquad
\gamma^{3}_{uu}=\frac{1}{2}( \gamma^{1}\gamma^{2}_{uu}-\gamma^{2}\gamma^{1}_{uu}).
\end{equation}
Hence we can compute
\begin{align*}
\frac{1}{|\gamma_{u}(u)|_{g}} T_u^3(u) &=
\frac{1}{|\gamma_{u}|_{g}} \left ( \frac{\gamma^{3}_{u}}{|\gamma_{u}|_{g}} \right )_{u} =
\frac{1}{|\gamma_{u}|_{g}} \left( \frac{\gamma^{3}_{uu}}{|\gamma_{u}|_{g}}-\frac{\gamma^{3}_{u}}{|\gamma_{u}|^{3}_{g}}(  \gamma_{u}^{2}\gamma^{2}_{uu} +\gamma^{1}_{u}\gamma^{1}_{uu}
 ) \right)\\
 &=\frac{1}{|\gamma_{u}|_{g}} \frac{1}{|\gamma_{u}|^{3}_{g}} \left(\gamma^{3}_{uu} ((\gamma^{1}_{u})^{2}+
 (\gamma^{2}_{u}))^{2} - \gamma^{3}_{u}(\gamma^{2}_{u}\gamma^{2}_{uu} +\gamma^{1}_{u}\gamma^{1}_{uu})\right)\\
 &=\frac{1}{|\gamma_{u}|_{g}} \frac{k(u)}{2} (\gamma^{2}\gamma^{2}_{u} +\gamma^{1}\gamma^{1}_{u}),
\end{align*}
where we have used \eqref{eq:star} and \eqref{eq:k} in the last equality.
\end{proof}

Note that
\begin{equation} \label{ONframenew}
g(T,T)=1, \quad g(N,N)=1, \quad g(T,N)=0, \quad g(T,X_{3})=0, \quad g(N,X_{3})=0,
\end{equation}
so that $\{ T,N,X_{3} \}$ provides another orthonormal frame with respect to the inner product $g(\cdot, \cdot)$.

Finally, let us recall how  time dependent motion of a curve in its $N$ direction can be corrected so that the curve $\gamma(t)$ stays horizontal along the flow (see \cite[Lemma~2]{Drugan}): 
\begin{lemma}\label{lemma-hproperty}
Let $\gamma_{0}: [0,1] \to \R^{3}$ be a horizontal regular curve. Suppose that the (sufficiently smooth) family of  curves $\gamma:[0,\hat{T}) \times [0,1] \to \R^{3}$, $\gamma=\gamma(t,u)$, satisfies $\gamma(0, \cdot)=\gamma_{0}(\cdot)$ and
\[
\gamma_{t}(t,u)=m(t,u) N(t,u) -\left(\int_{0}^{u}m(t, \xi) |\gamma_{u}(t, \xi)|_{g}d\xi \right)X_{3}(\gamma(t,u))
\]
for some sufficiently smooth $m:[0,\hat{T}) \times [0,1] \to \R$, with $N(t,\cdot)$ as in \eqref{eq:N}, and $|\gamma_{u}(t,u)|_{g}=\sqrt{ (\gamma^{1}_{u}(t,u))^{2} + (\gamma^{2}_{u}(t,u))^{2}} \neq 0$ for all $(t,u) \in [0,\hat{T}) \times [0,1]$.
Then $\gamma(t)$ is horizontal for any $t \in [0,\hat{T})$.
\footnote{The same result applies with $\gamma_{t}(t,u)=m(t,u) N(\gamma(t,u)) +\left(\int_{u}^{1}m(t ,\xi) |\gamma_{u}(t,\xi)|_{g}d\xi \right)X_{3}(\gamma(t,u))$.}
\end{lemma}
\begin{proof}
We compute
\begin{align*}
&
\left(\gamma^{3}_{u} +\frac{1}{2}\gamma^{2}\gamma^{1}_{u}- \frac{1}{2}\gamma^{1}\gamma^{2}_{u} \right)_{t} = (\gamma_{t}^{3})_{u} -\frac{1}{2} \gamma^{1}_{t}\gamma^{2}_{u} -\frac{1}{2} \gamma^{1}(\gamma^{2}_{t})_{u} + \frac{1}{2} \gamma^{2}_{t}\gamma^{1}_{u} +\frac{1}{2} \gamma^{2}(\gamma^{1}_{t})_{u}
\\ & \quad
=
\left( \frac{m}{|\gamma_{u}|_{g}} \frac{1}{2} (\gamma^{2}\gamma^{2}_{u} +\gamma^{1}\gamma^{1}_{u}
) \right)_{u} -m |\gamma_{u}|_{g}
 -\frac{1}{2} \gamma^{2}_{u}  \Big(-m\frac{\gamma^{2}_{u}}{|\gamma_{u}|_{g}}\Big) 
-\frac{1}{2} \gamma^{1} \Big(m\frac{\gamma^{1}_{u}}{|\gamma_{u}|_{g}}\Big)_{u}\\
& \qquad +\frac{1}{2} \gamma^{1}_{u}  \Big(m\frac{\gamma^{1}_{u}}{|\gamma_{u}|_{g}}\Big) 
+\frac{1}{2} \gamma^{2}  \Big(-m\frac{\gamma^{2}_{u}}{|\gamma_{u}|_{g}}\Big)_{u} =0.
\end{align*}
The claim then follows from the smoothness assumptions, the initial curve being horizontal and~\eqref{eq:iff}.
\end{proof}

We remark that motion in the $T$ direction maintains the horizontal property:
\begin{lemma}\label{lemma-hproperty2}
Let $\gamma_{0}: [0,1] \to \R^{3}$ be a horizontal regular curve. Suppose that the (sufficiently smooth) family of  curves $\gamma:[0,\hat{T}) \times [0,1] \to \R^{3}$, $\gamma=\gamma(t,u)$, satisfies $\gamma(0, \cdot)=\gamma_{0}(\cdot)$ and
\[
\gamma_{t}(t,u)=q(t,u) T(t,u)
\]
for some sufficiently smooth $q:[0,\hat{T}) \times [0,1] \to \R$, with $T(t,\cdot)$ as in \eqref{eq:T}, and $|\gamma_{u}(t,u)|_{g}=\sqrt{ (\gamma^{1}_{u}(t,u))^{2} + (\gamma^{2}_{u}(t,u))^{2}} \neq 0$ for all $(t,u) \in [0,\hat{T}) \times [0,1]$.
Then $\gamma(t)$ is horizontal for any $t \in [0,\hat{T})$.
\end{lemma}
\begin{proof}
As in the previous lemma we compute
\begin{align*}
\left(\gamma^{3}_{u} +\frac{1}{2}\gamma^{2}\gamma^{1}_{u}- \frac{1}{2}\gamma^{1}\gamma^{2}_{u} \right)_{t}&= (\gamma_{t}^{3})_{u} -\frac{1}{2} \gamma^{1}_{t}\gamma^{2}_{u} -\frac{1}{2} \gamma^{1}(\gamma^{2}_{t})_{u} + \frac{1}{2} \gamma^{2}_{t}\gamma^{1}_{u} +\frac{1}{2} \gamma^{2}(\gamma^{1}_{t})_{u}
\\
&=\left( \frac{q}{|\gamma_{u}|_{g}} ( -\frac{1}{2}\gamma^{2}\gamma^{1}_{u} +\frac{1}{2}\gamma^{1}\gamma^{2}_{u}
)\right)_{u}  -\frac{1}{2} \gamma_{u}^{2} q \frac{\gamma^{1}_{u}}{|\gamma_{u}|_{g}} -\frac{1}{2} \gamma^{1} \left( q\frac{\gamma^{2}_{u}}{|\gamma_{u}|_{g}}\right)_{u}\\
& \quad + \frac{1}{2} \gamma_{u}^{2} q \frac{\gamma^{1}_{u}}{|\gamma_{u}|_{g}} +\frac{1}{2} \gamma^{2} \left( q\frac{\gamma^{1}_{u}}{|\gamma_{u}|_{g}}\right)_{u} =\frac{1}{2} \gamma^{1}_{u}\gamma^{2}_{u} \frac{q}{|\gamma_{u}|_{g}} -\frac{1}{2} \gamma^{2}_{u}\gamma^{1}_{u} \frac{q}{|\gamma_{u}|_{g}}=0.
\end{align*}
The claim then follows as before from the smoothness assumptions and the initial curve being horizontal.
\end{proof}

\subsection{Horizontal curves of minimal length}
The existence of a length minimizing curve   connecting any two given points in $\mathcal{H}$ is a well studied topic. It is discussed for instance in \cite[\S1.6 and Appendix D.4]{Montgomery}, \cite{Calin}, and \cite[\S2.3]{CapognaDanielli}, but can also be derived anew with the techniques described below in Section~\ref{sec:Htriods}. 

In the following we derive the Euler--Lagrange equations for curves minimizing the length functional and investigate their solutions, since these computations will be useful later on in Section~\ref{sec:Htriods}.

\newcommand{\oldC}{\mathcal{C}_{single}}

\begin{rem}[Euler--Lagrange equations for length minimizing curves]\label{rem-EL}
Following the ideas presented in Remark~\ref{rem1.1}, to find  a shortest horizontal curve $\gamma:[0,1] \to \R^{3}$ connecting two given distinct points $P=(P^{1},P^{2},P^{3})^{t}$ and $Q=(Q^{1},Q^{2},Q^{3})^{t}\in \R^{3}$,  it is natural  to consider the class of projected curves
\[
\oldC:= \{ c \in C^{1}([0,1], \R^{2}) %
\,: \, c(0)=(P^{1},P^{2})^{t}, \ c(1)=(Q^{1},Q^{2})^{t} \}
\]
with Euclidean  length
\[
L_{E}(c):=\int_{0}^{1} |c_{u}(u)| du,
\]
where $|\cdot|$ denotes the Euclidean norm, and solve the problem
\begin{equation} \label{problem1}
\min \{ L_{E}(c) \,:\, c \in \oldC \text{ and } G(c)=Q^{3}-P^{3} \},
\end{equation}
where, similarly to \eqref{eq:G},  
\begin{equation} \label{funG}
G(c)=\int_{0}^{1} \left(-\frac{1}{2}c^{2}(u) c^{1}_{u}(u) +\frac{1}{2} c^{1}(u)c^{2}_{u}(u) \right) du.
\end{equation}
The horizontal lift, starting at $P$, of a minimizer for \eqref{problem1}, is then a horizontal length minimizing curve connecting $P$ and $Q$.

Note that for $\varphi\in C^{\infty}([0,1], \R^{2})$ the first variation of $G$ at $c$ in the direction $\varphi$ is given by

\begin{align}\label{varG}
\delta G(c)(\varphi)=\frac{d}{d \epsilon}\Big|_{\epsilon=0} G(c + \epsilon \varphi) 
& =\frac{1}{2} \int_{0}^{1}\varphi^{1}c^{2}_{u} + c^{1} \varphi^{2}_{u} -c^{2}\varphi^{1}_{u}-\varphi^{2}c^{1}_{u} du   \notag \\
&=\frac{1}{2}[c^{1}\varphi^{2}- c^{2}\varphi^{1}]_{0}^{1} + \int_{0}^{1}\varphi^{1}c^{2}_{u} -\varphi^{2}c^{1}_{u} du   \notag \\
&=\frac{1}{2}[c^{1}\varphi^{2}- c^{2}\varphi^{1}]_{0}^{1} + \int_{0}^{1} \binom{\varphi^{1}}{\varphi^{2}} \cdot \binom{c^{2}_{u}}{-c^{1}_{u}}
du.
\end{align}
In particular, we have that if $c\in \oldC$ is a regular curve (i.e.\ $|c'| \neq 0$ in $[0,1]$) then $\delta G(c)(\varphi)$ does not vanish for all $\varphi \in C_{0}^{\infty}((0,1), \R^{2})$. It then follows from \cite[Proposition 1.17, Proposition 1.2]{BGH} that a regular weak minimizer for problem \eqref{problem1} is such that there exists $\lambda \in \R$ such that
\[
\delta L_{E}(c)(\varphi)+ \lambda \delta G(c)(\varphi) =0 \qquad \forall \, \varphi  \in C_{0}^{\infty}((0,1), \R^{2}).
\]
Moreover, if $c \in C^{2}([0,1], \R^{2})$, then 
\begin{align*}
0&=-\int_{0}^{1}\frac{1}{|c'|}   \left( \frac{1}{|c'|} 
\binom{c^{1}_{u}}{c^{2}_{u}} \right)_{u}  \cdot 
\binom{\varphi^{1}}{\varphi^{2}} |c'| du
 + \lambda \int_{0}^{1} \binom{\varphi^{1}}{\varphi^{2}} \cdot \frac{1}{|c'|} 
\binom{c^{2}_{u}}{-c^{1}_{u}} |c'| du\\
 &=-\int_{0}^{1} \vec{k} \cdot \binom{\varphi^{1}}{\varphi^{2}} |c'| du
 -\lambda \int_{0}^{1} \binom{\varphi^{1}}{\varphi^{2}}
 \cdot \frac{1}{|c'|} \binom{-c^{2}_{u}}{c^{1}_{u}} |c'| du.
\end{align*}
for all test functions $\varphi  \in C_{0}^{\infty}((0,1), \R^{2})$.
It follows that $\lambda$ can be expressed through the following equation
\[
\int_{0}^{1} k(u) |c_u(u)| du + \lambda L_{E}(c)=0.
\]
Summing up, a horizontal length minimizing curve $\gamma:[0,1] \to \mathcal{H}$ connecting two given distinct points $P$ and $Q$ in $\R^{3}$ is a horizontal curve such that (in the notation of Lemma~\ref{keuclidea})
\begin{equation} \label{EL}
k(u)+\lambda =0 \text{ in } (0,1), \quad \text{ with }\quad  \lambda=- \frac{\int_{0}^{1} k(u) |\gamma_{u}(u)|_{g} du }{L(\gamma)}.
\end{equation}
\end{rem}

From \eqref{EL} we observe that length minimizing curves in $\mathcal{H}$ can be straight segments or smooth curves whose projection onto the plane are arcs of circles.
In Remark~\ref{rem1.10}, in the appendix, we derive an explicit solution formula and recognize the surprising fact that there is a unique length minimizing curve connecting the origin to a point $Q=(Q^{1},Q^{2},Q^{3})^t\neq 0$ if $(Q^{1},Q^{2})^t \neq 0$, whereas we have infinitely many length minimizing curves connecting the origin to $Q=(0,0,Q^{3})^t$.

\section{Horizontal triods} \label{sec:Htriods}
In the following we consider triods of horizontal curves.
Given are three distinct points $P_{\alpha}=(P^{1}_{\alpha}, P^{2}_{\alpha}, P^{3}_{\alpha})^{t} \in \R^{3}$, $\alpha=1,2,3$.
We consider three regular horizontal curves $\gamma_{\alpha}:[0,1] \to \R^{3}$, $\gamma_{\alpha}=\gamma_{\alpha}(u)$, connecting the points $P_{\alpha}$ to a common junction point, which we denote by $\Sigma$, i.e. 
\begin{align*}
\gamma_{\alpha}(1)&=P_{\alpha}, \qquad \alpha=1,2,3,\\
\gamma_{1}(0)&=\gamma_{2}(0)=\gamma_{3}(0)=:\Sigma \in \R^{3},
\end{align*}
with associated functional
\begin{equation} \label{Ltriods}
L(\Gamma)=\sum_{\alpha=1}^{3} L(\gamma_{\alpha})
\end{equation}
for $\Gamma=(\gamma_{1}, \gamma_{2}, \gamma_{3})$ the horizontal triod.
In the following we say that a horizontal triod $\Gamma$ is regular  if all its curves $\gamma_{\alpha}$ are regular (according to Definition~\ref{def:gammareg}).
Our goal is to study minimal configurations for $L$. 
First of all, we show existence in a proper class of functions.

\subsection{Existence of a minimal configuration}
As in the case of a single curve, we can reduce the problem to a constrained planar one.
We denote by
\[
\mathcal{C}:=\{\Gamma_{c}=(c_{1},c_{2},c_{3}) \,:\, c_{\alpha} \in H^{1,1}((0,1),\R^{2}), %
c_{\alpha}(1)=(P^{1}_{\alpha},P^{2}_{\alpha})^{t}, \alpha=1,2,3, \ c_{1}(0)=c_{2}(0)=c_{3}(0)\}
\]
the class of natural admissible functions with associated functional
\[
L_E(\Gamma_{c}):=\sum_{\alpha=1}^{3} L_{E}(c_{\alpha}).
\]
Note that because of the embedding of $H^{1,1}$ into the class of continuous maps, the pointwise boundary conditions are meaningful.
In the following we set $\overline{\Sigma}:=c_{1}(0)=c_{2}(0)=c_{3}(0)$ to be the planar junction point.
We want to solve the problem
\begin{equation*} %
\min \{L_E(\Gamma_{c}) \,:\,  \Gamma_{c} \in \mathcal{C}, \quad G(c_{1})-P_{1}^{3}=G(c_{2})-P_{2}^{3}=G(c_{3})-P_{1}^{3}\} ,
\end{equation*}
where $G$ is as in \eqref{funG} (recall also \eqref{eq:G} and \eqref{g3}). The horizontal lift of a minimal planar triod $\Gamma_c$, where each curve $c_\alpha$ is lifted to a horizontal curve $\gamma_\alpha$ ending at $P_\alpha$, is then a horizontal minimal triod.

Observe that the projected points $\overline{P}_{\alpha}:=(P^{1}_{\alpha},P^{2}_{\alpha})^{t}$, $\alpha=1,2,3$, need not be distinct, despite the fact that the original points $P_{\alpha}$ are. However, whenever two projected points coincide,  say $\overline{P}_{\alpha}=(P^{1}_{\alpha},P^{2}_{\alpha})^{t}=(P^{1}_{\beta},P^{2}_{\beta})^{t}=\overline{P}_{\beta}$ for some $\alpha, \beta \in \{ 1,2,3 \}$, then their third space component must be different, i.e.\ $P^{3}_{\alpha} \neq P^{3}_{\beta}$ if $\alpha\neq\beta$. On the other hand, if $P^{3}_{\alpha} = P^{3}_{\beta}$ for $\alpha\neq\beta$, then $\overline{P}_{\alpha} \neq \overline{P}_{\beta}$. This implies that for any triod $\Gamma_{c} \in \mathcal{C}$ satisfying the constraint $G(c_{1})-P_{1}^{3}=G(c_{2})-P_{2}^{3}=G(c_{3})-P_{1}^{3}$, at most one curve $c_{\alpha}$ can reduce to a constant map mapping to a point (in which case $G(c_{\alpha})=0$). The horizontal lift $\Gamma$ of such a configuration corresponds to a horizontal triod in which the junction point $\Sigma$ has ``collapsed'' to one of the points $P_{\alpha}$ and one curve of $\Gamma$ has reduced to a point: we refer to such a horizontal triod as a degenerate one.

Following the ideas presented in \cite[\S5.9]{BGH} we give the following definition. 
\begin{defi}\label{def:quasino}
We say that a curve $c\in H^{1,1}((0,1), \R^2)$, $c=c(u)$, is \emph{quasinormal} if $|c_{u}|=L_E(c)>0$ in $(0,1)$. A triod $\Gamma_{c} \in \mathcal{C}$ is quasinormal if all its components $c_{\alpha}$, $\alpha=1,2,3$, are quasinormal.
Note that the planar junction point $\overline{\Sigma}$ of a quasinormal triod might coincide with one of the  projected points $\overline{P}_{\alpha}$.

We say that a triod $\Gamma_{c} \in \mathcal{C}$ is \emph{degenerate} if exactly one of its components is a constant map. In other words, if the  junction point $\overline{\Sigma}$ coincides with one of the $\overline{P}_{\alpha}$ and the corresponding curve (starting and ending in the junction point) degenerates to one point. A degenerate triod is said to be quasinormal if its two remaining curves are quasinormal.
\end{defi}

Instead of working with the length functional, it is advantageous to exploit the properties of the Dirichlet functional, which is well known to ``dominate'' the length functional in a convenient way (see \cite[\S5.9]{BGH} and \cite{MH99} from which we have borrowed several ideas). Thus we study the following problem.
Let
\[
\mathcal{C}_{D}:=\{\Gamma_{c}=(c_{1},c_{2},c_{3}) \in \mathcal{C} \,:\, c_{\alpha} \in H^{1,2}((0,1),\R^{2})\}
\]
and 
\[
D(\Gamma_{c}):=\frac{1}{2}\sum_{\alpha=1}^{3} \int_{0}^{1}|(c_{\alpha})_{u}|^{2} du \geq 0.
\]
For $1\geq \epsilon >0$ let us consider the $\epsilon$-problem:
\begin{equation} \label{problem2Deps}
\min \{(L_E(\Gamma_{c}))^{2}+\epsilon D(\Gamma_{c}) \,:\,  \Gamma_{c} \in \mathcal{C}_{D}, \quad G(c_{1})-P_{1}^{3}=G(c_{2})-P_{2}^{3}=G(c_{3})-P_{1}^{3}\} .
\end{equation}
To prove existence of a minimizer, we apply standard direct methods of calculus of variations. More precisely, let 
$\Gamma_{c_{n}} \in \mathcal{C}_{D}$, $n \in \N$, be a minimizing sequence satisfying the integral constraints given in \eqref{problem2Deps}. 
Without loss of generality, we may assume that
\[
(L_E(\Gamma_{c_{n}}))^{2}+\epsilon D(\Gamma_{c_{n}}) \leq (L_E(\Gamma_{c_{1}}))^{2}+\epsilon D(\Gamma_{c_{1}}) \leq  (L_E(\Gamma_{c_{1}}))^{2}+ D(\Gamma_{c_{1}})=: C_{0} \qquad \forall n \in \N.
\]
It follows that
\[
\| \partial_{u}(c_{n,\alpha})\|_{L^{2}(0,1)} \leq C(\epsilon) \qquad \text{for all } n \in \N, \alpha=1,2,3,
\]
and, using the fact that the fundamental theorem of calculus holds for maps in $H^{1,m}(0,1)$, $m \geq 1$, and
\begin{align}\label{vediqui}
|c_{n,\alpha}(u)| &\leq |c_{n,\alpha}(u)-c_{n,\alpha}(1)| +|P_{\alpha}|
\leq  |P_{\alpha}|+ \int_{0}^{1} |\partial_{u}(c_{n,\alpha})(u)| du \notag \\
& \leq
|P_{\alpha}|+ \| \partial_{u}(c_{n,\alpha})\|_{L^{2}(0,1)} \leq C(\epsilon)  \text{ for all } u \in [0,1], 
\end{align}
we obtain
\[
\|c_{n,\alpha}\|_{C([0,1])} \leq C(\epsilon)  \text{ for any } n \in \N, \quad \alpha=1,2,3, 
\]
so that
\[
\|c_{n,\alpha}\|_{H^{1,2}(0,1)} \leq C(\epsilon)  \text{ for any } n \in \N, \quad \alpha=1,2,3, 
\]
yielding the existence of a subsequence converging weakly in $H^{1,2}(0,1)$.
Moreover, for any $u_{1},u_{2} \in [0,1]$ we have that
\[
|c_{n,\alpha}(u_{1}) -c_{n,\alpha}(u_{2})|\leq \|\partial_{u} (c_{n,\alpha})\|_{L^{2}(0,1)}|u_{1}-u_{2}|^{\frac{1}{2}} \leq C(\epsilon)  |u_{1}-u_{2}|^{\frac{1}{2}}  \text{ for all } n \in \N, \quad \alpha=1,2,3, 
\]
so that using the Theorem of Ascoli--Arzel\`a we infer the existence of $\Gamma_{c_{\epsilon}} \in \mathcal{C}_{D} \cap C^{0,\frac{1}{2}}([0,1])$ to which a subsequence of $\Gamma_{c_{n}}$ converges uniformly on $[0,1]$, as well as weakly in $H^{1,2}(0,1)$ (where the convergence is meant componentwise).
It is easy to verify that the triod $\Gamma_{c_{\epsilon}}$ satisfies also the constraints given in \eqref{problem2Deps}.

Using that the Dirichlet functional is sequentially lower semicontinuous with respect to weak convergence in $H^{1,2}$, and that the length functional $L_{E}$ is sequentially lower semicontinuous with respect to weak convergence in $H^{1,m}$ for any $m \geq 1$ (cf.\ for instance \cite[Theorem~1.3]{Dacorogna} or \cite[Theorem~3.5, Remark~2]{BGH}),
 we infer that $\Gamma_{c_{\epsilon}}$ is a minimum for the problem formulated in \eqref{problem2Deps}, satisfying
\[
(L_E(\Gamma_{c_{\epsilon}}))^{2} +\epsilon D(\Gamma_{c_{\epsilon}}) \leq C_{0}.
\]
 
Next, since  the functionals $L_E$ and $G$ are invariant under smooth reparametrizations that keep the orientation of the curve $c$, we infer by taking inner variations (see for instance \cite[Proposition~1.14, Remark~3]{BGH}) for any single component map of the triod $\Gamma_{c_{\epsilon}}$ that $|(c_{\epsilon, \alpha})_{u}|=constant$ for any $\alpha=1,2,3$, yielding that we have even more regularity, namely
\[
\Gamma_{c_{\epsilon}} \in \mathcal{C}_{D} \cap C^{0,1}([0,1]).
\]

It could happen that the junction point of the minimal triod  $\Gamma_{c_{\epsilon}}$ coincides with one of the points~$\overline{P}_{\alpha}$, and its corresponding component could also reduce to a point. In other words, it might be convenient in the minimization procedure for one curve to disappear entirely by mapping the interval $[0,1]$ to a single point. 
In any case, the minimal configuration (may it be degenerate or not) is a quasinormal configuration according to Definition~\ref{def:quasino}.

Summarizing our findings so far we can say that (for an arbitrary $0< \epsilon \leq 1$)
\begin{align*}
&
(L_E(\Gamma_{c_{\epsilon}}))^{2} +\epsilon D(\Gamma_{c_{\epsilon}})\\ & \quad
=\inf \{  (L_E(\Gamma_{{\tilde{c}}}))^{2}+\epsilon D(\Gamma_{\tilde{c}}) \,:\,  \Gamma_{\tilde{c}} \in \mathcal{C}_{D}, \quad G(\tilde{c}_{1})-P_{1}^{3}=G(\tilde{c}_{2})-P_{2}^{3}=G(\tilde{c}_{3})-P_{1}^{3}\} \\ &\quad
\leq \inf \{( L_E(\Gamma_{{\tilde{c}}}))^{2}+\epsilon D(\Gamma_{\tilde{c}}) \,:\,  \Gamma_{\tilde{c}} \in \mathcal{C}_{D}\cap C^{0,1},  \Gamma_{\tilde{c}} \text{ quasinormal}, \\ & \qquad\qquad\quad 
G(\tilde{c}_{1})-P_{1}^{3}=G(\tilde{c}_{2})-P_{2}^{3}=G(\tilde{c}_{3})-P_{1}^{3}\} \\ & \quad
\leq (L_E(\Gamma_{c_{\epsilon}}))^{2} +\epsilon D(\Gamma_{c_{\epsilon}})  \leq C_{0}. 
\end{align*}

In the next step we would like to send  the parameter $\epsilon $ to zero. To that end, using the quasinormality of the minima, we first observe that
\begin{align*}
(2+\epsilon)D(\Gamma_{c_{\epsilon}}) &=  (\sum_{\alpha=1}^{3} \int_{0}^{1} |(c_{\epsilon})_{u}|^{2} du) + \epsilon D(\Gamma_{c_{\epsilon}}) =(\sum_{\alpha=1}^{3}  |(c_{\epsilon})_{u}|^{2} ) + \epsilon D(\Gamma_{c_{\epsilon}})\\
& \leq (\sum_{\alpha=1}^{3}  |(c_{\epsilon})_{u}| )^{2} + \epsilon D(\Gamma_{c_{\epsilon}})=(\sum_{\alpha=1}^{3} \int_{0}^{1} |(c_{\epsilon})_{u}|  du)^{2} + \epsilon D(\Gamma_{c_{\epsilon}})\\
&=(L_E(\Gamma_{c_{\epsilon}}))^{2} +\epsilon D(\Gamma_{c_{\epsilon}})\leq
C_{0},
\end{align*}
from which we infer that
\begin{equation} \label{bounddiri}
D(\Gamma_{c_{\epsilon}}) \leq  C_{0}
\end{equation}
for all $0 < \epsilon \leq 1$.
Using the quasinormality we even have
\begin{equation} \label{boundderivative}
|(c_{\epsilon, \alpha})_{u}|^{2}=constant \leq 2C_{0} 
\end{equation}
for all $0 < \epsilon \leq 1$ and $\alpha=1,2,3$.
Repeating the same argument as above (cf.\ \eqref{vediqui}), replacing $\Gamma_{c_{n}}$ with $\Gamma_{c_{\epsilon}}$, we obtain
\[
\|c_{\epsilon,\alpha}\|_{H^{1,2}(0,1)} \leq C  \text{ for any } 0<\epsilon < 1, \,\alpha=1,2,3, 
\]
for a constant $C$ independent of $\epsilon$.
Moreover, for  any  $u_{1}, u_{2} \in [0,1]$, $u_{1} \leq u_{2}$, we have (using \eqref{boundderivative})
\[
|c_{\epsilon,\alpha}(u_{2}) -c_{\epsilon,\alpha}(u_{1})|\leq \int_{u_{1}}^{u_{2}} |\partial_{u} (c_{\epsilon,\alpha})| du = | (c_{\epsilon,\alpha})_{u}||u_{1}-u_{2}| \leq C  |u_{1}-u_{2}|   \quad \alpha=1,2,3, 
\]
 and for all $\epsilon \in (0,1)$.
Again using the Theorem of Ascoli--Arzel\`a we infer the existence of 
\[
\Gamma_{c} \in \mathcal{C}_{D} \cap C^{0,1}([0,1])
\]
to which a subsequence of $\Gamma_{c_{\epsilon}}$ converges uniformly on $[0,1]$ as well as weakly in $H^{1,2}(0,1)$.
The triod $\Gamma_{c}$ satisfies also the integral constraints given in \eqref{problem2Deps}.

Observe that \begin{align*}
0 \leq d:&=\inf \{(L_E( \Gamma_{\tilde{c}}))^{2} \,:\,  \Gamma_{\tilde{c}} \in \mathcal{C}_{D}, \quad G(\tilde{c}_{1})-P_{1}^{3}=G(\tilde{c}_{2})-P_{2}^{3}=G(\tilde{c}_{3})-P_{1}^{3}\} \\
&\leq\inf \{  (L_E(\Gamma_{{\tilde{c}}}))^{2}+\epsilon D(\Gamma_{\tilde{c}}) \,:\,  \Gamma_{\tilde{c}} \in \mathcal{C}_{D}, \quad G(\tilde{c}_{1})-P_{1}^{3}=G(\tilde{c}_{2})-P_{2}^{3}=G(\tilde{c}_{3})-P_{1}^{3}\}\\
&=(L_E( \Gamma_{c_{\epsilon}}))^{2} +\epsilon D(\Gamma_{c_{\epsilon}}) =:d(\epsilon).
\end{align*}
Since the map $[0,1] \ni \epsilon \mapsto d(\epsilon)$ is non decreasing, it has a limit as $\epsilon$ goes to zero and there holds
\[
d \leq \lim_{\epsilon \to 0} d(\epsilon).
\]
On the other hand, using \eqref{bounddiri}, we obtain for the limit of $d(\epsilon)$ that 
\begin{equation} \label{p1}
\lim_{\epsilon \to 0} d(\epsilon)=\lim_{\epsilon \to 0}\,  (L_E(\Gamma_{c_{\epsilon}}))^{2} +\epsilon D(\Gamma_{c_{\epsilon}})= \lim_{\epsilon \to 0} (L_E(\Gamma_{c_{\epsilon}}))^{2}. 
\end{equation}
Taking a sequence $\epsilon_{j} \to 0$ for which $\Gamma(c_{\epsilon_{j}})$ converges weakly to $\Gamma_{c}$ in $H^{1,2}(0,1)$ we infer using the sequentially weak lower semicontinuity of the length functional  that
\begin{equation} \label{p2}
d \leq (L_E(\Gamma_{c}))^{2} \leq \liminf_{\epsilon_{j} \to 0}  (L_E(\Gamma_{c_{\epsilon_{j}}}))^{2}
=\lim_{\epsilon \to 0} (L_E(\Gamma_{c_{\epsilon}}))^{2}. 
\end{equation}
On the other hand, using the minimality of $\Gamma_{c_{\epsilon}}$, we can write
\[
(L_E(\Gamma_{c_{\epsilon}}))^{2} +\epsilon D(\Gamma_{c_{\epsilon}}) \leq (L_E(\Gamma_{c}))^{2} +\epsilon D(\Gamma_{c}),
\]
so that 
\begin{equation} \label{p3}
\lim_{\epsilon \to 0} d(\epsilon) =\lim_{\epsilon \to 0} (L_E(\Gamma_{c_{\epsilon}}))^{2} +\epsilon D(\Gamma_{c_{\epsilon}}) \leq (L_E(\Gamma_{c}))^{2}.
\end{equation}
From \eqref{p2}, \eqref{p1}, \eqref{p3} we therefore obtain that
\[
d \leq (L_E(\Gamma_{c}))^{2} \leq \lim_{\epsilon \to 0} (L_E(\Gamma_{c_{\epsilon}}))^{2} =\lim_{\epsilon \to 0} d(\epsilon) \leq (L_E(\Gamma_{c}))^{2},
\]
so that
\begin{equation} \label{p4}
d \leq (L_E(\Gamma_{c}))^{2}=\lim_{\epsilon \to 0} d(\epsilon).
\end{equation}
On the other hand since $d(\epsilon) \leq (L_E(\Gamma_{\hat{c}}))^{2} + \epsilon D(\Gamma_{\hat{c}})$ for any admissible $\Gamma_{\hat{c}} \in \mathcal{C}_{D}$  satisfying the constraints given in \eqref{problem2Deps}, we also have that
\[
\lim_{\epsilon \to 0} d(\epsilon) \leq (L_E(\Gamma_{\hat{c}}))^{2} \quad \forall \, \Gamma_{\hat{c}} \in \mathcal{C}_{D} \text{ satisfying the constraints}.
\]
Hence
\begin{equation} \label{p5}
\lim_{\epsilon \to 0} d(\epsilon) \leq d =\inf \{(L_E(\Gamma_{\tilde{c}}))^{2} \,:\,  \Gamma_{\tilde{c}} \in \mathcal{C}_{D}, \quad G(\tilde{c}_{1})-P_{1}^{3}=G(\tilde{c}_{2})-P_{2}^{3}=G(\tilde{c}_{3})-P_{1}^{3}\} .
\end{equation}
From \eqref{p4} and \eqref{p5} we finally infer that
\[
 d= (L_E(\Gamma_{c}))^{2}=\lim_{\epsilon \to 0} d(\epsilon),
\]
in other words
\[
(L_E(\Gamma_{c}))^{2} =\inf \{(L_E(\Gamma_{\tilde{c}}))^{2} \,:\,  \Gamma_{\tilde{c}} \in \mathcal{C}_{D}, \quad G(\tilde{c}_{1})-P_{1}^{3}=G(\tilde{c}_{2})-P_{2}^{3}=G(\tilde{c}_{3})-P_{1}^{3}\}.
\]
It then follows that the triod $\Gamma_{c} \in \mathcal{C}_{D} \cap C^{0,1}([0,1])$ realizes also the following minimum problem
\begin{equation} \label{minimizers}
L_E( \Gamma_{c}) =\inf \{L_E( \Gamma_{\tilde{c}}) \,:\,  \Gamma_{\tilde{c}} \in \mathcal{C}_{D}, \quad G(\tilde{c}_{1})-P_{1}^{3}=G(\tilde{c}_{2})-P_{2}^{3}=G(\tilde{c}_{3})-P_{1}^{3}\},
\end{equation}
as well as
\begin{align*}
L_E( \Gamma_{c}) &=\inf \{L_E( \Gamma_{\tilde{c}}) \,:\,  \Gamma_{\tilde{c}} \in \mathcal{C}_{D}, \quad G(\tilde{c}_{1})-P_{1}^{3}=G(\tilde{c}_{2})-P_{2}^{3}=G(\tilde{c}_{3})-P_{1}^{3}\}\\
&=\inf \{L_E( \Gamma_{\tilde{c}}) \,:\,  \Gamma_{\tilde{c}} \in \mathcal{C}_{D}\cap C^{0,1}([0,1]), \quad G(\tilde{c}_{1})-P_{1}^{3}=G(\tilde{c}_{2})-P_{2}^{3}=G(\tilde{c}_{3})-P_{1}^{3}\}\\
&=\inf \{L_E( \Gamma_{\tilde{c}}) \,:\,  \Gamma_{\tilde{c}} \in \mathcal{C} \cap C^{0,1}([0,1]), \quad G(\tilde{c}_{1})-P_{1}^{3}=G(\tilde{c}_{2})-P_{2}^{3}=G(\tilde{c}_{3})-P_{1}^{3}\}.
\end{align*}

\subsection{Euler--Lagrange equations and regularity of minimizers}

For simplicity of exposition later on, let us first remark that quasinormal critical points for a linear combination of the functionals $L_{E}$ and $G$ are automatically smooth.
\begin{lemma}\label{regmin}
Let $\lambda \in \R$ and let 
$c \in C^{0,1}([0,1], \R^{2})$ be a quasinormal curve satisfying
\[
\frac{d}{d\epsilon}\Big|_{\epsilon=0} (L_{E}(c+\epsilon \phi)+ \lambda G(c+\epsilon \phi) ) =\delta L_{E}(c)\phi + \lambda \delta G(c)\phi=0 \qquad \forall \, \phi \in C_{0}^{\infty}((0,1), \R^{2}).
\]
Then $c \in C^{\infty}([0,1], \R^{2})$.
\end{lemma}
\begin{proof} Recalling \eqref{varG} we have that
\begin{equation} \label{ELweakunacurva}
0=\int_{0}^{1} \frac{c_{u}}{|c_{u}|} \cdot \phi_{u} du + \lambda \int_{0}^{1} \begin{pmatrix}
 \phi^{1}\\ \phi^{2}
\end{pmatrix}  \cdot \begin{pmatrix}
 c^{2}_{u}\\ -c^{1}_{u}
\end{pmatrix} 
du \qquad \forall \, \phi \in C_{0}^{\infty}((0,1), \R^{2}).
\end{equation}
Since $|c_{u}|=L_E(c)>0$ by assumption we obtain for each space component $j=1,2$ that
\[
\left| \int_{0}^{1} \partial_{u}c^{j} \, \partial_{u} \phi \, du \right | \leq C(c, \lambda) \|\phi \|_{L^{2}(0,1)}  \qquad \forall \, \phi \in C_{0}^{\infty}((0,1), \R),
\]
yielding $c \in H^{2,2}((0,1), \R^{2})$ by the characterization of Sobolev spaces (see for instance \cite[Theorem~2.9]{BGH}). Then, integrating by parts in \eqref{ELweakunacurva} and using the fundamental lemma of calculus of variation we infer (writing $ds= L_E(c) du$ for the length element) that
\[
\partial_{ss}c^{1} = \lambda c^{2}_{s}, \qquad \partial_{ss}c^{2} = -\lambda c^{1}_{s}
\]
and a boot-strap argument yields the claim.
\end{proof}

\begin{rem}Even without having calculated the Euler--Lagrange equation yet, we may have already some intuition about the shape of minimal triods.  Let  $\Gamma_{c} \in \mathcal{C}_{D} \cap C^{0,1}([0,1])$  be a  minimizer satisfying \eqref{minimizers} and let $\Gamma=(\gamma_{1},\gamma_{2},\gamma_{3})$ be its horizontal lift (with each curve $c_{\alpha}$ lifted to a horizontal curve $\gamma_{\alpha}$ ending at $P_{\alpha}$, $\alpha=1,2,3$).

Without loss of generality, i.e.\ after an appropriate left translation of $\Gamma$, we might assume that the junction point $\Sigma$ coincides with the origin (recall Remark~\ref{rem1.4}). Moreover, let us still denote by $P_{\alpha}$ the translated endpoints.

We immediately observe that if the points $P_{\alpha}$ lie in the $(x,y)$-plane, i.e.\ $P_{\alpha}^{3}=0$ for any $\alpha=1,2,3$,  then, since the shortest path between $P_{\alpha}$ and $\Sigma=0$ is a straight line, our minimal triod $\Gamma$ must be a solution of the classical Steiner problem. This means that if the largest angle of the triangle with vertices given by the $P_{\alpha}$ is larger than or equal to $120$ degrees, then $\Sigma$ coincides with the 
vertex of the largest angle. Otherwise $\Sigma$ is a point distinct from the points $P_{\alpha}$ and the three lines meet in $\Sigma$ at an angle of $120$ degrees (see for example \cite[\S5]{GP68}).
 
 On the other hand, if at least one of the points $P_{\alpha}$ does not lie in the $(x,y)$-plane, then the situation is not so clear anymore. But (in view of Remark~\ref{rem-EL} and Remark~\ref{rem1.10}) it is reasonable to believe that the projection of a minimal configuration $\Gamma_{c}$ will be composed of arcs of circles and/or straight segments.
\end{rem}
 
In the following let $\Gamma_{c} \in \mathcal{C}_{D} \cap C^{0,1}([0,1])$  be a  minimizer satisfying \eqref{minimizers} and let $\Gamma=(\gamma_{1},\gamma_{2},\gamma_{3})$ be its horizontal lift (with each curve $c_{\alpha}$ lifted to a horizontal curve $\gamma_{\alpha}$ ending at $P_{\alpha}$, $\alpha=1,2,3$).
For simplicity of exposition we assume here that the triod is not degenerate, i.e.\ no curves have collapsed to a point.
 
Upon a suitable reparametrization (\cite[Lemma~5.23]{BGH}), we may assume furthermore that the triod is quasinormal (according to Definition~\ref{def:quasino}).

To calculate the Euler--Lagrange equations denote by
\[
\mathcal{T}:=\{\Phi=(\phi_{1},\phi_{2},\phi_{3}) \,:\, \phi_{\alpha} \in C^{\infty}([0,1],\R^{2}), 
 \phi_{\alpha}(1)=(0,0)^{t}, \alpha=1,2,3, \quad \phi_{1}(0)=\phi_{2}(0)=\phi_{3}(0)\}
\]
the space of test functions, and by
\begin{align*}
\mathcal{G}_{1}(\Gamma_{\tilde{c}}):=G(\tilde{c}_{1})-G(\tilde{c}_{2})-P^{3}_{1}+P^{3}_{2},\\
\mathcal{G}_{2}(\Gamma_{\tilde{c}}):=G(\tilde{c}_{2})-G(\tilde{c}_{3})-P^{3}_{2}+P^{3}_{1},
\end{align*}
the two constraints given in \eqref{minimizers}. Choose $\varphi_{i} \in C^{\infty}_{0}((0,1), \R^{2})$, $i=1,3$,  such that
\[
0\neq \delta G(c_{i})(\varphi_{i})=\int_{0}^{1} \begin{pmatrix}
 \varphi_{i}^{1}\\ \varphi_{i}^{2}
\end{pmatrix}  \cdot \begin{pmatrix}
 (c^{2}_{i})_{u}\\ -(c^{1}_{i})_{u}
\end{pmatrix} 
du
\]
(which exists thanks to the regularity of the curves $c_{i}$) and set
\[
\Phi_{1}:=(\varphi_{1},0,0), \qquad \Phi_{3}:=(0,0,\varphi_{3}).
\]
Next, for $\Phi \in \mathcal{T}$, consider the map
\[
\R^{3 } \ni (\epsilon,t_{1},t_{2})^t \mapsto \mathcal{G}(\epsilon,t_{1},t_{2}):= \begin{pmatrix}
\mathcal{G}_{1}(\Gamma_{c}+ \epsilon \Phi + t_{1}\Phi_{1}+ t_{2}\Phi_{3})\\
\mathcal{G}_{2}(\Gamma_{c}+ \epsilon \Phi + t_{1}\Phi_{1}+ t_{2}\Phi_{3})
 \end{pmatrix}=
 \begin{pmatrix}
\mathcal{G}_{1}(\Gamma_{c}+ \epsilon \Phi + t_{1}\Phi_{1})\\
\mathcal{G}_{2}(\Gamma_{c}+ \epsilon \Phi + t_{2}\Phi_{3})
 \end{pmatrix} ,
\]
where in the last equality we have used the definition of $\mathcal{G}_{i}$, $i=1,2$. Note also that the triod $\Gamma_{c}+ \epsilon \Phi + t_{1}\Phi_{1}+ t_{2}\Phi_{3}$ belongs to $\mathcal{C}_{D}$ for any choice of $\epsilon,t_{1},t_{2}$.

We have $\mathcal{G}(0,0,0)=(0,0)^{t}$ by assumption and 
\[
\frac{\partial \mathcal{G}}{\partial t}(0):=\frac{\partial  \mathcal{G}}{\partial (t_{1},t_{2})} (0,0,0)= 
\begin{pmatrix}
\delta G(c_{1})(\varphi_{1}) & 0 \\
0& -\delta G(c_{3})(\varphi_{3})
\end{pmatrix}
\]
has maximal rank. By the implicit function theorem we can find $C^{1}$-maps $\sigma_{r}$, $r=1,2$, defined in a neighbourhood of zero such that $\sigma_{r}(0)=0$ and
\begin{equation} \label{eqGeps}
\mathcal{G}(\epsilon, \sigma_{1}(\epsilon), \sigma_{2}(\epsilon)) =(0,0)^{t}, \qquad \text{for } \epsilon \in (-\epsilon_{0}, \epsilon_{0}),
\end{equation}
which implies that $\Gamma_{c}+\epsilon \Phi+  \sigma_{1}(\epsilon)\Phi_{1}+ \sigma_{2}(\epsilon)\Phi_{3}$ is an admissible variation for any $\epsilon \in (-\epsilon_{0}, \epsilon_{0})$.
Differentiation of the above system \eqref{eqGeps} yields
\[
\begin{pmatrix}
\sigma_{1}'(0) \\
\sigma_{2}'(0)
\end{pmatrix}= - (\frac{\partial \mathcal{G}}{\partial t}(0))^{-1}
\begin{pmatrix}
\delta \mathcal{G}_{1}(\Gamma_{c}) \Phi  \\
\delta \mathcal{G}_{2} (\Gamma_{c}) \Phi
\end{pmatrix}.
\]
From the minimality of $\Gamma_{c}$ we infer
\[
0=\frac{d}{d \epsilon}\Big|_{\epsilon =0} L_E(\Gamma_{c}+\epsilon \Phi+  \sigma_{1}(\epsilon)\Phi_{1}+ \sigma_{2}(\epsilon)\Phi_{3})=\delta L_E(\Gamma_{c})\Phi + \sigma_{1}'(0)\delta L_E(\Gamma_{c})\Phi_{1} +\sigma_{2}'(0)\delta L_E(\Gamma_{c})\Phi_{3}
\]
and substituting the values for $\sigma_{r}'(0)$ we finally infer that 
there exists Lagrange multipliers $\lambda_{j} \in \R$, $j=1,2$, such that
\[
0=\frac{d}{d \epsilon}\Big|_{\epsilon =0}  \left( L_E(\Gamma_{c}+\epsilon \Phi) + \sum_{j=1}^{2}\lambda_{j}\mathcal{G}_{j}(\Gamma_{c}+\epsilon \Phi) \right)
\]
for any $\phi \in \mathcal{T}$.
In other words, for any $\Phi=(\phi_{1},\phi_{2},\phi_{3}) \in \mathcal{T}$  we have that
\begin{equation} \label{EL-triods}
0=\sum_{i=1}^{3}\delta L_{E}(c_{i})\phi_{i} + \lambda_{1} \delta G(c_{1})\phi_{1} - \lambda_{2}\delta G(c_{3})\phi_{3} + (\lambda_{2}-\lambda_{1}) \delta G(c_{2})\phi_{2}.
\end{equation}
This yields that
\begin{subequations} \label{eq:ELc}
\begin{align} 
0&=\delta L_{E}(c_{1})\phi + \lambda_{1} \delta G(c_{1})\phi,\\  
0&=\delta L_{E}(c_{2})\phi + (\lambda_{2}-\lambda_{1}) \delta G(c_{2})\phi,\\ 
0&=\delta L_{E}(c_{3})\phi - \lambda_{2} \delta G(c_{3})\phi,
\end{align}
\end{subequations}
for any $\varphi\in C_{0}^{\infty}((0,1),\R^{2})$, which in turns yields that the curves are smooth thanks to Lemma~\ref{regmin} and (recalling the computations performed in Remark~\ref{rem-EL}) that the curvatures $k_{\alpha}$ of the three curves $c_{\alpha}$ are constant and equal to:
\begin{align*}
&k_{1}(u)+\lambda_{1}=0 \quad \text{ for } u \in (0,1),\\
&k_{2}(u)+(\lambda_{2}-\lambda_{1})=0 \quad \text{ for } u \in (0,1),\\
&k_{3}(u)-\lambda_{2}=0 \quad \text{ for } u \in (0,1),
\end{align*}
in other words
\begin{equation} \label{sumk}
\sum_{\alpha=1}^{3}k_{\alpha}=0.
\end{equation}
Further we infer from \eqref{EL-triods} (considering now test functions $\Phi \in \mathcal{T}$ that do not vanish near the triple junction so that we may retrieve the natural boundary conditions, using \eqref{eq:ELc}, and recalling \eqref{varG}) that
\begin{align*}
0&=\sum_{\alpha=1}^{3} \left[ \begin{pmatrix}
\phi_\alpha^{1}\\ \phi_\alpha^{2}
\end{pmatrix}  \cdot \begin{pmatrix}
(c^{1}_\alpha)_{u}\\ (c^{2}_\alpha)_{u}
\end{pmatrix} \frac{1}{|(c_\alpha)_{u}|} \right ]_{0}^{1} \\
&+\frac{\lambda_{1}}{2}
\left[ \begin{pmatrix}
\phi_{1}^{1}\\ \phi_{1}^{2}
\end{pmatrix}  \cdot \begin{pmatrix}
-c^{2}_{1}\\ c^{1}_{1}
\end{pmatrix}  \right ]_{0}^{1}
+\frac{\lambda_{2}-\lambda_{1}}{2}
\left[ \begin{pmatrix}
\phi_{2}^{1}\\ \phi_{2}^{2}
\end{pmatrix}  \cdot \begin{pmatrix}
-c^{2}_{2}\\ c^{1}_{2}
\end{pmatrix}  \right ]_{0}^{1}
-\frac{\lambda_{2}}{2}
\left[ \begin{pmatrix}
\phi_{3}^{1}\\ \phi_{3}^{2}
\end{pmatrix}  \cdot \begin{pmatrix}
-c^{2}_{3}\\ c^{1}_{3}
\end{pmatrix}  \right ]_{0}^{1}\\
&=-\phi_{1}(0) \cdot \left\{ \left(\sum_{\alpha=1}^{3} \begin{pmatrix}
(c^{1}_\alpha)_{u}\\ (c^{2}_\alpha)_{u}
\end{pmatrix} \frac{1}{|(c_\alpha)_{u}|}  \right) +\frac{\lambda_{1}}{2} \begin{pmatrix}
-c^{2}_{1}\\ c^{1}_{1}
\end{pmatrix} 
+\frac{\lambda_{2}-\lambda_{1}}{2} \begin{pmatrix}
-c^{2}_{2}\\ c^{1}_{2}
\end{pmatrix} 
-\frac{\lambda_{2}}{2} \begin{pmatrix}
-c^{2}_{3}\\ c^{1}_{3}
\end{pmatrix} 
\right\}\Big|_{u=0} \\
&=-\phi_{1}(0) \cdot \left\{ \left(\sum_{\alpha=1}^{3} \begin{pmatrix}
(c^{1}_\alpha)_{u}\\ (c^{2}_\alpha)_{u}
\end{pmatrix} \frac{1}{|(c_\alpha)_{u}|}  \right) \right \}\Big|_{u=0} ,
\end{align*}
where  we have used the facts that $\phi_{1}(0)=\phi_{2}(0)=\phi_{3}(0)$  and $c_{1}(0)=c_{2}(0)=c_{3}(0)$.
In particular, we infer that the projected curves meet at the projected triple junction at angles of 120 degrees, i.e.
\begin{equation} \label{angle-cond}
0=\sum_{\alpha=1}^{3} \begin{pmatrix}
(c^{1}_\alpha)_{u}\\ (c^{2}_\alpha)_{u}
\end{pmatrix} \frac{1}{|(c_\alpha)_{u}|}\Big|_{u=0} =\sum_{\alpha=1}^{3} \frac{(c_\alpha)_{u}}{|(c_\alpha)_{u}|}\Big|_{u=0}  .
\end{equation}

\begin{rem}
If the minimal triod $\Gamma$ is degenerate, we can repeat the above arguments with $P_{2}=\Sigma$, $G(c_{2})=0$, i.e.\ with $c_{2}([0,1])=\overline{P}_{2}$. Then the integral constraints become 
\[ 
P^{3}_{2}=P^{3}_{1} -G(\tilde{c}_{1})= P^{3}_{3} -G(\tilde{c}_{3})
\]
(which are essentially decoupled), providing the existence of some $\lambda_{1}, \lambda_{2} \in \R$ for which
\begin{align*}
0&=\delta L_{E}(c_{1})\phi + \lambda_{1} \delta G(c_{1})\phi,\\  
0&=\delta L_{E}(c_{3})\phi - \lambda_{2} \delta G(c_{3})\phi.
\end{align*}
for any $\varphi\in C_{0}^{\infty}((0,1),\R^{2})$, so that 
\[
k_{1}+\lambda_{1}=0=k_{3}-\lambda_{2}  \text{ on }[0,1]. 
\]
\end{rem}

\begin{rem}\label{remEX}
From \eqref{sumk} we infer that for a non-degenerate minimal triod $\Gamma_c$ (or more generally a non-degenerate solution of the Euler--Lagrange equation \eqref{EL-triods}) the following situations can occur:
\begin{itemize}
\item[(i)] $k_{1}=k_{2}= 0 \Rightarrow k_{3}=0$, i.e.\ three straight segments meet at 120 degree angles;
\item[(ii)]  $k_{1}=0$, $k_{2}=-k_{3} \neq 0$, i.e.\ one straight segment and two arcs of curvature equal in absolute  value meet at a projected angles of 120 degrees; 
\item[(iii)] $k_{1}, k_{2}, k_{3} \neq 0$ and \eqref{sumk}, \eqref{angle-cond} hold. 
\end{itemize}

Again let $\Gamma=(\gamma_{1},\gamma_{2},\gamma_{3})$ be the horizontal lift of $\Gamma_{c}$  (with each curve $c_{\alpha}$ lifted to a horizontal curve $\gamma_{\alpha}$ ending at $P_{\alpha}$, $\alpha=1,2,3$).
A Steiner minimal configuration of type (i) is given for example by connecting through straight lines the junction point $\Sigma=0$ with  $P_{1}=(1,0,0)^{t}$, $P_{2}= (\cos (2\pi/3), \sin (2\pi/3),0)^{t}=(-\frac{\sqrt{3}}{2},\frac{1}{2}, 0)^{t}$, and  $P_{3}= (\cos ( -2\pi/3), \sin (-2\pi/3),0)^{t}=(-\frac{\sqrt{3}}{2},-\frac{1}{2}, 0)^{t}$ (cf.\ also Experiment~1 in Section~\ref{sec:nr}).

A critical configuration  of type (ii) can be constructed, for instance, by modifying the previous example as follows: we take the junction point  $\Sigma=0$  to be connected with a straight line with $P_{1}=(1,0,0)^{t}$ (here $\lambda_{1}=0$). 
The two remaining points $P_{2}$ and $P_{3}$ shall now not lie in the $(x,y)$-plane. To construct two horizontal geodesics meeting in $\Sigma$ at the right angle (and parametrized by arc length with  $L_{E}(c_{i})=s_{f}=1$ for simplicity) we use Remark~\ref{rem1.10}. For the curve $c_{2}$ connecting the origin to $P_{2}$, take for instance (using the notation adopted in Remark~\ref{rem1.10})  $c_{2}(0)=0$, $|B|=1$, $\alpha_{0}=\frac{2}{3}\pi$, so that by \eqref{eqQ3}, \eqref{eqQ3bis} with $\lambda_{2}=\pi$
\[
P_{2}^{3}= -\frac{1}{2\lambda}=-\frac{1}{2\pi},\quad  \qquad \gamma^{3}_{2}(s)=-\frac{1}{2\pi} s +\frac{1}{2\pi^{2}}\sin(\pi s), \quad s \in [0,1].  
\]
More precisely, the curve $c_{2}:[0,1] \to \R^{2}$ is parametrized by
\[
c_{2}(s)=\frac{1}{\pi} \begin{pmatrix} \sin(\frac{2\pi}{3})\\ -\cos (\frac{2\pi}{3})\end{pmatrix} +
\frac{1}{\pi} \begin{pmatrix} \sin(\pi s -\frac{2\pi}{3})\\ \cos (\pi s -\frac{2\pi}{3})\end{pmatrix}, \quad s \in [0,1] .
\]
For the curve $c_{3}$ connecting the origin $\Sigma$ to $P_{3}$, we obtain with $c_{3}(0)=0$, $|B|=1$, $\alpha_{0}=-\frac{2}{3}\pi$,$s_{f}=1$, $\lambda=-\pi$
\[
P_{3}^{3}= -\frac{1}{2\lambda}= +\frac{1}{2\pi}, \qquad \gamma^{3}_{3}(s)=\frac{1}{2\pi} s +\frac{1}{2\pi^{2}}\sin(-\pi s), \quad s \in [0,1]. 
\]
More precisely, the curve $c_{3}:[0,1] \to \R^{2}$ is parametrized by
\[
c_{3}(s)=\frac{1}{\pi} \begin{pmatrix} \sin(\frac{2\pi}{3})\\ \cos (\frac{2\pi}{3})\end{pmatrix} +
\frac{1}{\pi} \begin{pmatrix} \sin(\pi s -\frac{2\pi}{3})\\ -\cos (\pi s -\frac{2\pi}{3})\end{pmatrix} , \quad s \in [0,1].
\]

Further critical configurations of type (ii) are retrieved numerically for instance in Experiments~2 and 5 in Section~\ref{sec:nr}. Critical configurations of type (iii) are shown for instance in Experiments~11 and~13.
\end{rem}

\subsection{Formulation of horizontal curve shortening flow for triods}
Similarly to what is discussed in Appendix~\ref{sec:HCSF}, we can provide here evolution problems  whose stationary solutions yield critical points for the functional $L$ as given in \eqref{Ltriods}.
The following formulations are at the basis of the numerical schemes and simulations proposed in Section~\ref{sec:FEM} and~\ref{sec:nr}.
\begin{ass}\label{ass2}
Let $P_{1},P_{2},P_{3} \in \R^{3}$ be three given distinct points and let $\gamma_{0, \alpha}: [0,1] \to \R^{3}$, $\alpha=1,2,3$, be regular smooth horizontal curves such that 
\begin{align*}
\gamma_{0,\alpha}(1)&=P_{\alpha}, \qquad \alpha=1,2,3,\\
\gamma_{0,1}(0)&=\gamma_{0,2}(0)=\gamma_{0,3}(0). 
\end{align*}
 \end{ass}

Our aim now is to let such an initial  triod evolve in time in such a way that its flow converges to critical points of $L$. Since the junction point must be able to move in time, tangential movements must be allowed. Recalling \eqref{ONframenew} we formulate:
\begin{prob}\label{HCSF1triod}
Let Assumption~\ref{ass2} hold (together with appropriate compatibility conditions). 
Find  $\Gamma=(\gamma_{1},\gamma_{2},\gamma_{3})$ with
$\gamma_{\alpha}:[0,T) \times [0,1] \to \R^{3}$, $\gamma_{\alpha}=\gamma_{\alpha}(t,u)$, $\alpha=1,2,3$ (in an appropriate class of  functions), which  satisfies $\gamma_{\alpha}(0, \cdot)=\gamma_{0,\alpha}(\cdot)$, and
\begin{align*}
\partial_{t } \gamma_{1}^{\perp_{g}}(t,u)&=\vec{k}_{1,g}(t,u)+\lambda_{1}(t)N_{1}(t,u) +\Big(\int_{u}^{1}(k_{1}(t,\xi)+\lambda_{1}(t)) |\partial_{u}\gamma_{1}(t,\xi)|_{g}d\xi \Big)X_{3}(\gamma_{1}(t,u)), \\ 
\partial_{t } \gamma_{2}^{\perp_{g}}(t,u)&=\vec{k}_{2,g}(t,u)+(\lambda_{2}(t)-\lambda_{1}(t))N_{2}(t,u) \\ & \qquad +\Big(\int_{u}^{1}(k_{2}(t,\xi)+\lambda_{2}(t)-\lambda_{1}(t)) |\partial_{u}\gamma_{2}(t,\xi)|_{g}d\xi \Big)X_{3}(\gamma_{2}(t,u)),\\
\partial_{t } \gamma_{3}^{\perp_{g}}(t,u)&=\vec{k}_{3,g}(t,u)-\lambda_{2}(t)N_{3}(t,u) +\Big(\int_{u}^{1}(k_{3}(t,\xi)-\lambda_{2}(t)) |\partial_{u}\gamma_{3}(t,\xi)|_{g}d\xi \Big)X_{3}(\gamma_{3}(t,u)), \\
\gamma_{1}(t,0)&=\gamma_{2}(t,0)=\gamma_{3}(t,0)\quad \text{for all } t \in (0,T),\\
\gamma_{\alpha}(t,1)&=P_{\alpha} \quad \text{for all } t \in (0,T), \quad \alpha=1,2,3,\\
0&=\sum_{\alpha=1}^{3} \begin{pmatrix}
(\gamma^{1}_{\alpha})_{u}(t,u)\\ (\gamma^{2}_{\alpha})_{u}(t,u)
\end{pmatrix} \frac{1}{ |\partial_{u}\gamma_{\alpha}(t,u)|_{g}}\Big|_{u=0} \text{ for  all } t \in [0,T] ,
\end{align*}
where the Lagrange multipliers are given through \eqref{systemLM}, below, for the projected curves $c_{\alpha}$.
Here $\partial_{t } \gamma_{\alpha}^{\perp_{g}}$ is short for  $ \partial_{t } \gamma_{\alpha}^{\perp_{g}}:=\partial_{t } \gamma_{\alpha}-g(\partial_{t } \gamma_{\alpha},T_{\alpha})T_{\alpha}$.
\end{prob}

If $\Gamma=(\gamma_{1},\gamma_{2},\gamma_{3})$ is a solution to Problem~\ref{HCSF1triod}, then its projected triod $\Gamma_{c}=(c_{1},c_{2},c_{3})$, with $c_{\alpha}=(\gamma_{\alpha}^{1}, \gamma_{\alpha}^{2})^{t}$, $\alpha=1,2,3$,
 satisfies the following problem in the Euclidean setting.
\begin{prob}\label{HCSF1ptriod}
Let Assumption~\ref{ass2} hold (together with appropriate compatibility conditions).
Find $\Gamma_{c}=(c_{1},c_{2},c_{3})$, with  $c_{\alpha}:[0,T) \times [0,1] \to \R^{2}$, $c_{\alpha}=(\gamma^{1}_{\alpha}, \gamma^{2}_{\alpha})^{t}$, $\alpha=1,2,3$,  (in an appropriate class of  functions), which  satisfies $\gamma^{r}_{\alpha}(0, \cdot)=\gamma_{0,\alpha}^{r}$ for $r=1,2$, $\alpha=1,2,3$, and
\begin{subequations}\label{eq:strong}
\begin{align}\label{eq:strong-a}
(\partial_{t}c_{1}\cdot \vec{n}_{1} ) \vec{n}_{1}&=\vec{k}_{1}+\lambda_{1}(t)\vec{n}_{1}  \qquad &\text{ for } (t,u)\in (0,T) \times (0,1), \\ \label{eq:strong-b}
(\partial_{t}c_{2} \cdot \vec{n}_{2}) \vec{n}_{2}&=\vec{k}_{2}+(\lambda_{2}(t)-\lambda_{1}(t))\vec{n}_{2}  \qquad &\text{ for } (t,u)\in (0,T) \times (0,1),\\  \label{eq:strong-c}
(\partial_{t}c_{3} \cdot \vec{n}_{3}) \vec{n}_{3}&=\vec{k}_{3}-\lambda_{2}(t)\vec{n}_{3}  \qquad &\text{ for } (t,u)\in (0,T) \times (0,1),\\  \label{eq:strong-d}
c_{\alpha}(t,1)&=(P^{1}_{\alpha},P^{2}_{\alpha})^{t} \quad &\text{for all } t \in [0,T),\\  \label{eq:strong-e}
c_{1}(t,0)&= c_{2}(t,0)=c_{3}(t,0)\quad &\text{for all } t \in (0,T),\\  \label{eq:strong-f}
0&=\sum_{\alpha=1}^{3} \frac{\partial_{u} c_{\alpha}}{|\partial_{u} c_{\alpha}|}\Big|_{u=0} &\text{ for  all } t \in (0,T) ,  
\end{align}
\end{subequations}
where
\[
\vec{n}_{\alpha}:=\frac{1}{|\partial_{u} c_{\alpha}|}\begin{pmatrix}
-\partial_{u}\gamma^{2}_{\alpha}(t,u)\\ \partial_{u} \gamma^{1}_{\alpha}(t,u) 
\end{pmatrix}, \qquad 
I_{K}(c_{\alpha}(t)):=\int_{0}^{1} k_{\alpha}(t,u) |\partial_{u} c_{\alpha}| du, \, \quad \alpha=1,2,3, 
\]
and the Lagrange multipliers $\lambda_{1}$, $\lambda_{2}$ are solutions of
\begin{equation} \label{systemLM}
\begin{pmatrix}
L_{E}(c_{1}(t)) + L_{E}(c_{2}(t)) & -L_{E}(c_{2}(t)) \\
-L_{E}(c_{2}(t)) & L_{E}(c_{3}(t)) + L_{E}(c_{2}(t)) 
\end{pmatrix}
\begin{pmatrix}
\lambda_{1}(t)\\ \lambda_{2}(t)
\end{pmatrix}=
\begin{pmatrix}
I_{K}(c_{2}(t)) -I_{K}(c_{1}(t))\\ I_{K}(c_{3}(t)) -I_{K}(c_{1}(t))
\end{pmatrix}.
\end{equation}
\end{prob}
Note that the Lagrange multipliers are well defined as long as the lengths of at least two curves are not zero, and that if $L_{E}(c_{2})=0$ then the system decouples into two fully independent equations with Lagrange multipliers that are computed as in Problem~\ref{HCSF1}.
\begin{lemma}
For a (sufficiently smooth) solution to Problem~\ref{HCSF1ptriod} we have that
\[
\frac{d}{dt} \sum_{\alpha=1}^{3} L_{E}(c_{\alpha}(t)) \leq 0. 
\]
\end{lemma}
\begin{proof}
First of all recall that for a planar curve
$c=(c^{1},c^{2})^{t}=(\gamma^{1}, \gamma^{2})^{t}$ such that $c(t,1)$ is a fixed point and $c(t,0)$ a mobile one, we have (recall \eqref{funG})
\[
\frac{d}{dt} G(c) = -\int_{0}^{1} c_{t} \cdot \vec{n} |c_{u}| du -\frac{1}{2} [c_{t} \cdot c^{\perp}]_{u=0}.
\]
Note that for the lifted curve $\gamma $ (recall \eqref{g3})   with $\gamma(t,1)$ a fixed point and $\gamma(t,0)$ a mobile one, the above equation implies
\begin{equation} \label{veljunction}
-\partial_{t} \gamma^{3} (t,0) = \frac{d}{dt} ( \gamma^{3} (t,1)-\gamma^{3} (t,0) )  =\frac{d}{dt} G(c) = -\int_{0}^{1} c_{t} \cdot \vec{n} |c_{u}| du -\frac{1}{2} [c_{t} \cdot c^{\perp}]_{u=0}.
\end{equation}
Using  $c_{1}(t,0)= c_{2}(t,0)=c_{3}(t,0)$ it follows then that
\begin{equation} \label{oo1}
\frac{d}{dt} G(c_{\alpha})= \frac{d}{dt} G(c_{\beta})  \Longleftrightarrow \int_{0}^{1} \partial_{t}c_{\alpha} \cdot \vec{n}_{\alpha} |\partial_{u}c_{\alpha}| du =\int_{0}^{1} \partial_{t}c_{\beta} \cdot \vec{n}_{\beta} |\partial_{u}c_{\beta}| du \qquad \text{ for } \alpha, \beta \in \{1,2,3\}.
\end{equation}
A straightforward calculation yields then  that
\begin{equation} \label{oo2}
\int_{0}^{1} \partial_{t}c_{1} \cdot \vec{n}_{1} |\partial_{u}c_{1}| du =\int_{0}^{1} \partial_{t}c_{2} \cdot \vec{n}_{2} |\partial_{u}c_{2}| du =\int_{0}^{1} \partial_{t}c_{3} \cdot \vec{n}_{3} |\partial_{u}c_{3}| du  \quad \Longleftrightarrow \quad \text{ \eqref{systemLM} holds.}
\end{equation}
Now that the effect of the choice of the Lagrange multipliers is clarified, we can compute
\begin{align*}
\frac{d}{dt} \sum_{\alpha=1}^{3} L_{E}(c_{\alpha}(t)) &= \sum_{\alpha=1}^{3} \left[\frac{\partial_{u}c_{\alpha}}{|\partial_{u} c_{\alpha}|} \cdot \partial_{t} c_{\alpha} \right]_{0}^{1} -\sum_{\alpha=1}^{3} \int_{0}^{1} \vec{k}_{\alpha} \cdot \partial_{t} c_{\alpha} |\partial_{u} c_{\alpha}| du\\
& =-\sum_{\alpha=1}^{3} \int_{0}^{1} \vec{k}_{\alpha} \cdot (\partial_{t} c_{\alpha} \cdot \vec{n}_{\alpha}) \vec{n}_{\alpha} |\partial_{u} c_{\alpha}| du\\
& =-\left(\sum_{\alpha=1}^{3} \int_{0}^{1} | (\partial_{t} c_{\alpha} \cdot \vec{n}_{\alpha}) \vec{n}_{\alpha}|^{2} |\partial_{u} c_{\alpha}| du \right)
+\lambda_{1}(t) \int_{0}^{1} \partial_{t}c_{1} \cdot \vec{n}_{1} |\partial_{u}c_{1}| du \\
& \qquad + (\lambda_{2}(t)-\lambda_{1}(t))\int_{0}^{1} \partial_{t}c_{2} \cdot \vec{n}_{2} |\partial_{u}c_{2}| du  -\lambda_{2}(t) \int_{0}^{1} \partial_{t}c_{3} \cdot \vec{n}_{3} |\partial_{u}c_{3}| du \\
&=-\sum_{\alpha=1}^{3} \int_{0}^{1} |(\partial_{t} c_{\alpha} \cdot \vec{n}_{\alpha}) \vec{n}_{\alpha}|^{2} |\partial_{u} c_{\alpha}| du  \leq 0,
\end{align*}
where we have used the boundary conditions and \eqref{oo2}.
\end{proof}
From the above lemma we immediately infer that for a solution of Problem~\ref{HCSF1triod} the length decreases along the flow since
(recall Remark~\ref{rem1.1})
\begin{align*}
\frac{d}{dt} L(\Gamma(t)) = \frac{d}{dt} \left( \sum_{\alpha=1}^{3} L_{E}(c_{\alpha}(t)) \right) \leq 0.
\end{align*}
With similar considerations explained in Section~\ref{sec:HCSF}, and using \eqref{oo1}, \eqref{veljunction}, \eqref{g3}, and Assumption~\ref{ass2}, one can show that a horizontal lift of a solution $\Gamma_c$ to Problem~\ref{HCSF1ptriod}, where each curve $c_\alpha$ is lifted to a horizontal curve $\gamma_\alpha$ ending at $P_\alpha$, yields a solution to Problem~\ref{HCSF1triod}.
More precisely for a lifted curve $\gamma$ we  observe by differentiating in time \eqref{g3} that
\begin{align*}
\partial_{t}\gamma^{3}(t,u) &=\partial_{t}\gamma^{3}(t,0) +
\int_{0}^{u} \left(-\frac{1}{2}\gamma^{2}_{t}\gamma^{1}_{u} -\frac{1}{2}\gamma^{2}\gamma^{1}_{tu} +\frac{1}{2}\gamma^{1}_{t}\gamma^{2}_{u} +\frac{1}{2}\gamma^{1}\gamma^{2}_{tu}\right) d\xi \\
& = \partial_{t}\gamma^{3}(t,0) +\int_{0}^{u} \left(-\gamma^{2}_{t}\gamma^{1}_{u}  +\gamma^{1}_{t}\gamma^{2}_{u} \right) d\xi + [\frac{1}{2} \gamma^{1}\gamma^{2}_{t} -\frac{1}{2} \gamma^{2}\gamma^{1}_{t}]_{0}^{u}\\
&=\partial_{t}\gamma^{3}(t,0) - \int_{0}^{u} \partial_{t} c \cdot \vec{n} |c_{u}| d\xi +\frac{1}{2}[c_{t} \cdot c^{\perp}]_{0}^{u} \\
& =\partial_{t}\gamma^{3}(t,0) -\frac{1}{2} [c_{t} \cdot c^{\perp}]_{u=0} - \int_{0}^{1} \partial_{t} c \cdot \vec{n} |c_{u}| d\xi\\
& \quad + \int_{u}^{1} \partial_{t} c \cdot \vec{n} |c_{u}| d\xi  + \frac{1}{2} (c_{t} \cdot c^{\perp})(t,u)\\
& =\int_{u}^{1} \partial_{t} c \cdot \vec{n} |c_{u}| d\xi  + (c_{t} \cdot \vec{n})\frac{1}{2} (\vec{n} \cdot c^{\perp})(t,u) + (c_{t} \cdot \frac{c_{u}}{|c_{u}|})\frac{1}{2} (\frac{c_{u}}{|c_{u}|} \cdot c^{\perp})(t,u)\\
&=\int_{u}^{1} \partial_{t} c \cdot \vec{n} |c_{u}| d\xi  + (c_{t} \cdot \vec{n}) N^{3}(t,u) +(c_{t} \cdot \frac{c_{u}}{|c_{u}|}) T^{3}(t,u),
 \end{align*}
where we have used \eqref{veljunction}, \eqref{eq:T} and \eqref{eq:N}.

So far we have totally neglected the questions about the function class and the compatibility conditions.
Typically one works in the class of H\"older spaces, however, since the analysis of short and long-time existence and behaviours of the flow  are outside the scope of this paper we leave out here the details and refer for instance to \cite{Tortorelli}, \cite{KNP21} where such questions are thoroughly studied for similar problems. 

\subsubsection{Weak formulation used for numerics}
\newcommand{\kap}{\kappa} 
\newcommand{\ttau}{\Delta t} 

The numerical scheme we propose hinges on a suitable weak formulation of Problem~\ref{HCSF1ptriod}. To that end set $I=(0,1)$ and
\begin{align*}
V &:=\{(\eta_{1},\eta_{2},\eta_{3}) \,:\, \eta_{\alpha} \in H^{1,2}(I,\R^{2}), \,\eta_{\alpha}(1)=(P^{1}_{\alpha},P^{2}_{\alpha})^{t}, \alpha=1,2,3, \ \eta_{1}(0)=\eta_{2}(0)=\eta_{3}(0)\}, \\
V_0 &:=\{(\eta_{1},\eta_{2},\eta_{3}) \,:\, \eta_{\alpha} \in H^{1,2}(I,\R^{2}), \, \eta_{\alpha}(1)=0, \alpha=1,2,3, \ \eta_{1}(0)=\eta_{2}(0)=\eta_{3}(0)\}.
\end{align*}
Moreover, we let
\[
W :=\{(\chi_{1},\chi_{2},\chi_{3}) \,:\, \chi_{\alpha} \in H^{1,2}(I),\,  \alpha=1,2,3\}.
\]

Given  an initial triod configuration $(c_{1}(0), c_{2}(0), c_{3}(0))$, 
find $c \in H^1(0,T; V)$, $\kap \in L^2(0,T; W)$ and $\mu \in L^2(0,T; \R^3)$ 
such that for almost every $t\in (0,T]$ it holds that 
$\sum_{\alpha=1}^3 \mu_\alpha(t) =0$ and 
\begin{subequations} \label{eq:weak}
\begin{align}
\sum_{\alpha=1}^3
\int_I \partial_{t}c_{\alpha} \cdot \vec{n}_{\alpha} \chi_\alpha
|\partial_{u}c_{\alpha}| du 
- \sum_{\alpha=1}^3
\int_I (\kap_\alpha - \mu_\alpha) \chi_\alpha |\partial_{u}c_{\alpha}| du = 0 \qquad \forall\ \chi \in W, \label{eq:weaka} \\
\sum_{\alpha=1}^3
\int_I \kap_\alpha \vec{n}_{\alpha} \cdot \eta_\alpha
|\partial_{u}c_{\alpha}| du 
+ \sum_{\alpha=1}^3
\int_I \partial_{u}c_{\alpha} \cdot \partial_{u}\eta_{\alpha}
|\partial_{u}c_{\alpha}|^{-1} du = 0 \qquad \forall\ \eta \in V_0,\label{eq:weakb} \\
\int_I (\kap_1 - \mu_1) |\partial_{u}c_{1}| du =
\int_I (\kap_2 - \mu_2) |\partial_{u}c_{2}| du =
\int_I (\kap_3 - \mu_3) |\partial_{u}c_{3}| du.
\label{eq:weakc}
\end{align}
\end{subequations}
Upon recalling \eqref{oo2} and with the correspondence $\mu_{1}=-\lambda_{1}$,  $\mu_{2}= \lambda_{1} -\lambda_{2}$, $\mu_{3}=\lambda_{2}$, we observe that \eqref{eq:weak} is the natural weak formulation of 
\eqref{eq:strong}.
See also \cite[Remark~3.2]{clust3d}, where a corresponding strong
formulation for volume conserving mean curvature flow of a surface cluster 
with three surfaces and a triple junction line is given.

\section{Finite element approximation}\label{sec:FEM}
We consider a uniform partitioning of the reference domain as $\bar{I}=\bigcup_{j=1}^{J}I_j=\bigcup_{j=1}^J[u_{j-1},u_{j}]$ with $u_j= jh,\, h = 1/J$. We then introduce the finite element spaces
\begin{align*}
S^h&:=\left\{\chi\in C^0(\bar{I}): \chi\big|_{I_j}\quad\mbox{is affine}\quad j = 1,\ldots, J \right\},\\
V^h & := ([S^h]^2 \times [S^h]^2 \times [S^h]^2) \cap V, \quad 
V^h_0 := ([S^h]^2 \times [S^h]^2 \times [S^h]^2) \cap V_0, \\ 
W^h & := (S^h \times S^h \times S^h) \cap W.
\end{align*}
We also introduce the mass-lumped quadrature rule as
\begin{equation*}
\int_I^h \chi\, du = \frac{h}{2}\sum_{j=1}^J 
\left[\chi(u_j^-) + \chi(u_{j-1}^+)\right], 
\end{equation*}
for a function $\chi\in L^\infty(I)$ which is piecewise continuous 
with possible jumps at the nodes $\{u_j\}_{j=1}^J$, and 
$\chi(u_j^\pm)=\underset{\delta\searrow 0}{\lim}\ \chi(u_j\pm\delta)$.

We further divide the time interval uniformly by $[0,T]=\bigcup_{m=1}^{M}[t_{m-1}, t_{m}]$, where $t_m= m\ttau$ and  $\ttau = T/M$ is the time step size.  Let $c^m\in V^h$ be an approximation to $c(t_m)$ for $0\leq m\leq M$. Throughout this section we always assume that 
\[
|\partial_u c_\alpha^m|>0\quad\mbox{in}\quad I,\quad 
\alpha=1,2,3,\quad m = 0,1,\ldots, M,
\]
and define the unit normal on $c^m_\alpha(I)$ via
\[
\vec{n}^m_{\alpha} = \frac{(\partial_u c^m_\alpha)^\perp}{|\partial_u c^m_\alpha|}.
\]

Our fully discrete finite element approximation of \eqref{eq:weak} 
is given as follows.
Let $c^0 \in V^h$. Then, for $m\geq0$, given $c^m \in V^h$, find 
$c^{m+1} = c^m + \delta c^{m+1}\in V^h$ with 
$\delta c^{m+1} \in V^h_0$, $\kap^{m+1} \in W^h$ and
$\mu^{m+1} \in \R^3$ 
such that $\sum_{\alpha=1}^3 \mu^{m+1}_\alpha =0$ and 
\begin{subequations} \label{eq:fd}
\begin{align}
\sum_{\alpha=1}^3
\int_I^h \frac{\delta c^{m+1}_{\alpha}}{\ttau} \cdot \vec{n}^m_{\alpha} \chi_\alpha
|\partial_{u}c^m_{\alpha}| du 
- \sum_{\alpha=1}^3
\int_I^h (\kap^{m+1}_\alpha - \mu^{m+1}_\alpha) \chi_\alpha  |\partial_{u}c^m_{\alpha}| du = 0 \qquad \forall\ \chi \in W^h, \label{eq:fda} \\
\sum_{\alpha=1}^3
\int_I^h \kap^{m+1}_\alpha \vec{n}^m_{\alpha} \cdot \eta_\alpha
|\partial_{u}c^m_{\alpha}| du 
+ \sum_{\alpha=1}^3
\int_I \partial_{u}c^{m+1}_{\alpha} \cdot \partial_{u}\eta_{\alpha}
|\partial_{u}c^m_{\alpha}|^{-1} du = 0 \qquad \forall\ \eta \in V^h_0,\label{eq:fdb} \\
\int^h_I (\kap^{m+1}_1 - \mu^{m+1}_1) |\partial_{u}c^m_{1}| du =
\int^h_I (\kap^{m+1}_2 - \mu^{m+1}_2) |\partial_{u}c^m_{2}| du =
\int^h_I (\kap^{m+1}_3 - \mu^{m+1}_3) |\partial_{u}c^m_{3}| du.
\label{eq:fdc}
\end{align}
\end{subequations}
Observe that \eqref{eq:fd} is linear, and that it is the natural extension of the scheme (4.28) from
\cite[Remark~4.3]{clust3d}, to the special curve network studied here.

We make the following mild assumption for the existence and uniqueness proof.
\begin{ass} \label{ass:A}
${\rm span}\, \{ \int_I^h \xi \vec{n}^m_{\alpha} |\partial_{u}c^m_{\alpha}| du :
\xi \in S^h, \xi(0) = \xi(1) = 0 \} \not= \{0\}$, for $\alpha=1,2,3$.
\end{ass}

The above assumption is very mild. A rare case when it is violated is if
$J=2$ and $c_\alpha^m(0) = c_\alpha^m(1)$.
In fact, as the proof of Theorem~\ref{thm:stab} below will show,
Assumption~\ref{ass:A} can be weakened further to 
the span being nontrivial for just two out of the three curves. In any case,
the assumption is virtually always satisfied in practice, and is in line with
similar assumptions stated in the series of works
\cite{triplej,triplejMC,clust3d,ejam3d}. 

\begin{teo} \label{thm:stab}
Let Assumptions~\ref{ass:A} hold.
Then there exists a unique solution to \eqref{eq:fd}. 
Moreover, any solution to \eqref{eq:fd} is unconditionally stable. In particular
\begin{equation} \label{eq:stab}
L_E(c^{m+1}) + \ttau \sum_{\alpha=1}^3
\int_I^h (\kap^{m+1}_\alpha - \mu^{m+1}_\alpha)^2 |\partial_{u}c^m_{\alpha}| du
\leq L_E(c^{m}).
\end{equation}
\end{teo} 
\begin{proof}
As \eqref{eq:fd} is a linear system with the same number of equations as
unknowns, existence follows from uniqueness. It is hence sufficient to show
that the homogeneous system only has the trivial solution. Let
$\widehat{\delta c} \in V^h_0$, $\widehat\kap \in W^h$ and $\widehat\mu \in
\R^3$ be such that
\begin{subequations} \label{eq:homo}
\begin{align}
\sum_{\alpha=1}^3
\int_I^h \frac{\widehat{\delta c}_{\alpha}}{\ttau} \cdot \vec{n}^m_{\alpha} \chi_\alpha
|\partial_{u}c^m_{\alpha}| du 
- \sum_{\alpha=1}^3
\int_I^h (\widehat\kap_\alpha - \widehat\mu_\alpha) \chi_\alpha  |\partial_{u}c^m_{\alpha}| du = 0 \qquad \forall\ \chi \in W^h, \label{eq:homoa} \\
\sum_{\alpha=1}^3
\int_I^h \widehat\kap_\alpha \vec{n}^m_{\alpha} \cdot \eta_\alpha
|\partial_{u}c^m_{\alpha}| du 
+ \sum_{\alpha=1}^3
\int_I \partial_{u} \widehat{\delta c}_{\alpha} \cdot \partial_{u}\eta_{\alpha}
|\partial_{u}c^m_{\alpha}|^{-1} du = 0 \qquad \forall\ \eta \in V^h_0,\label{eq:homob} \\
\int^h_I (\widehat\kap_1 - \widehat\mu_1) |\partial_{u}c^m_{1}| du =
\int^h_I (\widehat\kap_2 - \widehat\mu_2) |\partial_{u}c^m_{2}| du =
\int^h_I (\widehat\kap_3 - \widehat\mu_3) |\partial_{u}c^m_{3}| du.
\label{eq:homoc}
\end{align}
\end{subequations}
Choosing $\chi = \widehat\kap \in W^h$ in \eqref{eq:homoa} and
$\eta = \widehat{\delta c} \in V^h_0$ in \eqref{eq:homob} yields that
\begin{align} \label{eq:homo1}
0 & =
\sum_{\alpha=1}^3
\int_I |\partial_{u} \widehat{\delta c}_{\alpha}|^2 
|\partial_{u}c^m_{\alpha}|^{-1} du
+ \ttau \sum_{\alpha=1}^3
\int_I^h (\widehat\kap_\alpha - \widehat\mu_\alpha) \widehat\kap_\alpha |\partial_{u}c^m_{\alpha}| du \nonumber \\ & 
= \sum_{\alpha=1}^3
\int_I |\partial_{u} \widehat{\delta c}_{\alpha}|^2 
|\partial_{u}c^m_{\alpha}|^{-1} du
+ \ttau \sum_{\alpha=1}^3
\int_I^h (\widehat\kap_\alpha - \widehat\mu_\alpha)^2 |\partial_{u}c^m_{\alpha}| du \,,
\end{align}
where we have noticed from \eqref{eq:homoc} that
\begin{equation} \label{eq:homo2}
\sum_{\alpha=1}^3 \widehat\mu_\alpha 
\int_I^h (\widehat\kap_\alpha - \widehat\mu_\alpha) |\partial_{u}c^m_{\alpha}| du 
= \int_I^h (\widehat\kap_1 - \widehat\mu_1) |\partial_{u}c^m_1| du
\sum_{\alpha=1}^3 \widehat\mu_\alpha
 = 0.
\end{equation}
It follows immediately from \eqref{eq:homo1} 
that $\widehat{\delta c}$ is a constant, and hence
$\widehat{\delta c}=0$ because of the boundary conditions in $V^h_0$. Moreover,
we have that $\widehat\kap_\alpha = \widehat\mu_\alpha$, $\alpha=1,2,3$.
Combining these results with \eqref{eq:homob} implies that
\begin{equation} \label{eq:homo3}
\sum_{\alpha=1}^3 \widehat\mu_\alpha 
\int_I^h \vec{n}^m_{\alpha} \cdot \eta_\alpha
|\partial_{u}c^m_{\alpha}| du = 0 \quad \forall\ \eta \in V^h_0.
\end{equation}
Assumption~\ref{ass:A} now yields that $\widehat\mu_\alpha = 0$,
$\alpha=1,2,3$. For example, for $\alpha=1$ we choose a nonzero vector
$v = \int_I^h \xi \vec{n}^m_{1} |\partial_{u}c^m_{1}| du$ for
$\xi \in S^h$ with $\xi(0)= \xi(1) = 0$. Hence we can choose
$\eta = (\widehat\mu_1 \xi v, 0, 0) \in V^h_0$ in \eqref{eq:homo3} to obtain 
that
\[
0 = \widehat\mu_1
\int_I^h \vec{n}^m_1 \cdot (\widehat\mu_1 \xi v)
|\partial_{u}c^m_1| du 
= (\widehat\mu_1)^2 v \cdot \int_I^h \vec{n}^m_1 \xi 
|\partial_{u}c^m_1| du 
= (\widehat\mu_1)^2 |v|^2,
\]
which implies that $\widehat\mu_1 = 0$. We can argue similarly to obtain
$\widehat\mu_2 = 0$, and then $\widehat\mu_3 = - \widehat\mu_1 - \widehat\mu_2
= 0$. Overall we have that $(\widehat{\delta c}, \widehat\kap, \widehat\mu) =
(0,0,0)$. This proves the existence of a unique solution to \eqref{eq:fd}. 
 
In order to prove \eqref{eq:stab}, we choose $\chi = \kap^{m+1} \in W^h$ 
in \eqref{eq:fda} and $\eta = \delta c^{m+1} \in V^h_0$ in \eqref{eq:fdb} 
to obtain
\begin{align} \label{eq:stab1}
0 & =
\sum_{\alpha=1}^3
\int_I \partial_{u} c^{m+1}_{\alpha} \cdot \partial_{u}\delta c^{m+1}_{\alpha}
|\partial_{u}c^m_{\alpha}|^{-1} du
+ \ttau \sum_{\alpha=1}^3
\int_I^h (\kap^{m+1}_\alpha - \mu^{m+1}_\alpha) \kap^{m+1}_\alpha |\partial_{u}c^m_{\alpha}| du \nonumber \\ & 
= \sum_{\alpha=1}^3
\int_I \partial_{u} c^{m+1}_{\alpha} \cdot \partial_{u}(c^{m+1}_{\alpha} -
c^m_\alpha) 
|\partial_{u}c^m_{\alpha}|^{-1} du
+ \ttau \sum_{\alpha=1}^3
\int_I^h (\kap^{m+1}_\alpha - \mu^{m+1}_\alpha)^2 |\partial_{u}c^m_{\alpha}| du \,,
\end{align}
where we have noticed from \eqref{eq:fdc}, similarly to \eqref{eq:homo2}, that
\begin{equation*}
\sum_{\alpha=1}^3 \mu^{m+1}_\alpha 
\int_I^h (\kap^{m+1}_\alpha - \mu^{m+1}_\alpha) |\partial_{u}c^m_{\alpha}| du 
= \int_I^h (\kap^{m+1}_1 - \mu^{m+1}_1) |\partial_{u}c^m_1| du
\sum_{\alpha=1}^3 \mu^{m+1}_\alpha
 = 0.
\end{equation*}
In addition, we recall from \cite[Lemma~57]{bgnreview} that
\begin{equation} \label{eq:lemma57}
\int_I \partial_{u} c^{m+1}_{\alpha} \cdot \partial_{u}(c^{m+1}_{\alpha} -
c^m_\alpha) 
|\partial_{u}c^m_{\alpha}|^{-1} du \geq L_E(c^{m+1}_\alpha) 
- L_E(c^{m}_\alpha), \quad \alpha=1,2,3.
\end{equation}
Combining \eqref{eq:stab1} and \eqref{eq:lemma57} yields the desired result
\eqref{eq:stab}. 
\end{proof}

\section{Numerical results} \label{sec:nr}

\newcommand{\proj}[1]{\overline{#1}}

We implemented the scheme \eqref{eq:fd}  
within the finite element toolbox Alberta, \cite{Alberta}, using
the sparse factorization package UMFPACK, see \cite{Davis04},
for the solution of the linear systems of equations arising at each time level.
Unless otherwise stated, for the discretization parameters we
always use $J=100$ and $\Delta t = 10^{-4}$.

In order to describe the horizontal initial data, we will always state the junction point $\Sigma \in \R^3$
and then either state the coordinates of the fixed points
$P_\alpha \in \R^3$, $\alpha=1,2,3$, directly, or at first only specify the
coordinates of their projections into the two-dimensional Euclidean plane,
$\proj{P_\alpha} \in \R^2$, $\alpha=1,2,3$, from which the missing coordinates can be 
computed with the help of the chosen initial planar curves via \eqref{g3}.
 
For the colours in our plots we follow the following convention: if the triod
(or its projection) is shown at a single time $t=t_m$, 
then we use the colours olive,
purple and gold for the curves $c_1^m$, $c_2^m$ and $c_3^m$, respectively.
If a plot shows the evolution of a triod (or its projection) at several times, 
then we use the colours blue, black and read for the initial data,
the intermediate times and the final time, respectively.

\subsection{Planar straight line initial data}

In this subsection we always choose $\Sigma=(0,0,0)^t$ and the initial curves to be straight lines
in the two-dimensional Euclidean plane. In general this means that the 120 
degree angle condition at the projected triple junction $\proj\Sigma$ will
not be satisfied at time $t=0$.

\begin{figure}
\center
\includegraphics[angle=-0,width=0.28\textwidth]{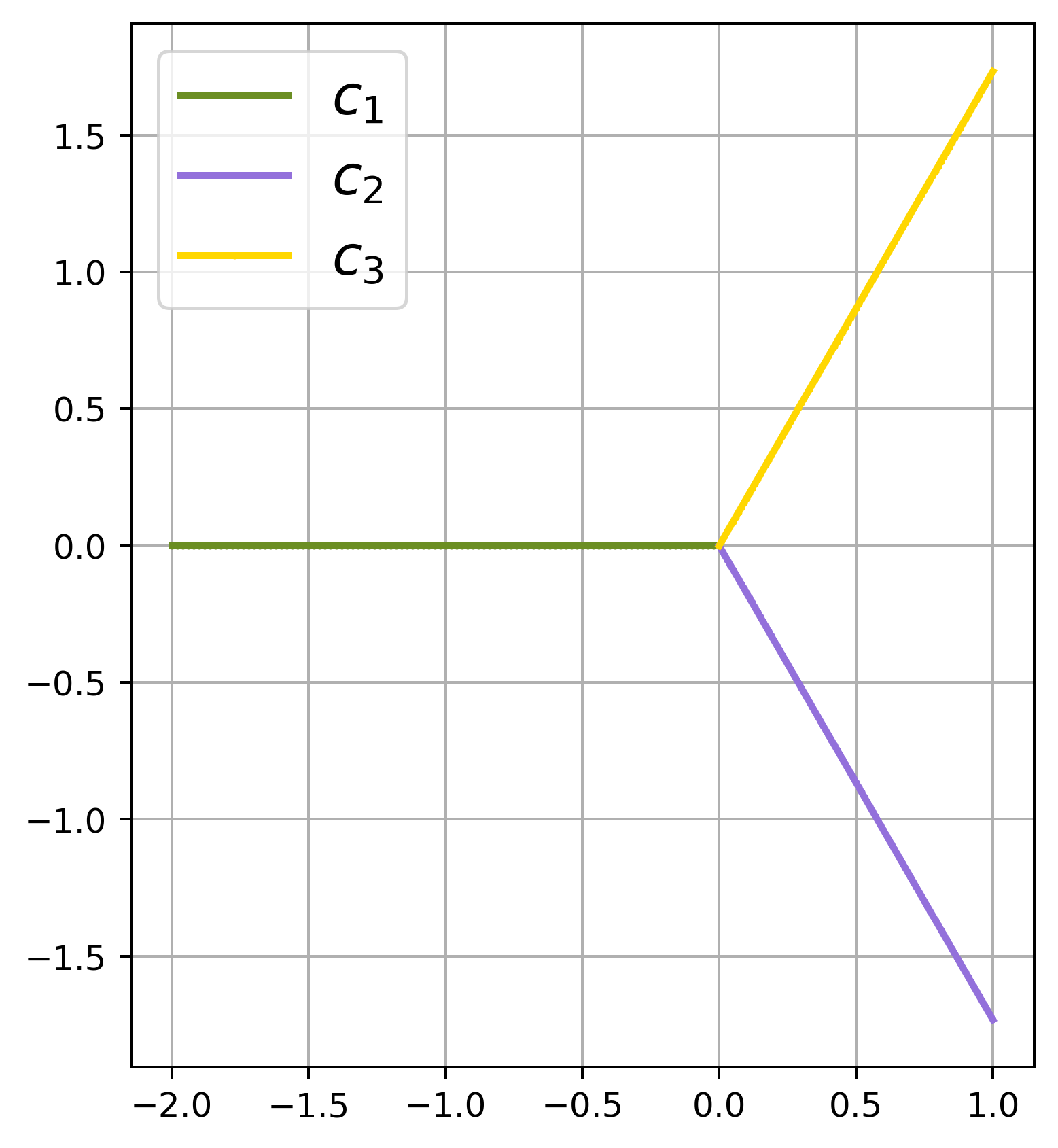}
\includegraphics[angle=-0,width=0.28\textwidth]{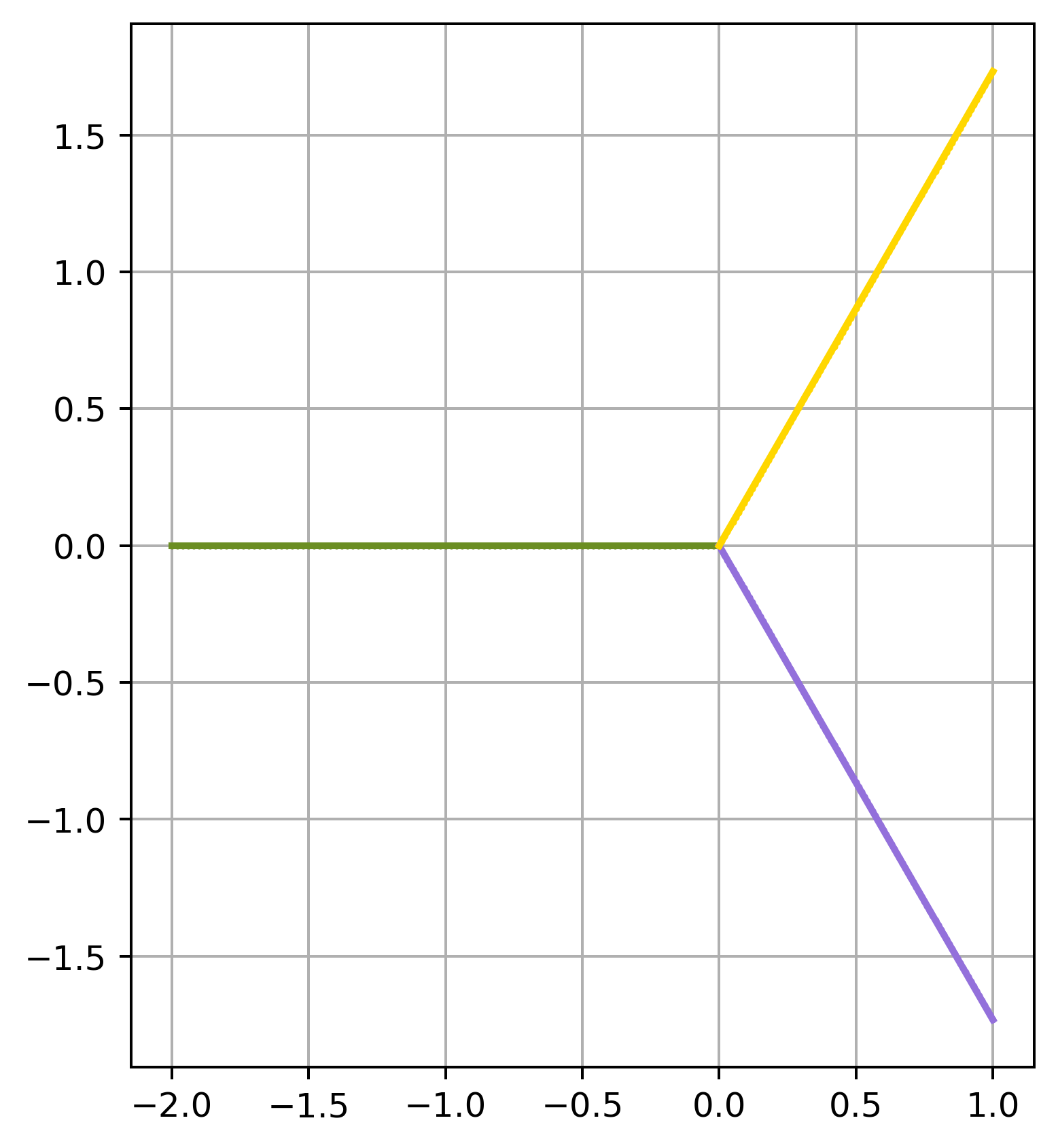}
\includegraphics[angle=-0,width=0.4\textwidth]{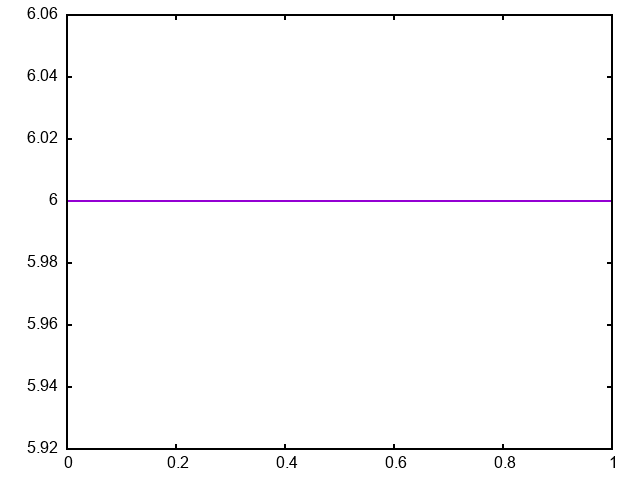}
\caption{The triod at time $t=0$ and at time $t=1$.
We use the colours olive, purple and gold for the curves $c_1^m$, $c_2^m$ and 
$c_3^m$, respectively.
On the right we show a plot of the discrete energy $L_E(c^m)$ over time.
}
\label{fig:steinergood}
\end{figure}

{\bf Experiment 1}: We begin with a classical Steiner configuration, that is
three straight line segments that meet at the origin at $120^\circ$ angles.
In particular $\Sigma = (0,0,0)^t$,
$P_1=(-2, 0, 0)^{t}$, $P_2 = (1, -\sqrt{3}, 0)^{t}$ and $P_3 = (1, \sqrt{3}, 0)^{t}$.
This setup is a steady state for the continuous problem, and also our 
finite element approximation does not change this initial data when we
integrate from $t=0$ to $t=1$. In Figure~\ref{fig:steinergood} we show the
initial and final discrete solution, together with a plot of the (constant)
discrete energy over time.

\begin{figure}
\center
\includegraphics[angle=-0,width=0.1\textwidth]{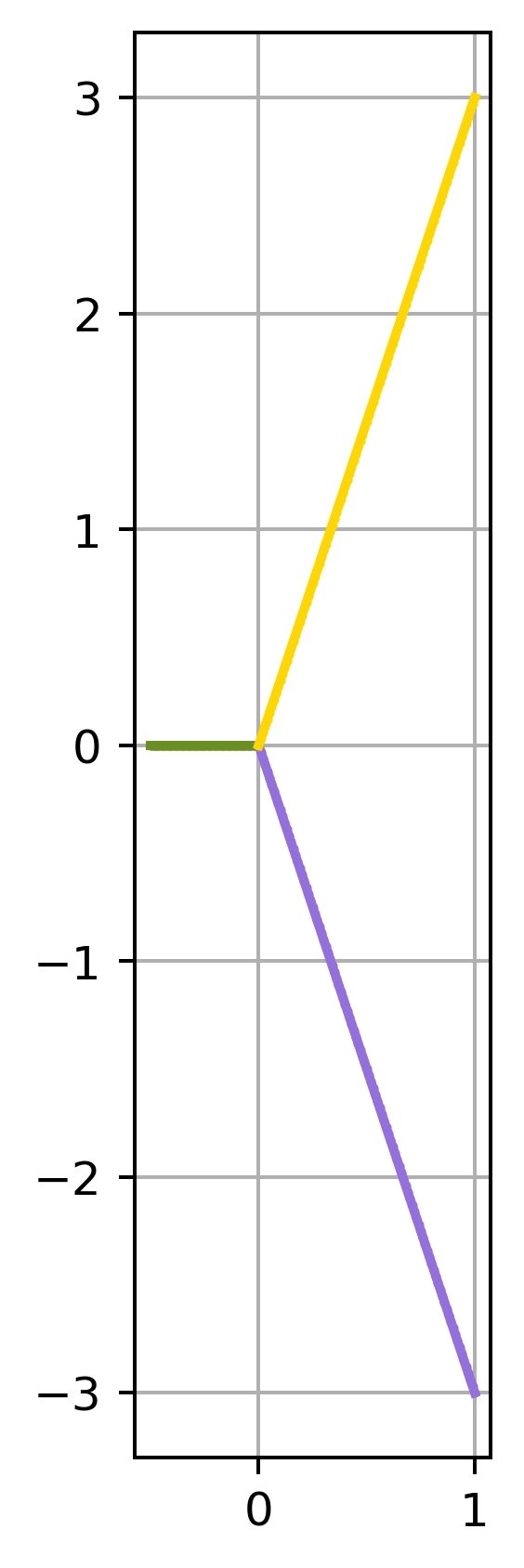}
\includegraphics[angle=-0,width=0.1\textwidth]{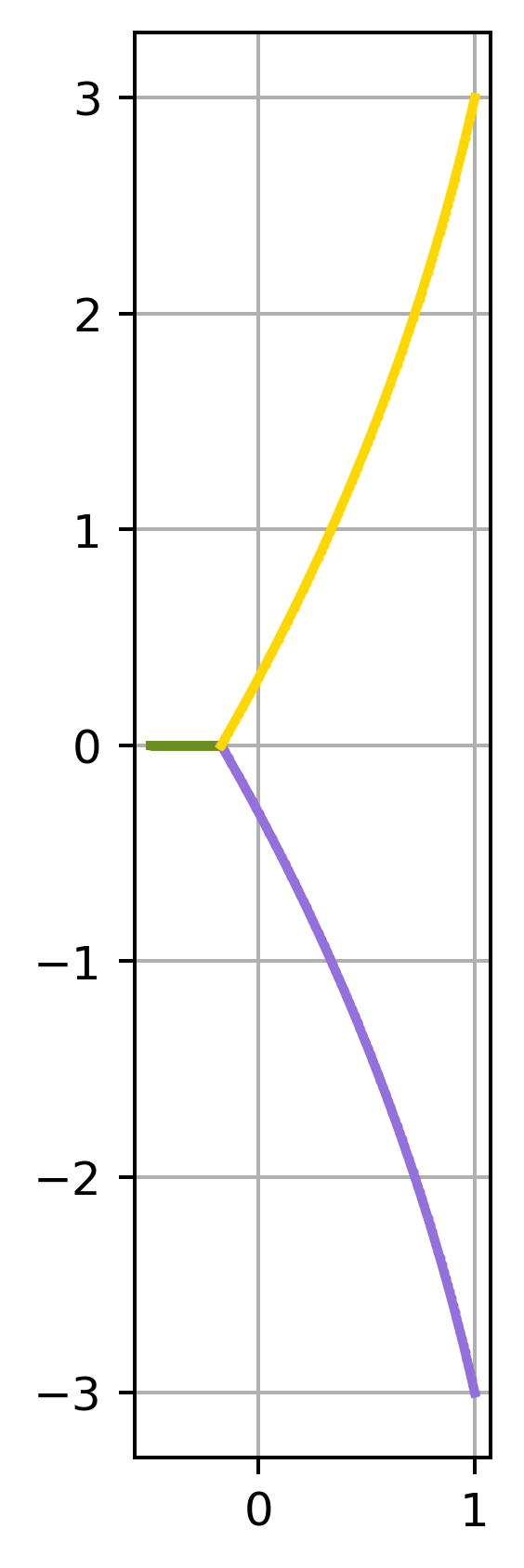}
\includegraphics[angle=-0,width=0.25\textwidth]{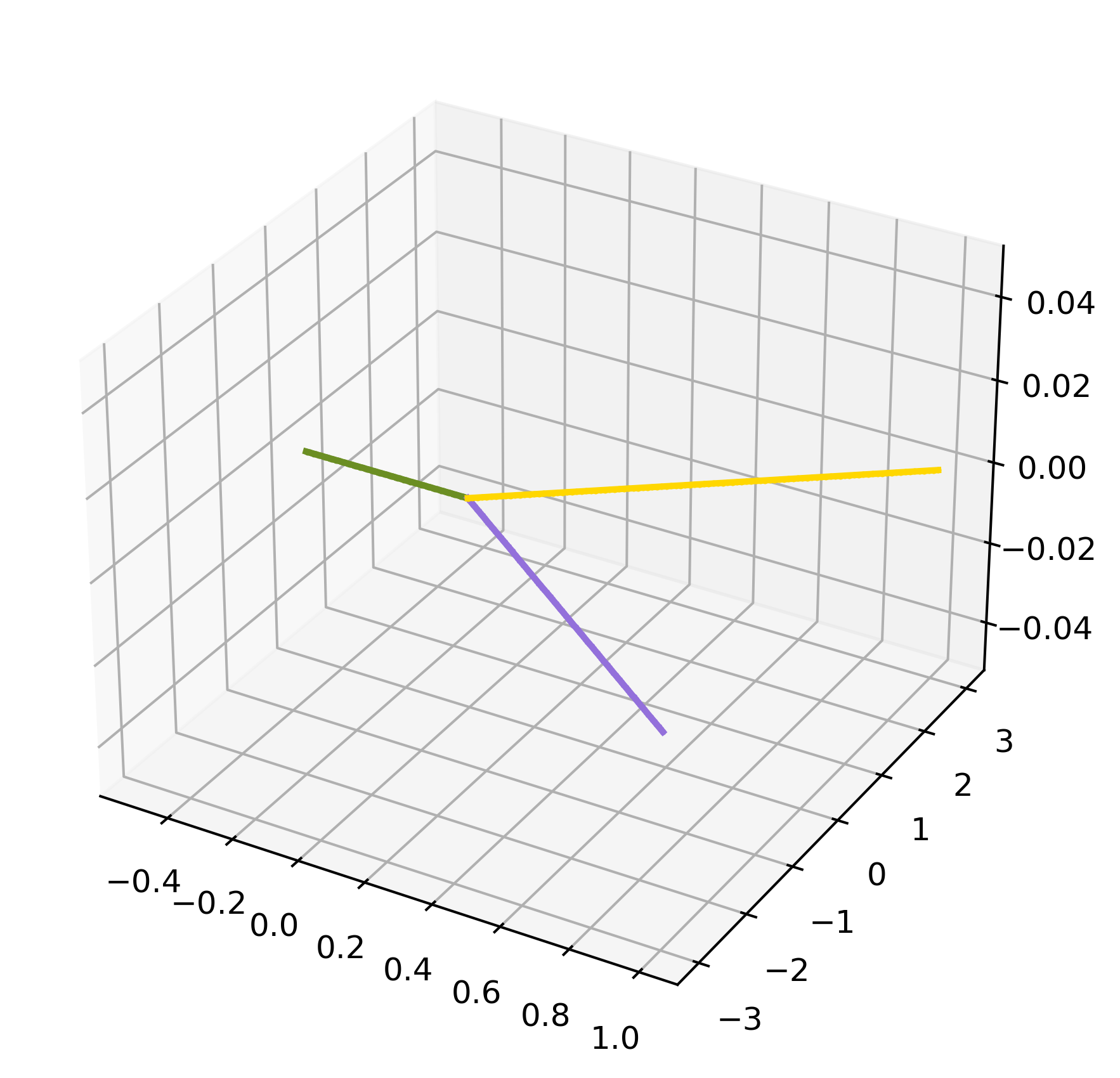}
\includegraphics[angle=-0,width=0.25\textwidth]{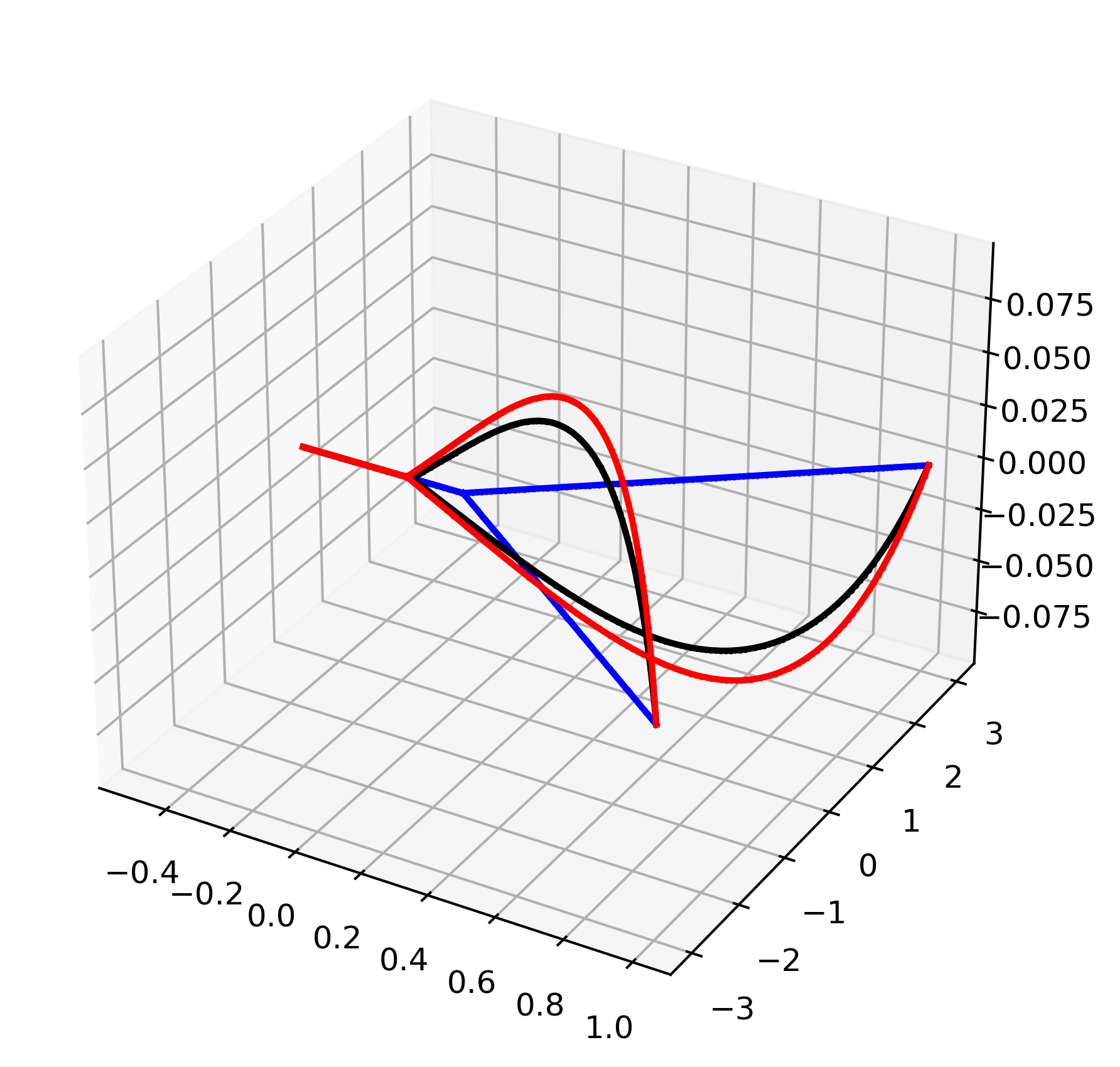}
\includegraphics[angle=-0,width=0.25\textwidth]{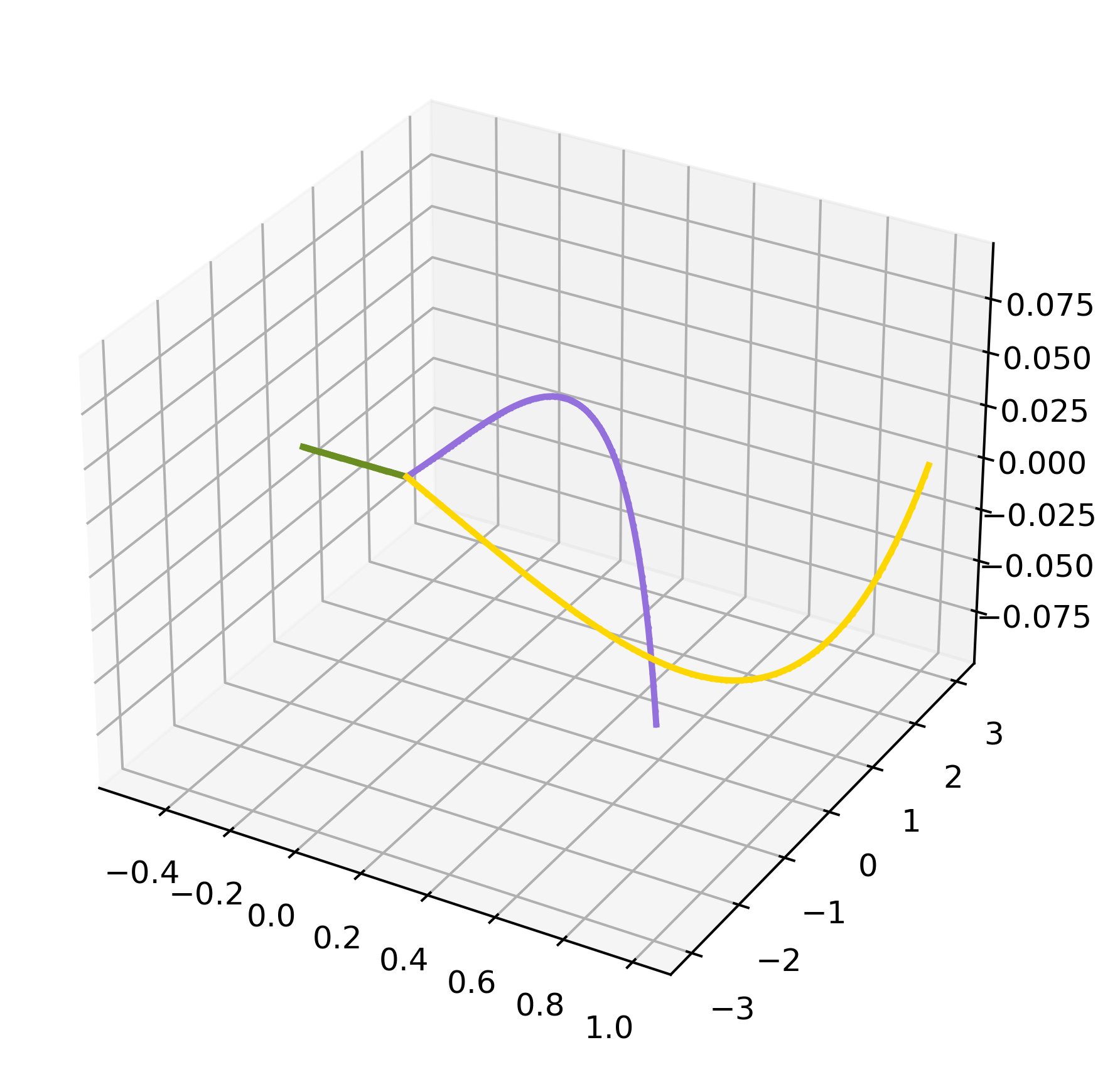}
\includegraphics[angle=-0,width=0.4\textwidth]{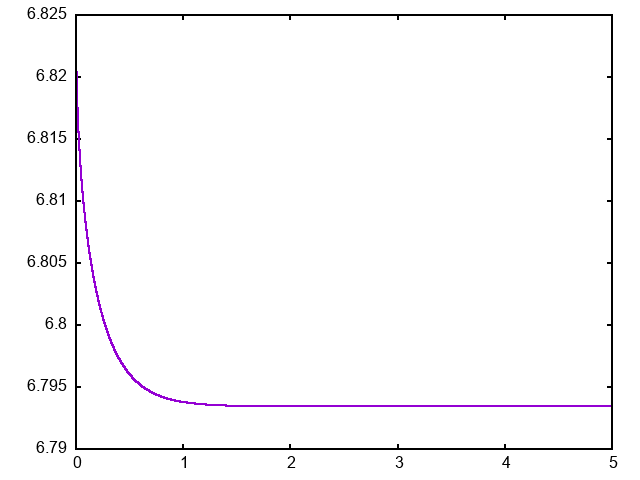}
\includegraphics[angle=-0,width=0.4\textwidth]{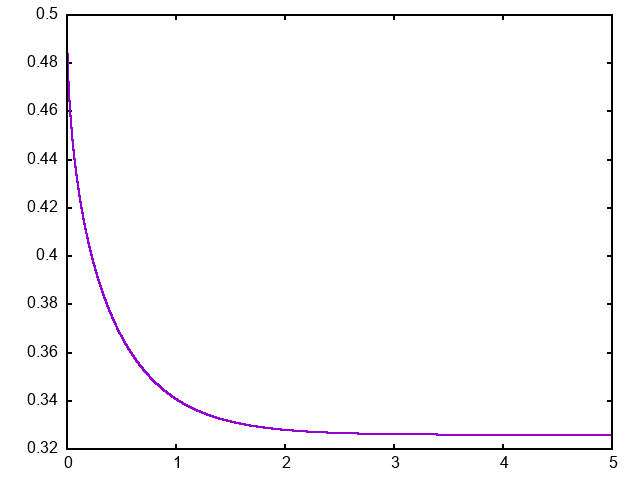}
\caption{The projected triod at time $t=0$ and at time $t=5$. 
The triod in $\R^3$ is shown at time $t=0$, at times $t=0,1,5$, and at 
time $t=5$.
Below we show a plot of the discrete energies $L_E(c^m)$ and $L_E(c^m_1)$
over time.
}
\label{fig:steinerbadstraight}
\end{figure}

{\bf Experiment 2}:
We would like to investigate what happens if we change the positions of the
points $P_\alpha$, $\alpha=1,2,3$, in such a way, that they form
a triangle with an angle larger than 120$^\circ$ in $P_1$. This is
motivated by the fact that for classical curve shortening flow
a singularity is expected with the triple junction $\Sigma$ moving
towards $P_1$ so that the curve $c_1$ eventually vanishes (\cite{Tortorelli}).
We begin with a moderate setting of $P_1=(-0.5,0,0)^{t}$,
$P_2 = (1,-3, 0)^{t}$, $P_3 = (1,3,0)^{t}$, meaning that the angle in
$P_1$ is $126.87^\circ$. Since we use planar straight lines for the
initial curve setup, the $120^\circ$ angle condition at the triple junction is
violated at time $t=0$.
The evolution shown in Figure~\ref{fig:steinerbadstraight} demonstrates that
the triod quickly relaxes to a configuration that satisfies the $120^\circ$
triple junction angle condition. Eventually a steady state is approached, 
as can be seen from the discrete energy plots. At the final time the 
shortest curve has reached a length of about $0.326$.

\begin{figure}
\center
\includegraphics[angle=-0,width=0.07\textwidth]{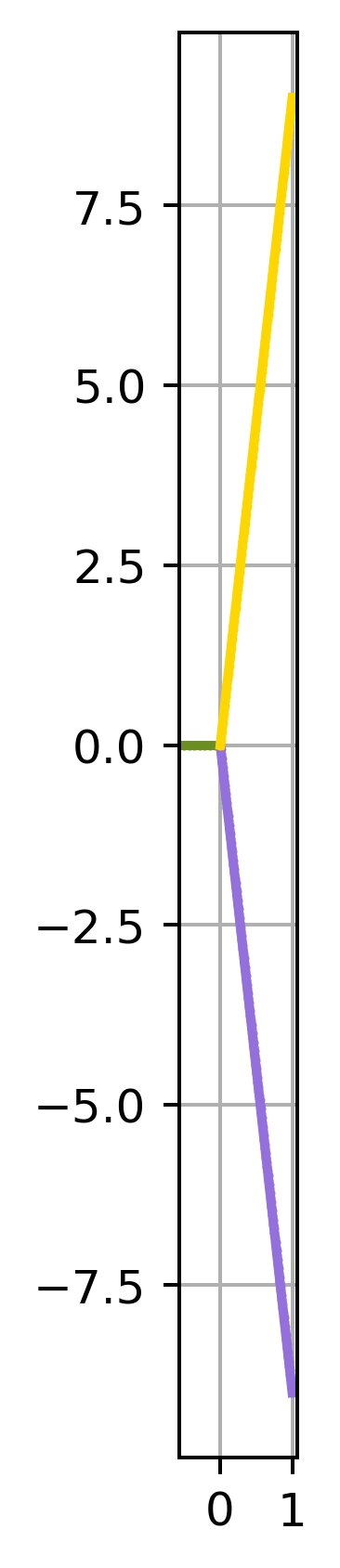}
\includegraphics[angle=-0,width=0.07\textwidth]{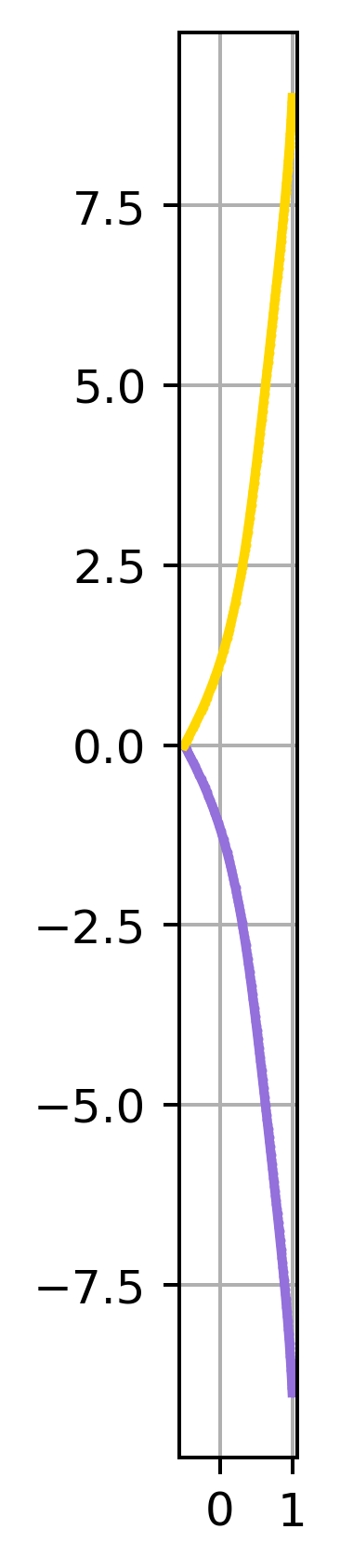}
\includegraphics[angle=-0,width=0.25\textwidth]{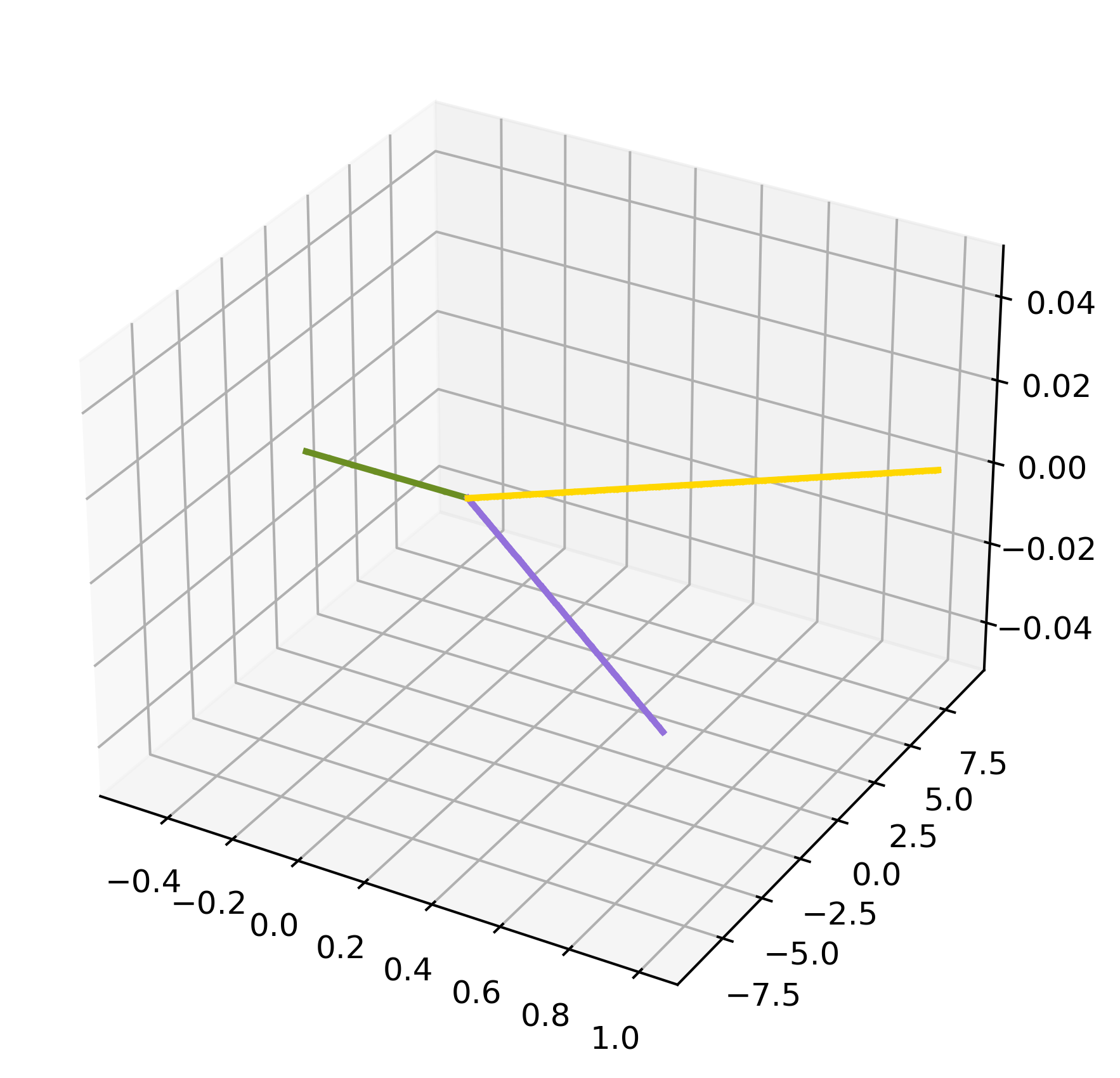}
\includegraphics[angle=-0,width=0.25\textwidth]{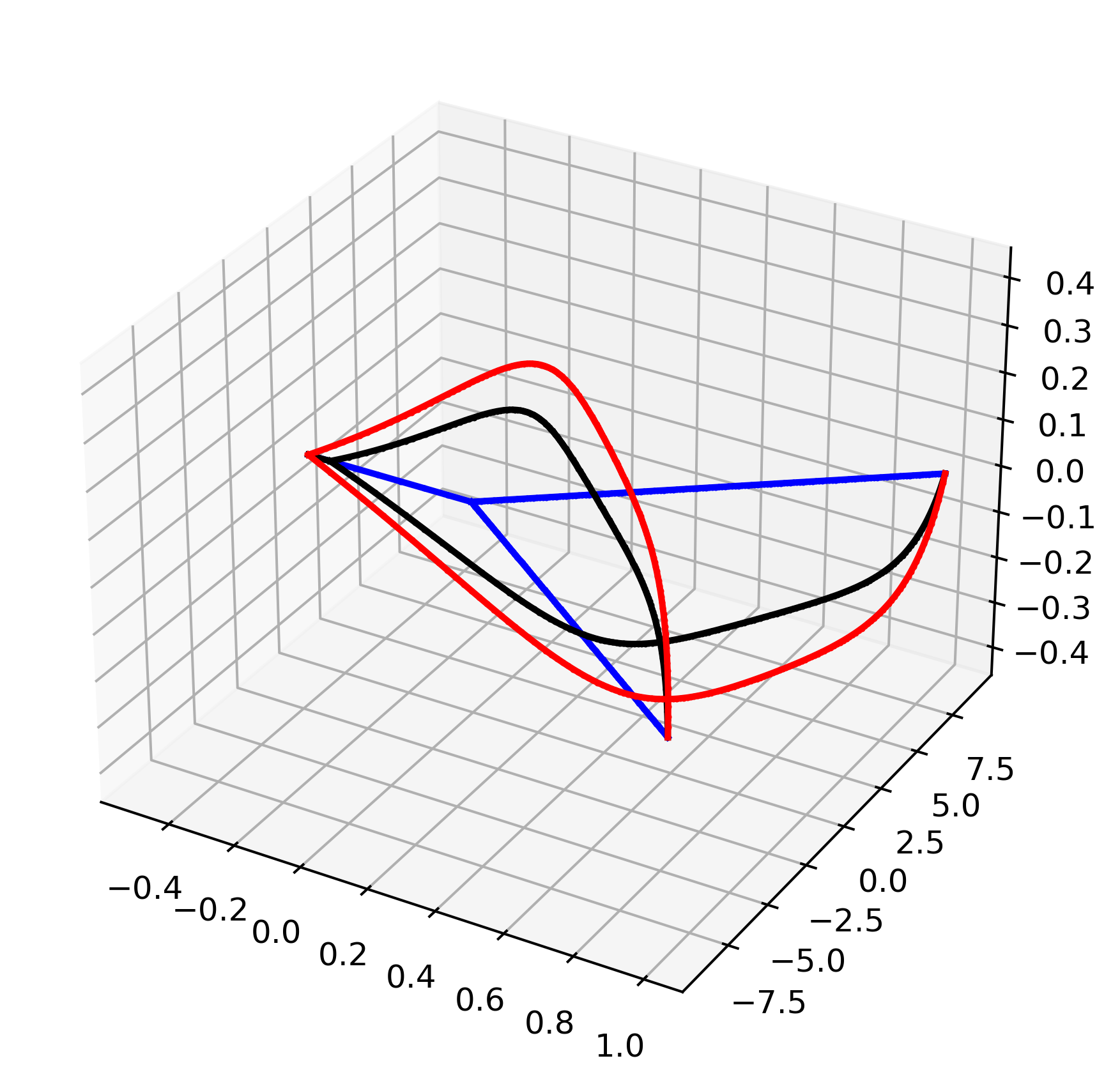}
\includegraphics[angle=-0,width=0.25\textwidth]{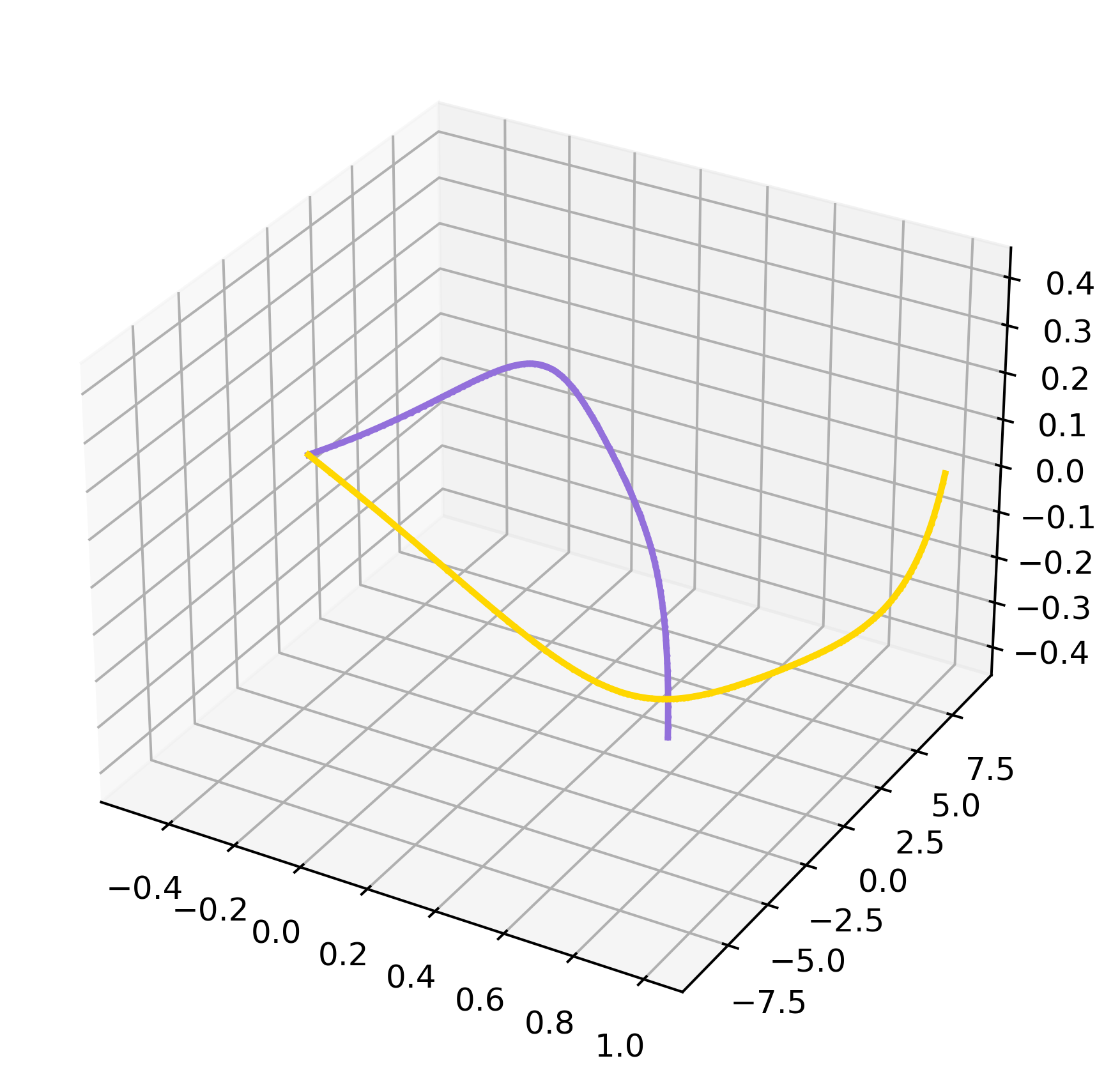}
\includegraphics[angle=-0,width=0.4\textwidth]{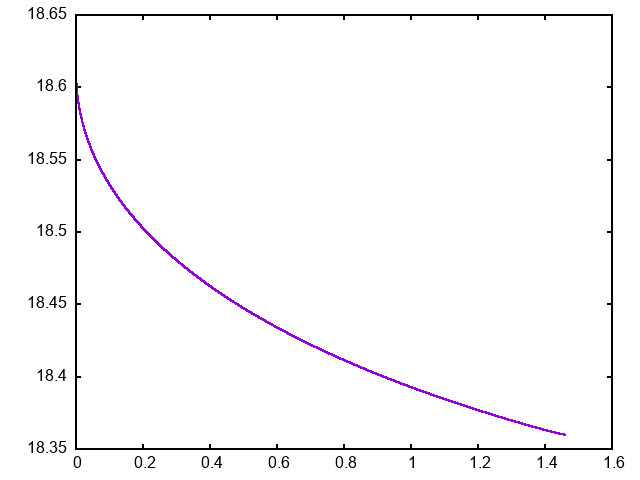}
\includegraphics[angle=-0,width=0.4\textwidth]{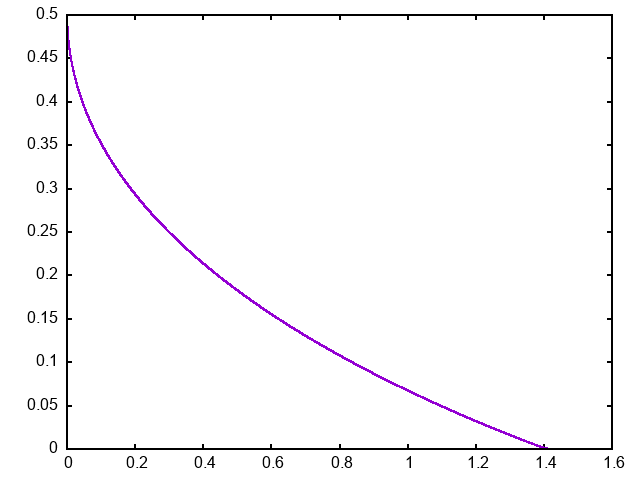}
\caption{The projected triod at time $t=0$ and at time $t=1.46$.
The triod in $\R^3$ is shown at time $t=0$, at times $t=0,1,1.46$, and at 
time $t=1.46$.
Below we show a plot of the discrete energies $L_E(c^m)$ and $L_E(c^m_1)$
over time.
}
\label{fig:steinerbadastraight}
\end{figure}

{\bf Experiment 3}:
In order to try to trigger an evolution towards a singularity, we make the
triangle formed by $P_\alpha$, $\alpha=1,2,3$, even more obtuse. That is,
we let $P_1=(-0.5,0,0)^{t}$,
$P_2 = (1,-9,0)^{t}$, $P_3 = (1,9,0)^{t}$, meaning that the angle in
$P_1$ is $161.06^\circ$.
The evolution shown in Figure~\ref{fig:steinerbadastraight} indicates that now
the flow does indeed encounter a singularity, as the triple junction $\Sigma$
moves towards $P_1$, so that the curve $c_1^m$ eventually vanishes. 
Note that the
curve vanishing represents a singularity for our parametric formulation, and so
our simulation breaks down. Observe from the combined energy plot that the
remaining curves $c_2^m$ and $c_3^m$ are not geodesics yet.

\subsection{Straight line initial data}

In this subsection we choose the initial curves to be straight lines in $\R^3$.

\begin{figure}
\center
\mbox{
\includegraphics[angle=-0,width=0.2\textwidth]{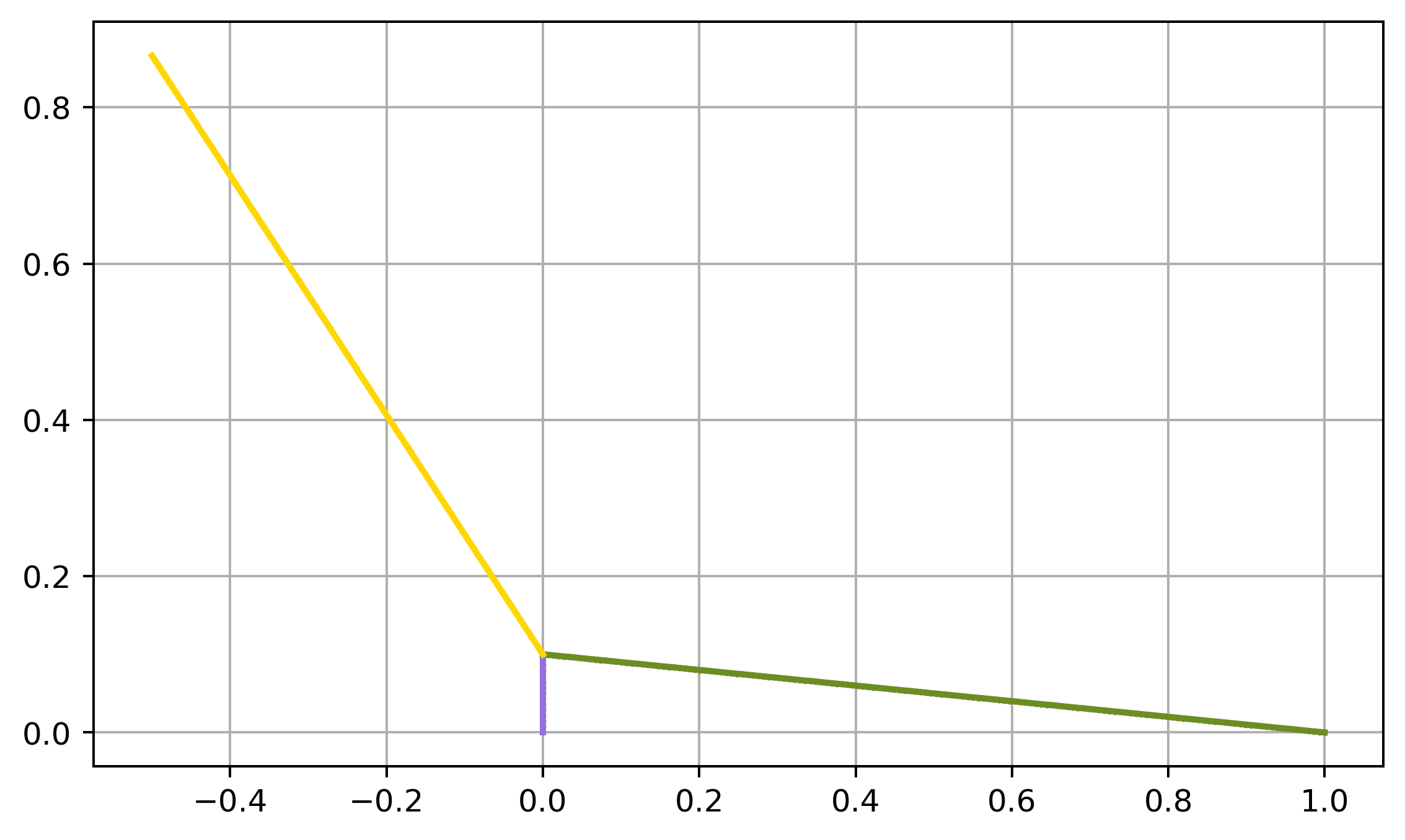}
\includegraphics[angle=-0,width=0.2\textwidth]{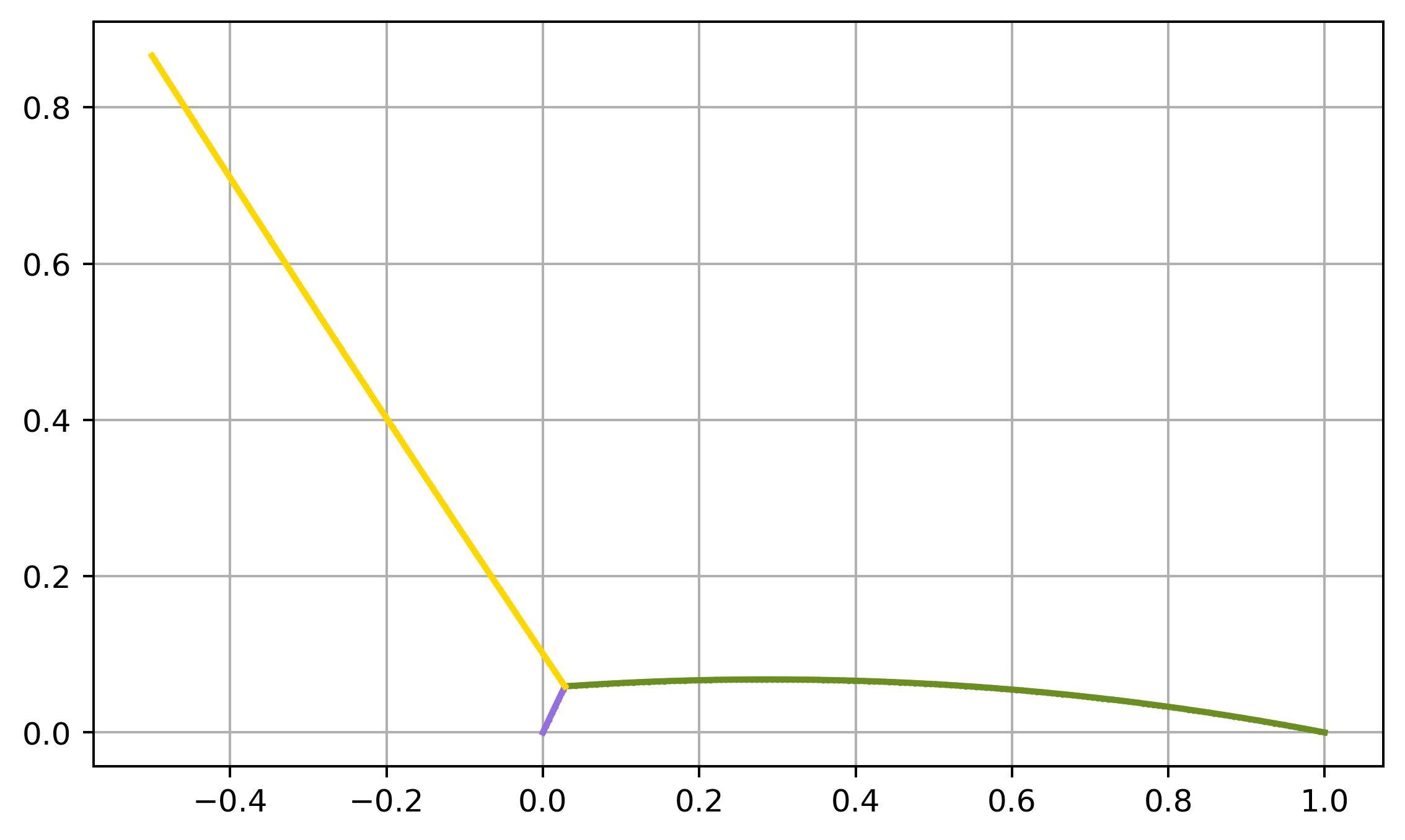}
\includegraphics[angle=-0,width=0.3\textwidth]{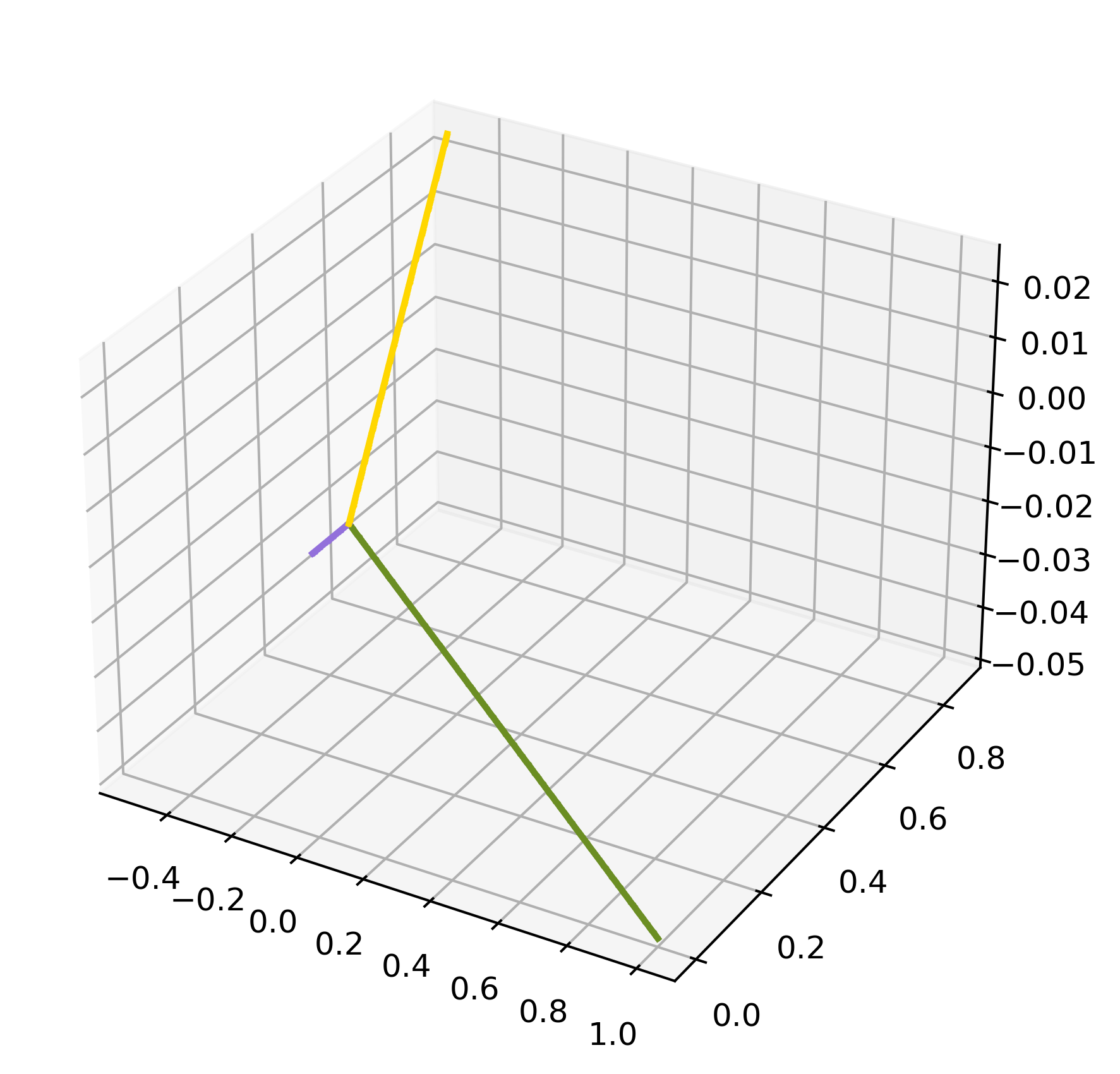}
\includegraphics[angle=-0,width=0.3\textwidth]{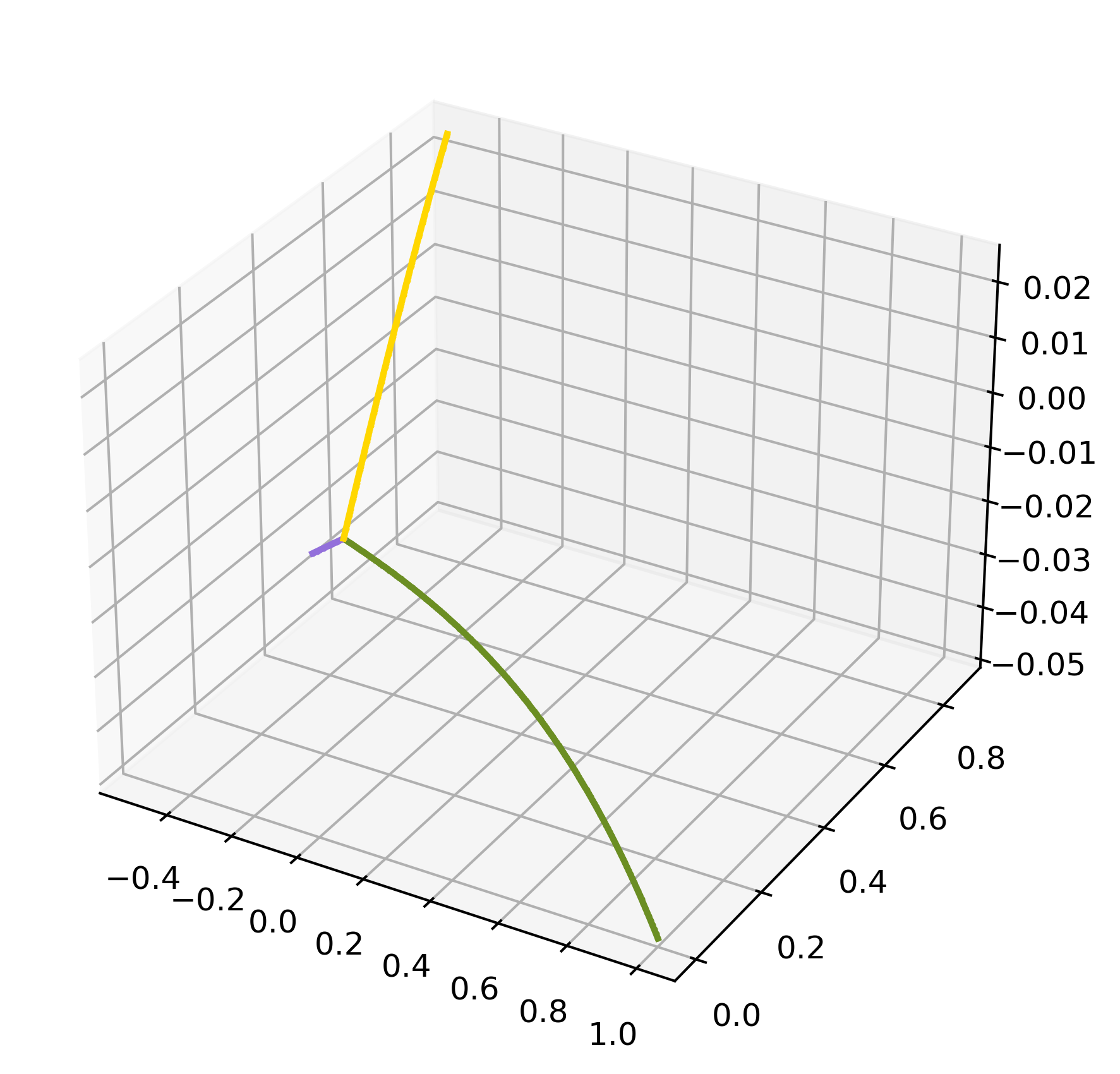}
}
\includegraphics[angle=-0,width=0.4\textwidth]{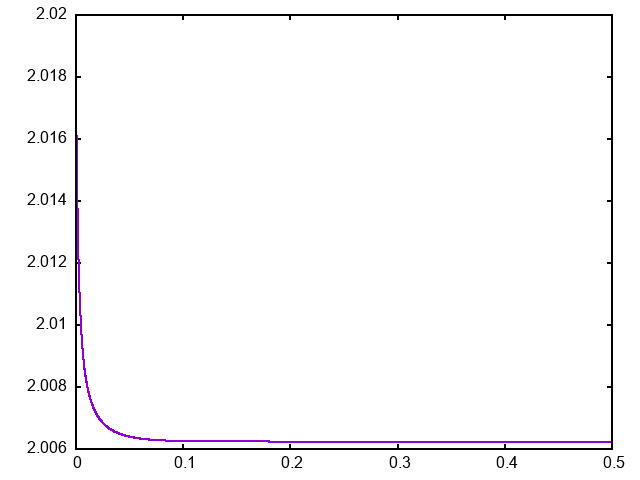}
\includegraphics[angle=-0,width=0.4\textwidth]{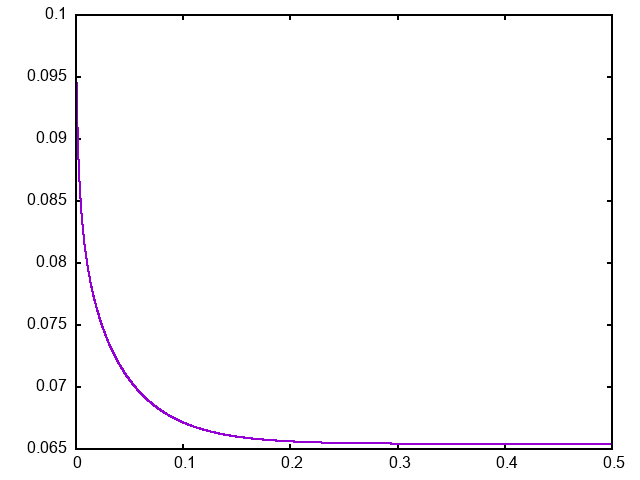}
\caption{The (projected) triod at time $t=0$ 
and at time $t=0.5$.
Below we show a plot of the discrete energies $L_E(c^m)$ and $L_E(c^m_2)$
over time.
}
\label{fig:newstraightz}
\end{figure}

{\bf Experiment 4}:
For this experiment we choose the location of the initial triple junction
to be different from the origin. In fact, we let $\Sigma = (0, 0.1, 0)^{t}$
as well as $P_1 = (1, 0, -0.05)^{t}$, $P_2=(0,0,0)^{t}$ and
$P_3=(\cos(\frac23\pi), \sin(\frac23\pi), 0.05)^{t} = 
(-0.5, \frac12\sqrt{3}, 0.05)^{t}$.
As a consequence, the three points $\proj{P_\alpha}$, $\alpha=1,2,3$, 
form a triangle with an angle equal to 120$^\circ$ in $\proj{P_2}$. 
Once again this choice is motivated by the fact that for standard curve
shortening flow in the plane the triple junction $\proj{\Sigma}$ would move 
towards $\proj{P_2}$, and the curve $c_2$ would vanish.
We show the results of our scheme in Figure~\ref{fig:newstraightz}.
Similarly to Figure~\ref{fig:steinerbadstraight}, for the flow considered in
this paper no singularity occurs. In fact, 
eventually a steady state is approached, with the shortest curve having 
reached a length of about $0.065$.

\subsection{Compatible curved initial data}

In this subsection we choose the projected initial data in such a way, that the 
$120^\circ$ angle condition at the triple junction is satisfied. 
This is achieved with the help of cubic B{\'e}zier curves, and so the initial 
data is in general curved.
Unless otherwise stated it is $\Sigma=(0,0,0)^{t}$.

\begin{figure}
\center
\includegraphics[angle=-0,width=0.1\textwidth]{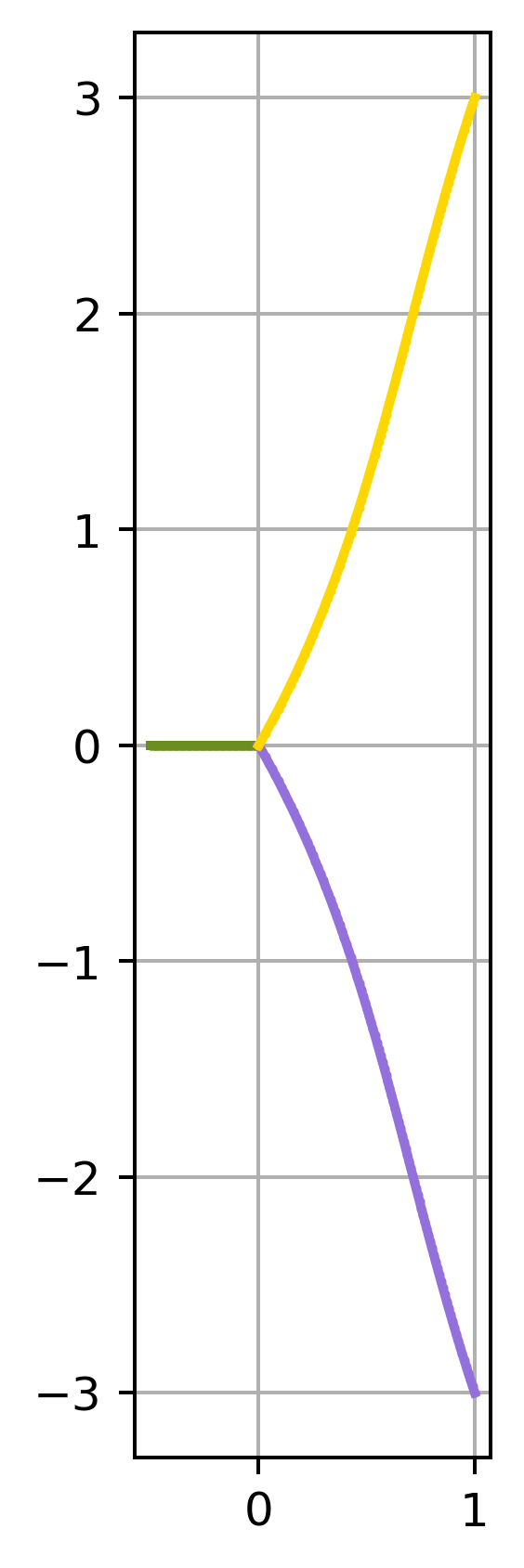}
\includegraphics[angle=-0,width=0.1\textwidth]{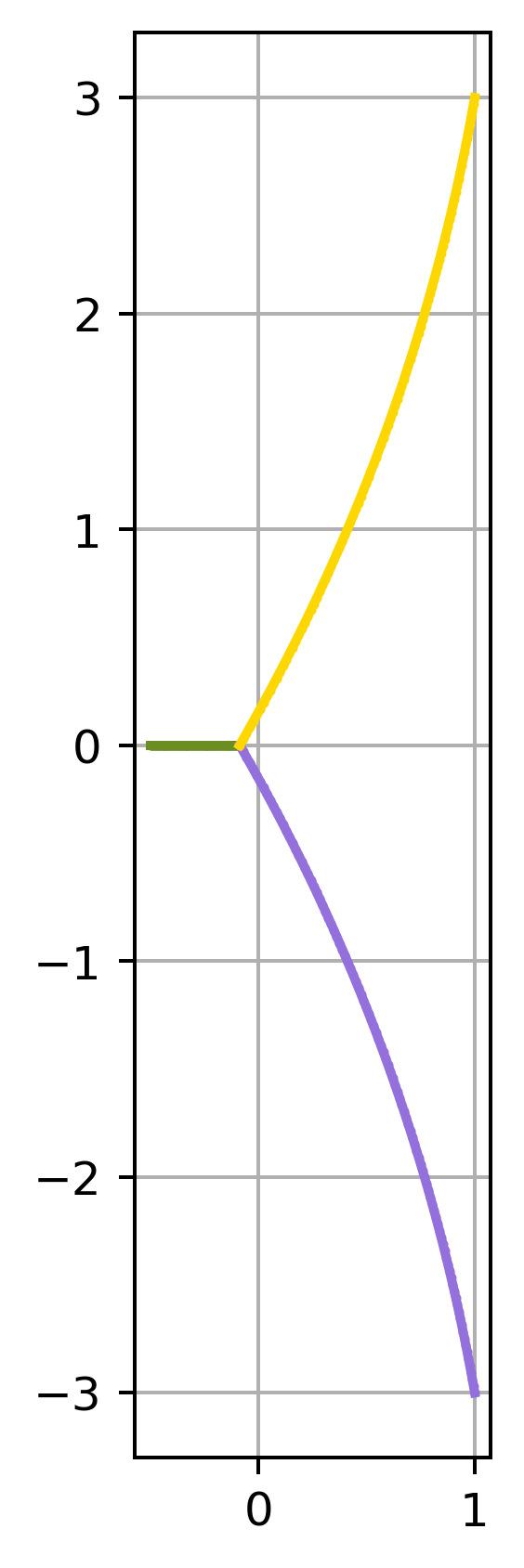}
\includegraphics[angle=-0,width=0.25\textwidth]{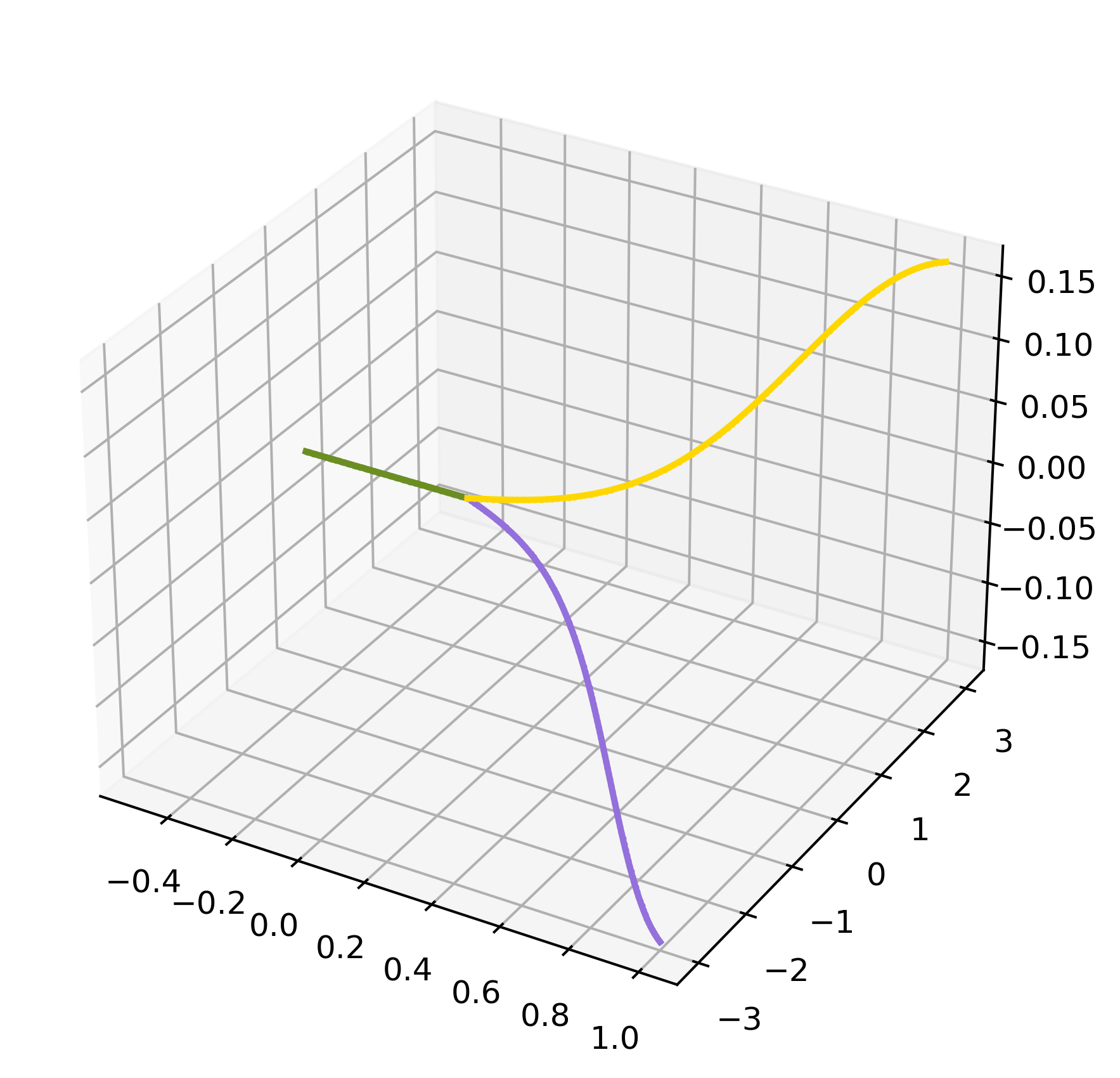}
\includegraphics[angle=-0,width=0.25\textwidth]{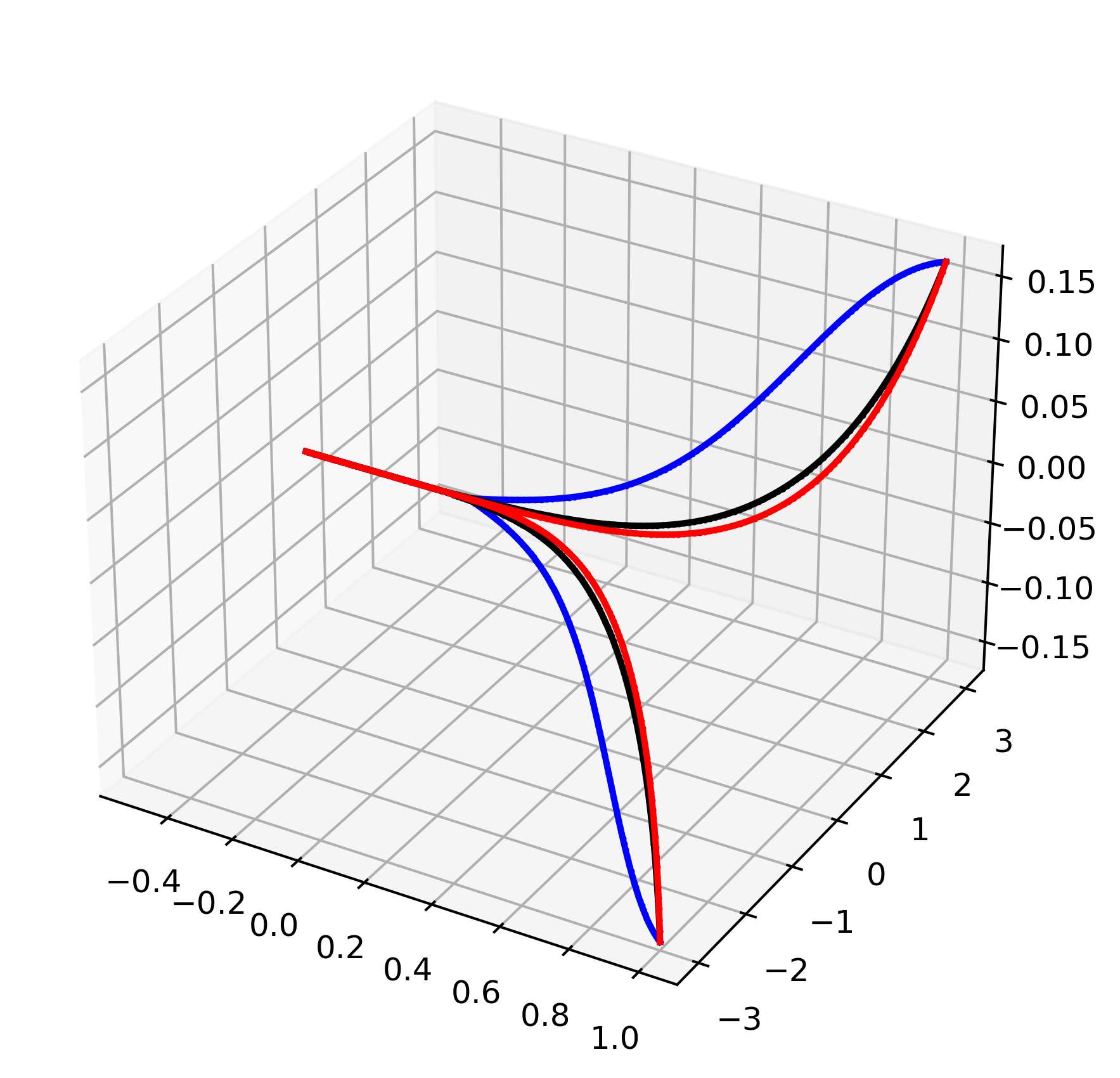}
\includegraphics[angle=-0,width=0.25\textwidth]{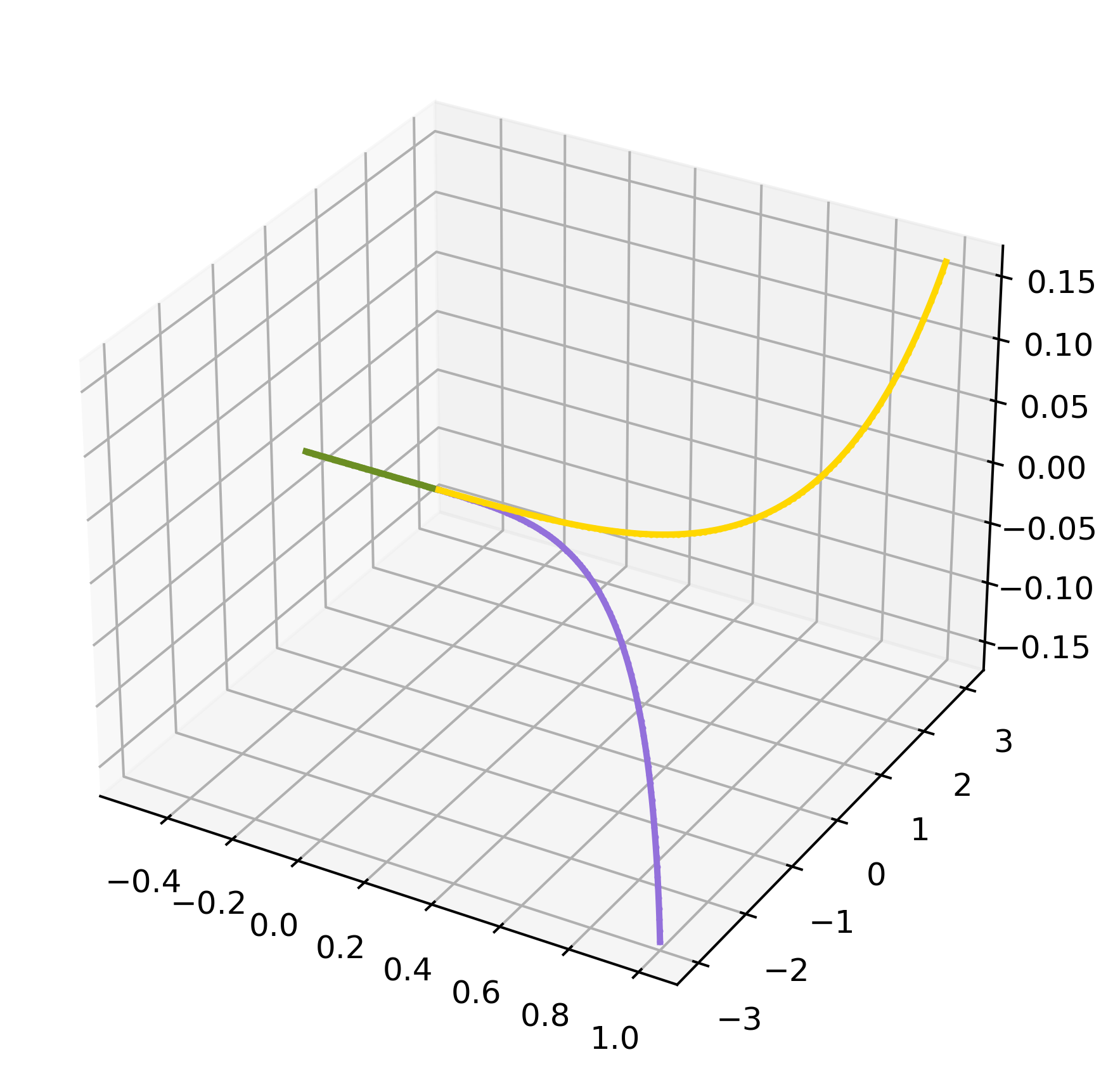}
\includegraphics[angle=-0,width=0.4\textwidth]{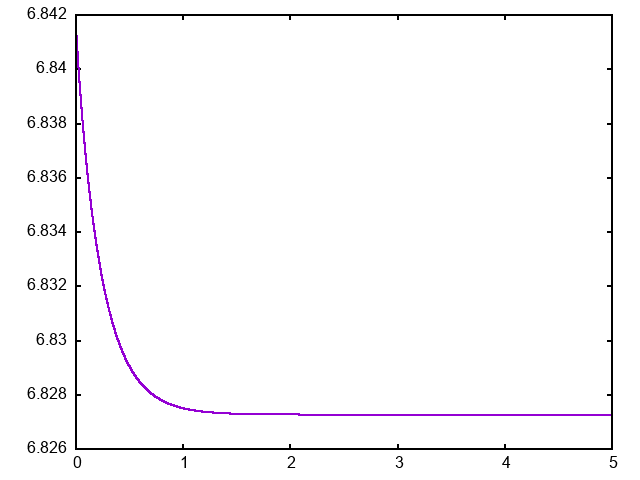}
\includegraphics[angle=-0,width=0.4\textwidth]{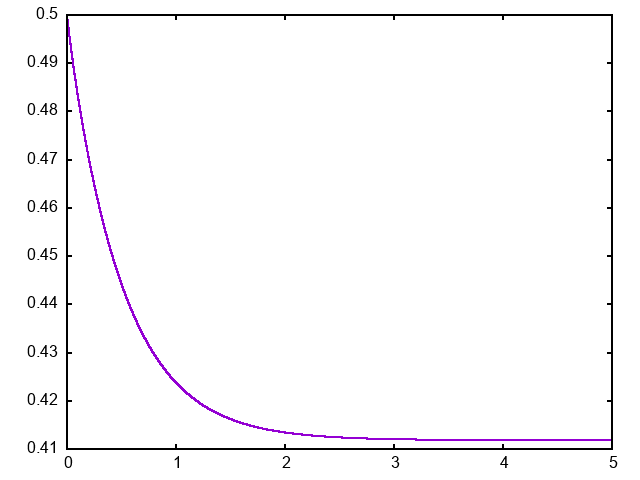}
\caption{The projected triod at time $t=0$ and at time $t=5$. 
The triod in $\R^3$ is shown at time $t=0$, at times $t=0,1,5$, and at 
time $t=5$.
Below we show a plot of the discrete energies $L_E(c^m)$ and $L_E(c^m_1)$
over time.
}
\label{fig:steinerbad}
\end{figure}

{\bf Experiment 5}:
We begin with an analogue of the simulation shown in
Figure~\ref{fig:steinerbadstraight}, but now with planar initial data that
satisfies the $120^\circ$ triple junction angle condition. Note that despite
the projections $\proj{P_\alpha}$, $\alpha=1,2,3$, being the same as before, due to the
curved nature of the initial curves the points $P_\alpha$, $\alpha=1,2,3$, here do not
lie in the two dimensional Euclidean plane. In fact, here we have
$P_1=(-0.5,0,0)^{t}$, $P_2 = (1,-3, -0.167)^{t}$, $P_3 = (1,3,0.167)^{t}$.
We show the results from our finite element approximation in 
Figure~\ref{fig:steinerbad}, where we see once more that for the flow
considered in this paper, no singularity arises from the chosen initial data.

\begin{figure}
\center
\includegraphics[angle=-0,width=0.07\textwidth]{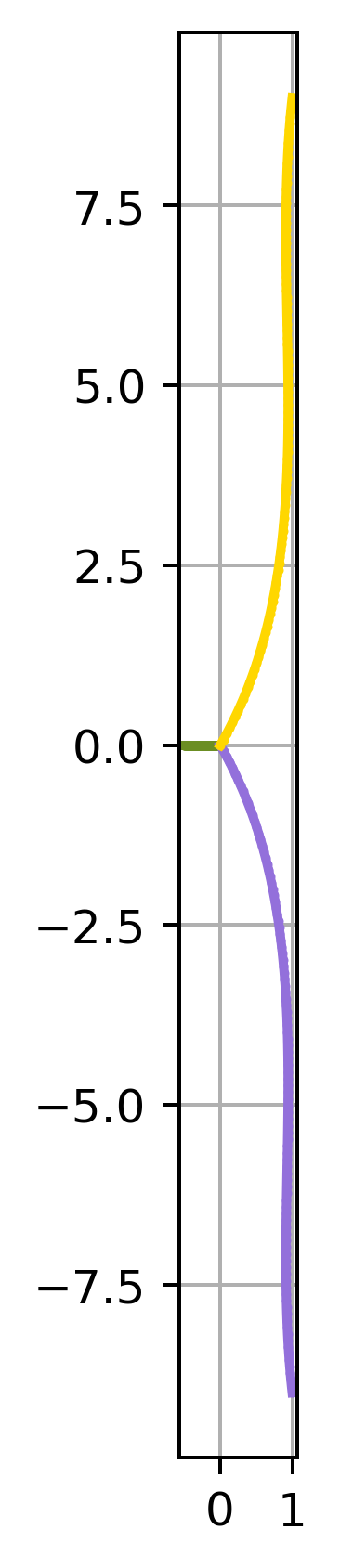}
\includegraphics[angle=-0,width=0.07\textwidth]{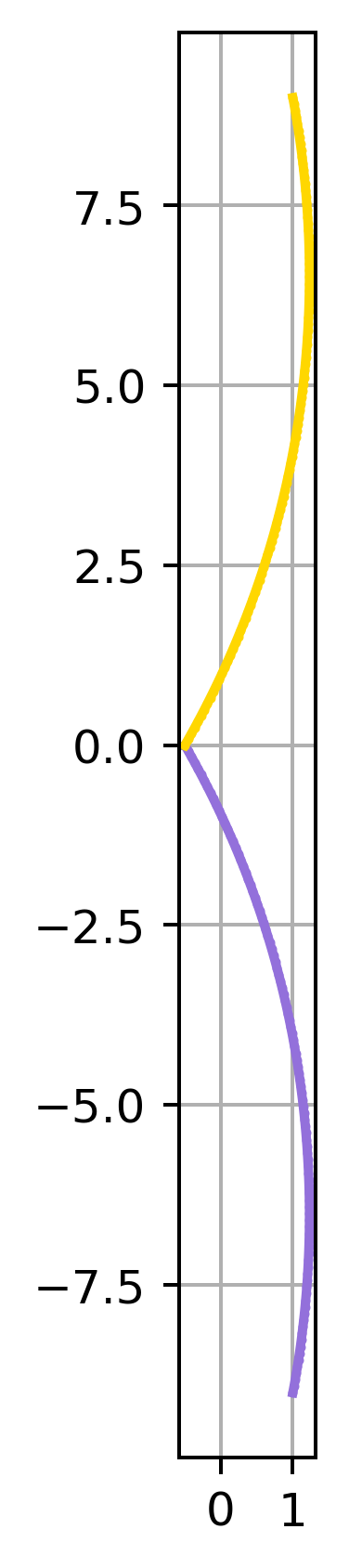}
\includegraphics[angle=-0,width=0.25\textwidth]{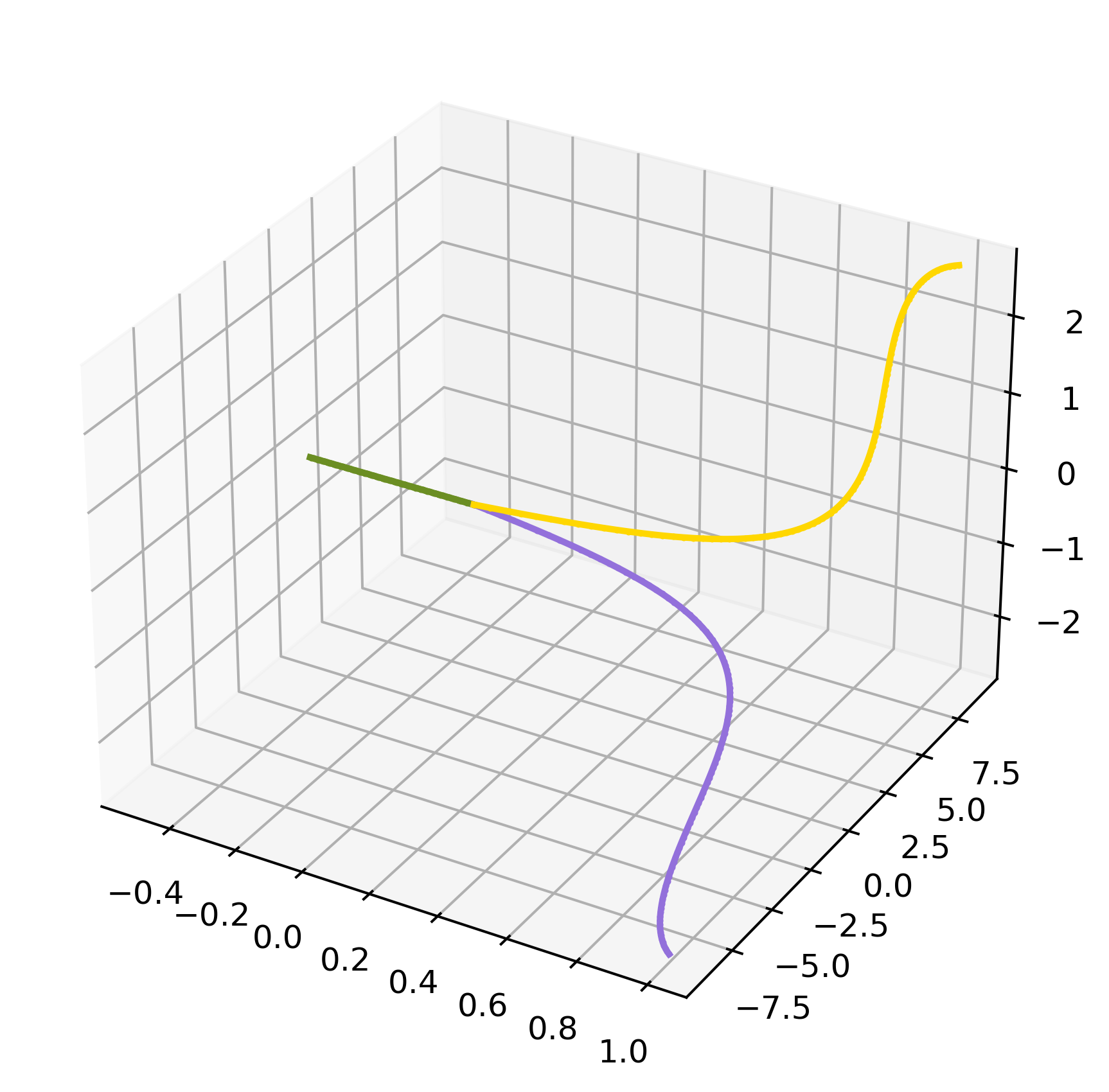}
\includegraphics[angle=-0,width=0.25\textwidth]{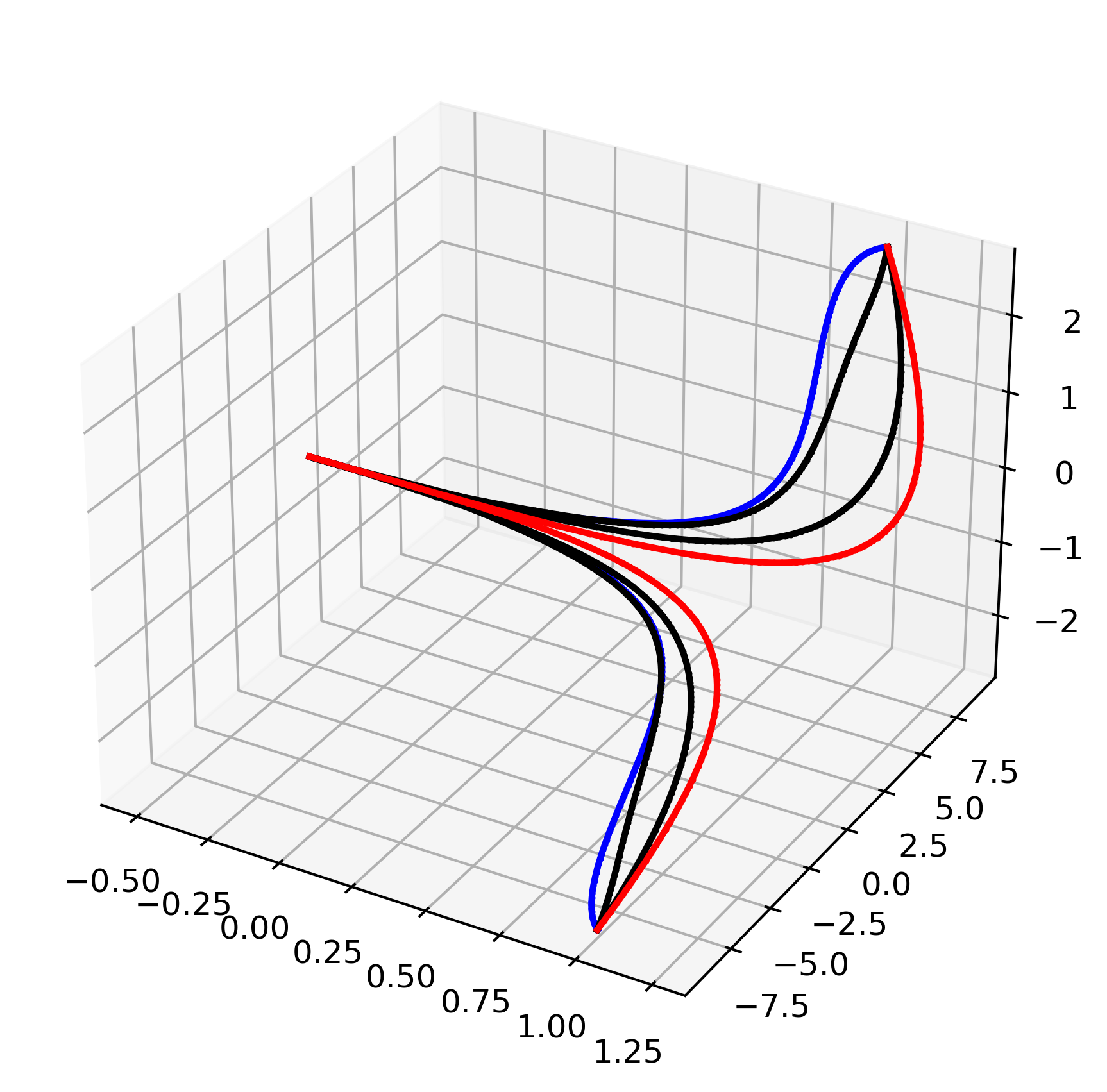}
\includegraphics[angle=-0,width=0.25\textwidth]{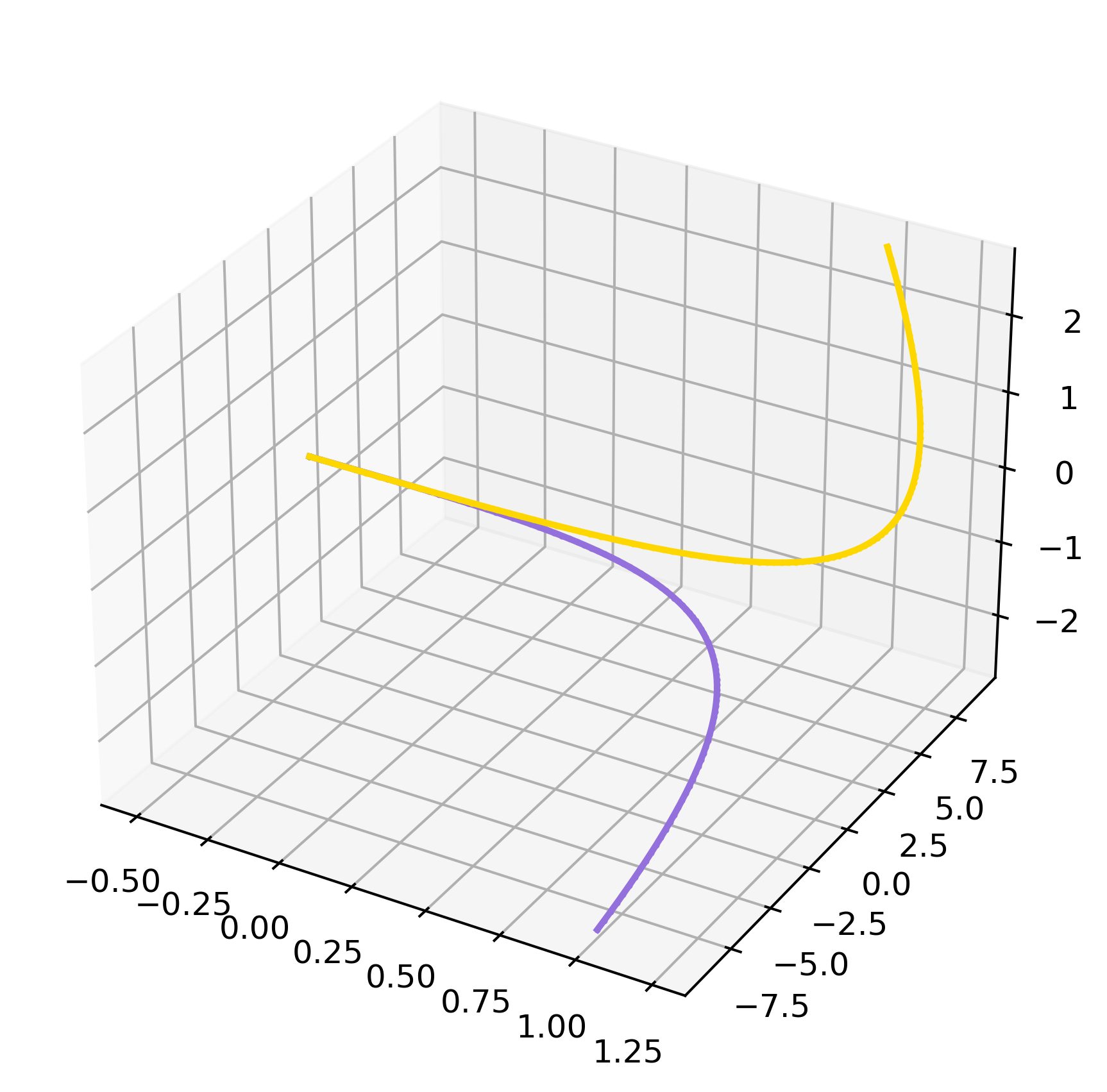}
\includegraphics[angle=-0,width=0.4\textwidth]{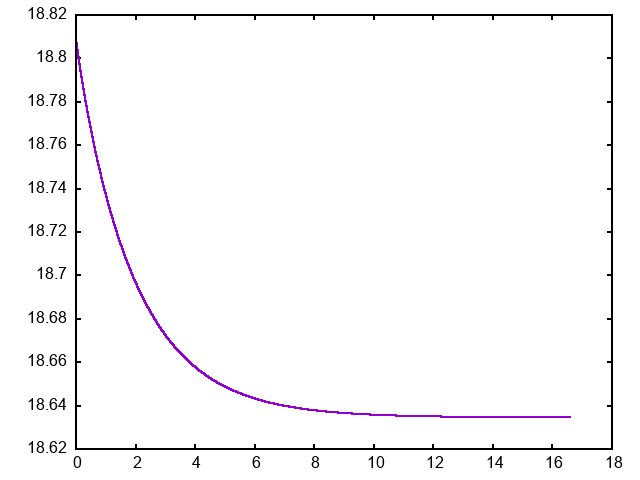}
\includegraphics[angle=-0,width=0.4\textwidth]{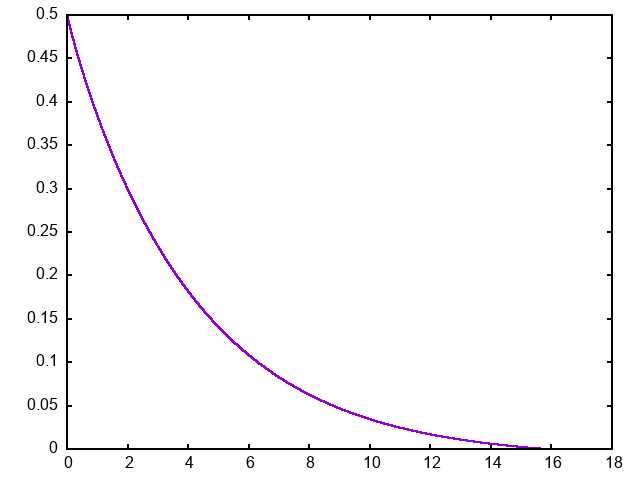}
\caption{The projected triod at time $t=0$ and at time $t=16.6$.
The triod in $\R^3$ is shown at time $t=0$, at times $t=0,1,5,16.6$ and at
time $t=16.6$.
Below we show a plot of the discrete energies $L_E(c^m)$ and $L_E(c^m_1)$
over time.
}
\label{fig:steinerbada}
\end{figure}

{\bf Experiment 6:}
Correspondingly, we show an analogue of Figure~\ref{fig:steinerbadastraight}
in Figure~\ref{fig:steinerbada}. That is, we let $P_1=(-0.5,0,0)^{t}$,
$P_2 = (1,-9,-2.74)^{t}$, $P_3 = (1,9,2.74)^{t}$.
Similarly to before, the simulation encounters
a singularity, when the shortest curve tries to vanish.

\begin{figure}
\center
\mbox{
\includegraphics[angle=-0,width=0.3\textwidth]{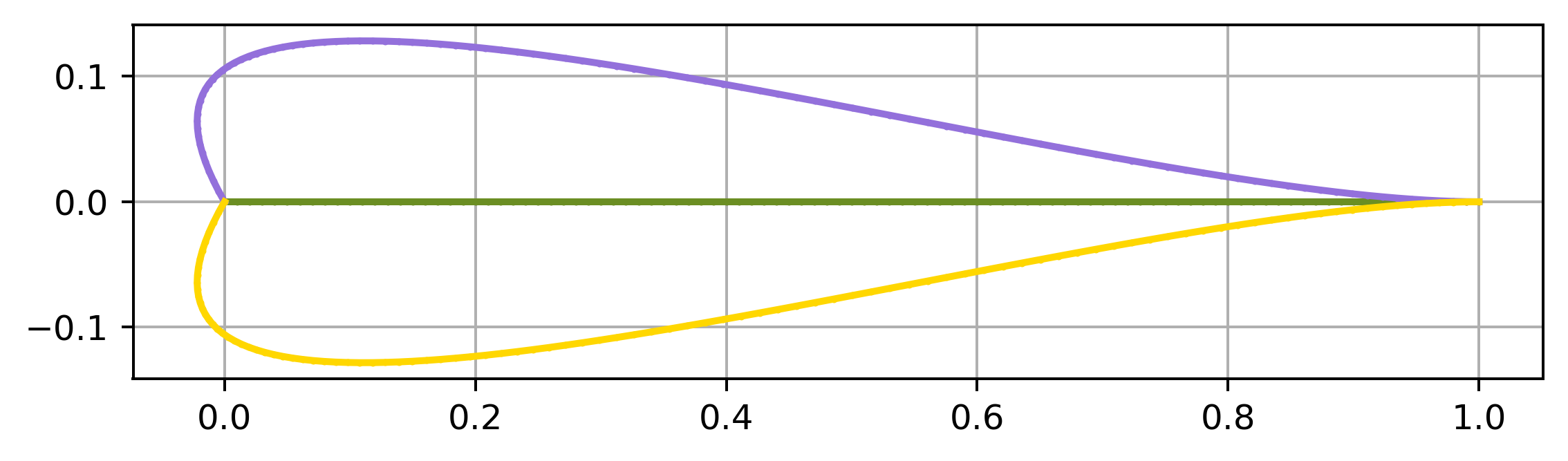}
\includegraphics[angle=-0,width=0.25\textwidth]{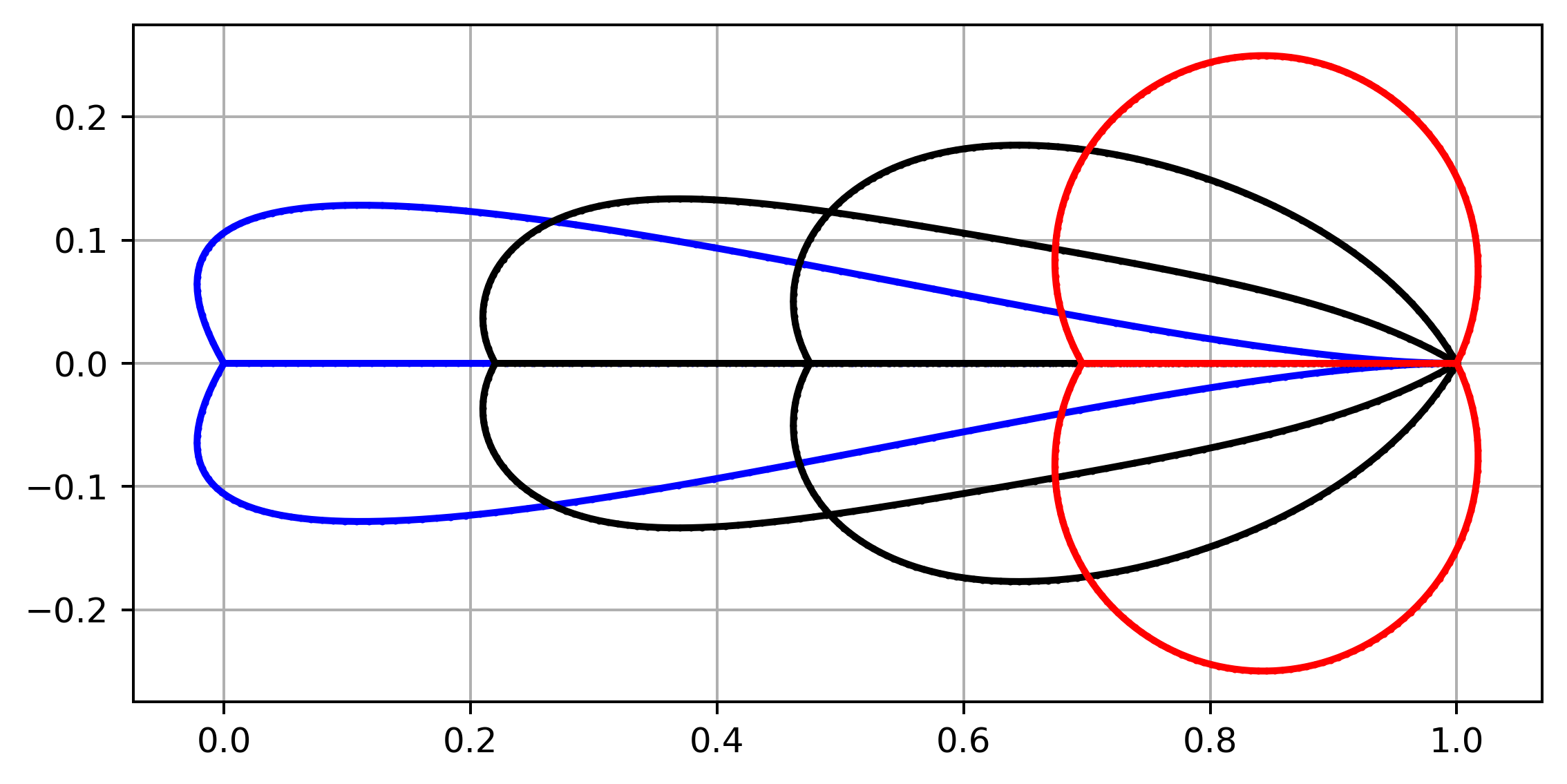}
\includegraphics[angle=-0,width=0.13\textwidth]{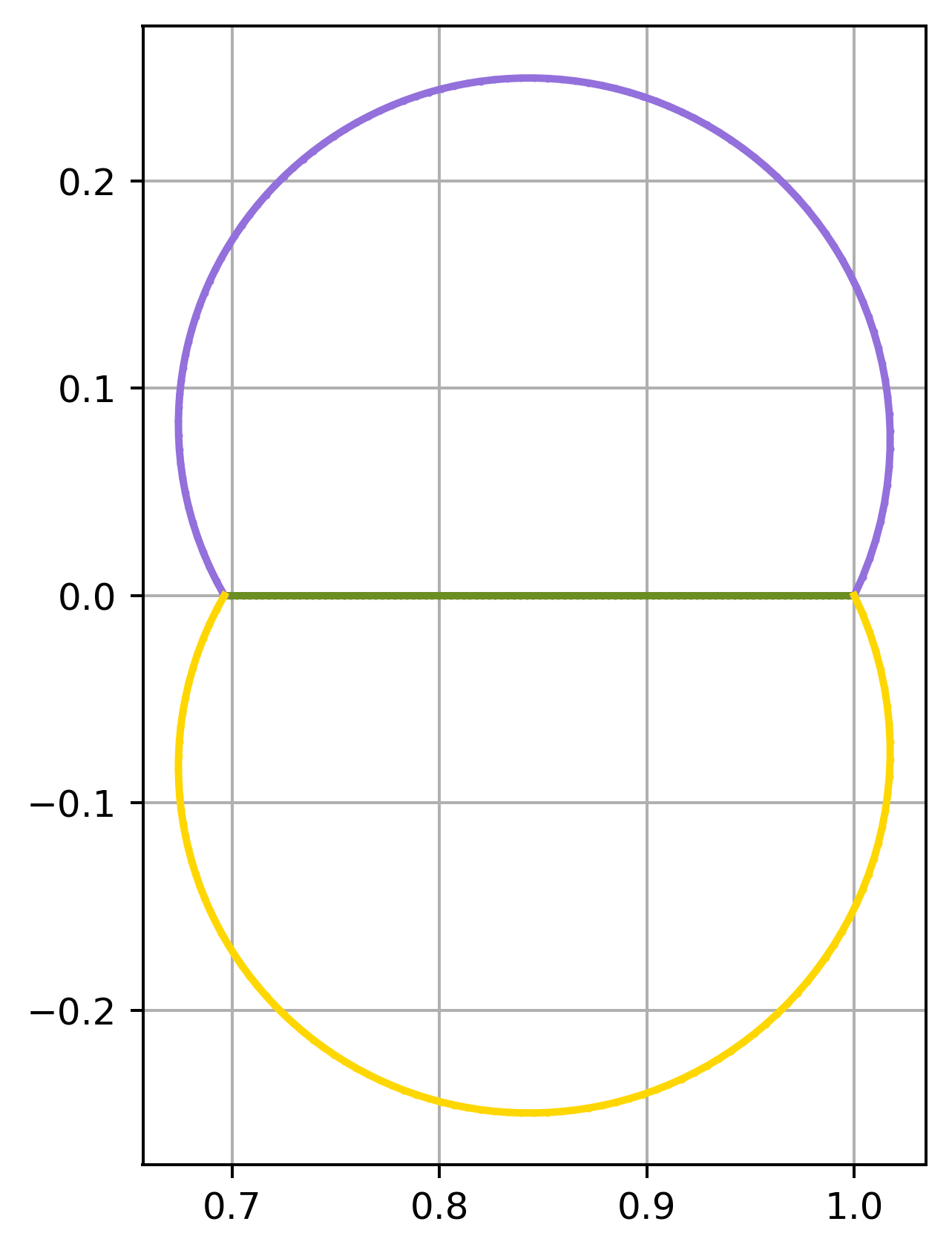}
\includegraphics[angle=-0,width=0.3\textwidth]{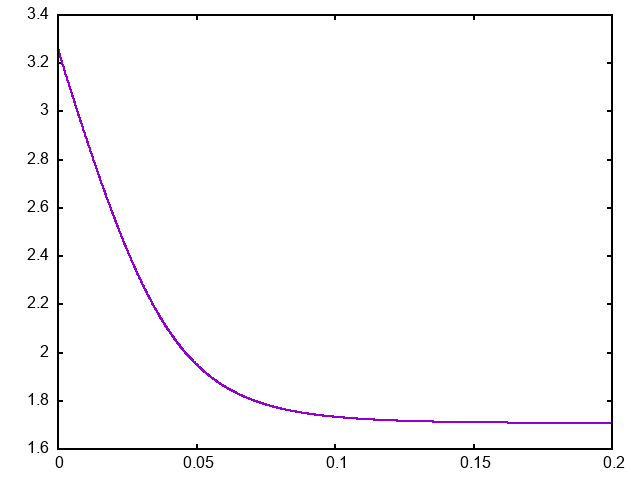}
}
\includegraphics[angle=-0,width=0.3\textwidth]{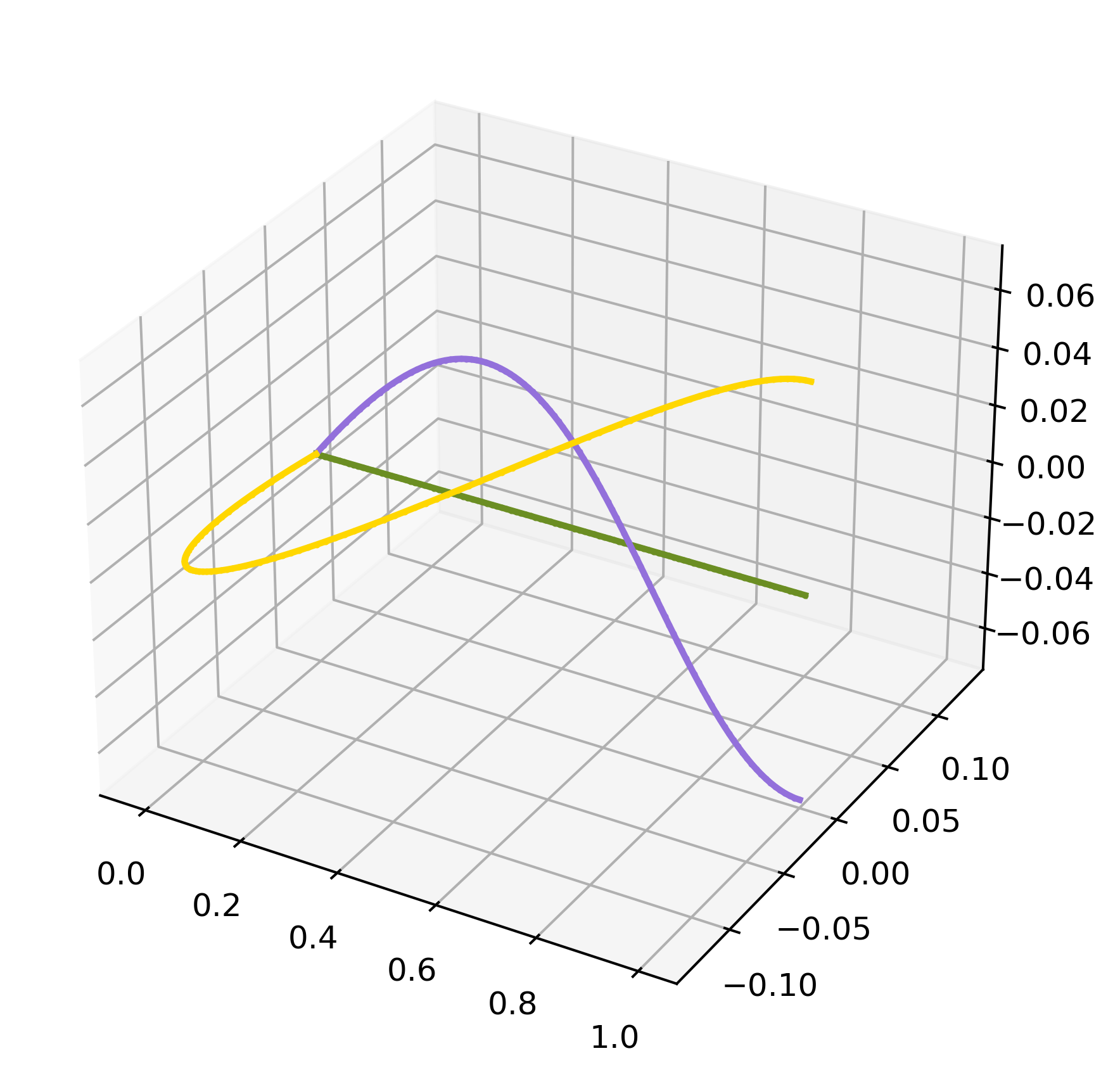}
\includegraphics[angle=-0,width=0.3\textwidth]{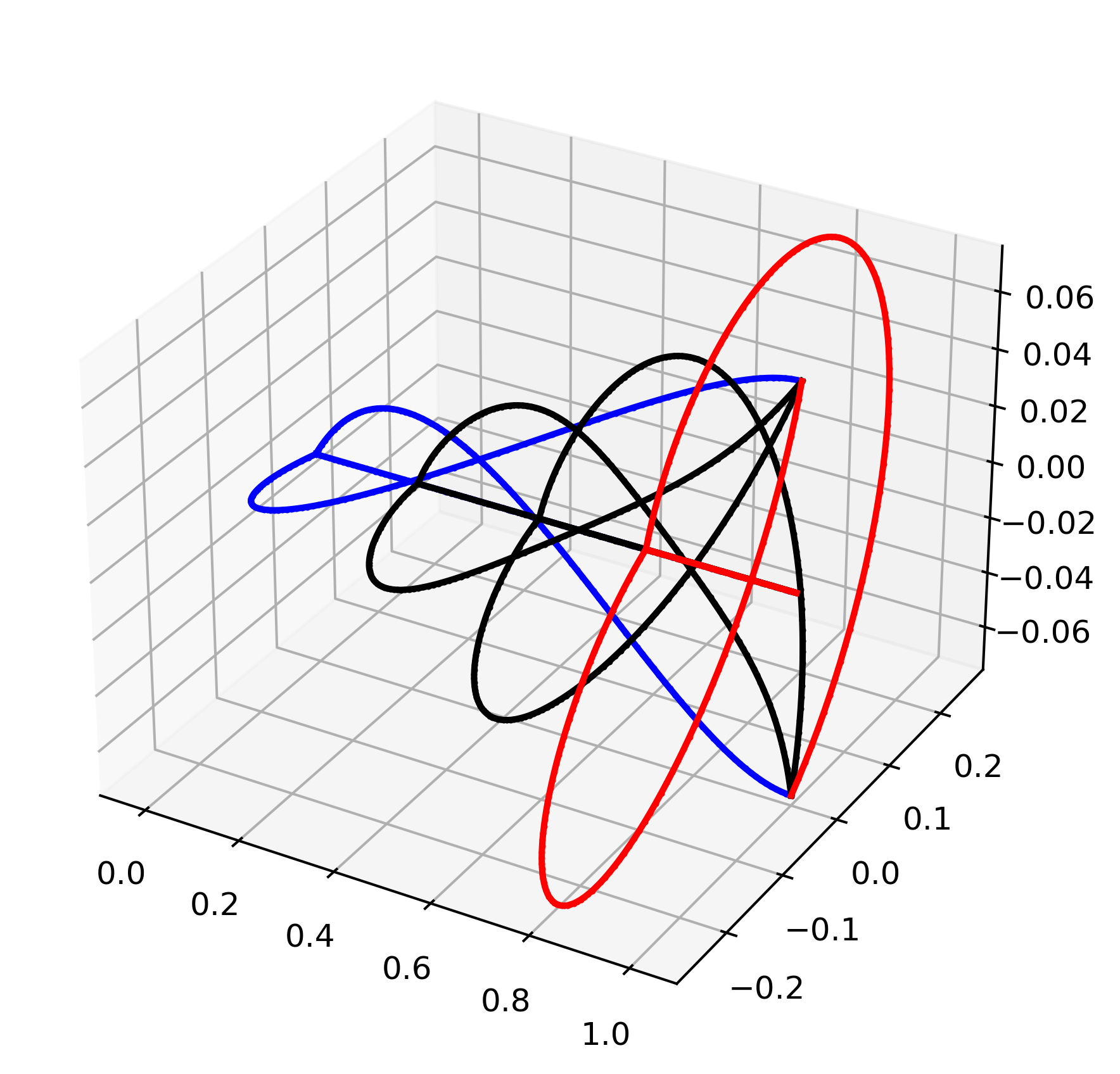}
\includegraphics[angle=-0,width=0.3\textwidth]{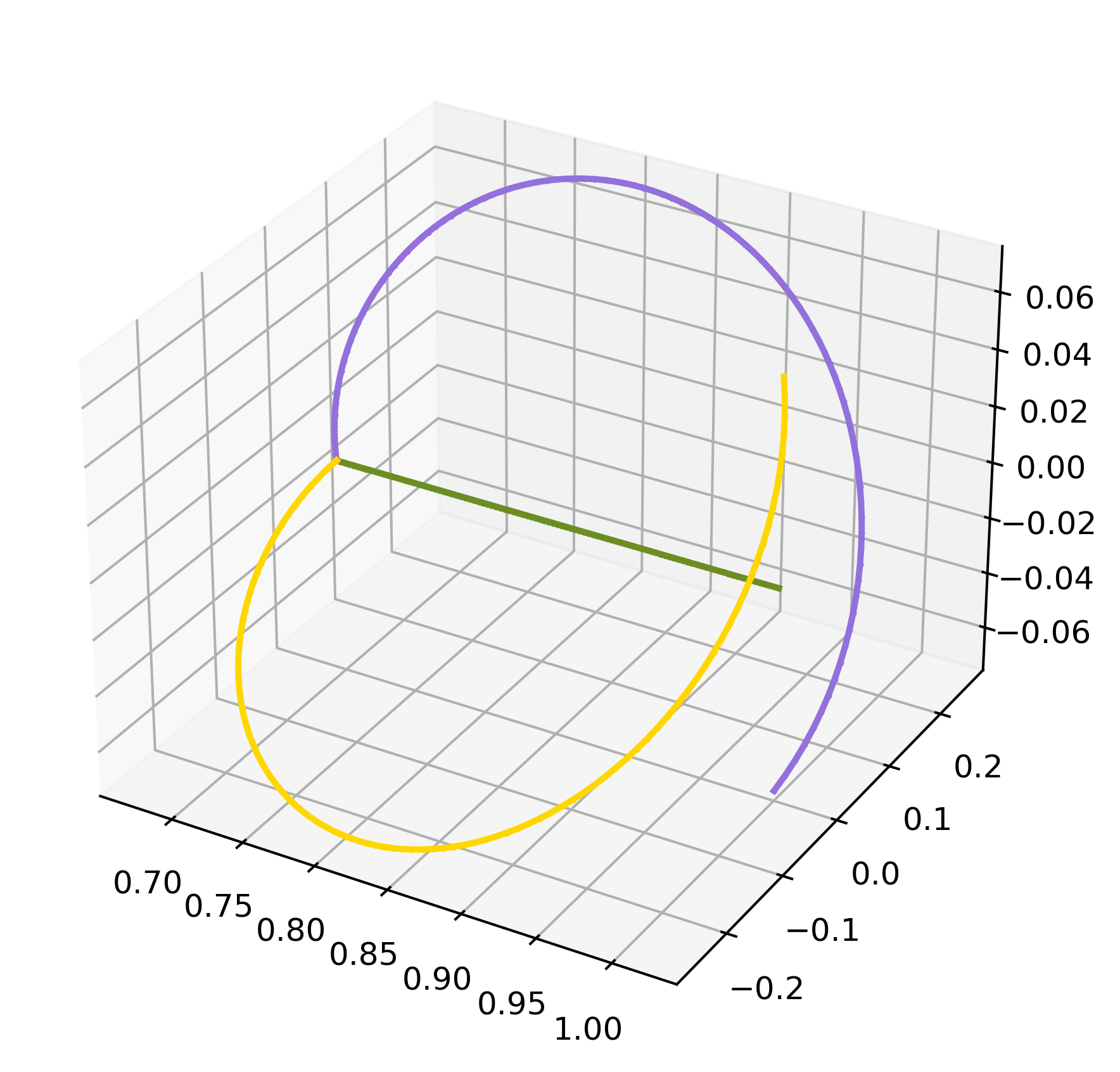}
\caption{The (projected) triod at time $t=0$, at times 
$t=0,0.02,0.05,0.2$ and at time $t=0.2$.
At the top right we show a plot of the discrete energy $L_E(c^m)$ over time.
}
\label{fig:threesame}
\end{figure}

\begin{figure}
\center
\includegraphics[angle=-0,width=0.3\textwidth]{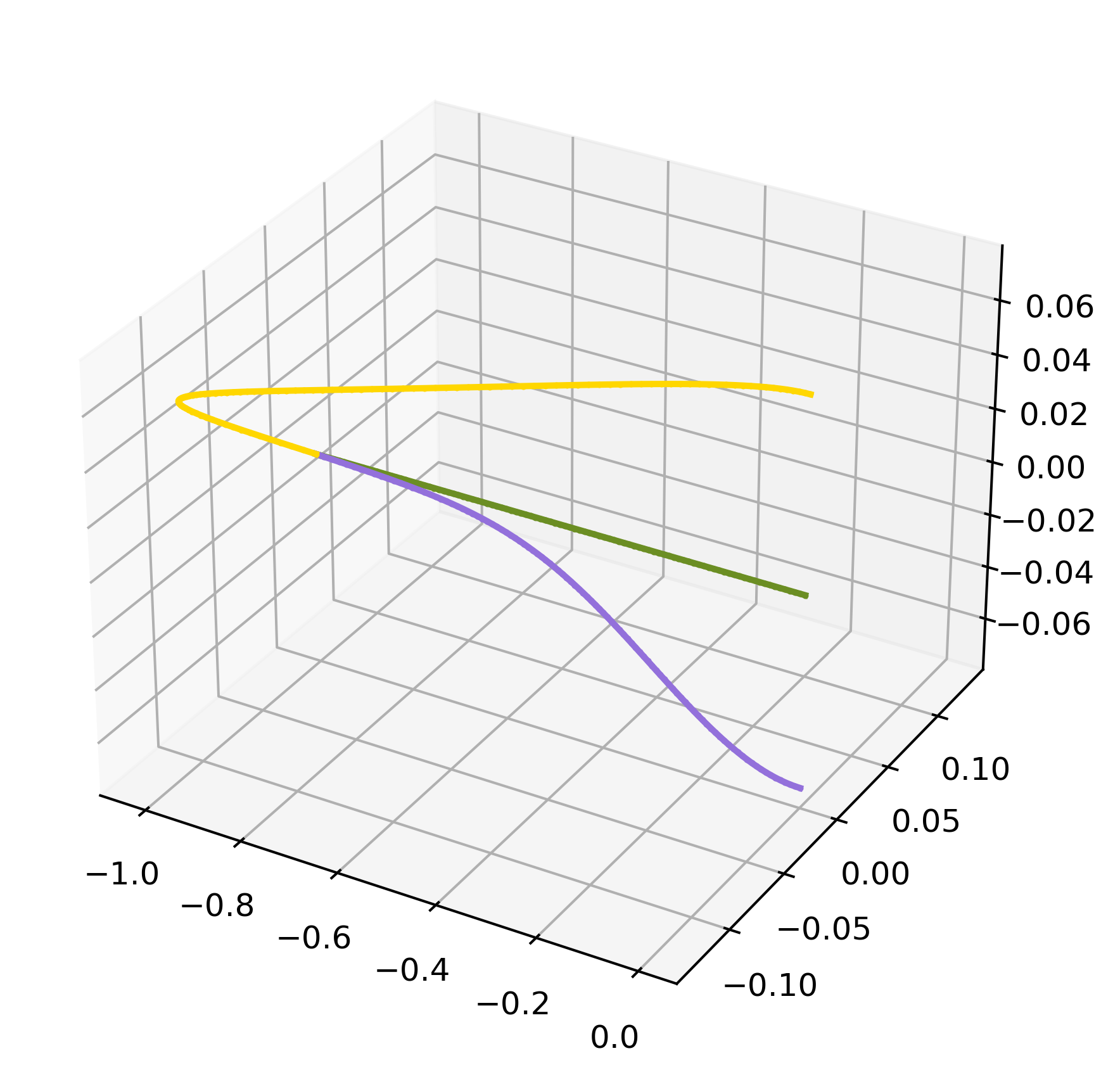}
\includegraphics[angle=-0,width=0.3\textwidth]{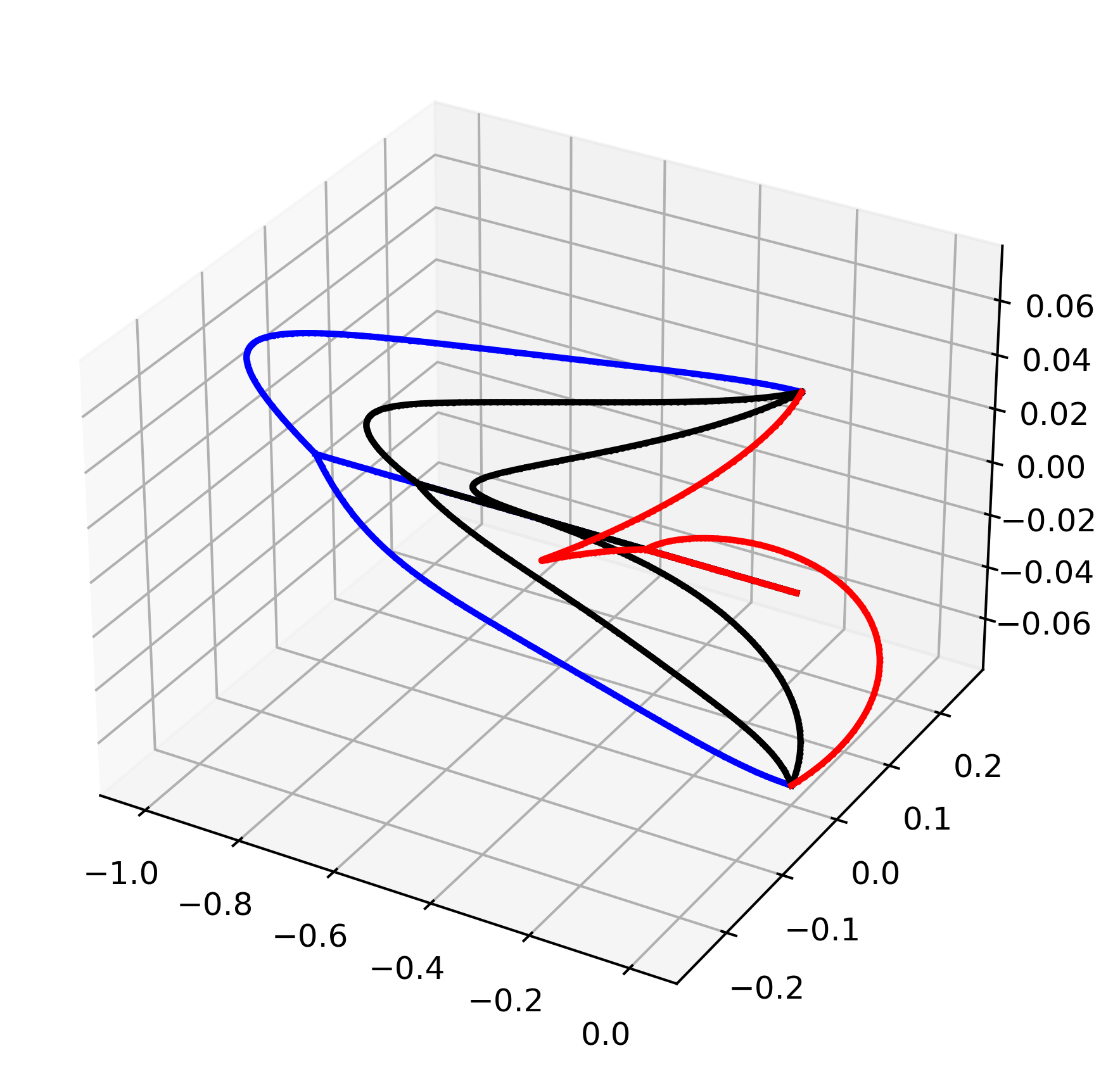}
\includegraphics[angle=-0,width=0.3\textwidth]{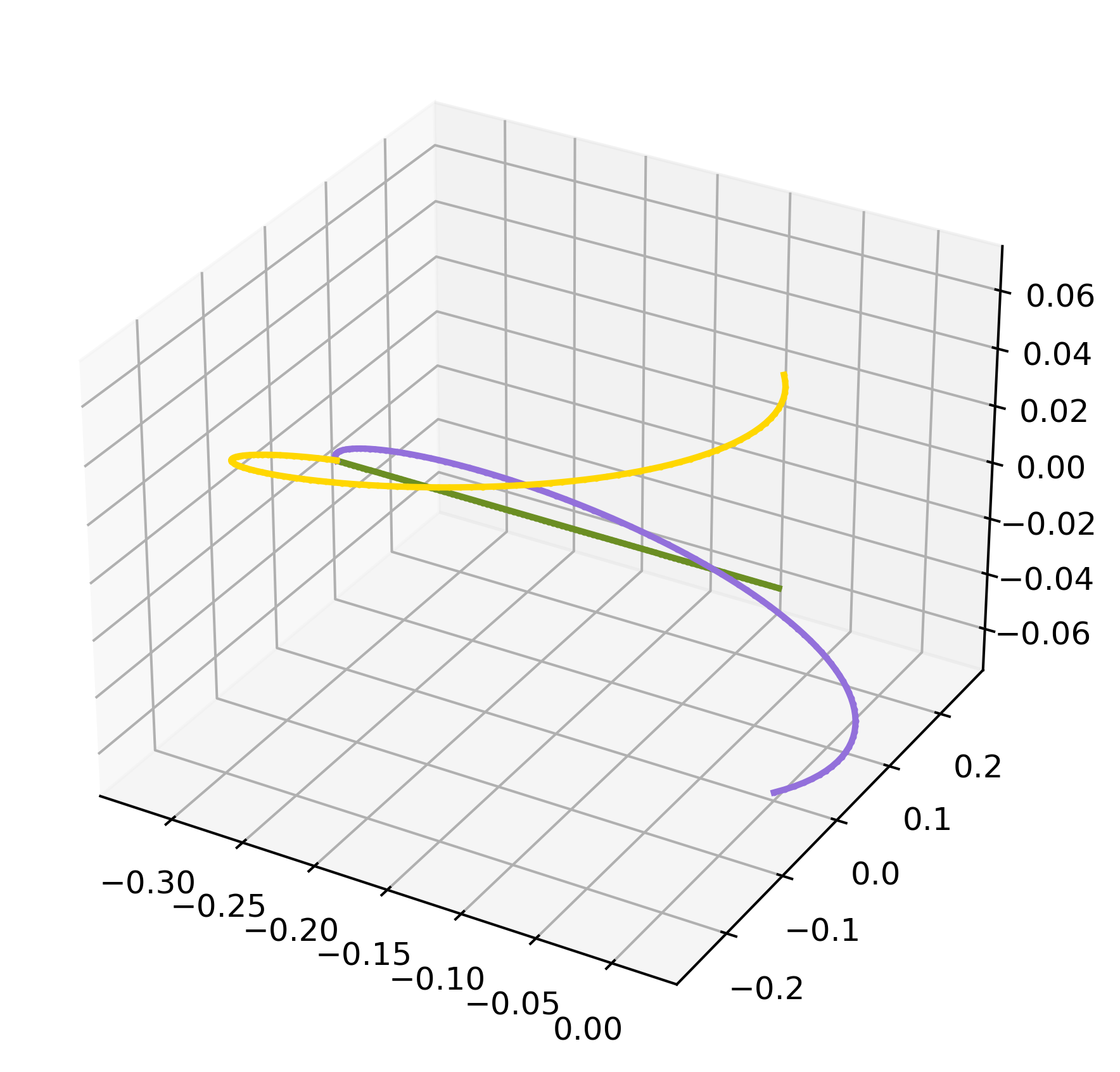}
\caption{The triod at time $t=0$, at times 
$t=0,0.02,0.05,0.2$ and at time $t=0.2$.
}
\label{fig:threesameb}
\end{figure}

{\bf Experiment 7}:
Here we choose the 2d projections of $P_\alpha$, $\alpha=1,2,3$, to be the same, i.e.\
all equal to $(1,0)^{t}$, with $\Sigma = (0,0,0)^{t}$. In fact, we let
$P_1=(1,0,0)^{t}$, $P_2=(1,0,-0.07)^{t}$ and $P_3=(1,0,0.07)^{t}$.
The evolution of the planar
triod then converges towards a standard double bubble in the plane,
see \cite{FoisyABHZ93}. That is because the planar geodesics are given by
straight lines and circle segments, and as they meet at $120^\circ$ at
$\proj\Sigma$, due to symmetry the same applies to the fixed planar triple
junction point $(1,0)^{t}$. We show the results in Figure~\ref{fig:threesame}.
In a near repeat of the previous simulation we shift the planar initial data by
one unit to the left. That is, $\Sigma = (-1, 0, 0)^{t}$ and
$\proj{P_\alpha} = (0,0)^{t}$, $\alpha=1,2,3$. In fact, we have
$P_1=(0,0,0)^{t}$, $P_2=(0,0,-0.07)^{t}$ and $P_3=(0,0,0.07)^{t}$.
Since the projections will evolve exactly as
shown in Figure~\ref{fig:threesame}, we only display the 3d curves for this
experiment in Figure~\ref{fig:threesameb}.

\begin{figure}
\center
\mbox{
\includegraphics[angle=-0,width=0.25\textwidth]{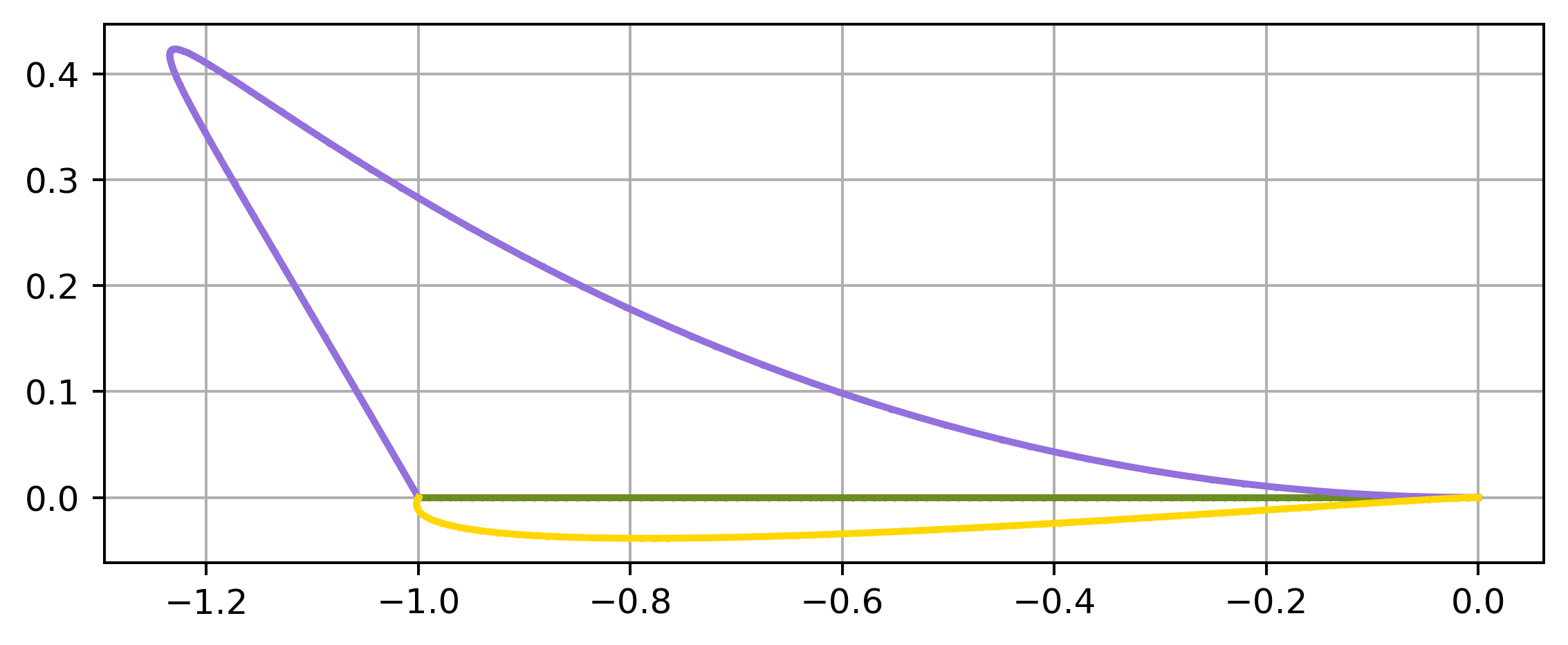}
\includegraphics[angle=-0,width=0.25\textwidth]{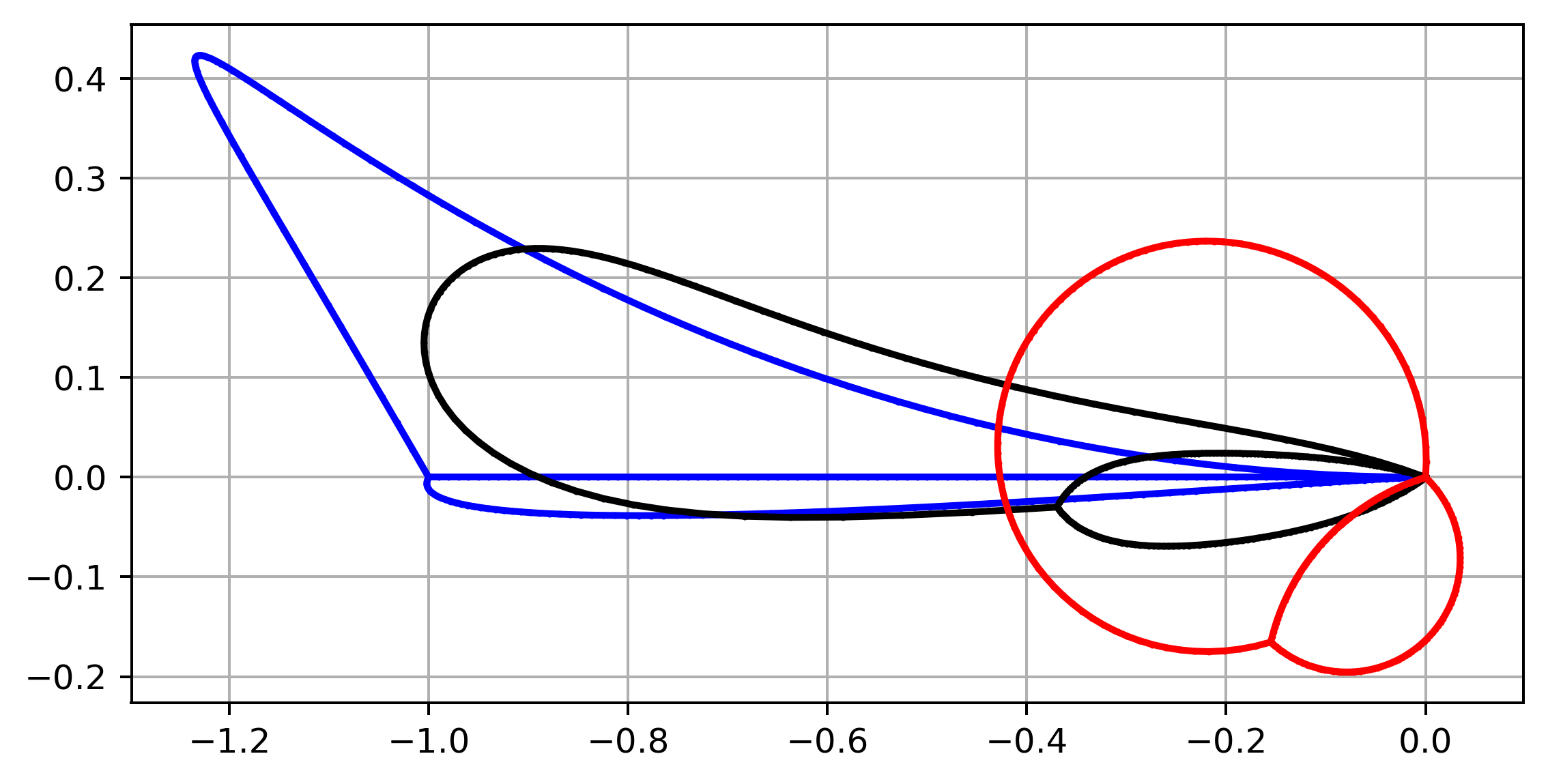}
\includegraphics[angle=-0,width=0.15\textwidth]{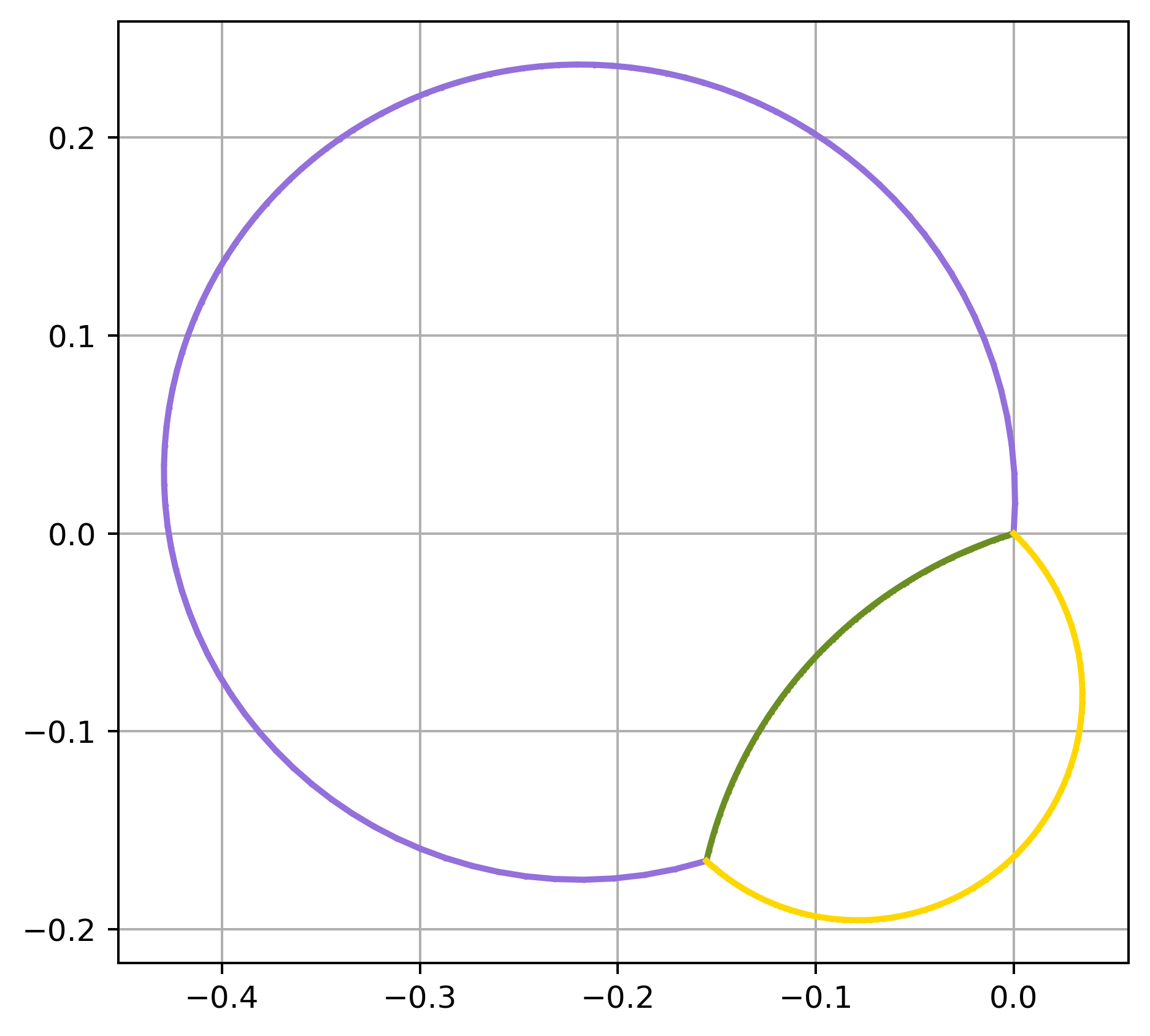}
\includegraphics[angle=-0,width=0.3\textwidth]{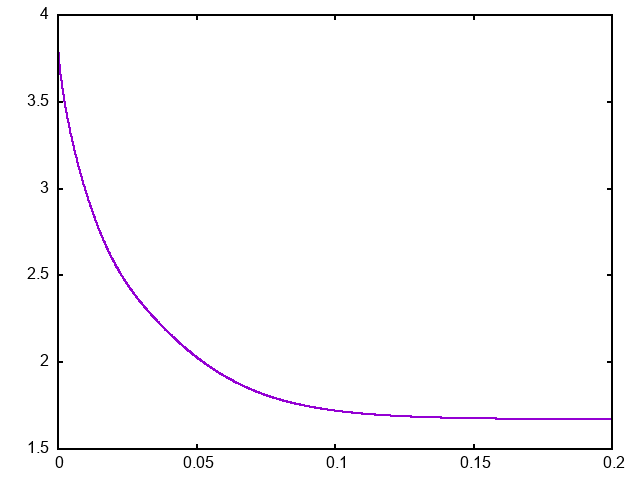}
}
\includegraphics[angle=-0,width=0.3\textwidth]{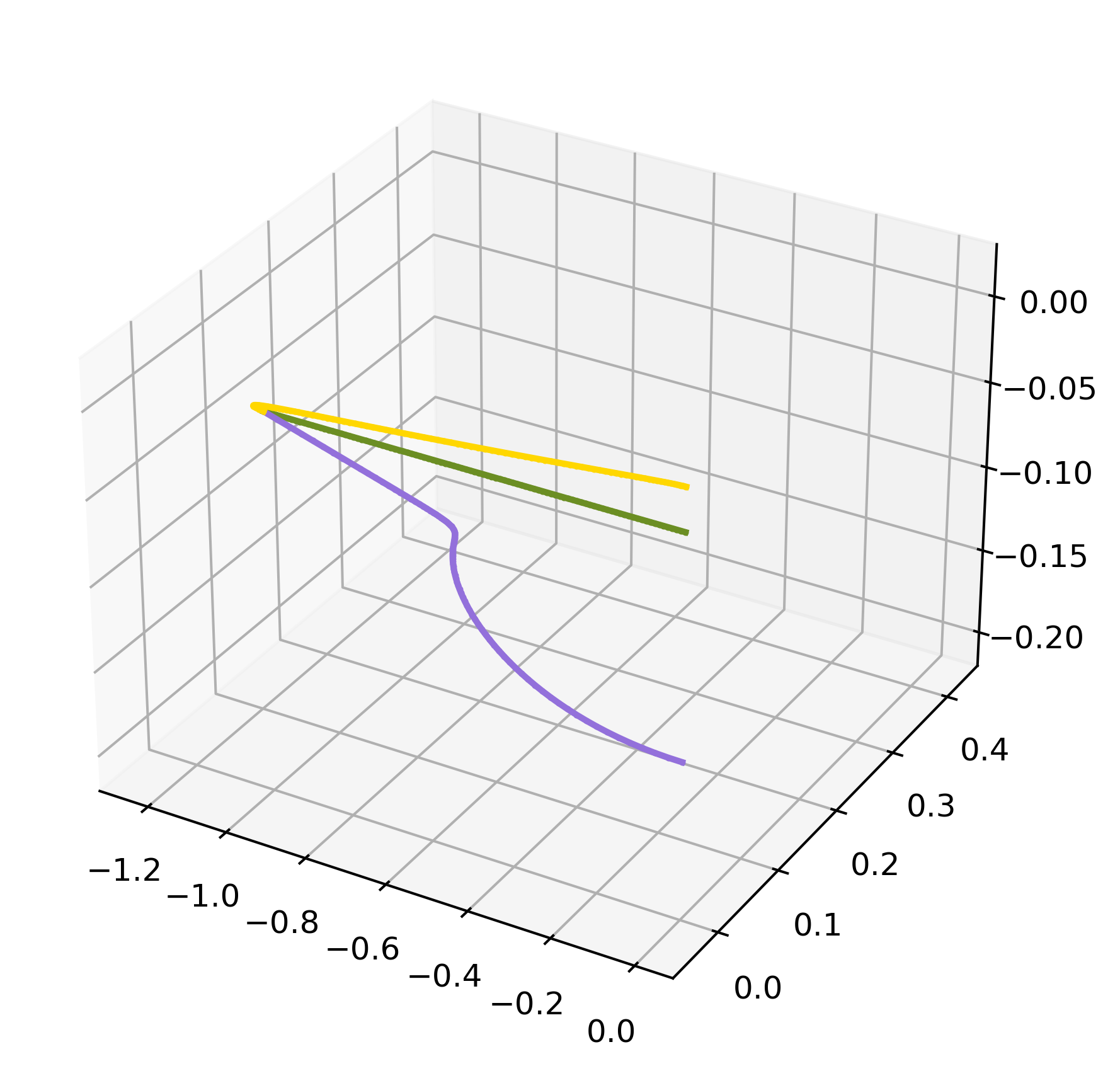}
\includegraphics[angle=-0,width=0.3\textwidth]{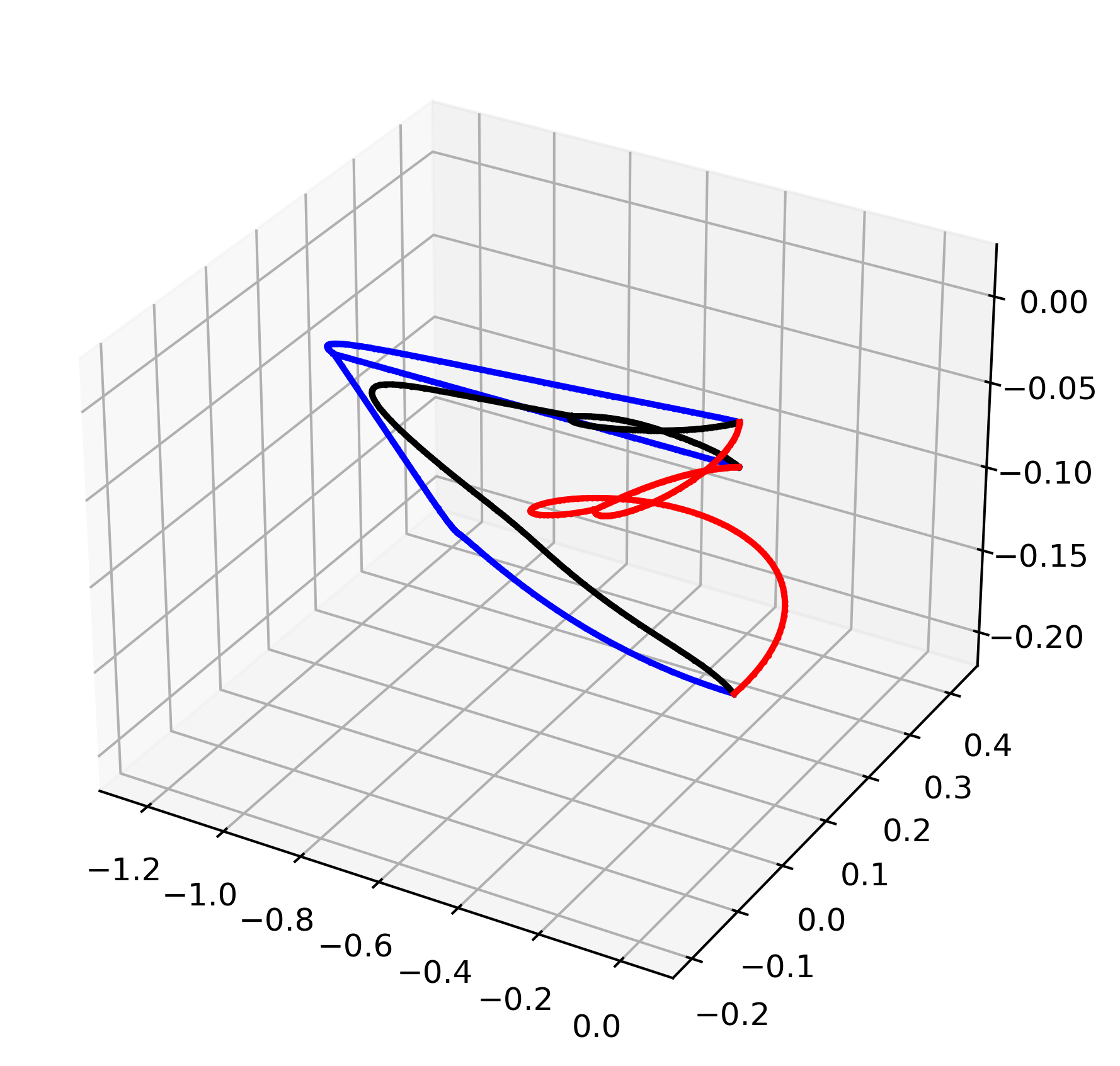}
\includegraphics[angle=-0,width=0.3\textwidth]{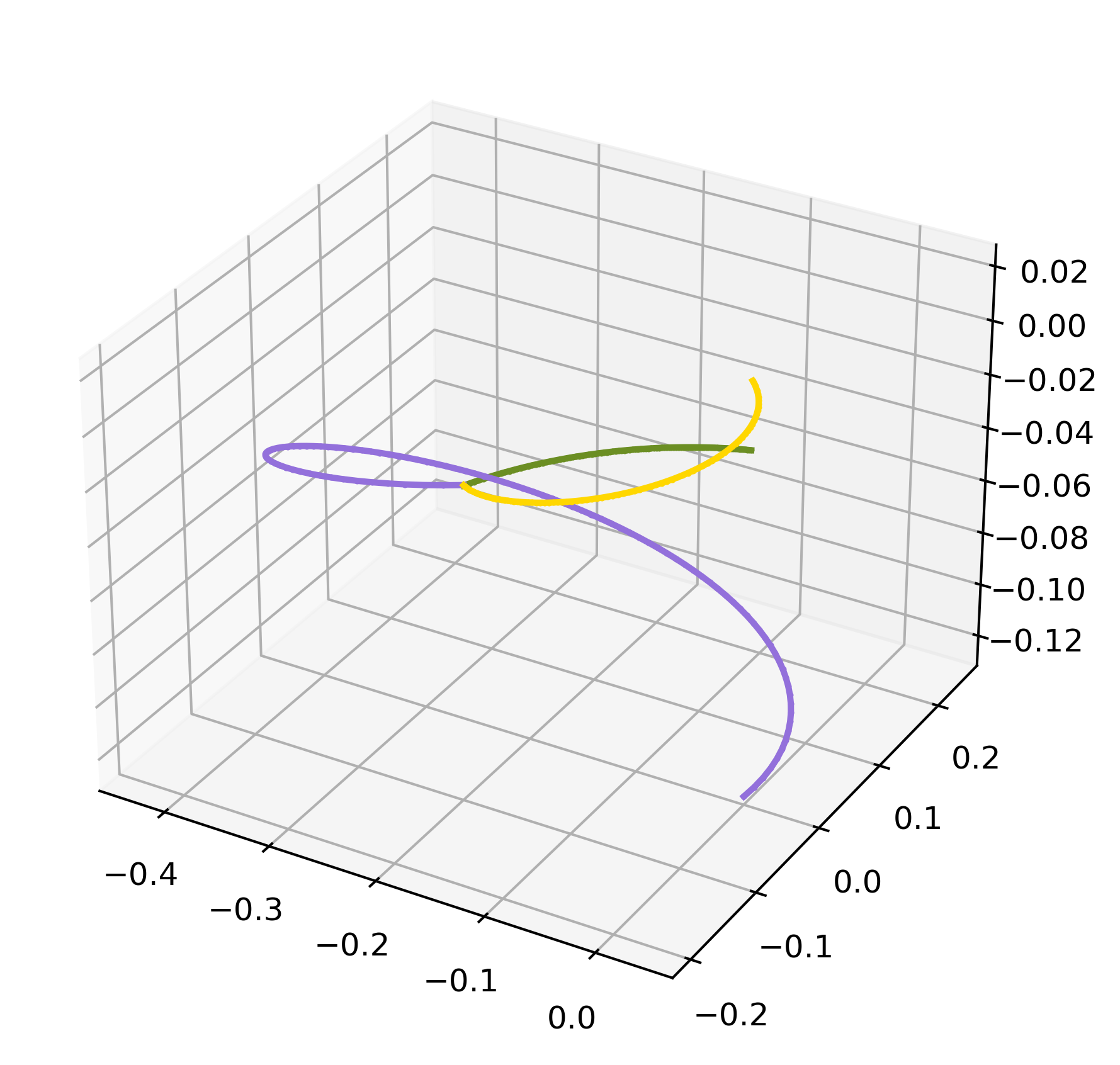}
\caption{The (projected) triod at time $t=0$, at times 
$t=0,0.02,0.2$ and at time $t=0.2$.
At the top right we show a plot of the discrete energy $L_E(c^m)$ over time.
}
\label{fig:threesamec}
\end{figure}

{\bf Experiment 8}:
In this experiment we choose the initial position of $\Sigma$ to be $(-1,0,0)^{t}$
and let $\proj{P_\alpha} = (0,0)^{t}$, $\alpha=1,2,3$. In order to break the symmetry between
the planar regions enclosed by $c_1^0$ and $c_2^0$, and by $c_1^0$ and 
$c_3^0$, we choose $c_2^0$ much longer than $c_3^0$. 
As a result, the planar triod now
evolves towards a standard double bubble with non-equal enclosed areas.
We show the results from our numerical method in Figure~\ref{fig:threesamec}.
Notice that the choice of planar initial data means that
$P_1=(0,0,0)^{t}$, $P_2=(0,0,-0.13)^{t}$ and $P_3=(0,0,0.02)^{t}$.

\begin{figure}
\center
\mbox{
\includegraphics[angle=-0,width=0.25\textwidth]{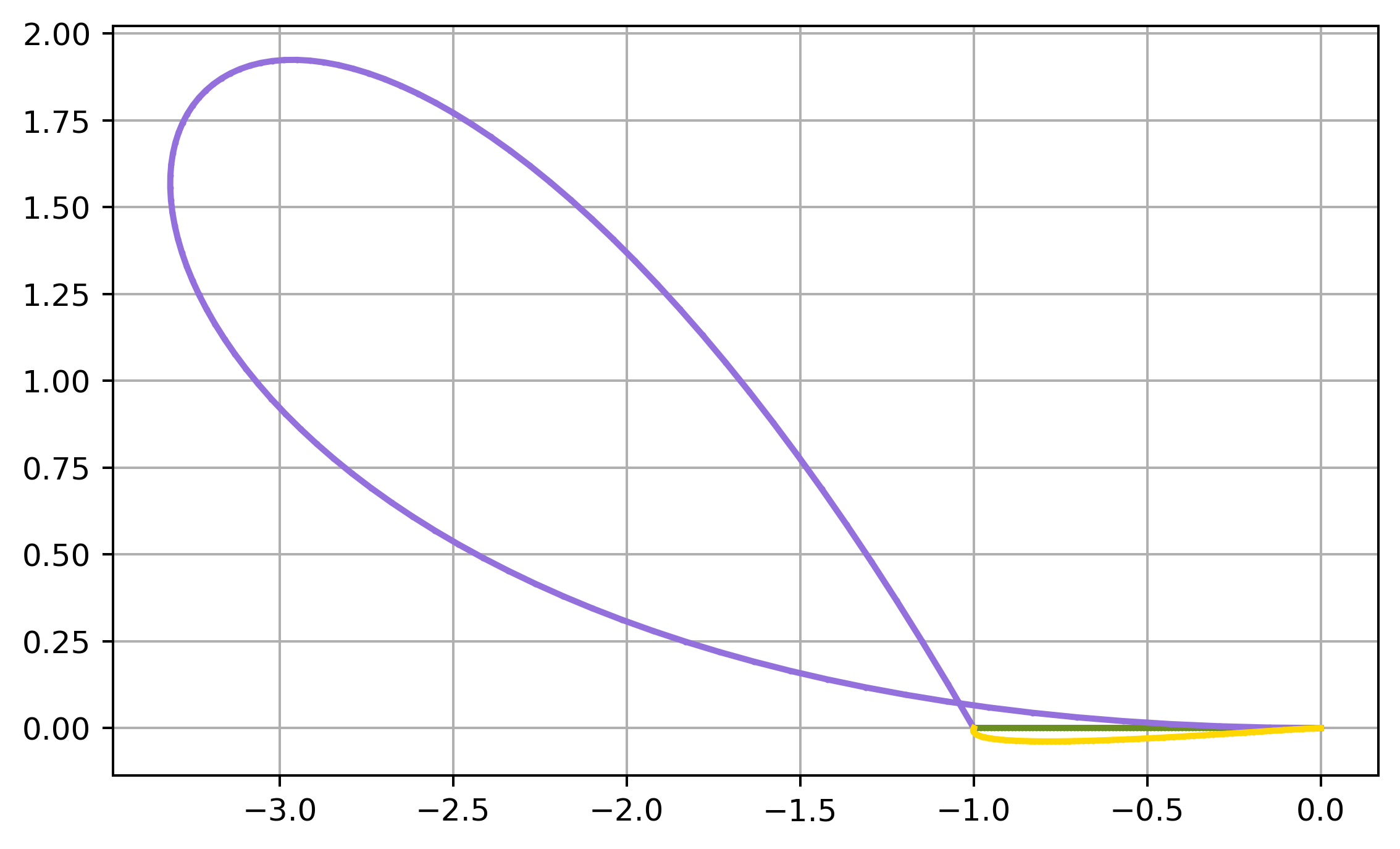}
\includegraphics[angle=-0,width=0.23\textwidth]{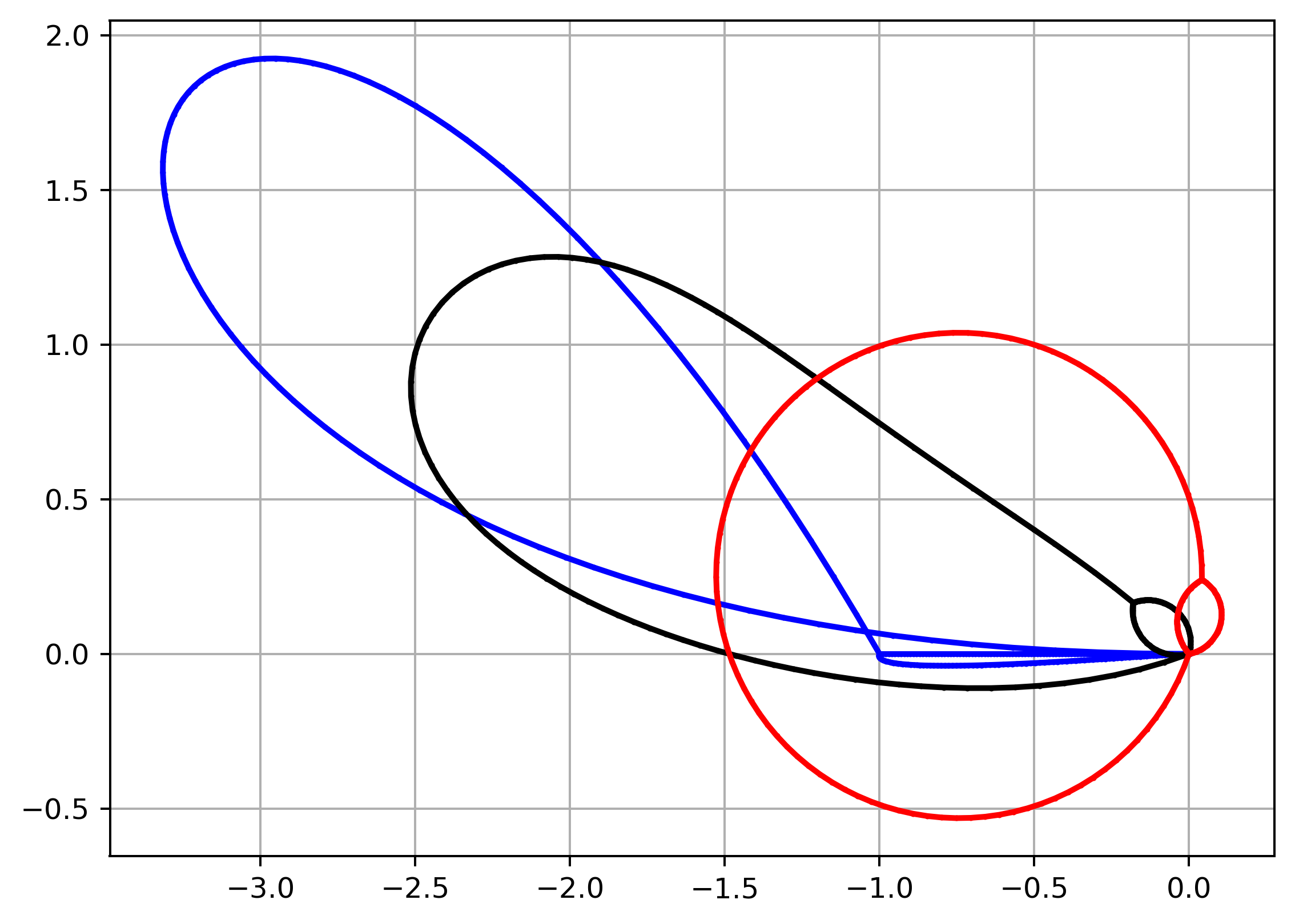}
\includegraphics[angle=-0,width=0.2\textwidth]{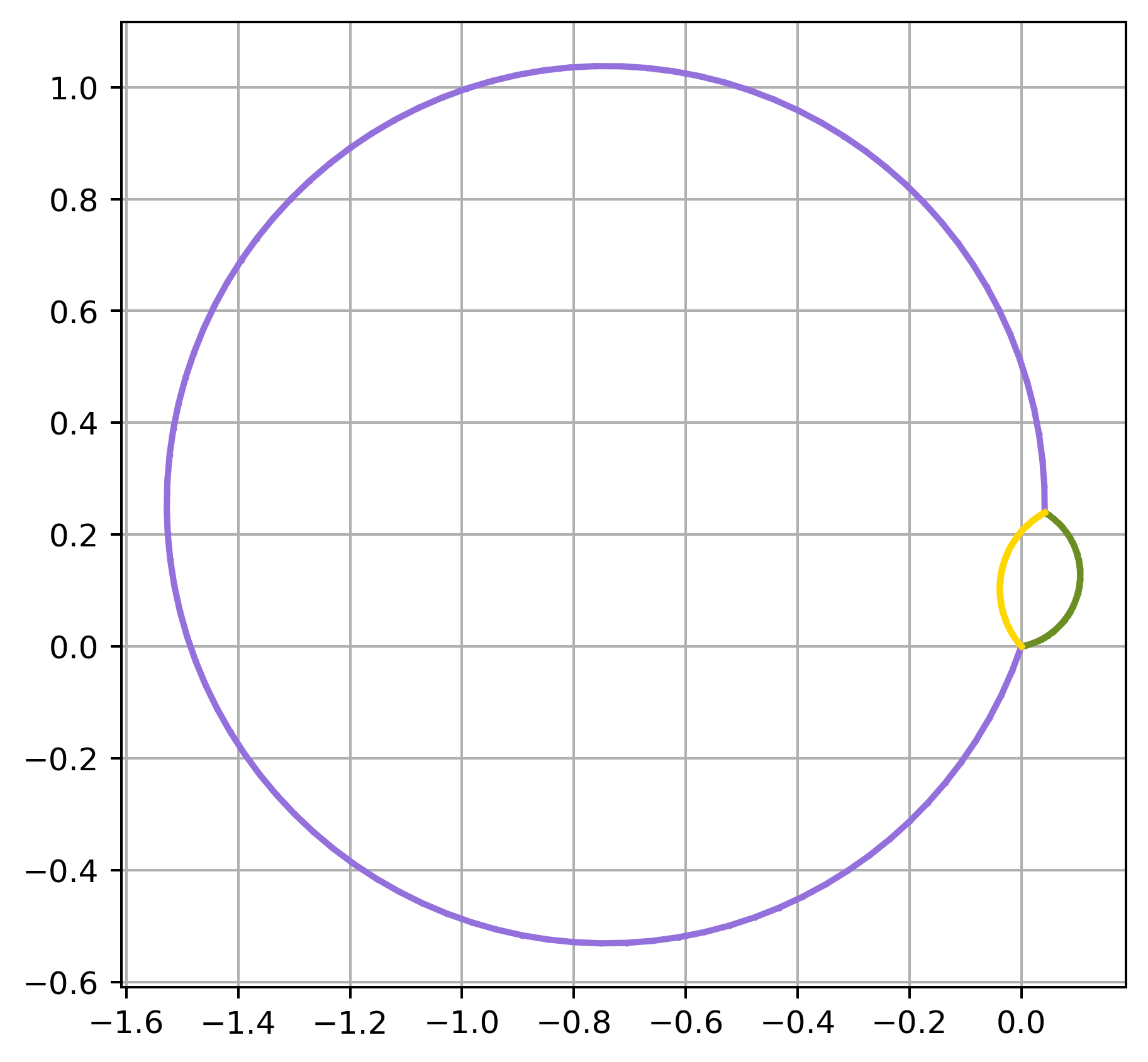}
\includegraphics[angle=-0,width=0.3\textwidth]{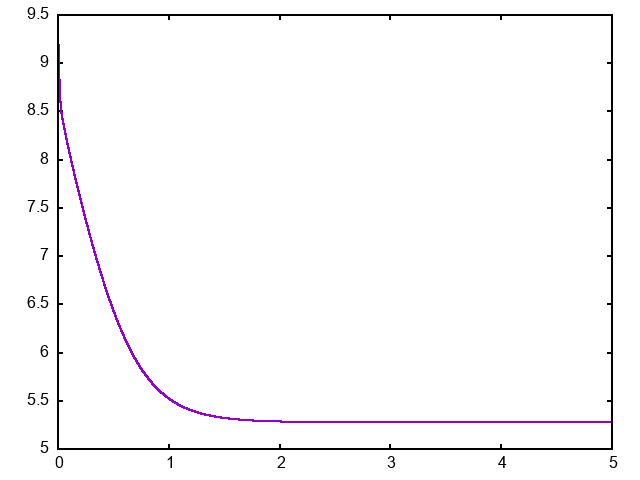}
}
\includegraphics[angle=-0,width=0.3\textwidth]{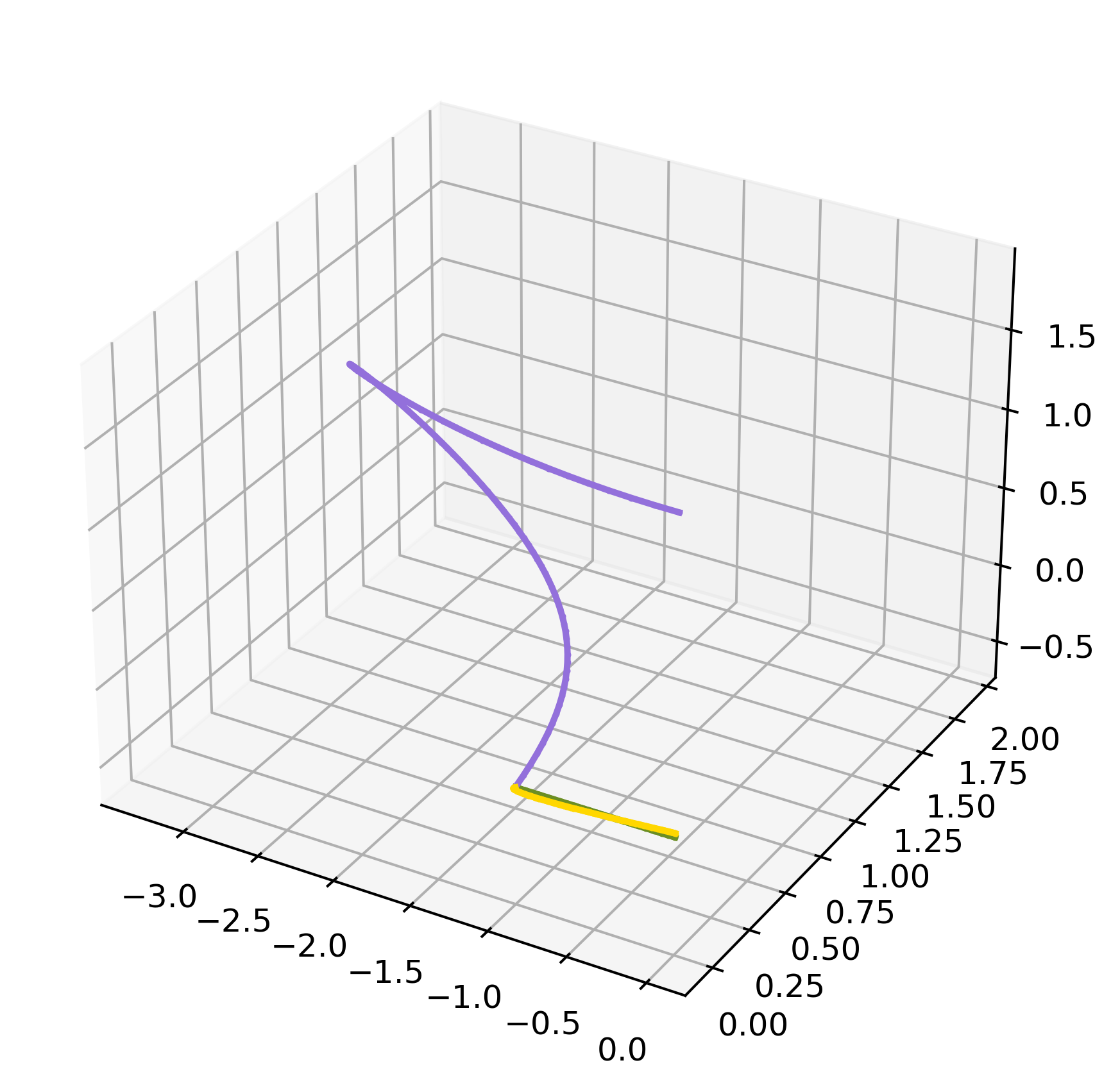}
\includegraphics[angle=-0,width=0.3\textwidth]{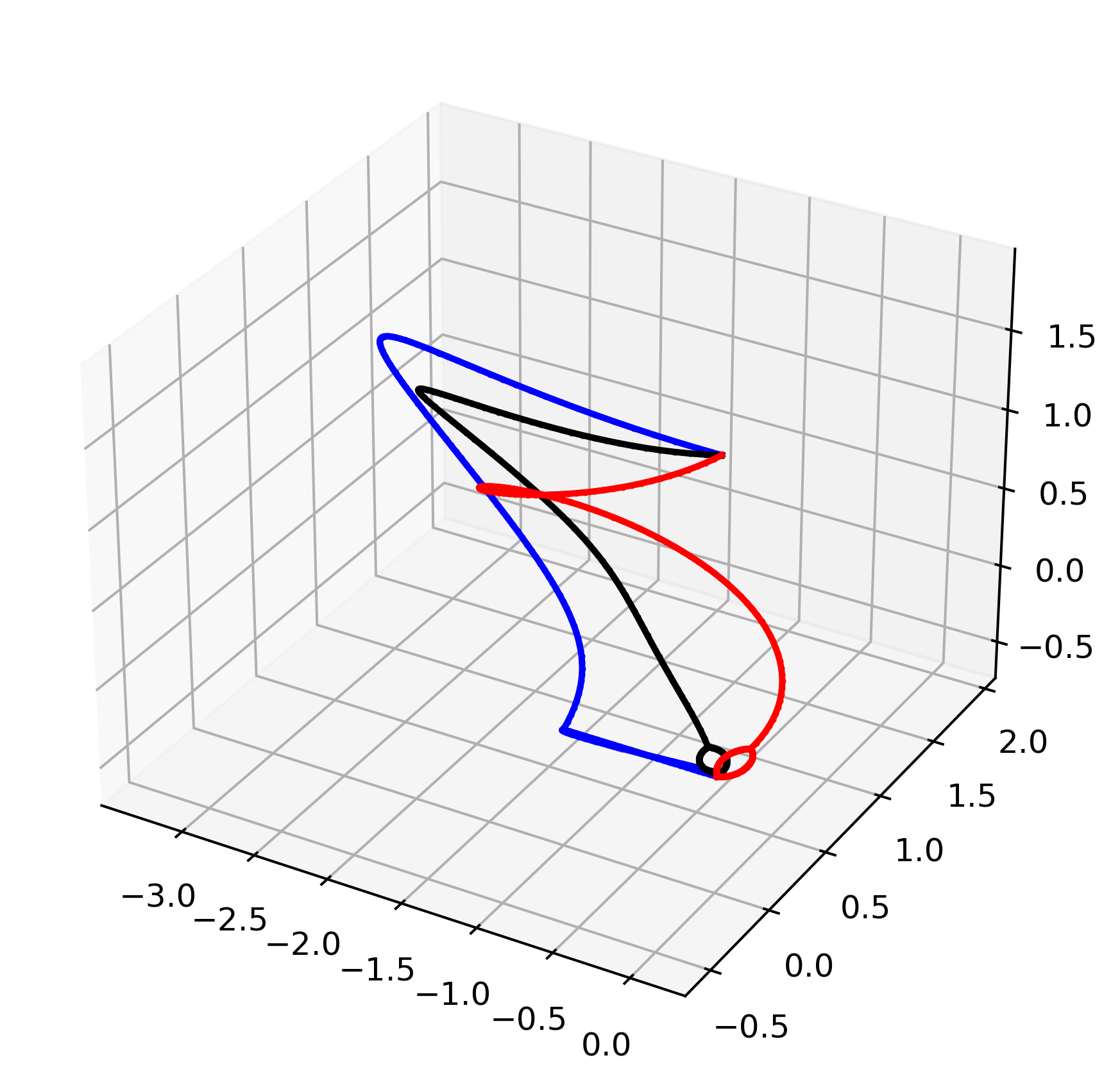}
\includegraphics[angle=-0,width=0.3\textwidth]{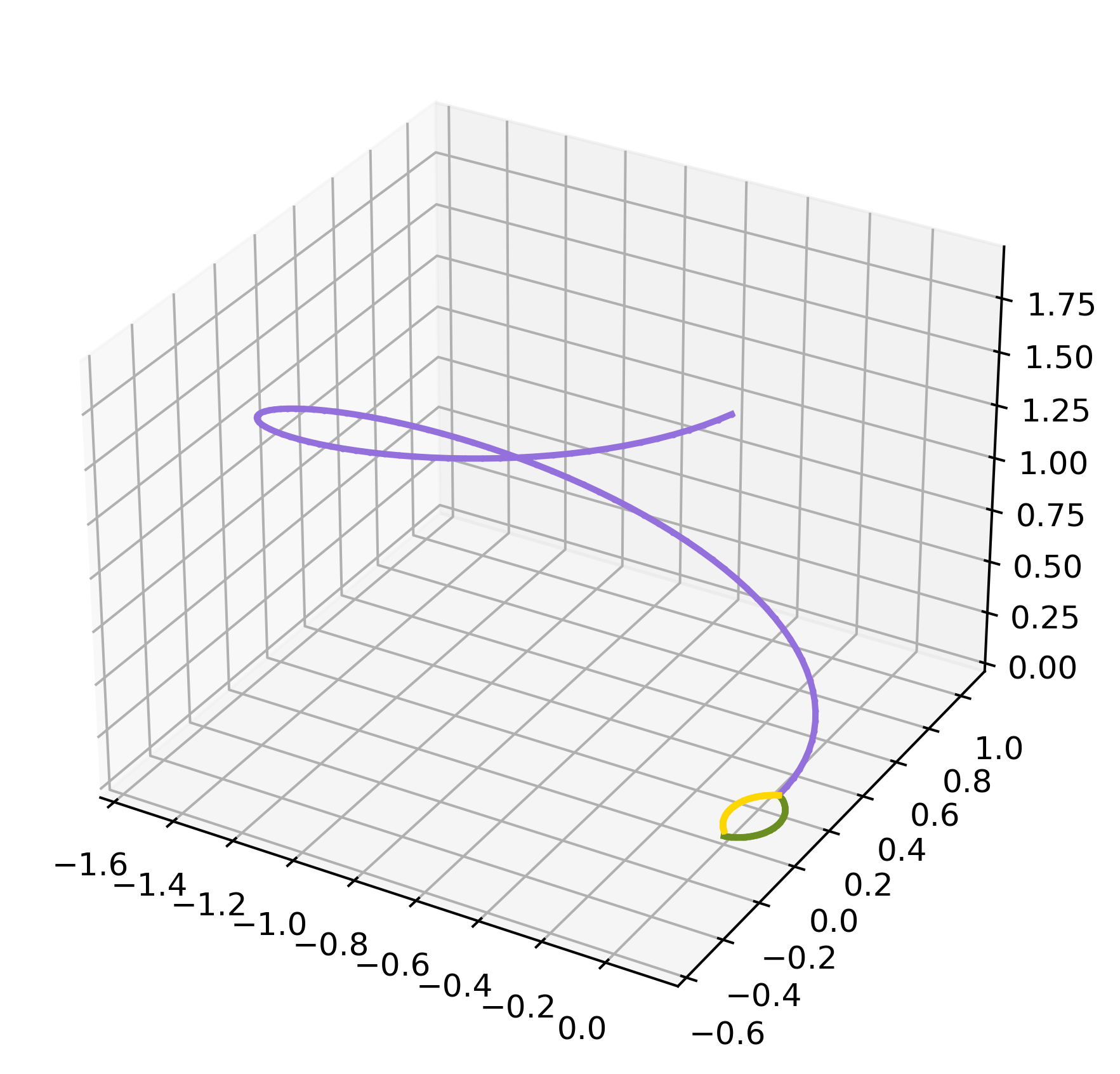}
\caption{The (projected) triod at time $t=0$, at times $t=0,0.5,5$
and at time $t=20$.
At the top right we show a plot of the discrete energy $L_E(c^m)$ over time.
}
\label{fig:threesamed}
\end{figure}

{\bf Experiment 9}:
An even bigger imbalance between the two enclosed areas is shown in
Figure~\ref{fig:threesamed}. This time we let $\Sigma=(-1,0,0)^{t}$ and
$\proj{P_\alpha}= (0,0)^{t}$, $\alpha=1,2,3$. Observe that the initial data is such that one
of the projected curves has a self-intersection.
Moreover, it holds that
$P_1=(0,0,0)^{t}$, $P_2=(0,0,1.95)^{t}$ and $P_3=(0,0,0.02)^{t}$.

\begin{figure}
\center
\includegraphics[angle=-0,width=0.25\textwidth]{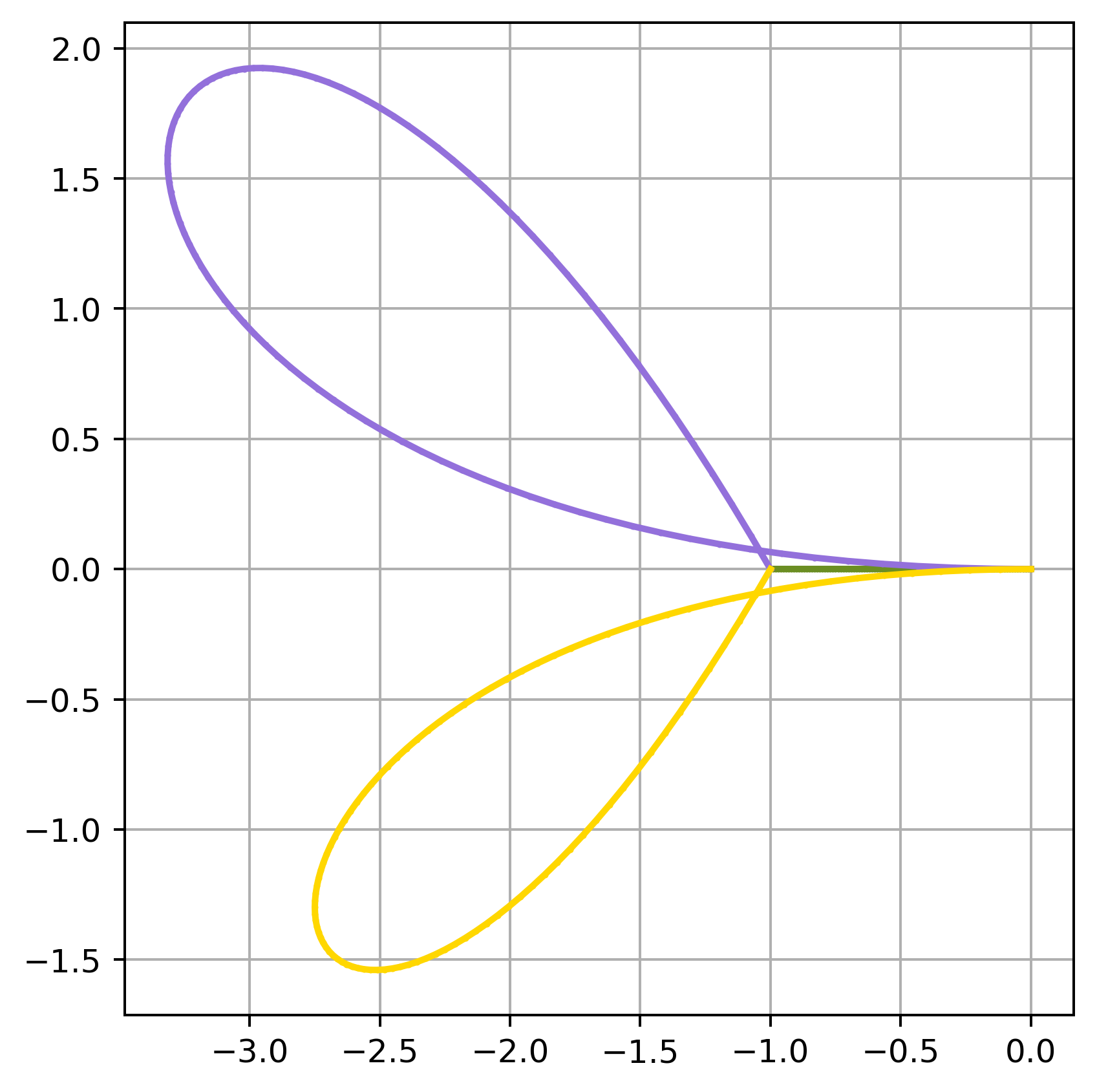}
\includegraphics[angle=-0,width=0.25\textwidth]{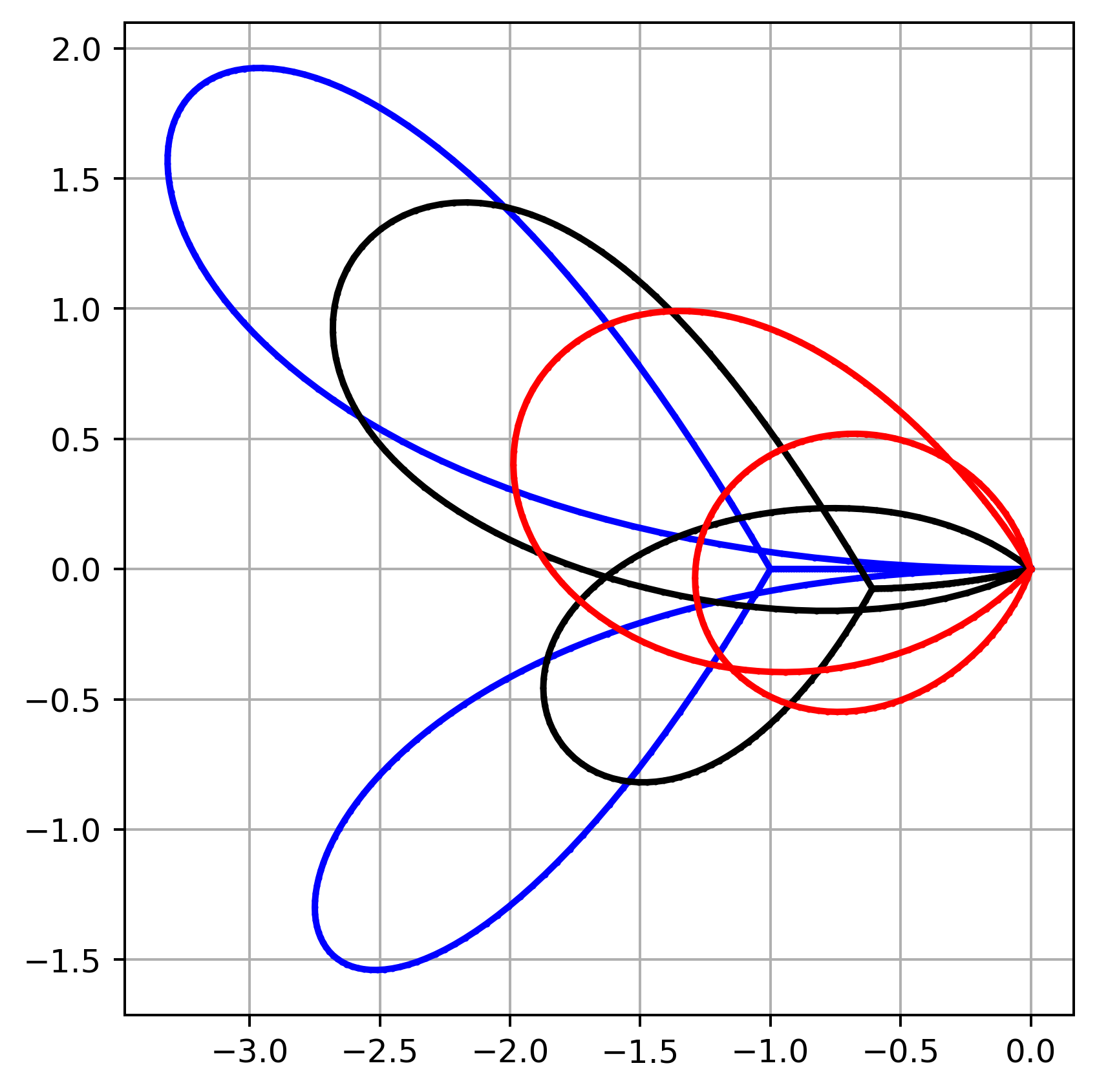}
\includegraphics[angle=-0,width=0.3\textwidth]{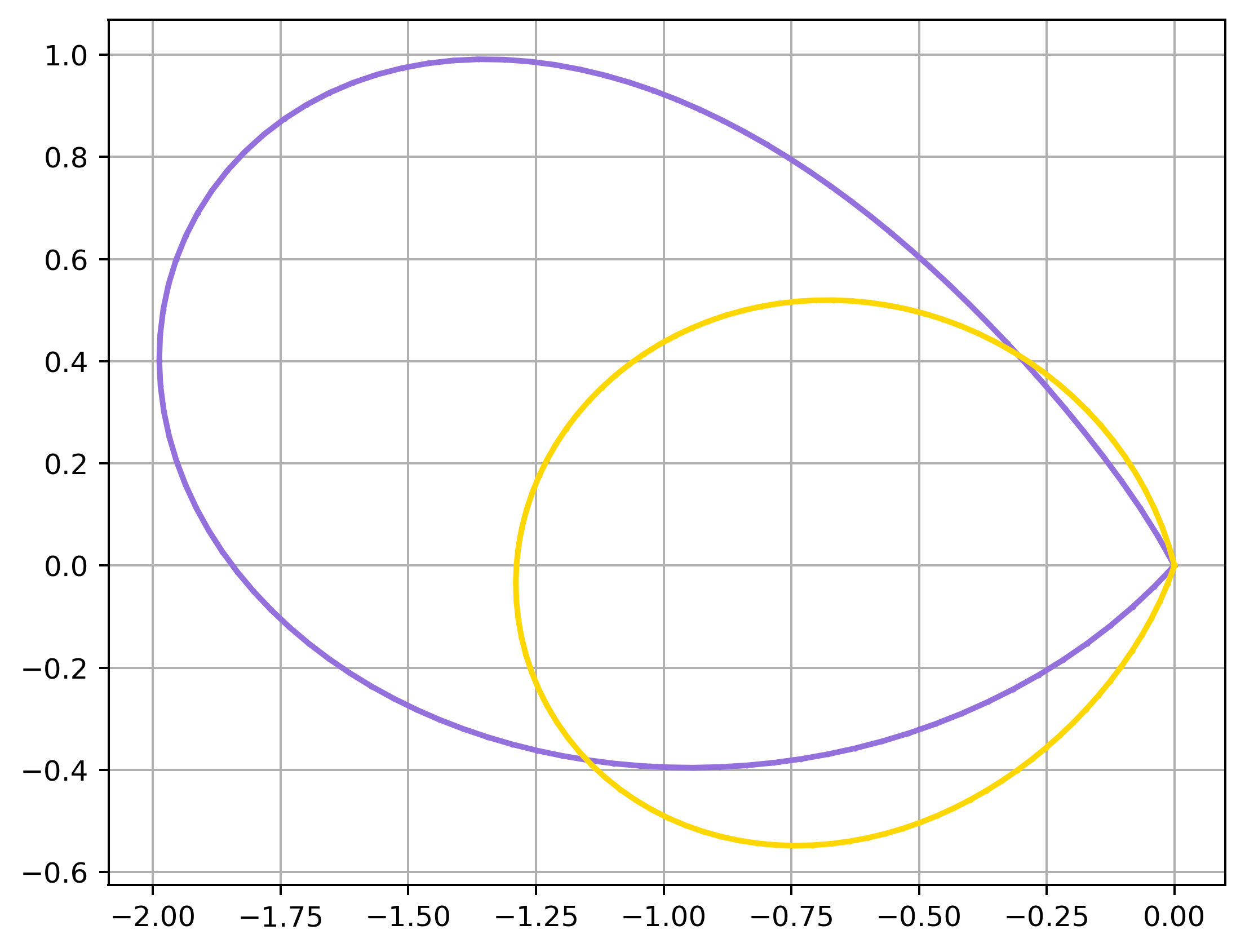}
\includegraphics[angle=-0,width=0.3\textwidth]{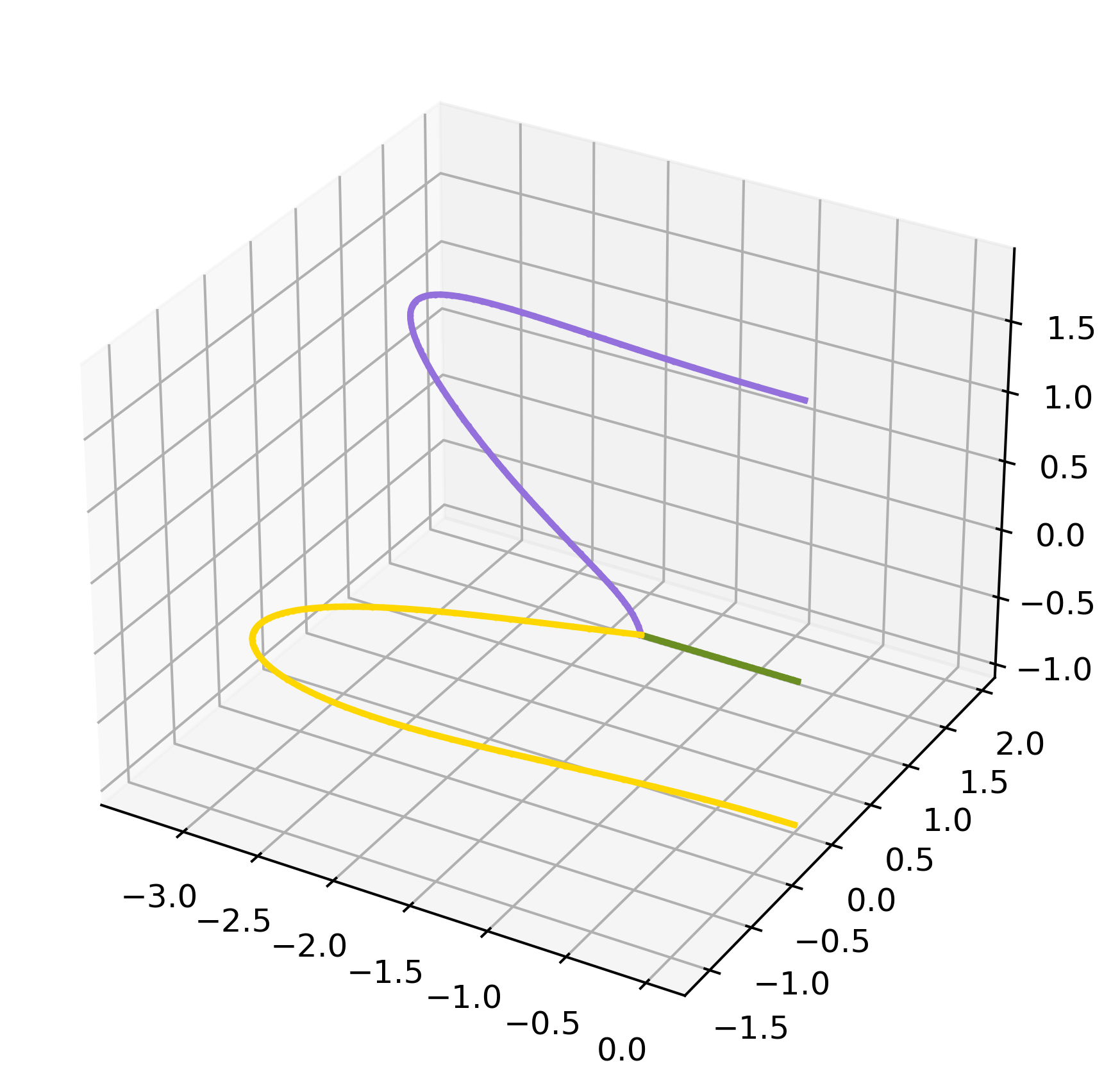}
\includegraphics[angle=-0,width=0.3\textwidth]{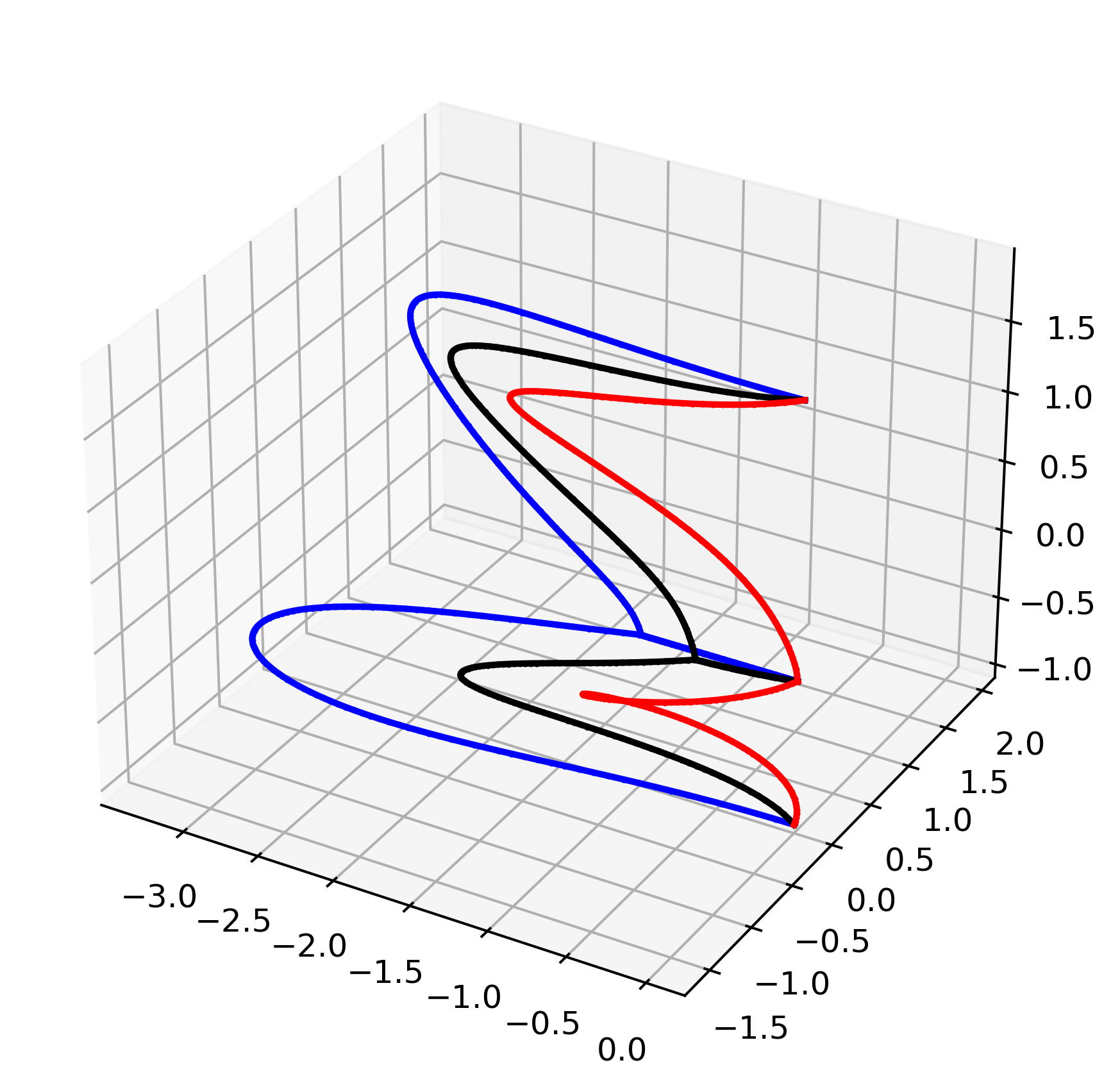}
\includegraphics[angle=-0,width=0.3\textwidth]{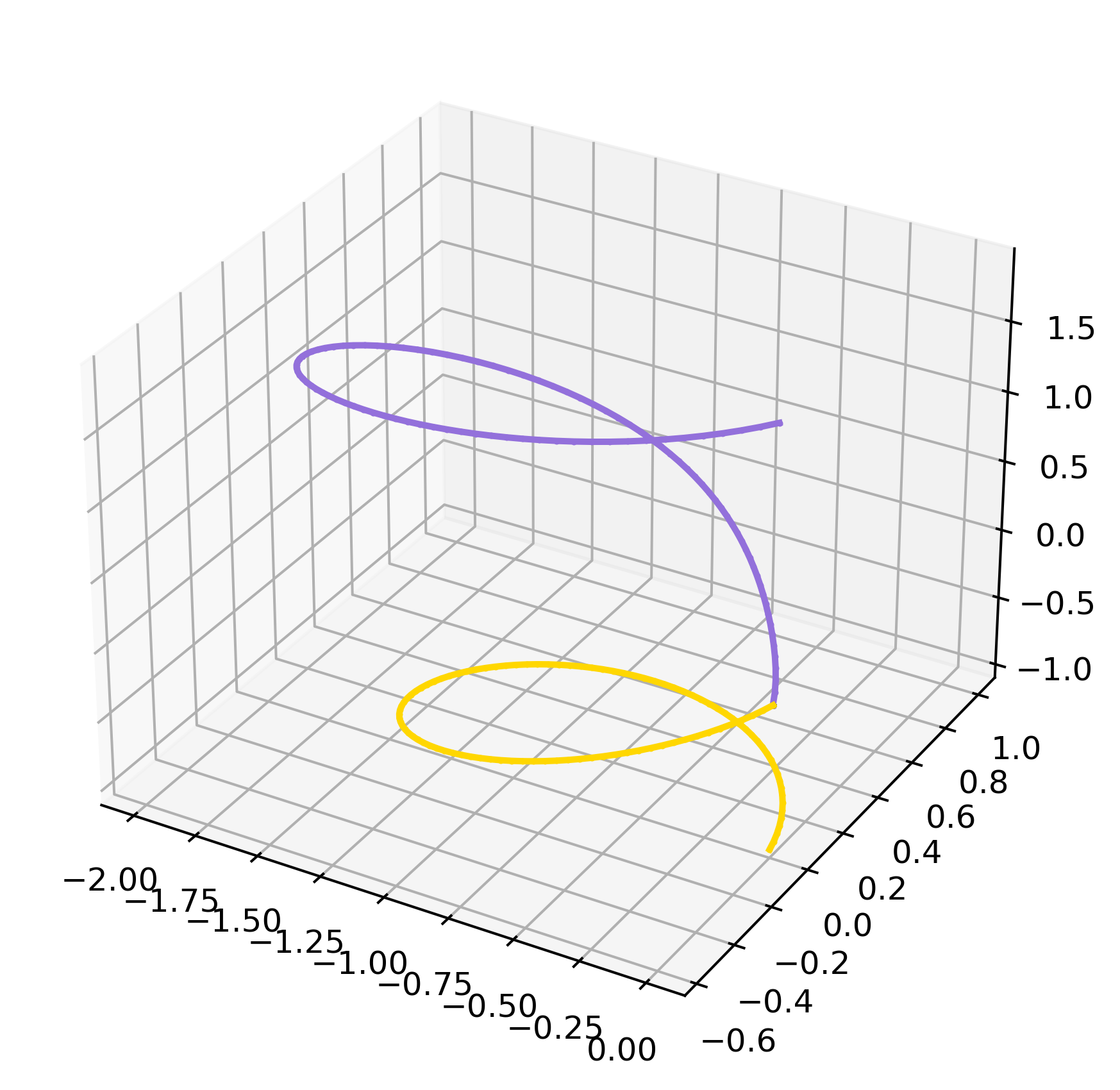}
\includegraphics[angle=-0,width=0.4\textwidth]{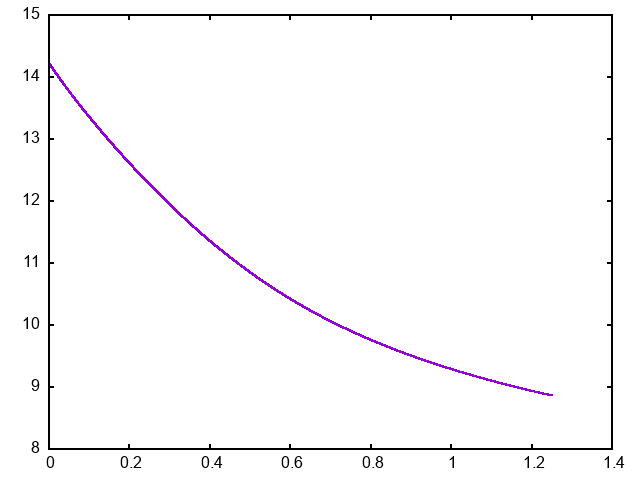}
\includegraphics[angle=-0,width=0.4\textwidth]{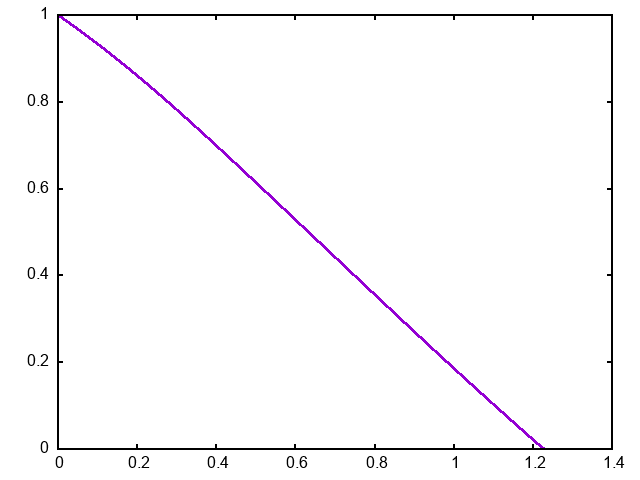}
\caption{The (projected) triod at time $t=0$, at times $t=0,0.5,1.2$
and at time $t=1.2$.
Below we show a plot of the discrete energies $L_E(c^m)$ and $L_E(c^m_1)$
over time.
}
\label{fig:threesamee}
\end{figure}

{\bf Experiment 10}:
By increasing the relative lengths of the second and third curves,
we manage to create an initial data from which a singularity seems to 
arise. That is, we let
$P_1=(0,0,0)^{t}$, $P_2=(0,0,1.95)^{t}$, $P_3=(0,0,-1.04)^{t}$ and $\Sigma=(-1,0,0)^{t}$.
In fact, for the evolution shown in Figure~\ref{fig:threesamee} we
observe that eventually the shortest curve vanishes, because the (projected)
triple junction has moved to the origin. Once again we observe that the planar
initial data features self-intersections. Moreover, at the time our numerical
breaks down due to the singularity, the remaining two curves do not represent
geodesics. This can be seen from the negative slope of the graph of the
discrete energy $L_E(c^m)$ at the final time.

\begin{figure}
\center
\mbox{
\includegraphics[angle=-0,width=0.2\textwidth]{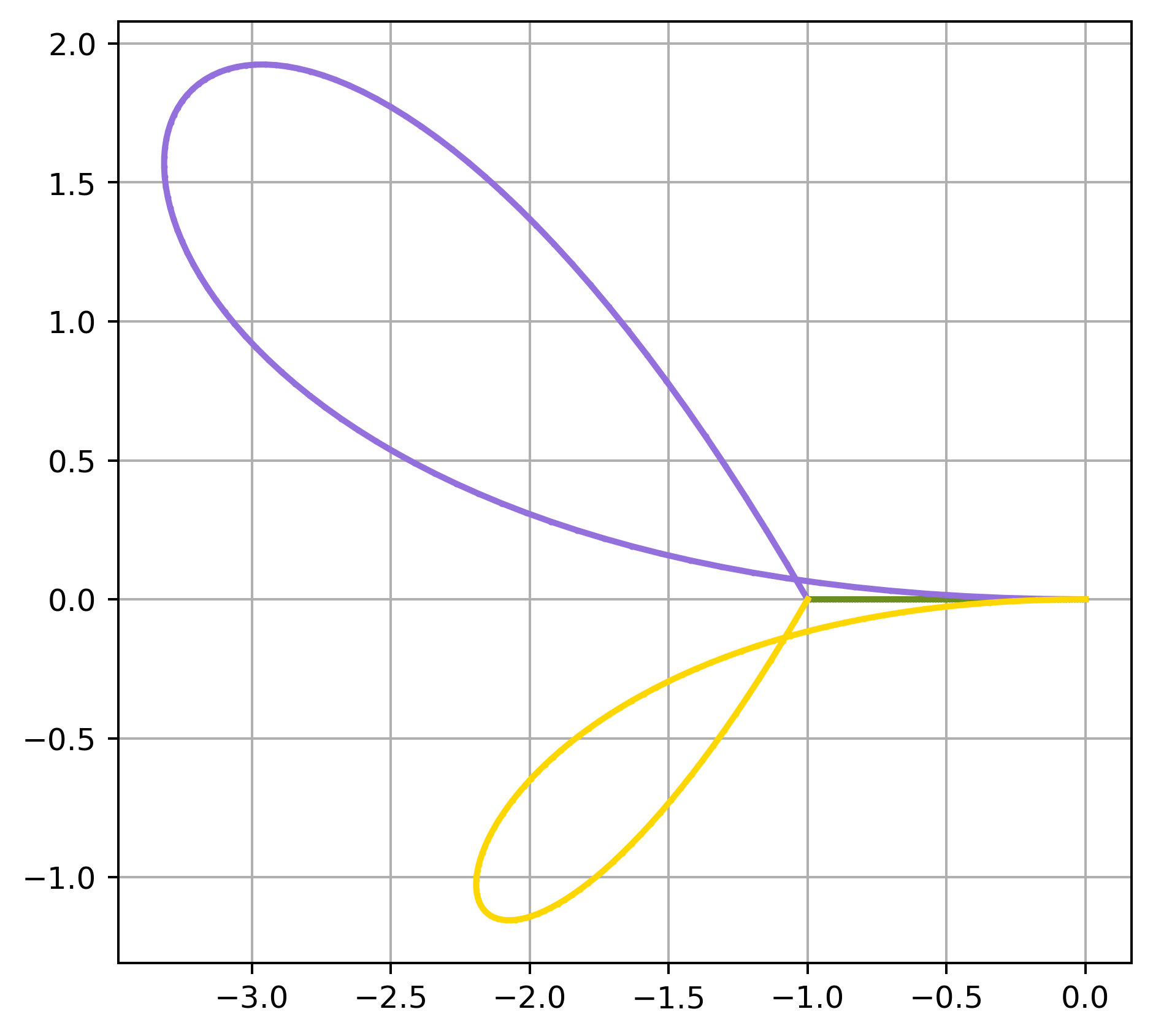}
\includegraphics[angle=-0,width=0.22\textwidth]{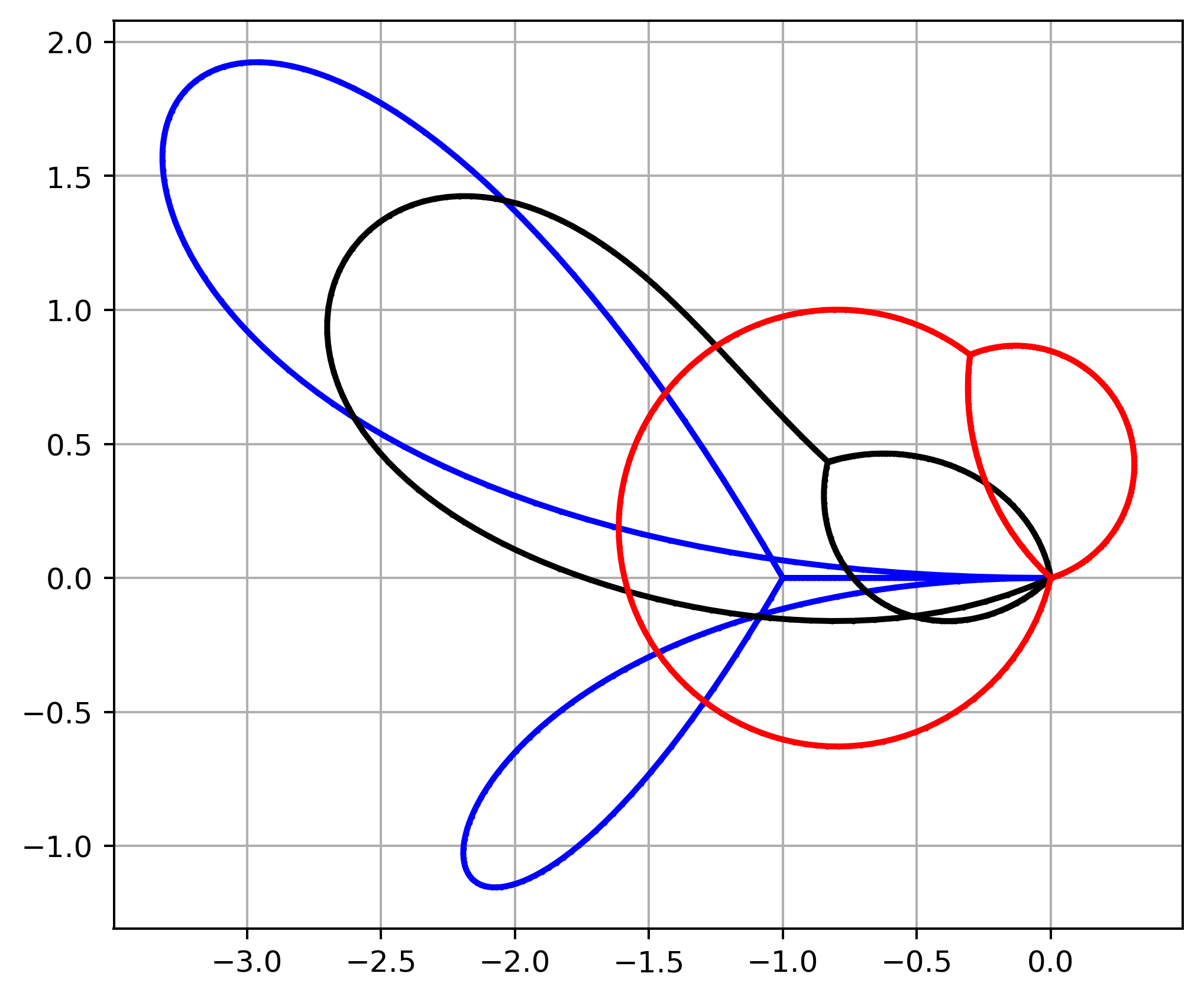}
\includegraphics[angle=-0,width=0.25\textwidth]{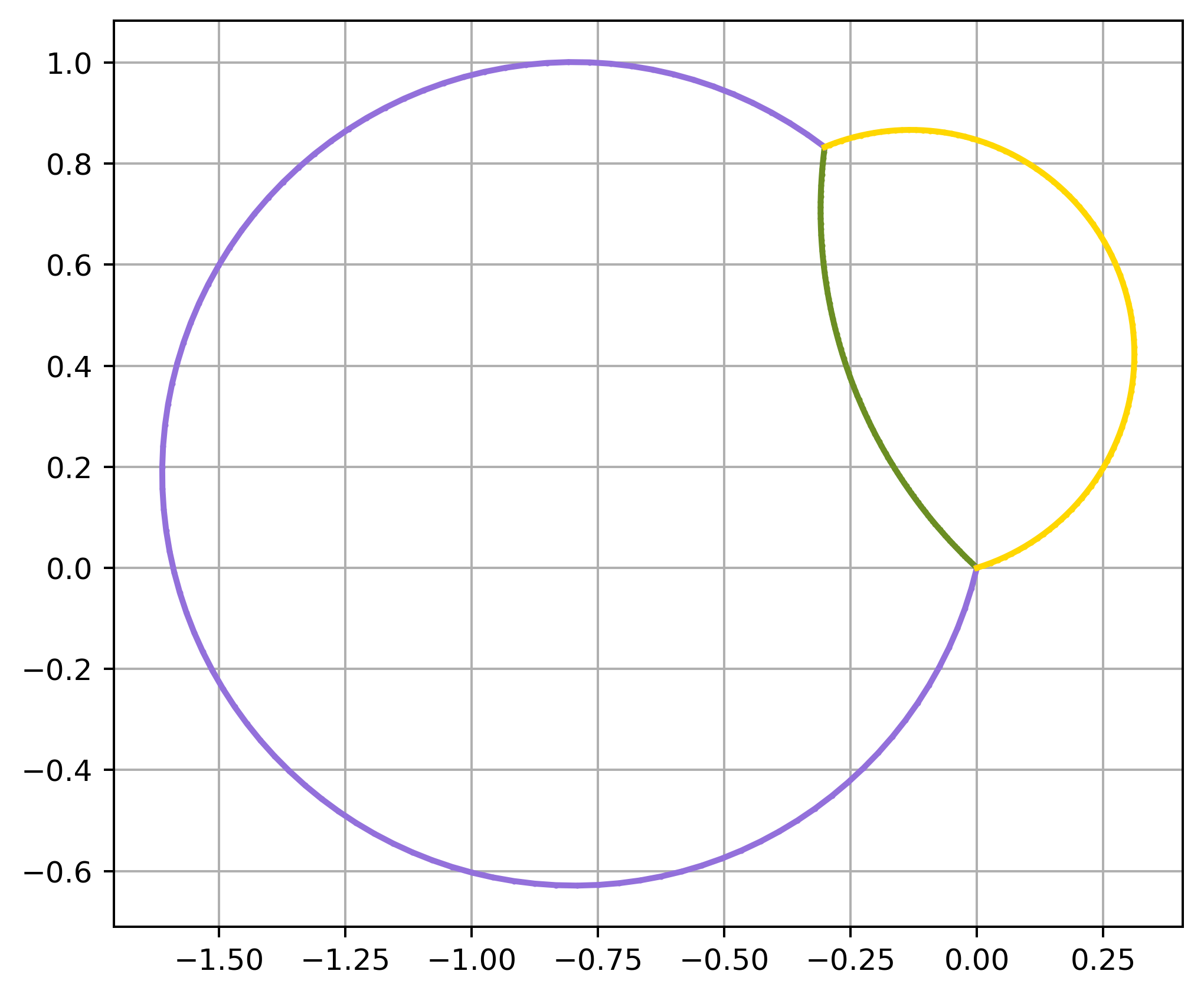}
\includegraphics[angle=-0,width=0.3\textwidth]{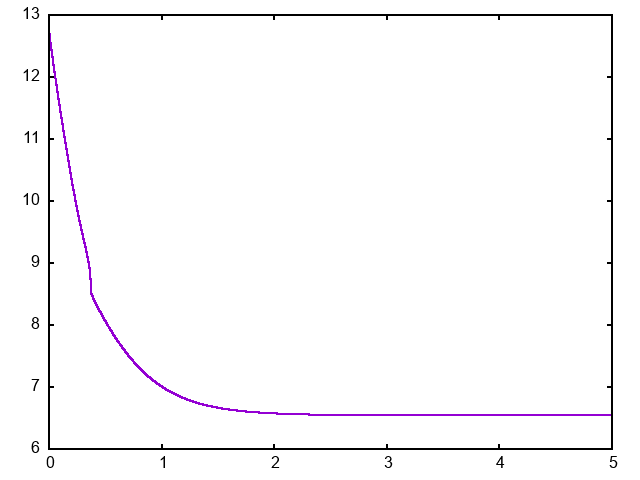}
}
\includegraphics[angle=-0,width=0.3\textwidth]{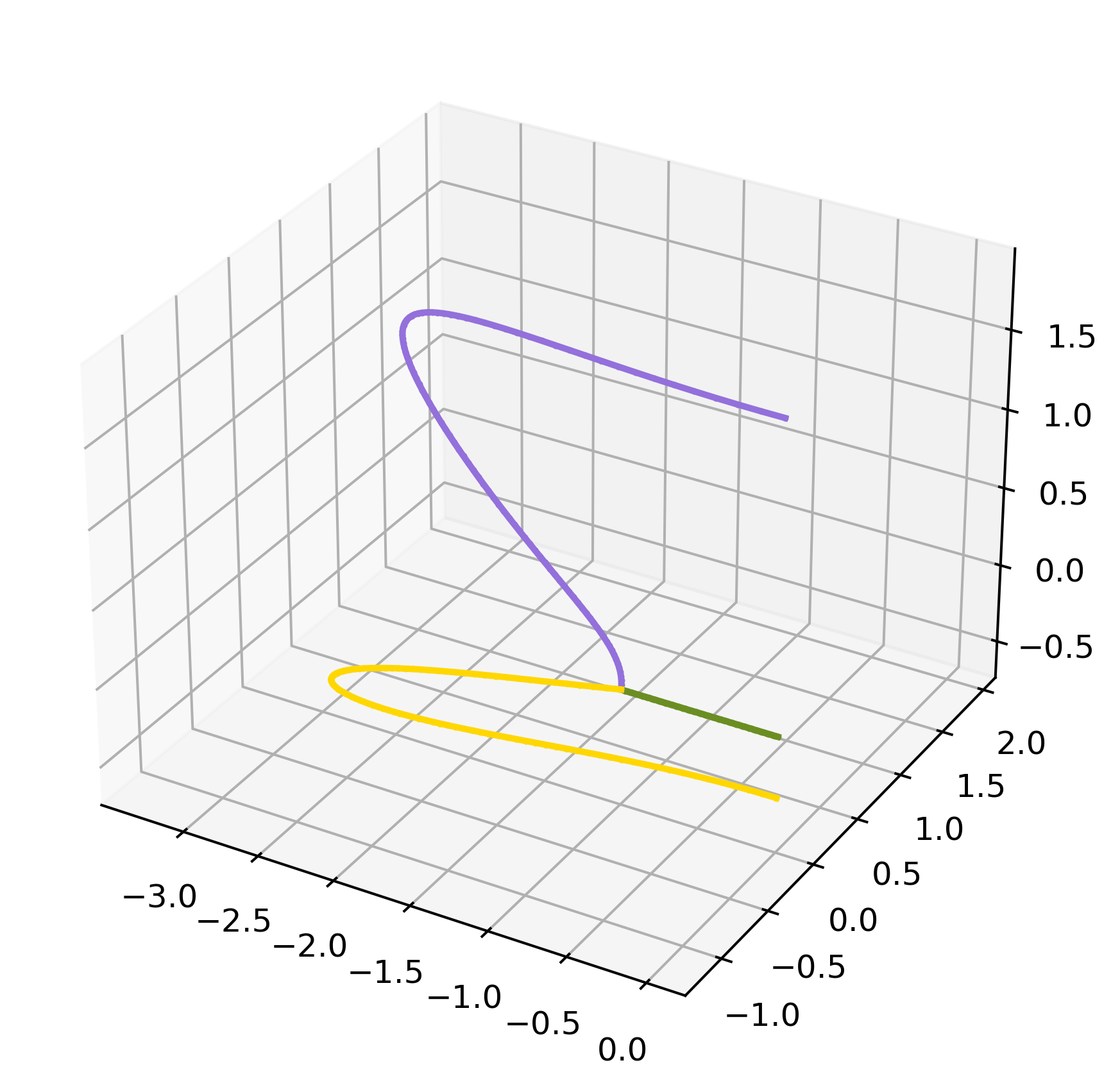}
\includegraphics[angle=-0,width=0.3\textwidth]{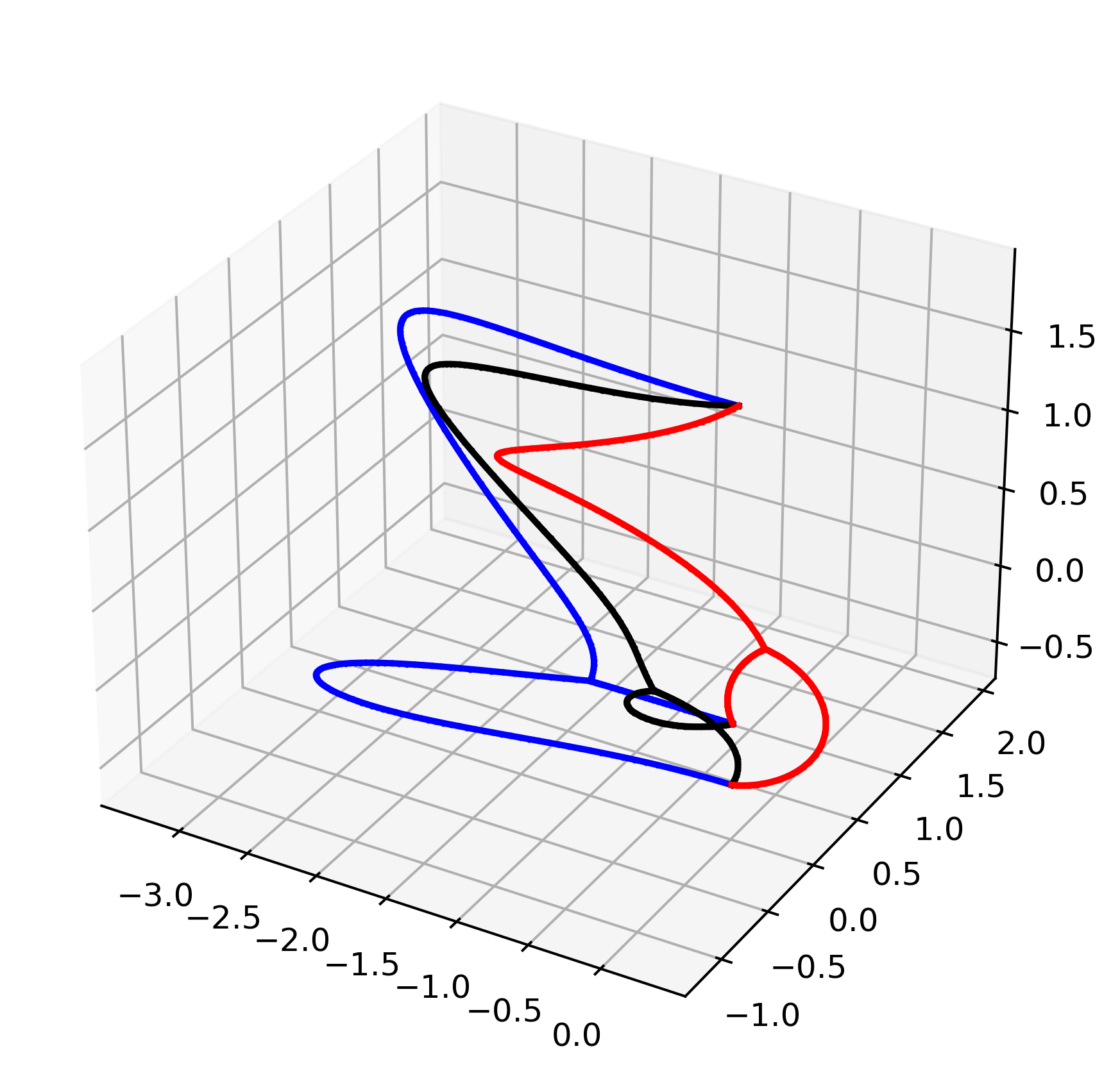}
\includegraphics[angle=-0,width=0.3\textwidth]{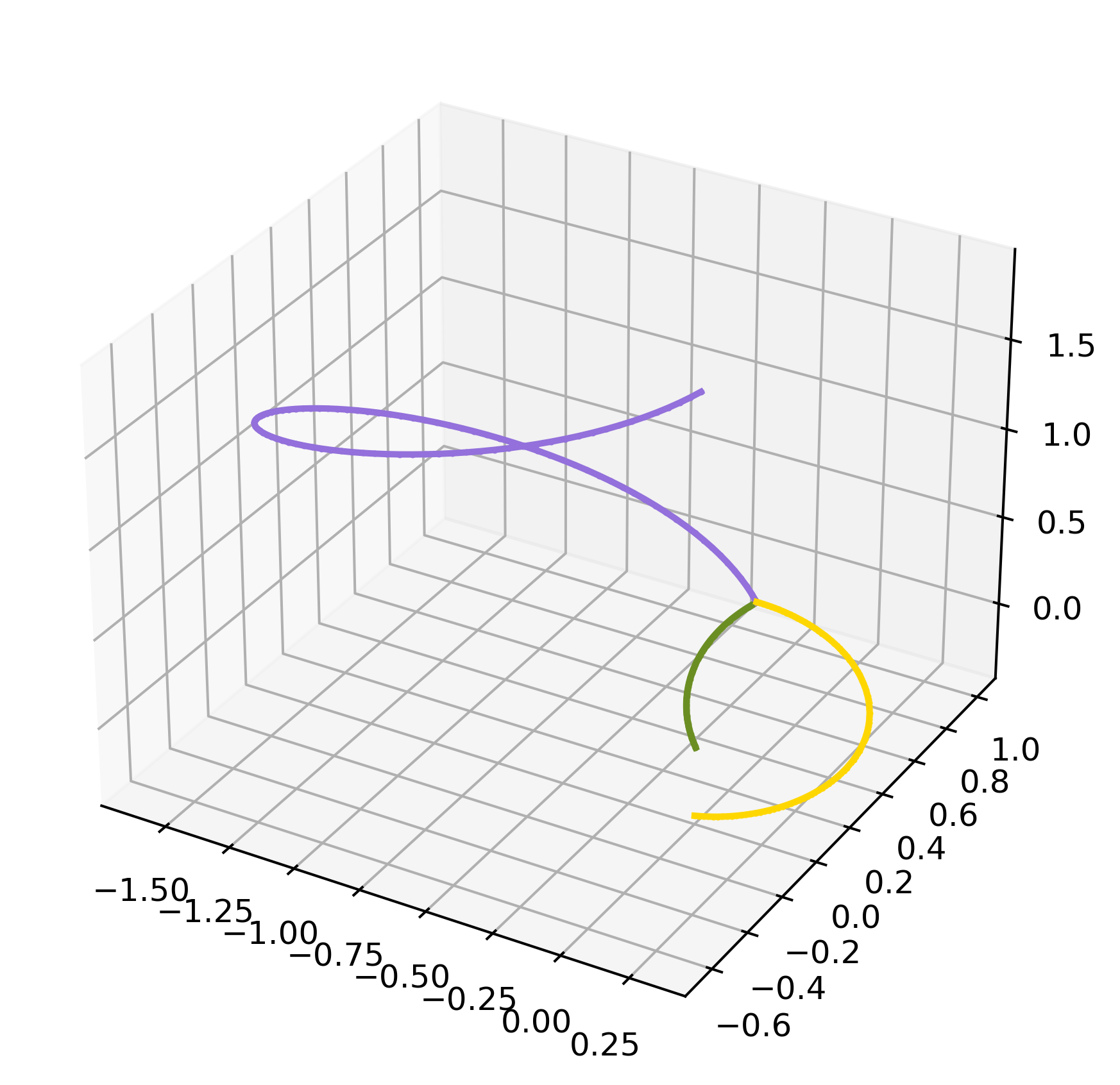}
\includegraphics[angle=-0,width=0.3\textwidth]{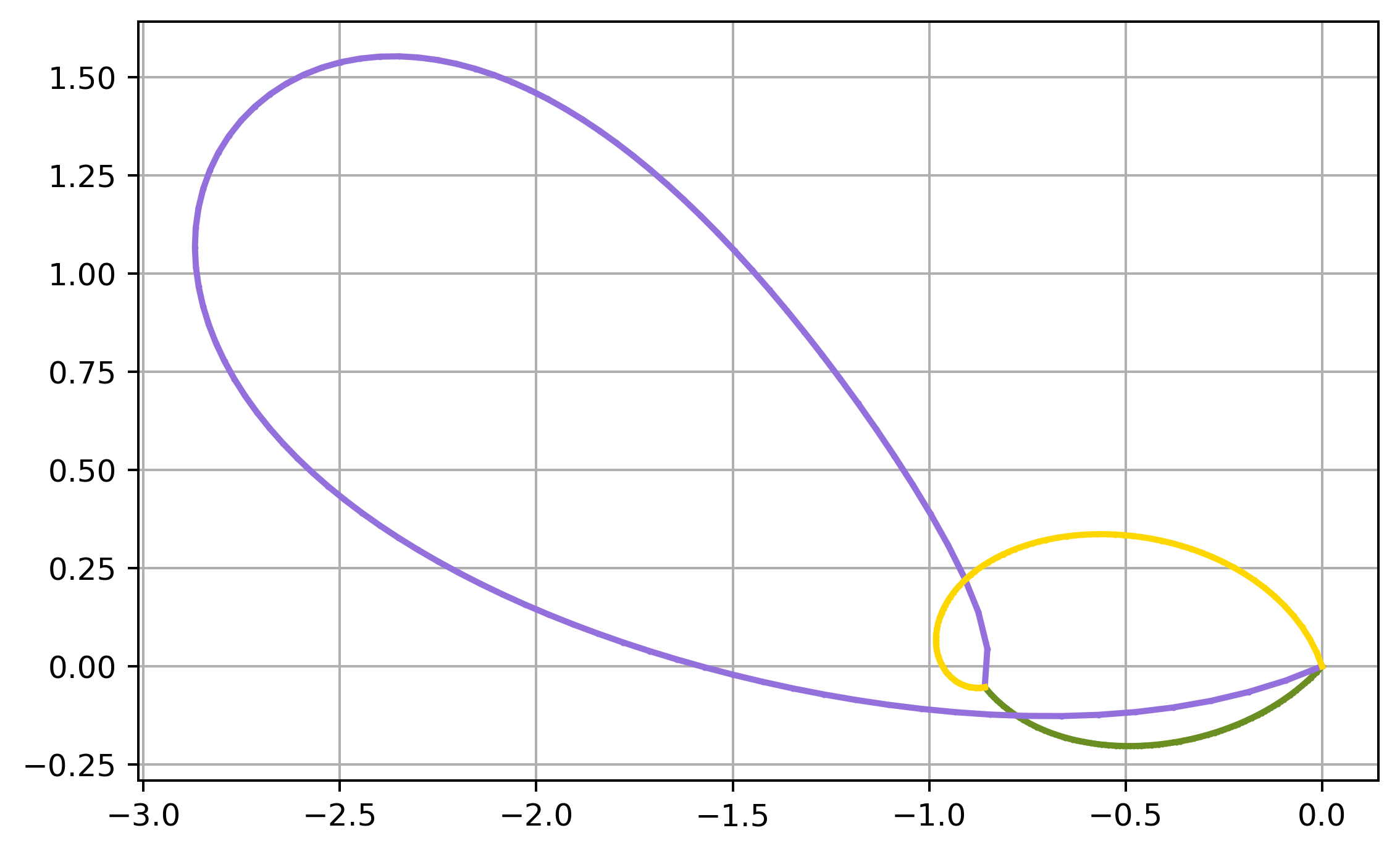}
\includegraphics[angle=-0,width=0.3\textwidth]{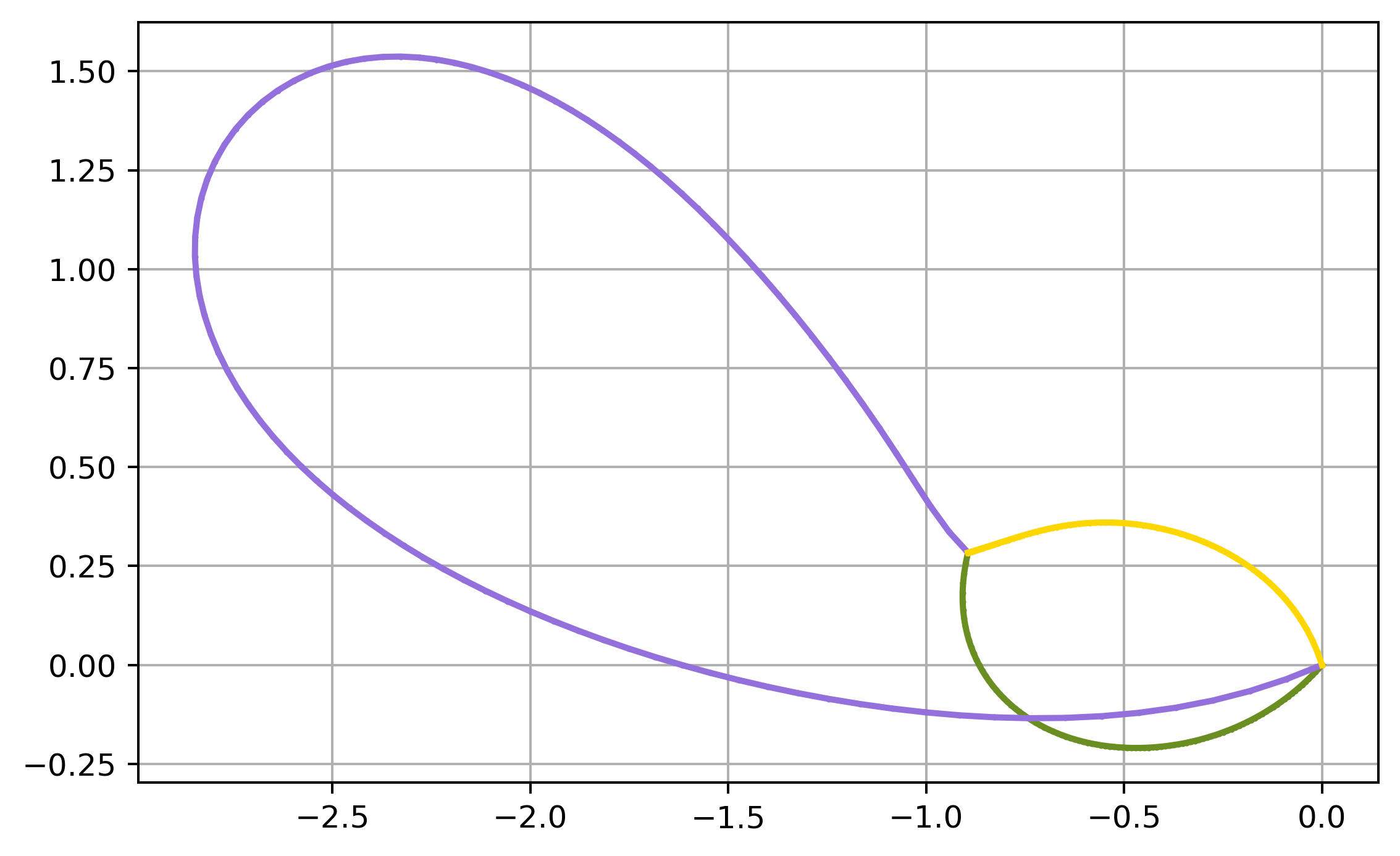}

\caption{The (projected) triod at time $t=0$, at times $t=0,0.5,5$
and at time $t=5$.
At the top right we show a plot of the discrete energy $L_E(c^m)$ over time.
In the bottom row we show the projected triod at time $t=0.36$ and at time
$t=0.38$.
}
\label{fig:threesamef}
\end{figure}

{\bf Experiment 11}:
Another simulation with a planar setup very similar to the previous experiment
is shown in Figure~\ref{fig:threesamef}. The positions of $\Sigma = (-1,0,0)^{t}$ 
and $\proj{P_\alpha} = (0,0)^{t}$, $\alpha=1,2,3$, are as before. But this time the area
enclosed by the curves $c_1$ and $c_2$ is significantly larger than the area
enclosed by $c_1$ and $c_3$. In fact, we now have that
$P_1=(0,0,0)^{t}$, $P_2=(0,0,1.95)^{t}$ and $P_3=(0,0,-0.39)^{t}$.
As a consequence, the observed evolution is very
different to the previous experiment. In particular, we do not observe any
singularity this time. In fact, the planar evolution appears to approach a non-equal
area standard double bubble. However, we notice a very rapid decrease in energy
at around time $t=0.37$, followed by a more moderate decrease in energy. To
visualize what is causing this behaviour, we also show the projected triod
shortly before and after this time at the bottom of
Figure~\ref{fig:threesamef}. What we see is that the position of the triple
junction $\proj\Sigma$ relative to the enclosed regions changes, which we
conjecture is behind the change in the speed of the energy decrease.

\begin{figure}
\center
\includegraphics[angle=-0,width=0.3\textwidth]{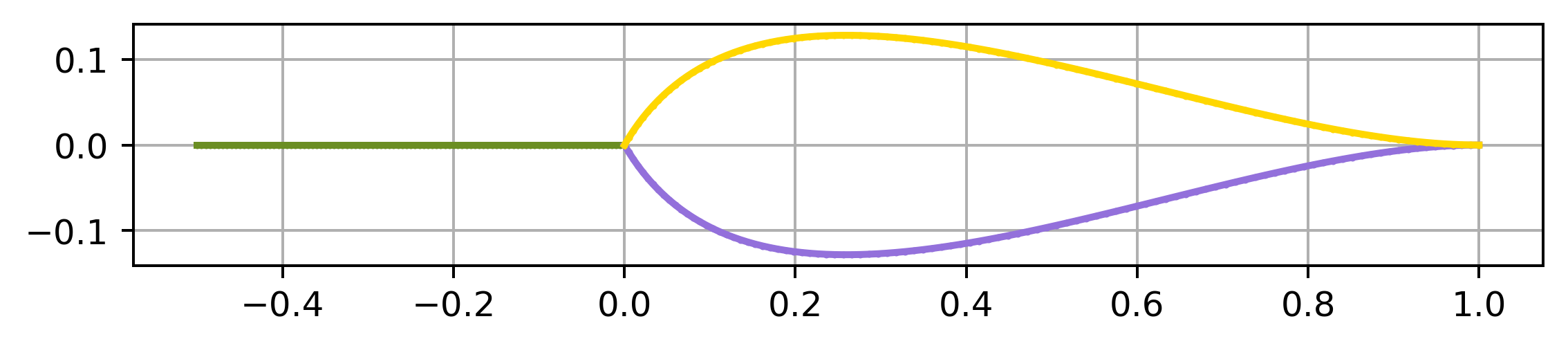}
\includegraphics[angle=-0,width=0.3\textwidth]{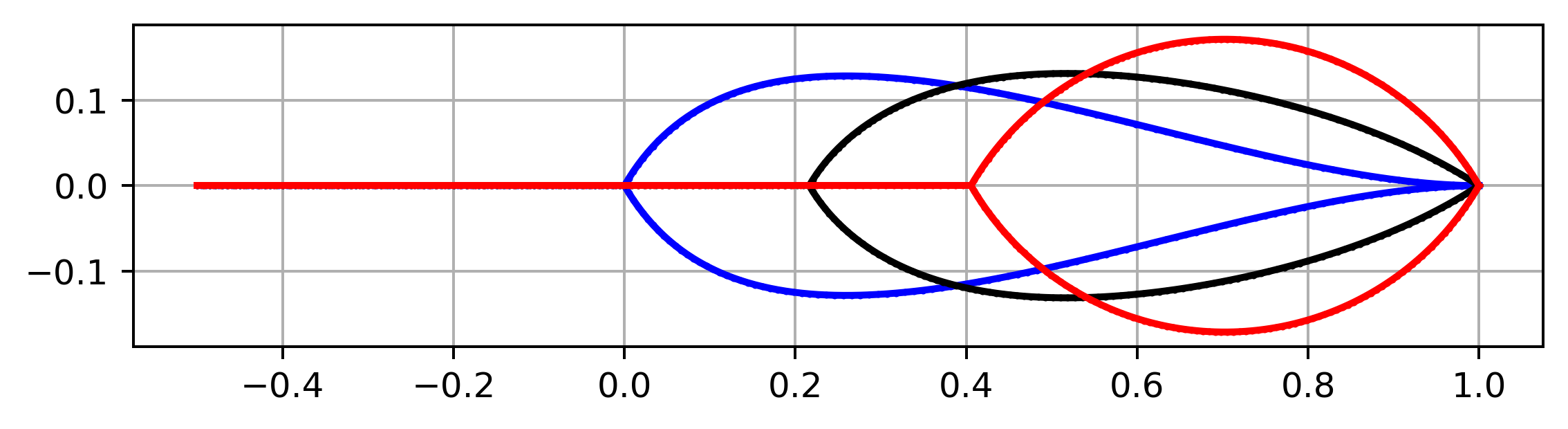}
\includegraphics[angle=-0,width=0.3\textwidth]{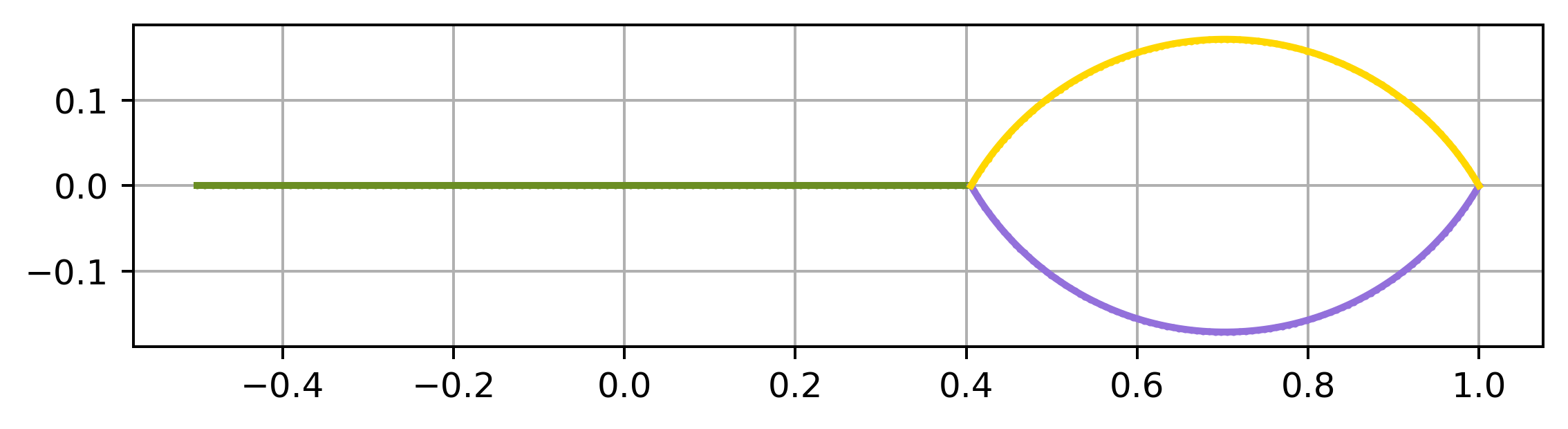}
\includegraphics[angle=-0,width=0.3\textwidth]{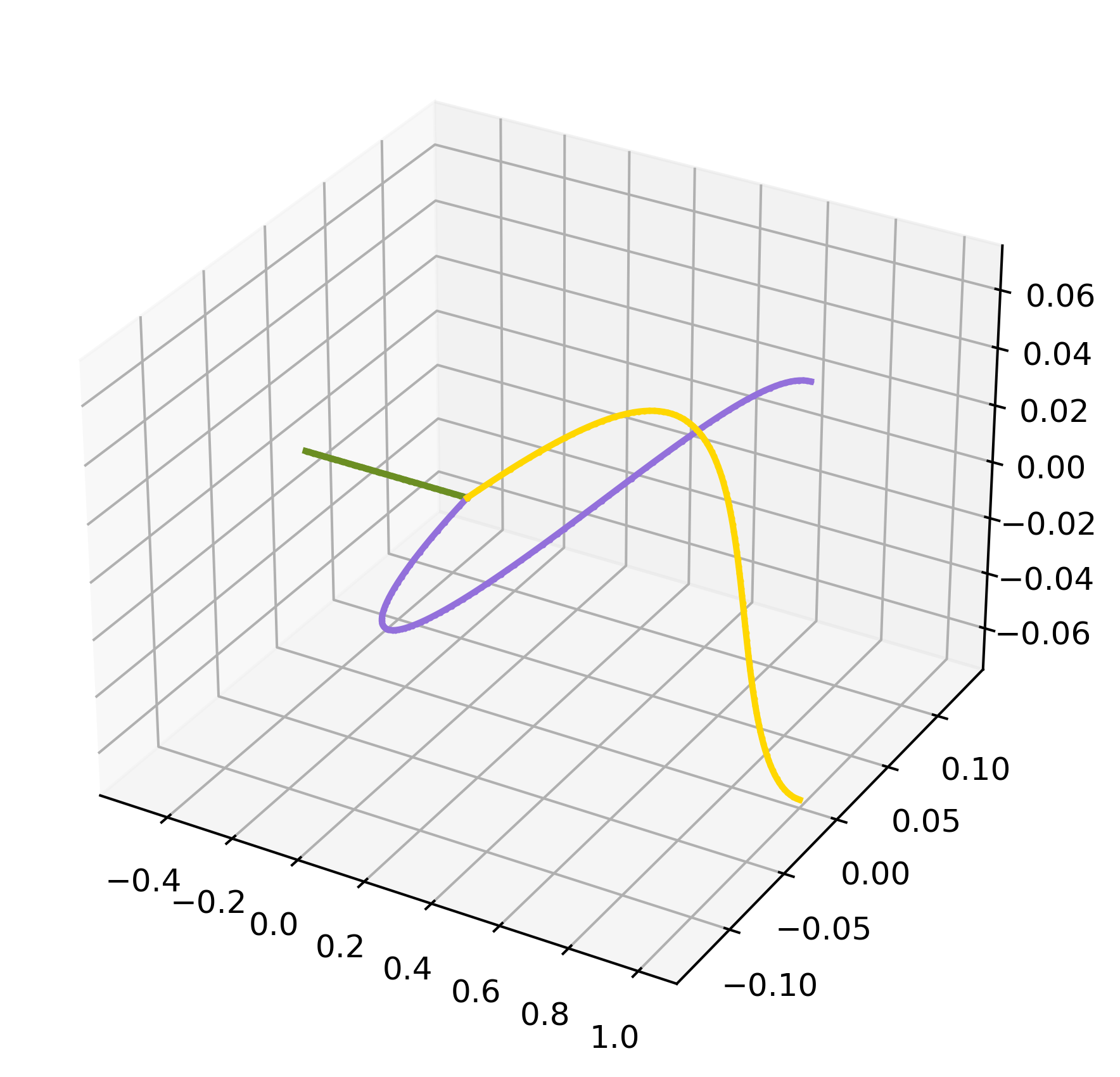}
\includegraphics[angle=-0,width=0.3\textwidth]{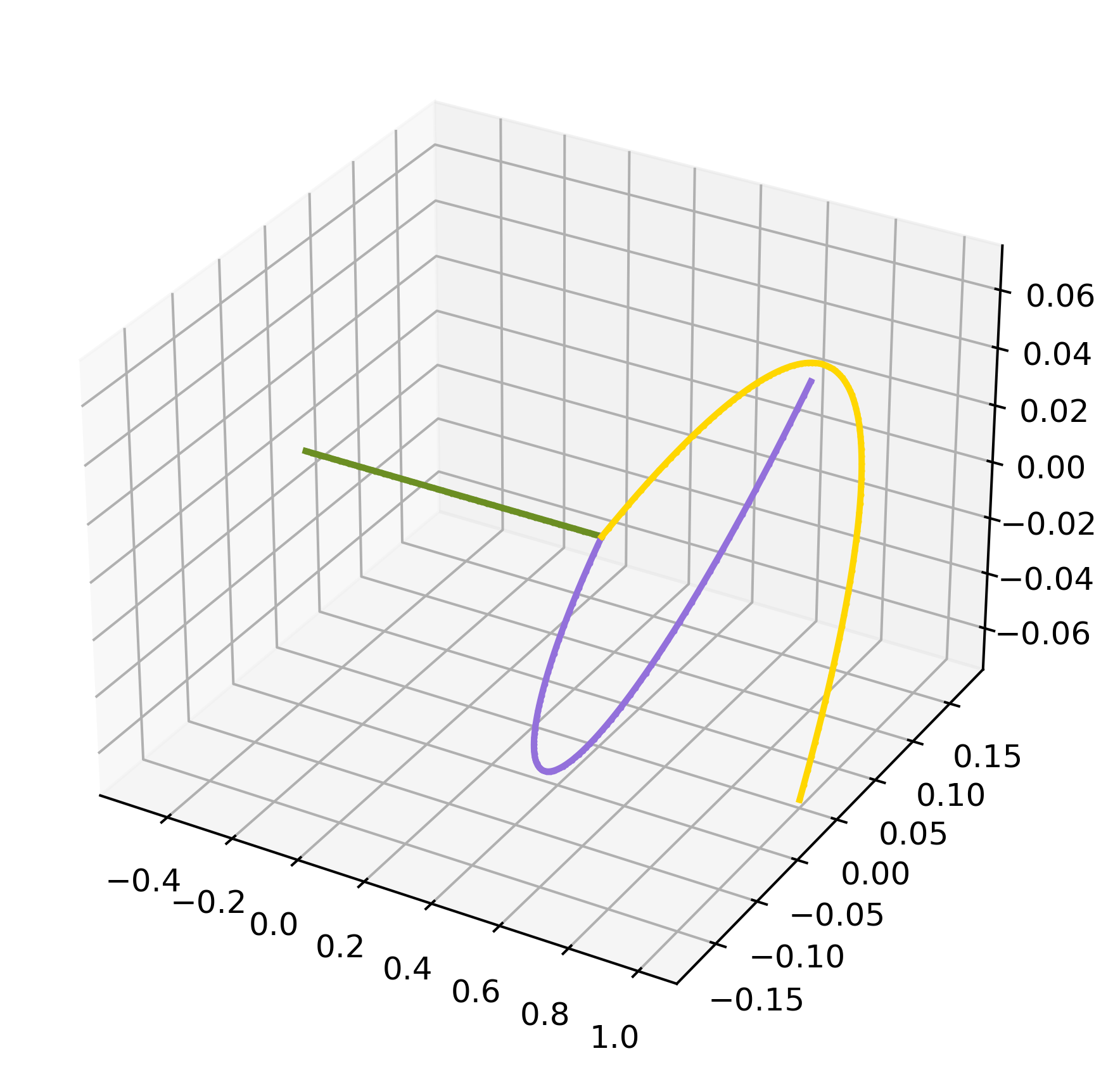}
\includegraphics[angle=-0,width=0.3\textwidth]{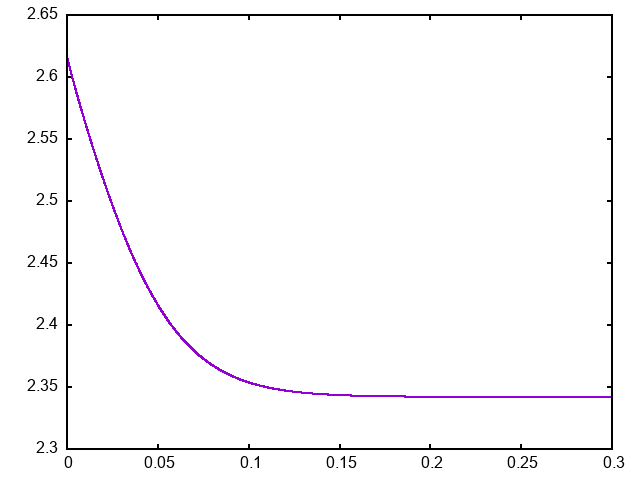}
\caption{The projected triod at time $t=0$, at times 
$t=0,0.05,0.3$ and at time $t=0.3$.
Below we show the triod in $\R^3$ at time $t=0$ and at time $t=0.3$,
as well as a plot of the discrete energy $L_E(c^m)$ over time.
}
\label{fig:twoandone}
\end{figure}

{\bf Experiment 12}:
For the initial setup we let $\Sigma=(0,0,0)^{t}$, as well as
$\proj{P_1} = (-1,0)^{t}$ and $\proj{P_\alpha} = (1,0)^{t}$, $\alpha=2,3$. The initial curves
imply that
$P_1=(-1,0,0)^{t}$, $P_2=(1,0,0.07)^{t}$ and $P_3=(1,0,-0.07)^{t}$.
In Figure~\ref{fig:twoandone} we observe that
the planar triod evolves towards a lens steady state, as it is known from
simulations for grain boundary motions or multiphase fluids. See e.g.\
\cite[Fig.\ 5]{ejam3d} and \cite[Fig.\ 13]{fluidfbptj}.

\subsection{Incompatible curved initial data}

In this subsection we once again choose initial data for which the 120 degree
angle condition at the planar triple junction is not satisfied. The initial
curves will in general be curved. 

\begin{figure}
\center
\mbox{
\includegraphics[angle=-0,width=0.23\textwidth]{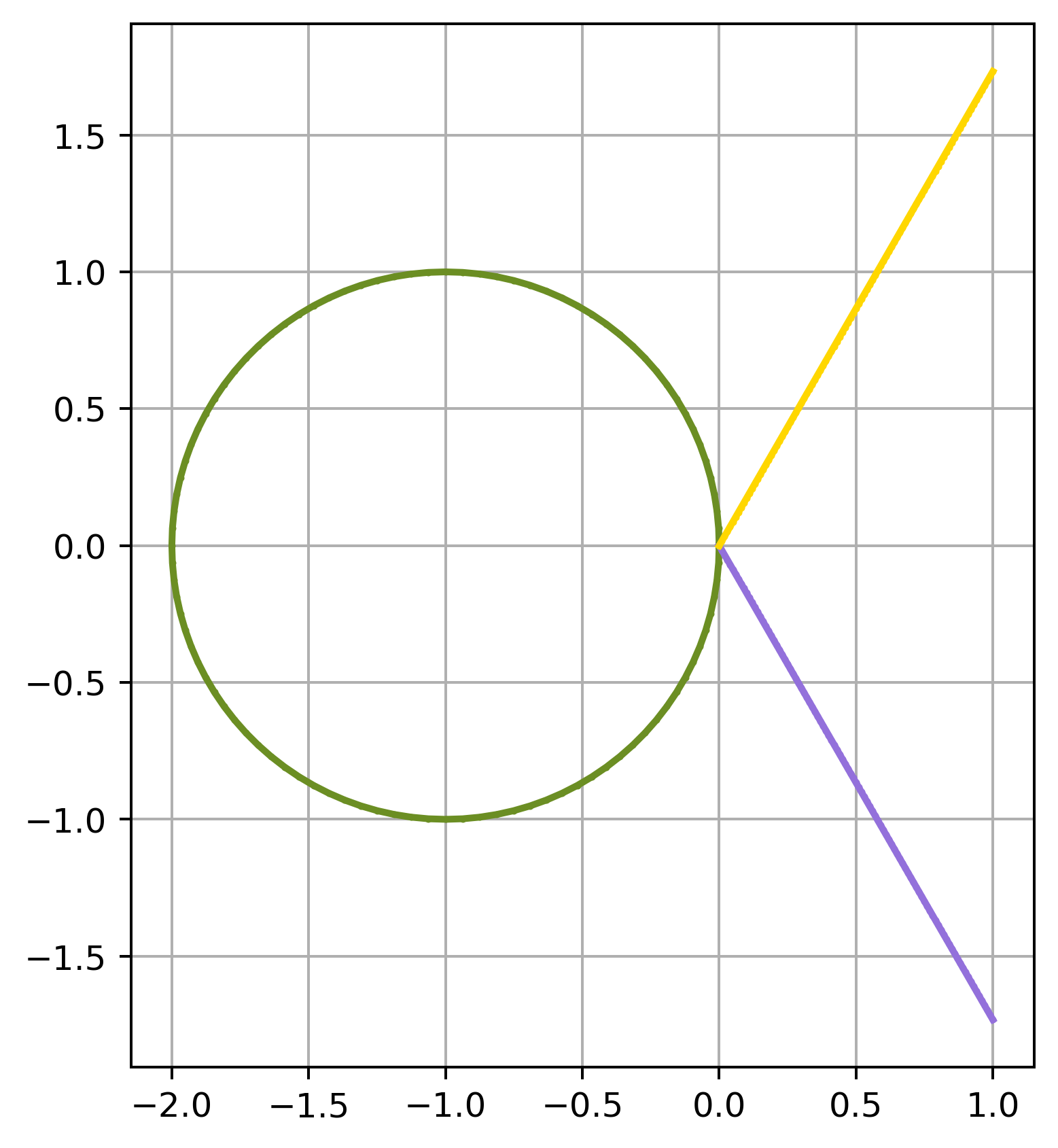}
\includegraphics[angle=-0,width=0.22\textwidth]{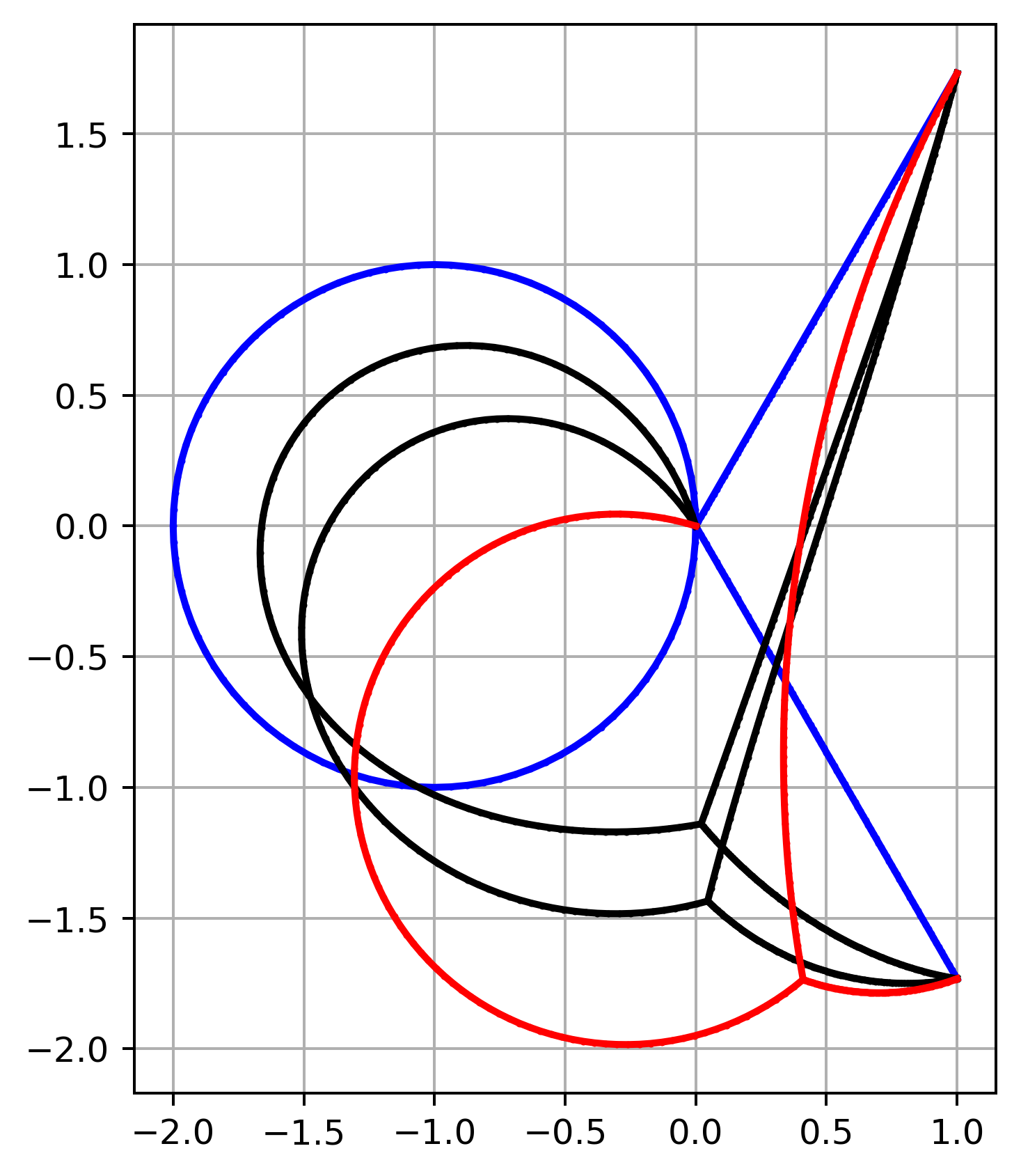}
\includegraphics[angle=-0,width=0.2\textwidth]{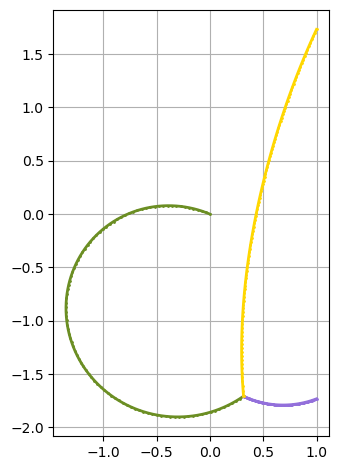}
\includegraphics[angle=-0,width=0.3\textwidth]{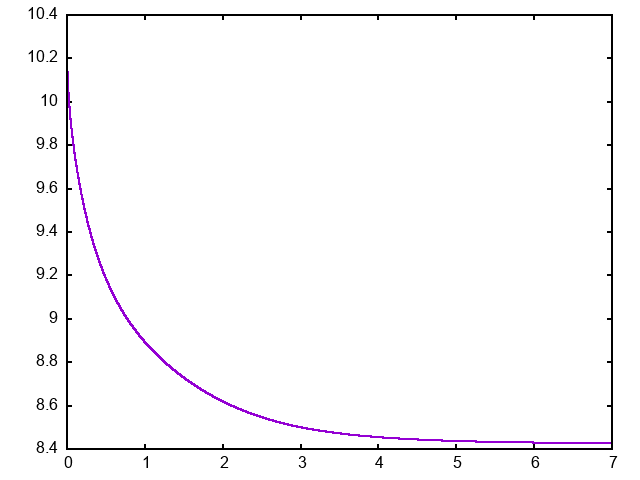}
}
\includegraphics[angle=-0,width=0.3\textwidth]{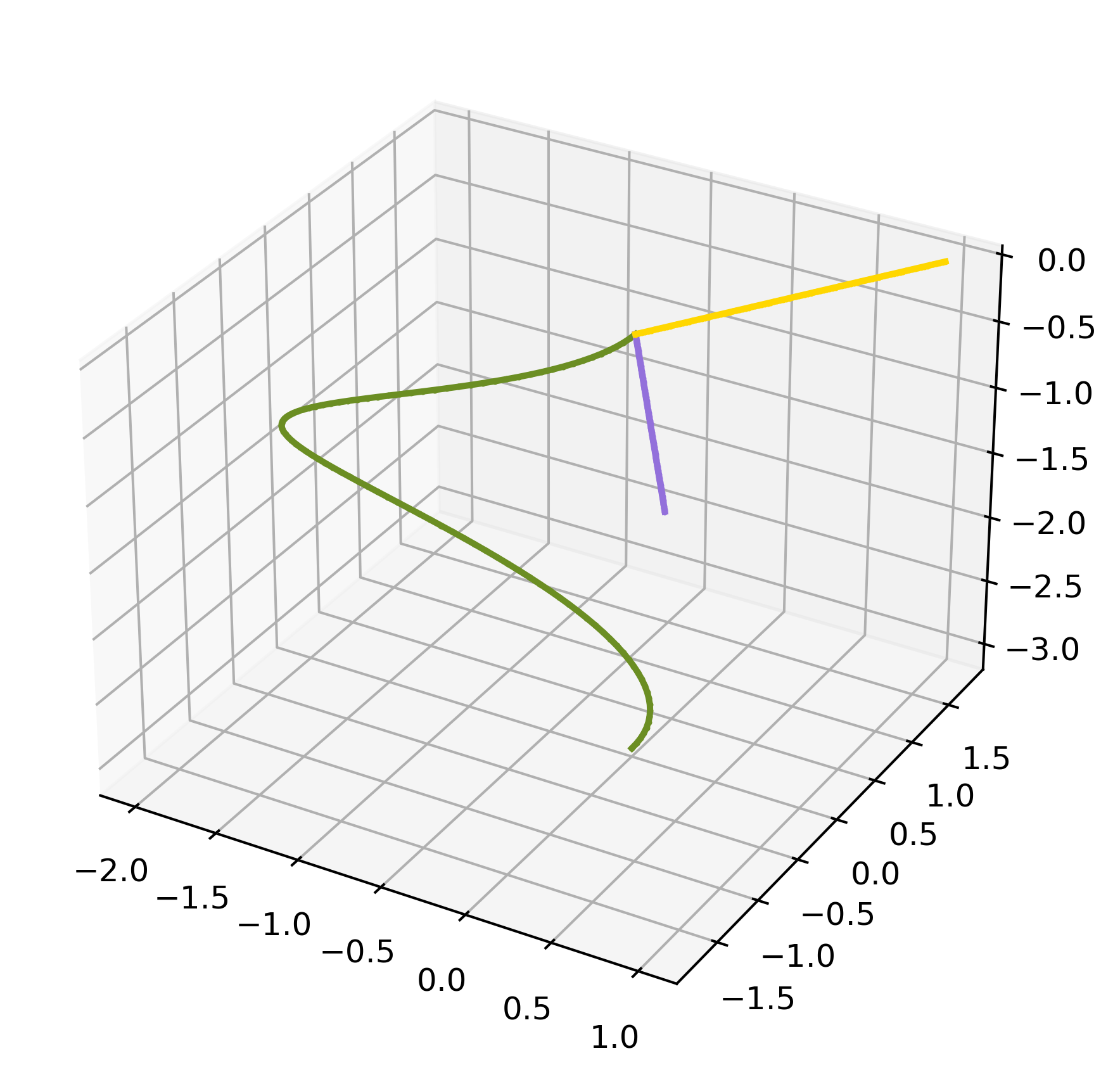}
\includegraphics[angle=-0,width=0.3\textwidth]{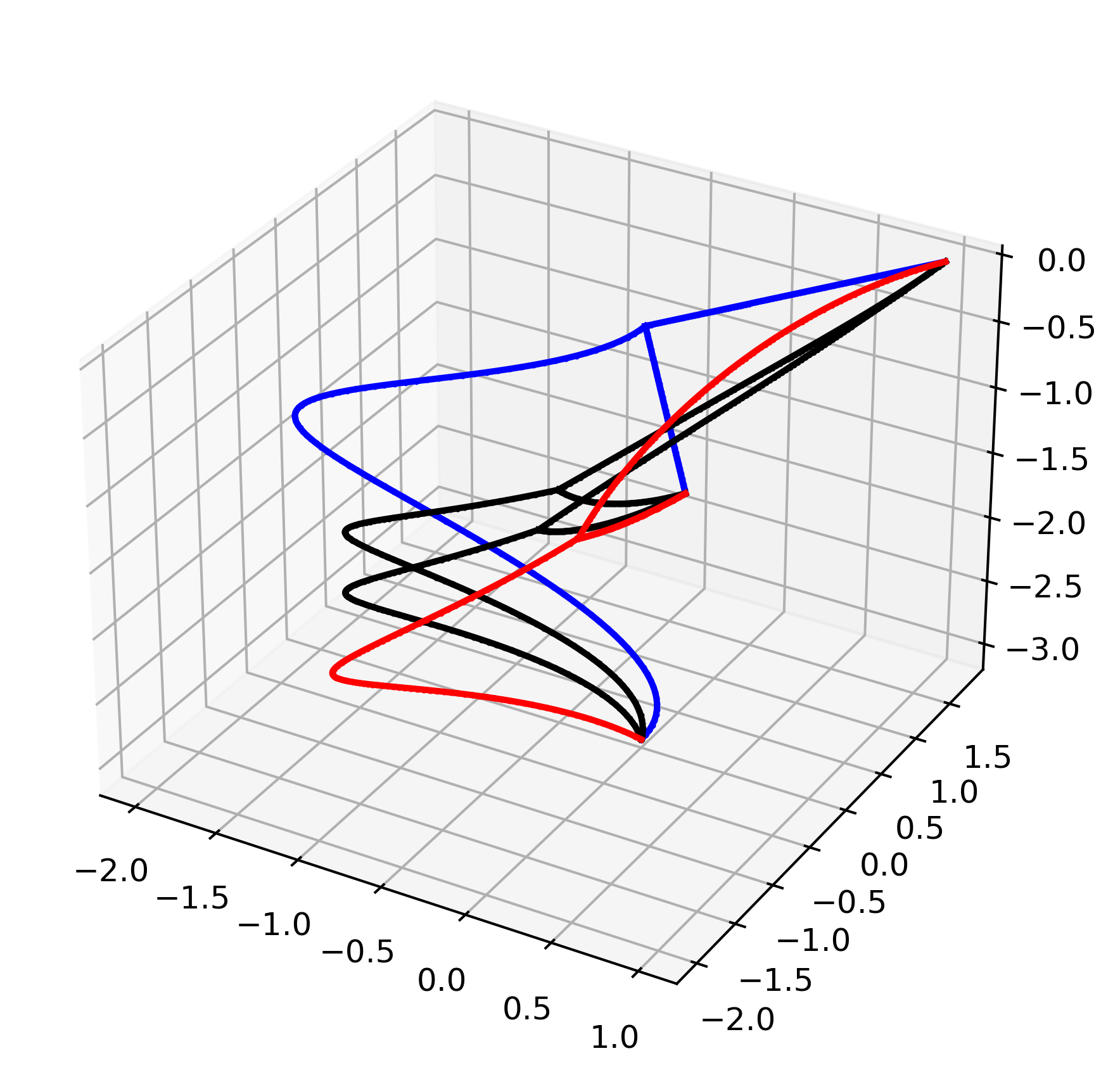}
\includegraphics[angle=-0,width=0.3\textwidth]{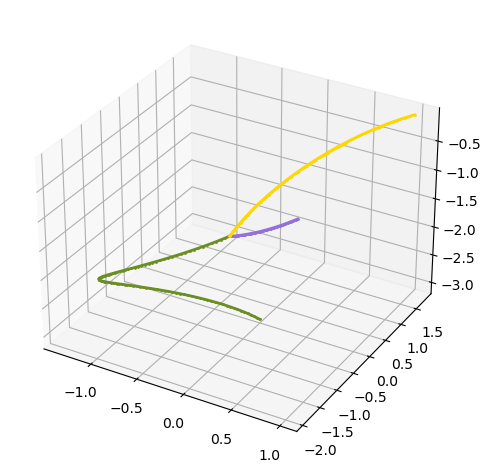}
\caption{The (projected) triod at time $t=0$, at times $t=0,1,2,7$ and 
at time $t=7$.
At the top right we show a plot of the discrete energy $L_E(c^m)$ over time.
}
\label{fig:vcircle}
\end{figure}

{\bf Experiment 13}:
We choose $\Sigma = (0,0,0)^{t}$ and let $P_1 = (0,0,-\pi)^{t}$,
$P_2 = (1, -\sqrt{3}, 0)^{t}$, $P_3 = (1, \sqrt{3}, 0)^{t}$.
As a consequence, the 2d projections of $\Sigma$ and $P_1$ are both equal to 
the origin. We display the evolution in Figure~\ref{fig:vcircle}, where we
notice that the projected triple junction $\proj\Sigma$ moves away from the
origin and that the projected triod eventually reaches a steady state made up
of three circle segments.

\begin{figure}
\center
\includegraphics[angle=-0,width=0.11\textwidth]{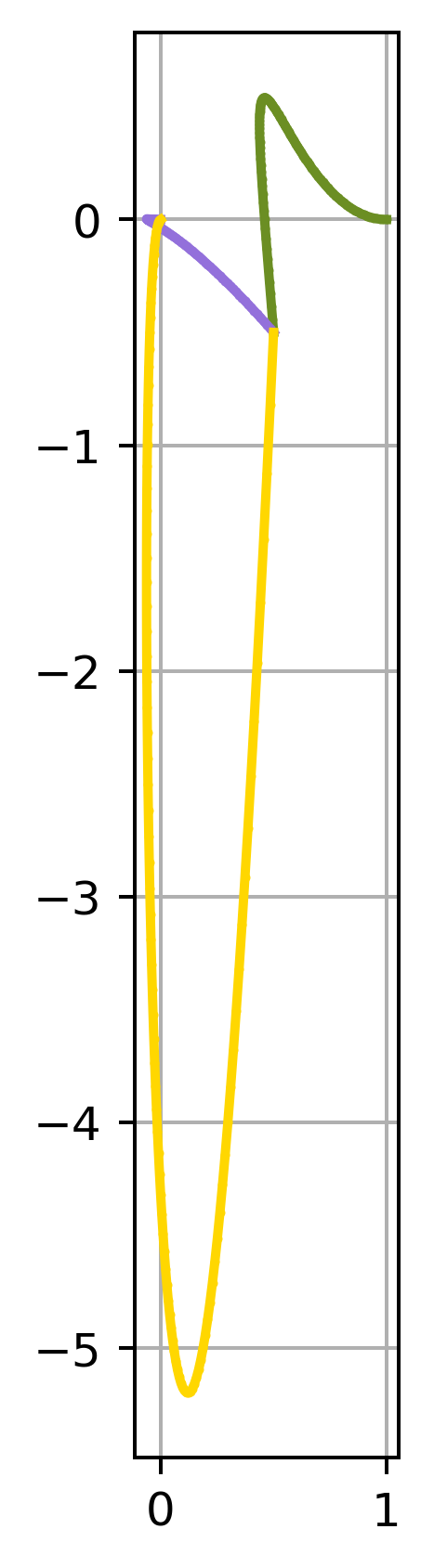}
\includegraphics[angle=-0,width=0.14\textwidth]{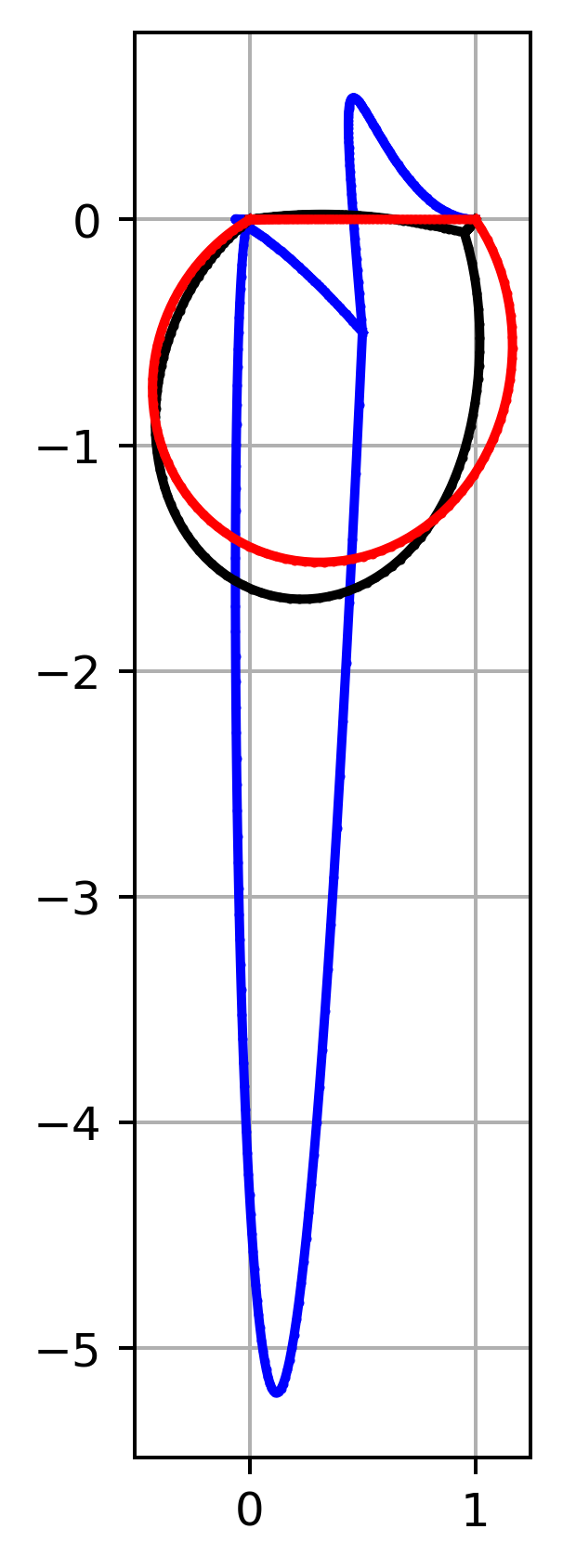}
\includegraphics[angle=-0,width=0.3\textwidth]{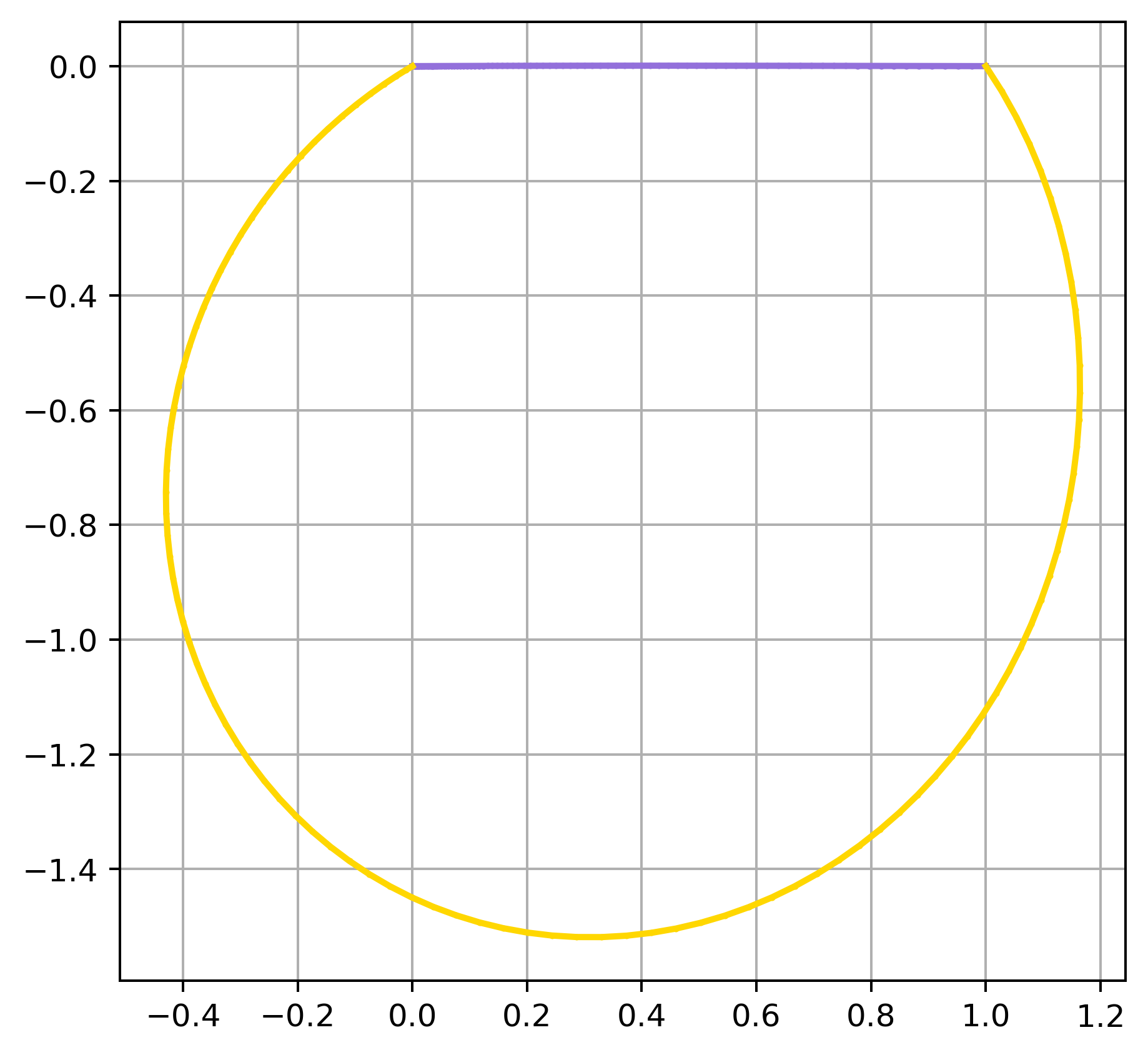}\\
\includegraphics[angle=-0,width=0.3\textwidth]{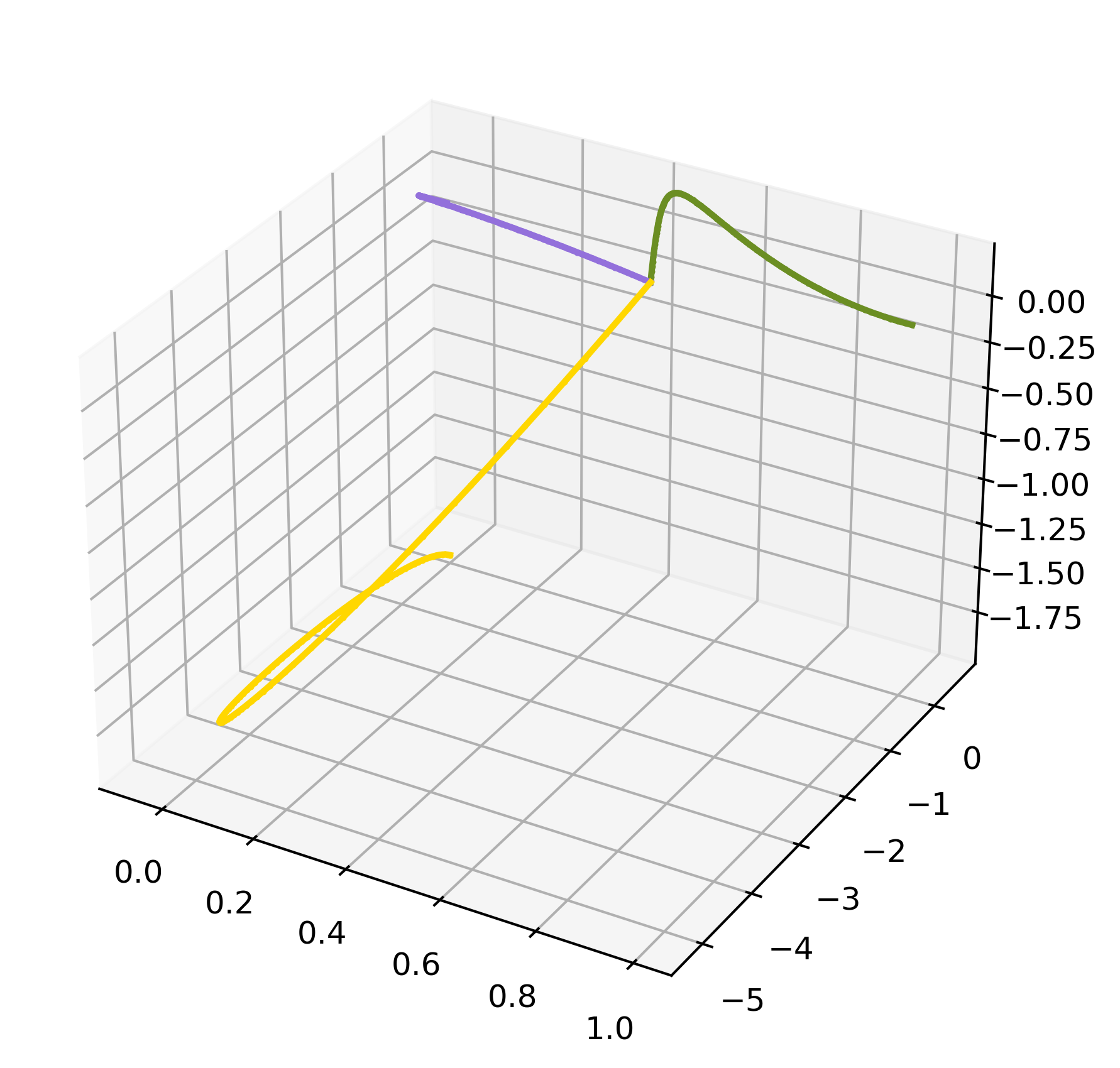}
\includegraphics[angle=-0,width=0.3\textwidth]{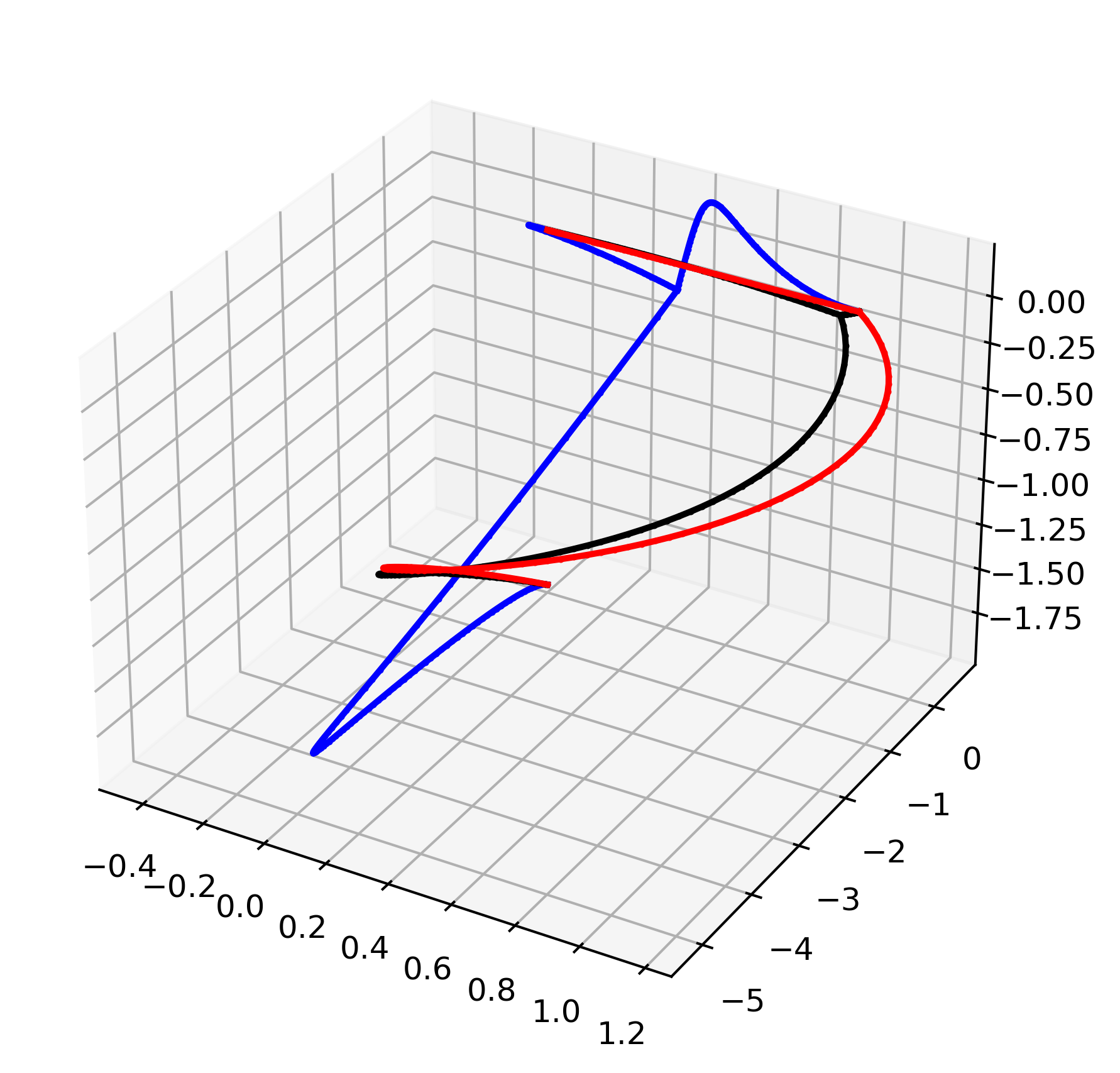}
\includegraphics[angle=-0,width=0.3\textwidth]{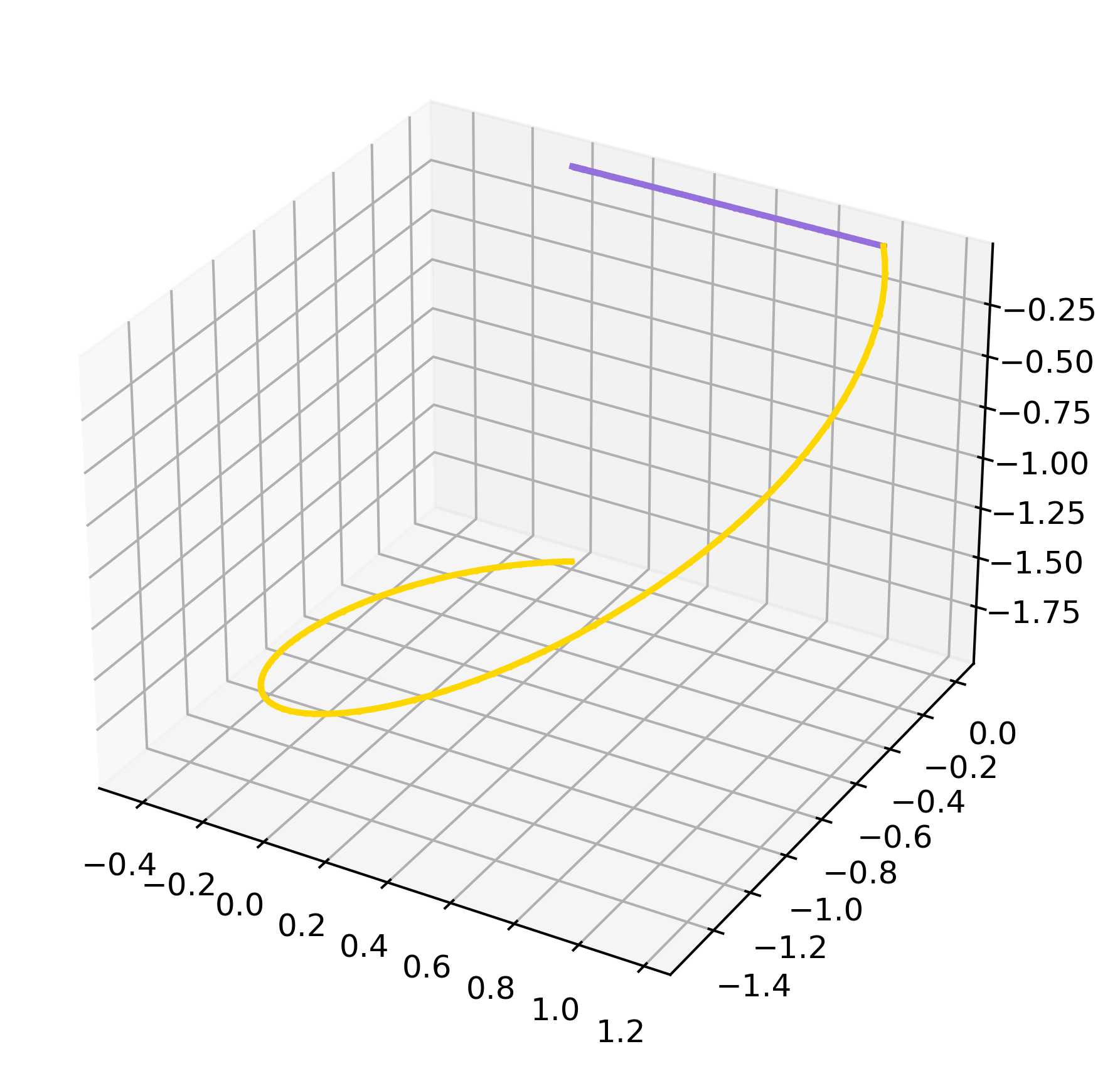}\\
\includegraphics[angle=-0,width=0.4\textwidth]{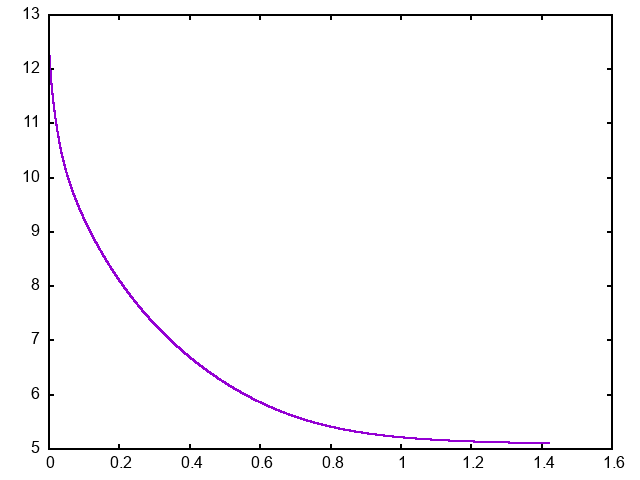}
\includegraphics[angle=-0,width=0.4\textwidth]{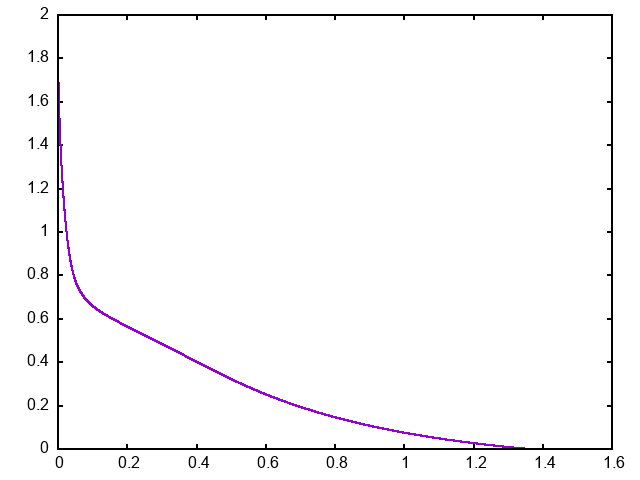}
\caption{The (projected) triod at time $t=0$, at times $t=0,1,1.42$ and 
at time $t=1.42$.
Below we show a plot of the discrete energies $L_E(c^m)$ and $L_E(c^m_1)$
over time.
}
\label{fig:ex27a}
\end{figure}

{\bf Experiment 14}:
We let $\Sigma = (0.5, -0.5, 0)^{t}$ and 
$P_1 = (1, 0, 0)^{t}$, $P_2=(0,0,0)^{t}$, $P_3 = (0, 0, -2)^{t}$, and compute the
initial horizontal curves connecting $\Sigma$ with $P_\alpha$, $\alpha=1,2,3$,
with the help of \eqref{curvacostruita} from Example~\ref{ex2.7}, where we use the
parameter $b=1$ for each curve. The ensuing evolution is shown in 
Figure~\ref{fig:ex27a}. It can be seen that the triple junction moves towards
$P_1$, leading to a singularity when the curve $c_1^m$ vanishes.

\begin{figure}
\center
\includegraphics[angle=-0,width=0.3\textwidth]{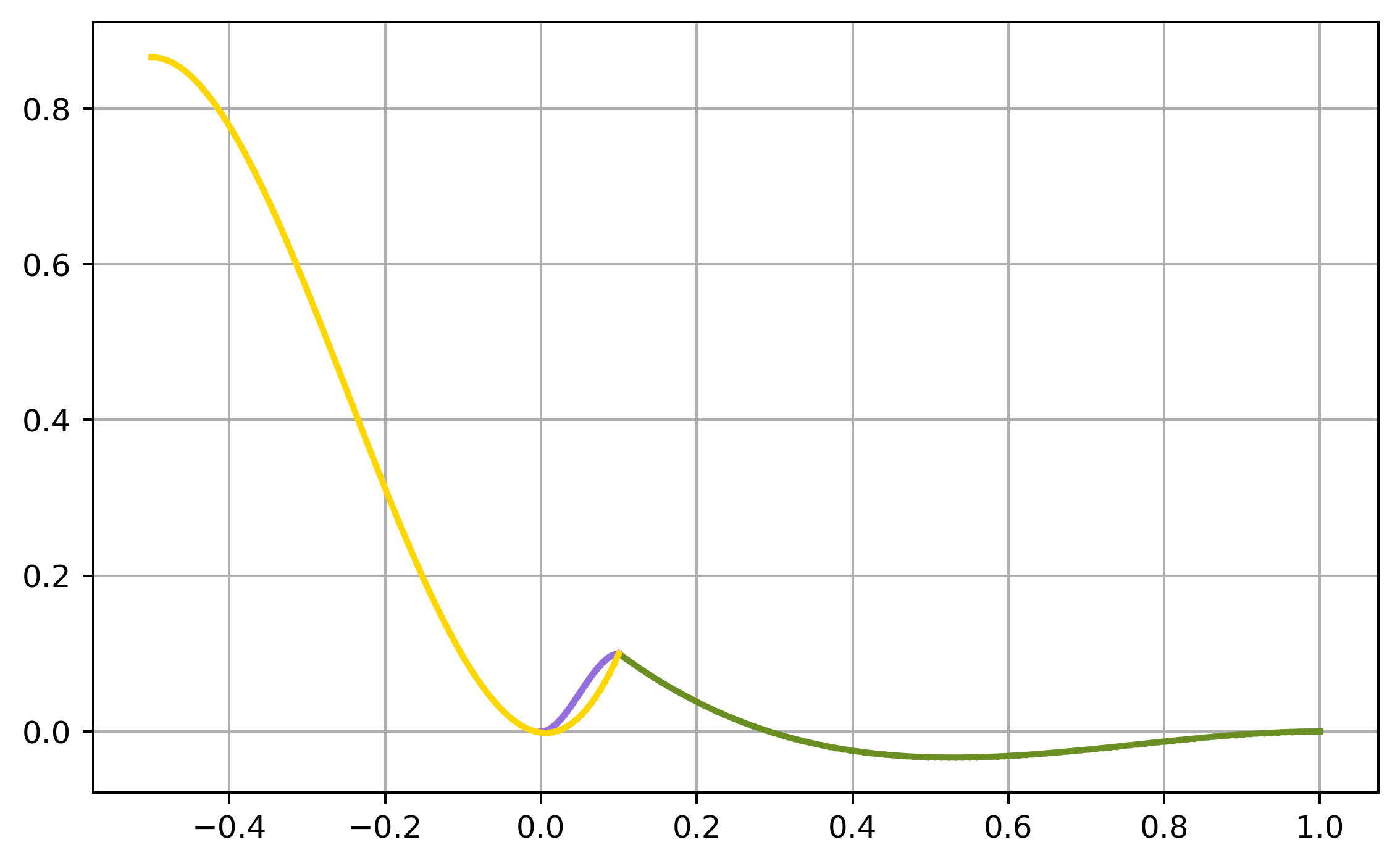}
\includegraphics[angle=-0,width=0.3\textwidth]{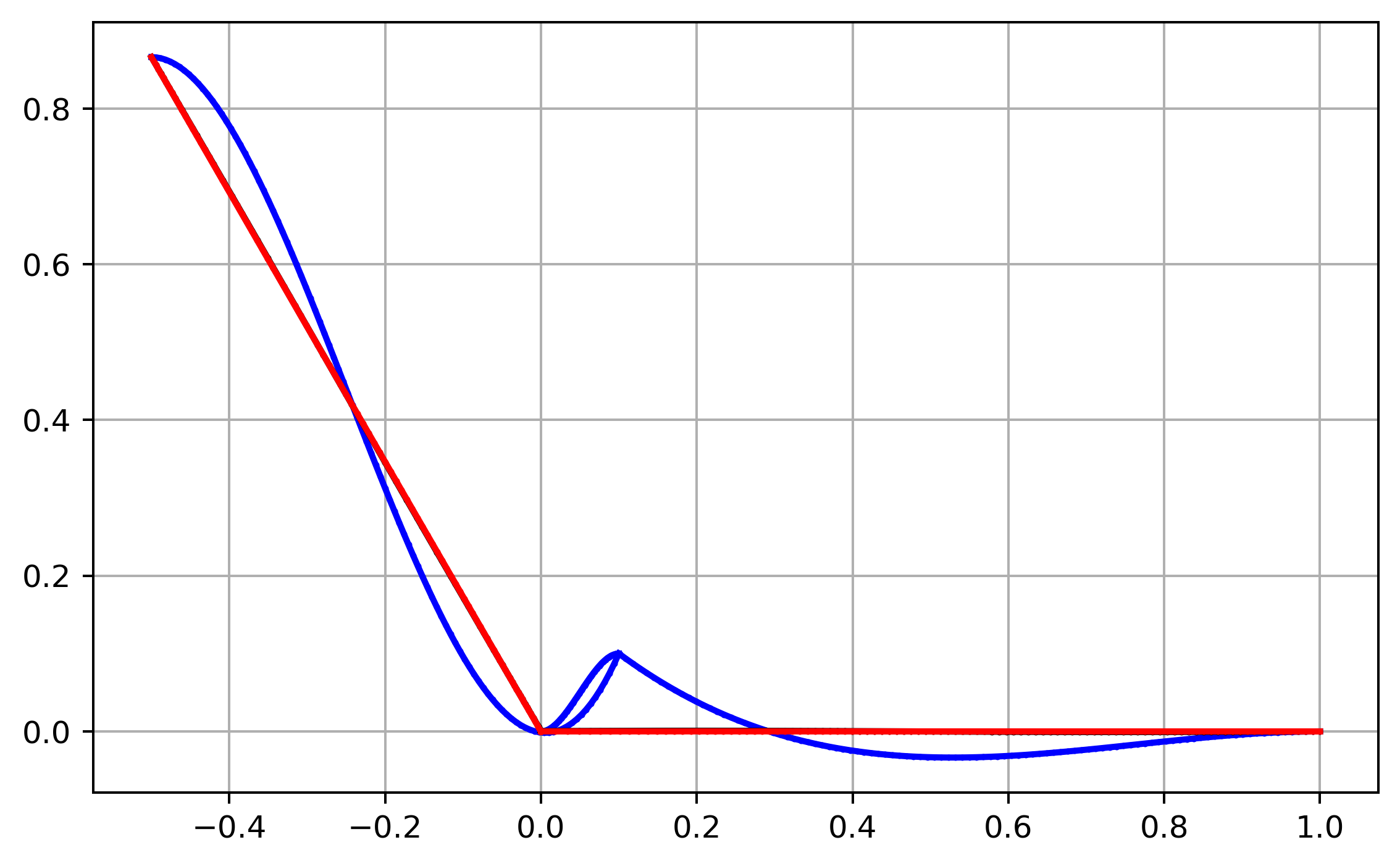}
\includegraphics[angle=-0,width=0.3\textwidth]{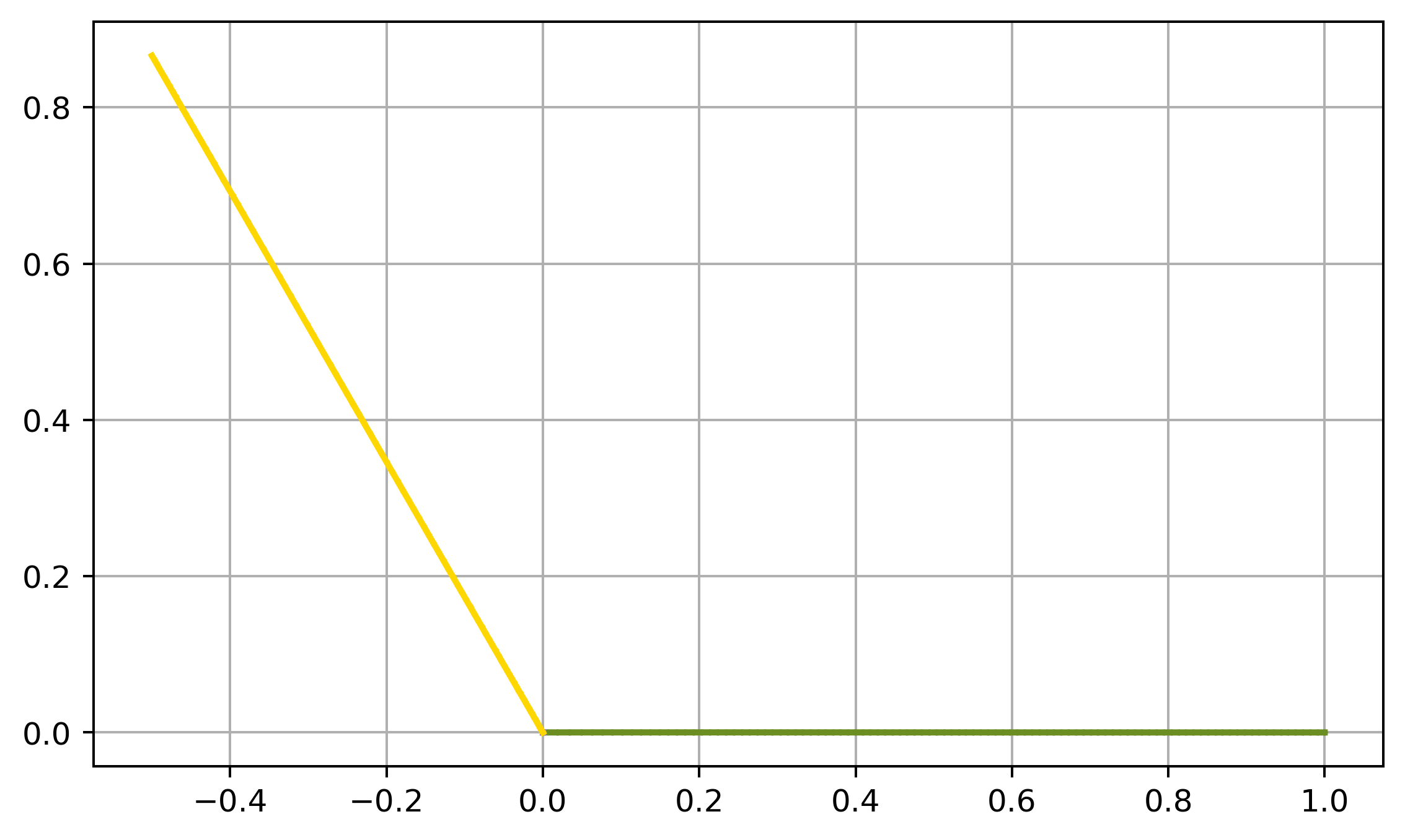} \\
\includegraphics[angle=-0,width=0.3\textwidth]{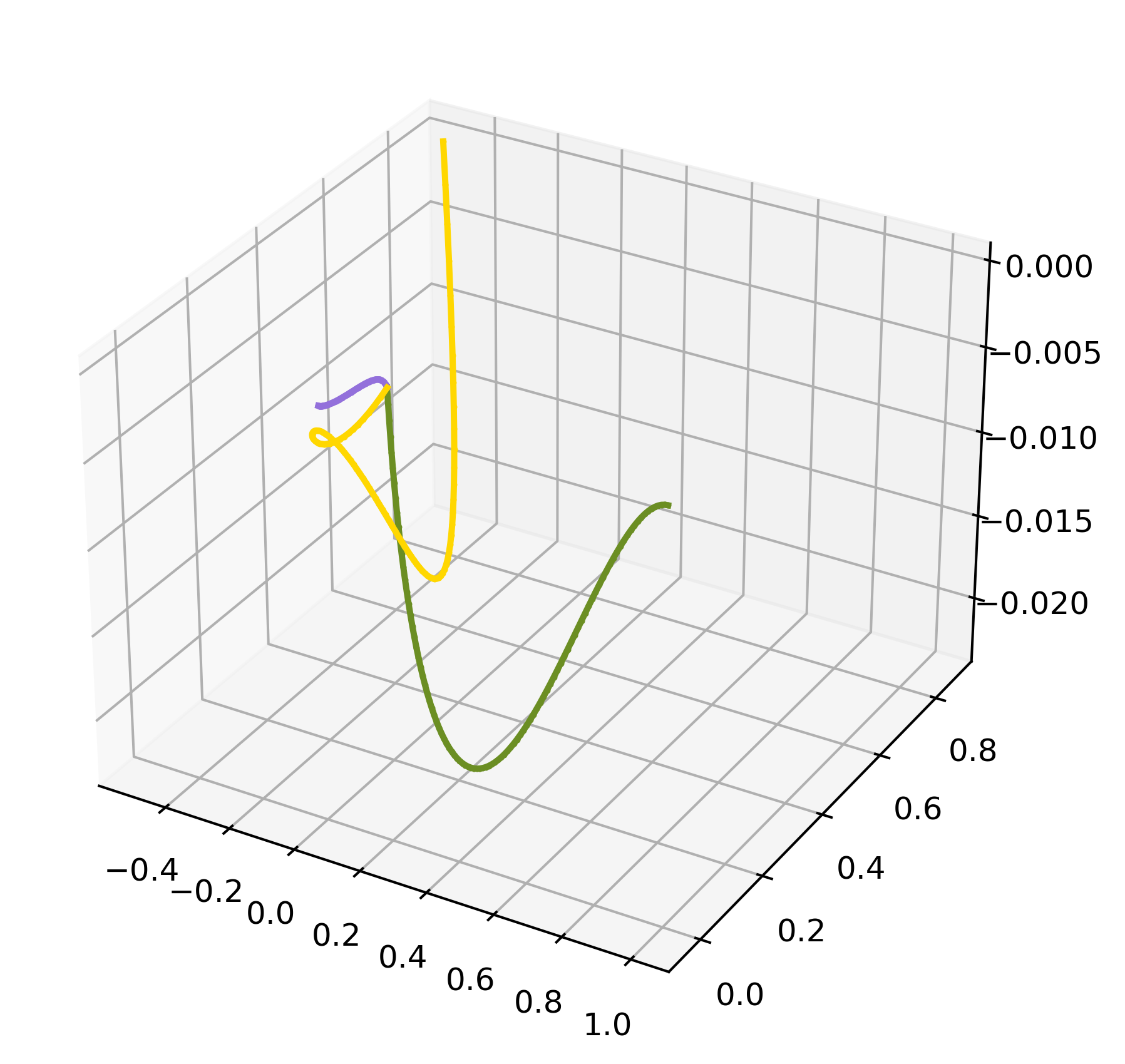}
\includegraphics[angle=-0,width=0.3\textwidth]{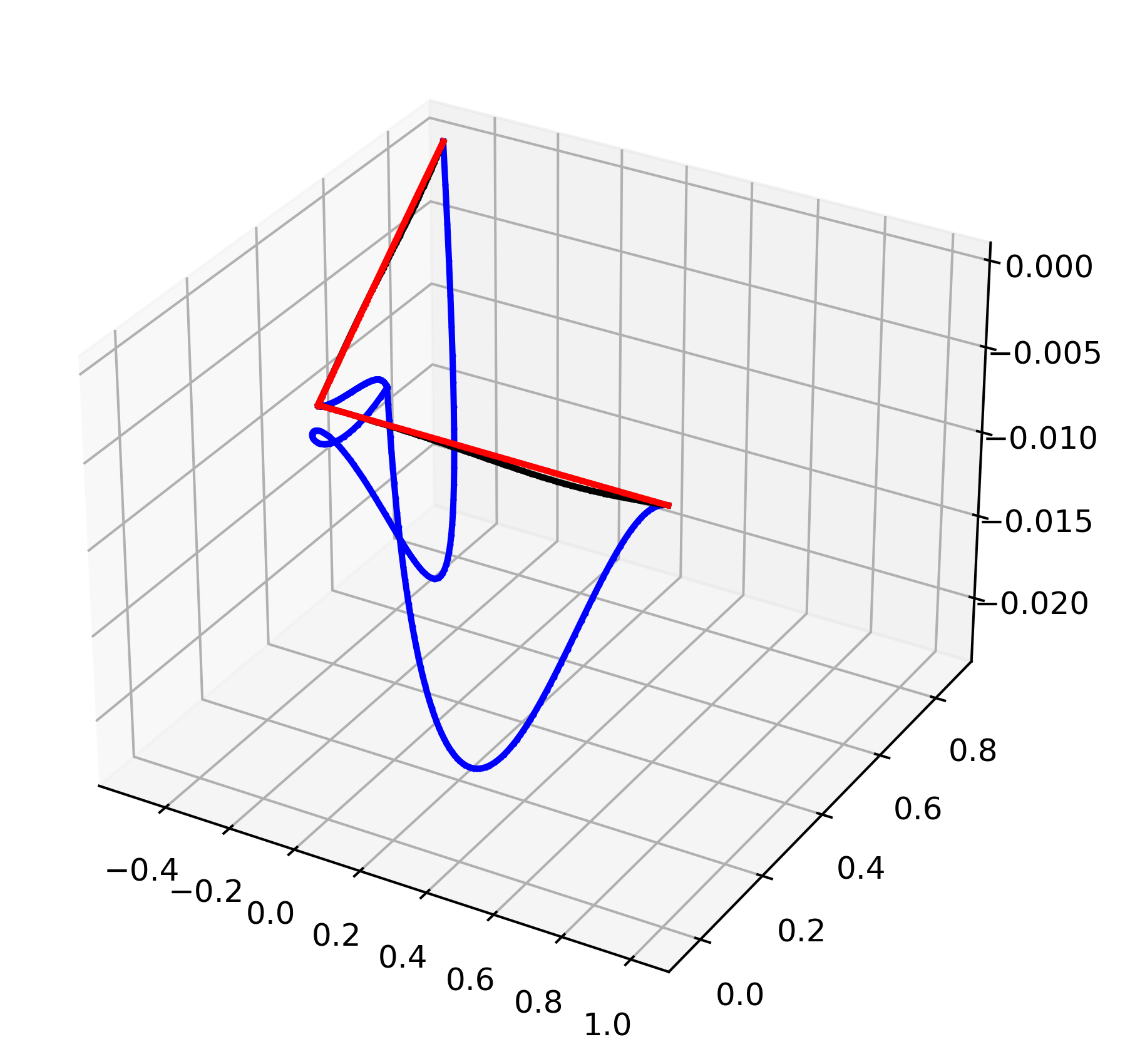}
\includegraphics[angle=-0,width=0.3\textwidth]{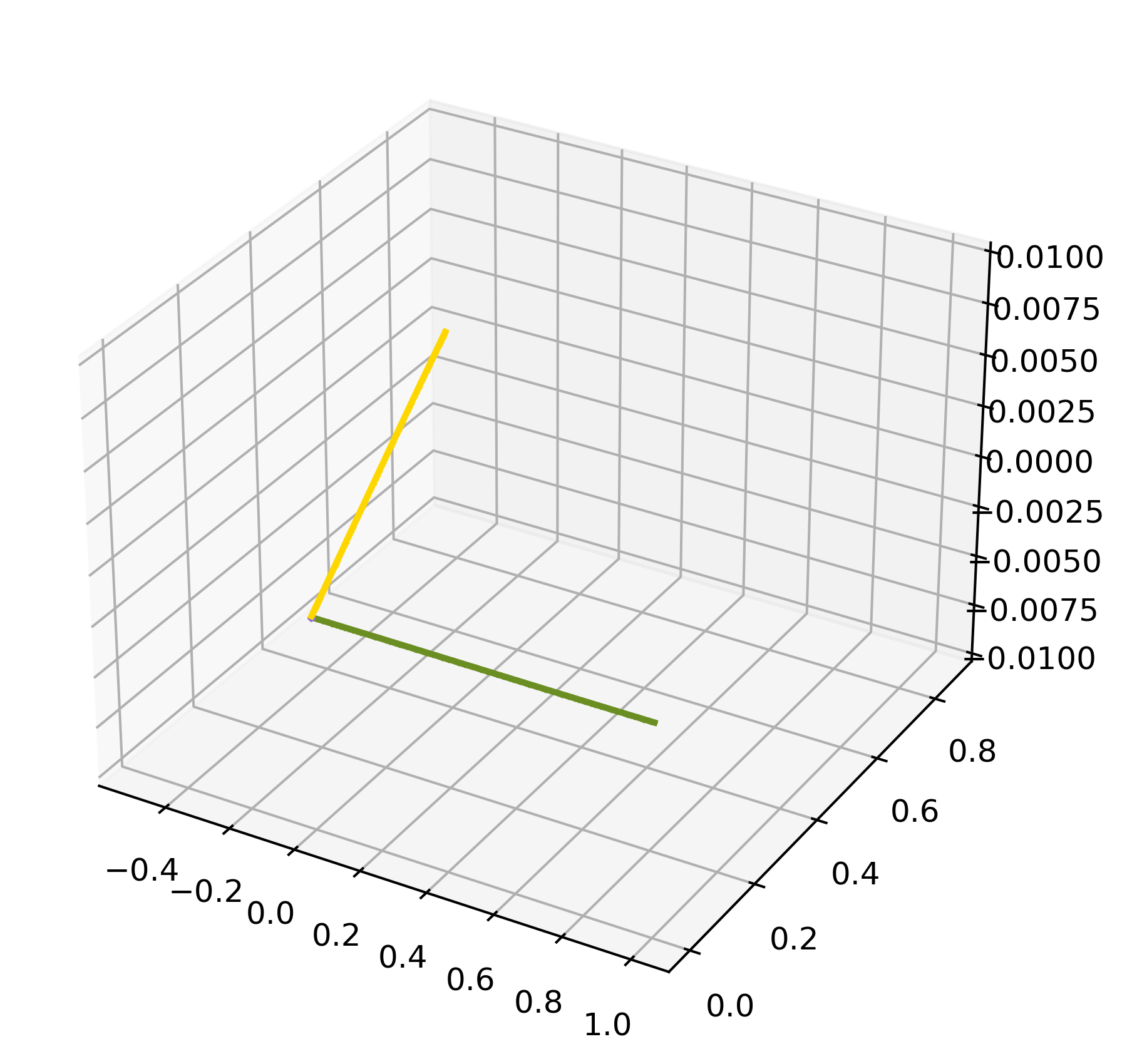} \\
\includegraphics[angle=-0,width=0.4\textwidth]{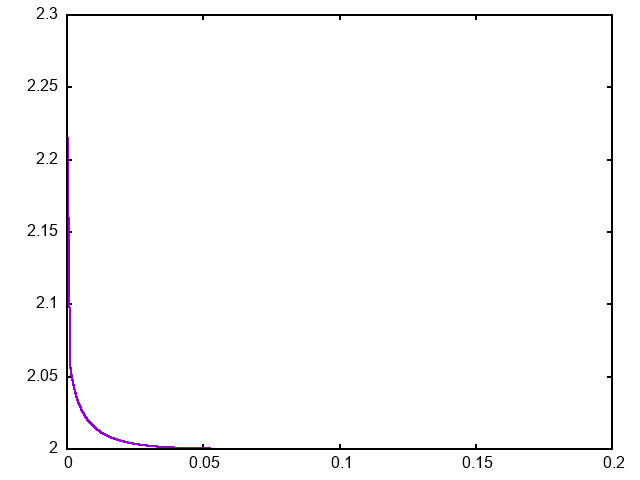}
\includegraphics[angle=-0,width=0.4\textwidth]{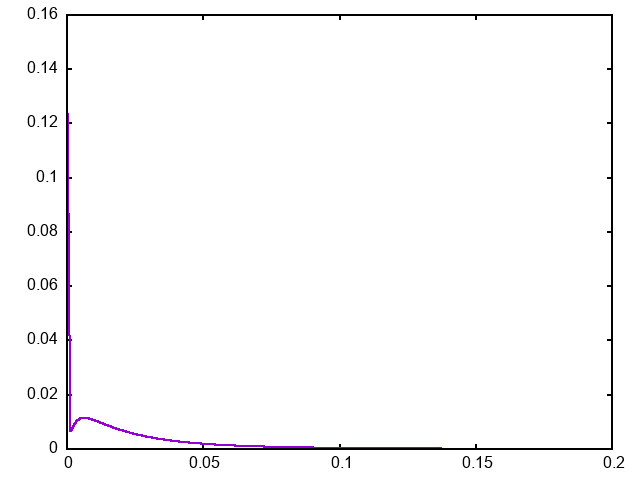}
\caption{The (projected) triod at time $t=0$, at times $t=0,0.1,0.2$ and 
at time $t=0.2$.
Below we show a plot of the discrete energies $L_E(c^m)$ and $L_E(c^m_2)$
over time.
}
\label{fig:ex27d}
\end{figure}

{\bf Experiment 15}:
We let $\Sigma = (0.1, 0.1, 0)^{t}$ and 
$P_1 = (1, 0, 0)^{t}$, $P_2=(0,0,0)^{t}$, \linebreak
$P_3 = (\cos(\frac23\pi), \sin(\frac23\pi),0)^{t} = (-0.5, \frac12\sqrt{3}, 0)^{t}$,
and once again compute the
initial horizontal curves connecting $\Sigma$ with $P_\alpha$, $\alpha=1,2,3$,
with the help of \eqref{curvacostruita}, where this time we always choose $b=0$.
Observe that the points $P_{\alpha}$ form a triangle with an angle of 120 degrees in $P_{2}$.
The evolution is shown in Figure~\ref{fig:ex27d}. 
We observe that the projected triple junction $\proj\Sigma$ moves towards
the origin, meaning that the curve $c^m_2$ vanishes.

\renewcommand{\thesection}{}

\appendix\renewcommand{\thesection}{\Alph{section}}
\setcounter{equation}{0}
\renewcommand{\theequation}{\Alph{section}.\arabic{equation}}

\section{Supporting materials}\label{secA}

For numerical simulations it is advantageous to be able to construct horizontal curves between given points in $\R^{3}$.
Following \cite[Lemma~1.13]{Calin} (or by direct computation) we give the following example.
\begin{ex}\label{ex2.7}
Let $\gamma:[0,1] \to \R^{3}$, $\gamma=\gamma(u)$, be given through
\[
\gamma(u)=\begin{pmatrix} au+bu^{2} \\ \alpha u^{2}+\beta u^{3}\\
\alpha a \frac{u^{3}}{6} + a\beta \frac{u^{4}}{4} + b \beta \frac{u^{5}}{10}
\end{pmatrix}.
\]
The curve $\gamma$ is horizontal and passes through the origin.
Let $M=(M^{1},M^{2},M^{3})^{t} \in \R^{3}$. Choosing
\begin{equation} \label{valoripara}
\alpha=M^{2}-\beta, \qquad a=M^{1}- b, \qquad \beta =\frac{60 M^{3} -10 M^{1}M^{2} + 10b M^{2}}{5M^{1}+b}
\end{equation}
we obtain a one-parameter family of horizontal curves (depending on $b$) such that $\gamma(0)=0$ and $\gamma(1)=M$.

Let now $P,Q \in \R^{3}$ be given. Using left translation we obtain that the curves (depending on $b$)
\begin{equation}\label{curvacostruita}
 \hat{\gamma}(u): =\tau_{P}(\gamma(u))= \begin{pmatrix}
P^{1}+ (au+bu^{2}) \\  P^{2}+(\alpha u^{2}+\beta u^{3})\\
P^{3}+ (\alpha a \frac{u^{3}}{6} + a\beta \frac{u^{4}}{4} + b \beta \frac{u^{5}}{10})
+\frac{1}{2} P^{1}(\alpha u^{2}+\beta u^{3}) -\frac{1}{2} P^{2}(au+bu^{2})
\end{pmatrix}
\end{equation}
are horizontal curves such that $\hat{\gamma}(0)=P$ and $\hat{\gamma}(1)=Q$, provided
we set
\[
M^{1}:=Q^{1}-P^{1}, \quad M^{2}:=Q^{2}-P^{2}, \quad M^{3}:=Q^{3}-P^{3}-\frac{1}{2}(P^{1}(Q^{2}-P^{2})-P^{2} (Q^{1}-P^{1}) ) 
\]
in \eqref{valoripara}.
\end{ex}

\begin{rem}[Solutions of the Euler--Lagrange equations \eqref{EL}]\label{rem1.10}
To derive an explicit solution for curves $\gamma$ satisfying \eqref{EL},
 we distinguish between two possible cases. 
In the following we let $c$ always denote the projection of $\gamma$.

\underline{First case:} if $\lambda=0$, i.e.\ $k(u)=0$ in $[0,1]$, then $c$ is a straight line connecting $(P^{1},P^{2})^{t}$ and $(Q^{1},Q^{2})^{t}$, that is
\[
\binom{\gamma^{1}(u)}{\gamma^{2}(u)}
=\binom{P^{1}}{P^{2}}
+u \binom{Q^{1}-P^{1}}{Q^{2}-P^{2}}
, \quad  \qquad u \in [0,1].
\]
From \eqref{g3} a straightforward calculation yields
\[
\gamma^{3}(u)= P^{3}+u\left(-\frac{1}{2}P^{2}Q^{1}+ \frac{1}{2}P^{1}Q^{2}\right)  \quad  \qquad u \in [0,1],
\]
and $Q^{3}$ must satisfy the equation
$Q^{3}=P^{3}-\frac{1}{2}P^{2}Q^{1}+ \frac{1}{2}P^{1}Q^{2}$.
We can write the straight line also as the horizontal curve
\[
\gamma(u)=P +u(Q^{1}-P^{1})X_{1}(P)+u(Q^{2}-P^{2})X_{2}(P)
\]
in the horizontal plane through $P$.
Observe that if $P$ is the origin (which we can always attain after an appropriate left translation, recall Lemma~\ref{leftT}), then the line lies in the $(x,y)$-plane.

\underline{Second case:} if $\lambda \neq 0$ then the Euler--Lagrange equations yield
\begin{align*}
\frac{1}{|c'|}   \left( \frac{1}{|c'|}  
\gamma^{1}_{u} \right)_{u} &=\lambda \frac{\gamma^{2}_{u}}{|c'|} \quad  \text{ in } (0,1),\\
\frac{1}{|c'|}   \left( \frac{1}{|c'|}  
\gamma^{2}_{u} \right)_{u} &=-\lambda \frac{\gamma^{1}_{u}}{|c'|} \quad  \text{ in } (0,1).
\end{align*}
Reparametrizing $c$ by arc length (i.e.\ $|c'(s)|=1$, $ds=|c'(u)|du$), the above system reads
\[
\gamma^{1}_{ss}=\lambda \gamma^{2}_{s}, \qquad \gamma^{2}_{ss}=-\lambda \gamma^{1}_{s} \qquad \text{ for } s \in (0, L_{E}(c)),
\]
i.e.
\[
c_{ss}=\lambda J c_{s}, \qquad \text{ with } J=\begin{pmatrix}
0& 1\\ -1 & 0 \end{pmatrix},
\]
with solution $c_{s}(s)= e^{\lambda J s}c_{s}(0)= \begin{pmatrix}
\cos(\lambda s)& \sin(\lambda s)\\ -\sin(\lambda s) & \cos(\lambda s)
\end{pmatrix} c_{s}(0)$ (cf.\ \cite[Lemma 1.16]{Calin}). Integration gives
\[
c(s)=c(0)+ \int_{0}^{s}e^{\lambda J \xi}c_{s}(0) d\xi = c(0)+\left[\frac{1}{\lambda}J^{-1}e^{\lambda J \xi}c_{s}(0) \right]_{0}^{s} =c(0)- \frac{1}{\lambda}J^{-1} c_{s}(0)+\frac{1}{\lambda}J^{-1}e^{\lambda J s}c_{s}(0). 
\]
Therefore  a general solution is of type
\[
c(s)=A+\frac{1}{\lambda}\begin{pmatrix} 0& -1\\ 1 & 0 \end{pmatrix} 
\begin{pmatrix}
\cos(\lambda s)& \sin(\lambda s)\\ -\sin(\lambda s) & \cos(\lambda s)
\end{pmatrix} B = A+\frac{1}{\lambda} \begin{pmatrix}
\sin(\lambda s)& -\cos(\lambda s)\\ \cos(\lambda s) & \sin(\lambda s)
\end{pmatrix} B,
\]
for vectors $A,B \in \R^{2}$,  with $B=c_{s}(0)$ and $A=c(0)- \frac{1}{\lambda}J^{-1} c_{s}(0)$. In particular, we see that the projected curve $c$ is an arc of a circle of radius $|B|/|\lambda|= \frac{1}{|\lambda|}$ (since $|B|=|c_{s}(0)|=1$ due to the arc length parametrization). 
Writing
\[
c_{s}(0)=B=(\cos(\alpha_{0}), \sin(\alpha_{0}))^{t}
\]
for some $\alpha_{0} \in [0, 2\pi)$, we obtain
\begin{equation} \label{eq-c}
c(s)= A+\frac{1}{\lambda} \begin{pmatrix}
 \sin(\lambda s)& -\cos(\lambda s)\\ \cos(\lambda s) & \sin(\lambda s)
 \end{pmatrix} \begin{pmatrix}
 \cos (\alpha_{0})\\ \sin(\alpha_{0})
 \end{pmatrix} =A + \frac{1}{\lambda} \begin{pmatrix}
 \sin (\lambda s-\alpha_{0})\\ \cos(\lambda s-\alpha_{0})
 \end{pmatrix}.
\end{equation}
The constants $A$, $B$ are determined through the equations
\[
\begin{pmatrix}
 P^{1}\\ P^{2}
 \end{pmatrix}= c(0)=A+\frac{1}{\lambda} \begin{pmatrix}
 0& -1\\ 1 & 0
 \end{pmatrix} B
\]
and, with $s_f:=L_{E}(c)=L(\gamma)$,
\[
\begin{pmatrix}
Q^{1}\\ Q^{2}
\end{pmatrix}= c(s_{f})=A+\frac{1}{\lambda} \begin{pmatrix}
\sin(\lambda s_{f})& -\cos(\lambda s_{f})\\ \cos(\lambda s_{f}) & \sin(\lambda s_{f})
\end{pmatrix} B.
\]
This gives for $ \lambda s_{f} \in (0, 2\pi)+2\pi \mathbb{Z}$:
\begin{align}
B&=\frac{\lambda}{\sin^{2}(\lambda s_{f})+ (1-\cos(\lambda s_{f}))^{2}}\begin{pmatrix}
\sin(\lambda s_{f})& \cos(\lambda s_{f})-1\\ 1-\cos(\lambda s_{f}) & \sin(\lambda s_{f})
\end{pmatrix}\begin{pmatrix}
Q^{1}-P^{1}\\Q^{2} -P^{2}
\end{pmatrix} \notag \\ \label{eq-B}
&=\frac{\lambda}{2(1-\cos(\lambda s_{f}))}\begin{pmatrix}
\sin(\lambda s_{f})& \cos(\lambda s_{f})-1\\ 1-\cos(\lambda s_{f}) & \sin(\lambda s_{f})
\end{pmatrix}\begin{pmatrix}
Q^{1}-P^{1}\\Q^{2} -P^{2}
\end{pmatrix}, \\  \label{eq-A}
A&=\begin{pmatrix}
P^{1}\\ P^{2}
\end{pmatrix} -\frac{1}{\lambda} \begin{pmatrix}
0& -1\\ 1 & 0
\end{pmatrix} B .
\end{align}
We deduce in particular that $(P^{1},P^{2}) \neq (Q^{1},Q^{2})$.
After computing the expressions \eqref{eq-B}, \eqref{eq-A}, we obtain the curve $c$ by expression \eqref{eq-c}.
The third component of the horizontal lift of $c$ starting at $P$ is finally computed via \eqref{g3}, and hence $Q^{3}$ must satisfy
\begin{align*}
Q^{3} &= P^{3} + \int_{0}^{s_{f}} \left(-\frac{1}{2}\gamma^{2}(s)\gamma^{1}_{s}(s) +\frac{1}{2}\gamma^{1}(s)\gamma^{2}_{s}(s) \right) ds\\
&=P^{3} +\int_{0}^{s_{f}}\Big ( -\frac{1}{2}A^{2}[B^{1} \cos(\lambda s) + B^{2} \sin(\lambda s)]
+ \frac{1}{2}A^{1}[B^{2} \cos(\lambda s) - B^{1} \sin(\lambda s)] -\frac{|B|^{2}}{2\lambda} \Big ) ds.
\end{align*}
In particular if $P=0$ (which we can always attain after appropriate left translation) then
$A^{1}=\frac{1}{\lambda}B^{2}$, $A^{2}=-\frac{1}{\lambda}B^{1}$ and 
\begin{equation} \label{eqQ3bis}
\gamma^{3}(s)=-\frac{1}{2\lambda} s +\frac{1}{2\lambda^{2}}\sin(\lambda s)= \frac{\sin(\lambda s) -\lambda s }{2\lambda^{2}} \qquad \text{for }s \in [0,s_{f}],
\end{equation}
so that 
\begin{equation} \label{eqQ3}
Q^{3}
= \frac{\sin(\lambda s_{f}) -\lambda s_{f} }{2\lambda^{2}}  .
\end{equation}

If $ \lambda s_{f} =2\pi k$, for some $ k \in  \mathbb{Z}$, then $(P^{1},P^{2}) =(Q^{1},Q^{2})=c(0)$, i.e.\ the projected curve $c$ is a ($k$-multiply covered) circle, whereby we have freedom in choosing the vector $B=c_{s}(0)=(\cos (\alpha_{0}), \sin (\alpha_{0}))^{t}$ for $\alpha_{0} \in [0,2\pi)$. For this one-parameter family of solutions we obtain  by \eqref{g3} and \eqref{eq-c} that   
\[
Q^{3}=P^{3}+ \int_{0}^{\lambda s_{f}}-\frac{1}{2}\gamma^{2}(s)\gamma^{1}_{s}(s) + \frac{1}{2}\gamma^{1}(s)\gamma^{2}_{s}(s) ds = P^{3}- (2\pi k)\frac{|B|^{2}}{2\lambda^{2}}.
\]
In particular 
\[
|Q_{3}-P_{3}|=\frac{(2\pi |k|)}{2\lambda^{2}}=\frac{\pi |k|}{\lambda^{2}}.
\]
Note that unlike the Euclidean setting, we have here a one-parameter family of solution curves having all the same length.
\end{rem}

\section{Horizontal curve shortening flow} \label{sec:HCSF}
Although Remark~\ref{rem1.10}  (see also \cite[\S1.6]{Calin} for further information) fully describes horizontal curves of minimal length, an alternative way to obtain the critical points \eqref{EL}, is to let a given horizontal curve flow into a critical configuration. This approach is interesting in its own right. In the following we provide an informal definition of such a flow and call it horizontal curve shortening flow. 
\begin{ass}\label{ass1}
Let $P\neq Q$ be two given points in $\R^{3}$ and let $\gamma_{0}:[0,1] \to \R^{3}$ be a regular smooth horizontal curve connecting $P$ and $Q$.
\end{ass}
Taking into consideration Remark~\ref{rem-EL}, \eqref{kg}, and  Lemma~\ref{lemma-hproperty}, it is natural to consider the following flow, which we call \emph{horizontal curve shortening flow}: 
\begin{prob}\label{HCSF1}
Let Assumption~\ref{ass1} hold. 
Find 
$\gamma:[0,T) \times [0,1] \to \R^{3}$, $\gamma=\gamma(t,u)$, which  satisfies $\gamma(0, \cdot)=\gamma_{0}$, and
\begin{align*}
\gamma_{t}(t,u)&=\vec{k}_{g}(t,u)+\lambda(t)N(t,u) -\Big(\int_{0}^{u}(k(t,\xi)+\lambda(t)) |\gamma_{u}(t,\xi)|_{g}d\xi \Big)X_{3}(\gamma(t,u)) \\ & \hspace{8cm}
\qquad \text{ for all } (t,u)\in (0,T) \times (0,1), \\
\gamma(t,0)&=P \quad \text{for all } t \in (0,T),\\
\gamma(t,1)&=Q \quad \text{for all } t \in (0,T),
\end{align*}
where
\[
\lambda(t)=- \frac{\int_{0}^{1} k(t,u) |\gamma_{u}(t,u)|_{g} du
}{L(\gamma(t))}.
\]
\end{prob}

\begin{lemma}\label{lemma1.2}
If $\gamma:[0,T) \times [0,1] \to \R^{3}$ is a solution to Problem~\ref{HCSF1}, then the length $L(\gamma(t))$ of the horizontal curve $\gamma(t)$ decreases along the flow.
\end{lemma}
\begin{proof}
A direct computation gives
\begin{align*}
\frac{d}{dt}L(\gamma(t))&= \frac{d}{dt} \int_{0}^{1} \sqrt{(\gamma^{1}_{u})^{2}+ (\gamma^{2}_{u})^{2}} du\\
&=-\int_{0}^{1}\frac{1}{|\gamma_{u}|_{g}}   \left( \frac{1}{|\gamma_{u}|_{g}}  \begin{pmatrix}
\gamma^{1}_{u}\\\gamma^{2}_{u} 
\end{pmatrix} \right)_{u}  \cdot \begin{pmatrix}
 \gamma_{t}^{1}\\ \gamma_{t}^{2}
\end{pmatrix} |\gamma_{u}|_{g} du
 + \left [ \frac{1}{|\gamma_{u}|_{g}}  \begin{pmatrix}
\gamma^{1}_{u}\\\gamma^{2}_{u} 
\end{pmatrix}   \cdot \begin{pmatrix}
 \gamma_{t}^{1}\\ \gamma_{t}^{2}
\end{pmatrix}  \right ]_{0}^{1}\\
&=-\int_{0}^{1}k(k+\lambda) |\gamma_{u}|_{g} du =- \int_{0}^{1}k^{2} |\gamma_{u}|_{g} du -\lambda\int_{0}^{1}k |\gamma_{u}|_{g} du \leq0,
\end{align*}
where we have used that $\gamma_{t}=0$ at the boundary, the definition of $\lambda$, and a H\"older inequality.
\end{proof}

Note that if $\gamma$ solves the above horizontal curve shortening flow (Problem~\ref{HCSF1}), then its projected curve
$c=(\gamma^{1},\gamma^{2})^{t}$ satisfies the following problem in the Euclidean setting:
\begin{prob}\label{HCSF1p}
Let Assumption~\ref{ass1} hold.
Find $c:[0,T) \times [0,1] \to \R^{2}$, $c=(\gamma^{1}, \gamma^{2})^{t}$, which  satisfies $\gamma^{r}(0, \cdot)=\gamma_{0}^{r}$ for $r=1,2$, and
\begin{align}\label{ACCSF1}
c_{t}(t,u)&=\vec{k}(t,u)+\lambda(t)\vec{n}(t,u)  \qquad \text{ for all } (t,u)\in (0,T) \times (0,1), \notag \\
c(t,0)&=(P^{1},P^{2})^{t} \quad \text{for all } t \in (0,T), \notag \\
c(t,1)&=(Q^{1},Q^{2})^{t} \quad \text{for all } t \in (0,T),  
\end{align}
where
\begin{align} \label{ACCSF2}
\vec{n}(t,u)&:=\frac{1}{|c_{u}(t,u)|}\begin{pmatrix}
-\gamma^{2}_{u}(t,u)\\\gamma^{1}_{u}(t,u) 
\end{pmatrix}, \notag \\
\lambda(t)&=- \frac{\int_{0}^{1} k(t,u) |c_{u}(t,u)| du }{L_{E}(c(t))} 
=- \frac{\int_{0}^{1}\vec{k}(t,u) \cdot \vec{n}(t,u) |c_{u}(t,u)| du }
{L_{E}(c(t))} .
\end{align}
\end{prob}
On the other hand, if $c=(\gamma^{1},\gamma^{2})^{t}$ is a solution of Problem~\ref{HCSF1p} then its horizontal lift $\gamma:[0,T) \times [0,1] \to \R^{3}$ with
\[
\gamma^{3}(t,u) = P^{3} + \int_{0}^{u} \left(-\frac{1}{2}\gamma^{2}(t,\xi)\gamma^{1}_{u}(t,\xi) +\frac{1}{2}\gamma^{1}(t,\xi)\gamma^{2}_{u}(t,\xi) \right) d\xi
\]
yields (using integration by parts and the fact that $\gamma^{1}_{t}=\gamma^{2}_{t}=0$ at the boundary)
\begin{align*}
\gamma^{3}_{t}(t,u)&= \int_{0}^{u} \left(-\frac{1}{2}\gamma^{2}_{t}\gamma^{1}_{u} -\frac{1}{2}\gamma^{2}\gamma^{1}_{tu} +\frac{1}{2}\gamma^{1}_{t}\gamma^{2}_{u} +\frac{1}{2}\gamma^{1}\gamma^{2}_{tu}\right) d\xi \\
& =\int_{0}^{u} \left(-\gamma^{2}_{t}\gamma^{1}_{u}  +\gamma^{1}_{t}\gamma^{2}_{u} \right) d\xi +\frac{1}{2} \gamma^{1}(u)\gamma^{2}_{t}(u) -\frac{1}{2} \gamma^{2}(u)\gamma^{1}_{t}(u)\\
&=\int_{0}^{u} \left(-\gamma^{2}_{t}\gamma^{1}_{u}  +\gamma^{1}_{t}\gamma^{2}_{u} \right) d\xi 
+\frac{1}{2}\gamma^{1}(u)(k(t,u) +\lambda(t)) \frac{\gamma^{1}_{u}}{|c_{u}|} -
\frac{1}{2}\gamma^{2}(u)(k(t,u) +\lambda(t))(- \frac{\gamma^{2}_{u}}{|c_{u}|})
\\
&= \vec{k}_{g}^{3}(t,u)+\lambda(t)N^{3}(t,u)+\int_{0}^{u} -(k+\lambda)n^{2} n^{2}|c_{u}| +(k+\lambda)n^{1}(- n^{1}|c_{u}|) d\xi\\
&=
\vec{k}_{g}^{3}(t,u) +\lambda(t)N^{3}(t,u)-\int_{0}^{u}(k+\lambda(t))|c_{u}| d\xi.
\end{align*}
By assumption $\gamma^{3}_{0}$ fulfills  the boundary conditions and has the structure given in \eqref{g3}. By Lemma~\ref{lemma-hproperty} the component $\gamma^{3}(t, \cdot)$ maintains this structure and fulfills the boundary and initial condition. Hence we see that by taking the horizontal lift of a solution to Problem~\ref{HCSF1p}, we recover a solution to Problem~\ref{HCSF1}.
 
\begin{rem}
(i) Note that Assumption~\ref{ass1} restricts the class of initial curves for which Problem~\ref{HCSF1p} is solved (recall Remark \ref{rem1.2}).

(ii) Recall that $\int_{0}^{1}k |c_{u}|du = \theta(1)-\theta(0)$, where $\theta$ is the angle that the unit tangent of the planar curve $c$ makes  with the $x$-axis. Even if $(P^{1},P^{2})^t=(Q^{1},Q^{2})^t$, the quantity $(\theta(1)-\theta(0))$ need not be a multiple of $2\pi$, since no angle condition is imposed at the boundary. In other words, the projected curve might be closed, in the sense that $c(0)=c(1)$, but not smooth across $c(0)=c(1)$.
\end{rem}

\begin{rem}
As in Lemma~\ref{lemma1.2}, we see that $\frac{d}{dt} L_{E}(c(t)) \leq 0$ for a smooth solution to Problem~\ref{HCSF1p}. We give here a different way to derive the same result, as it provides further geometrical insight into the role of the Lagrange multiplier.
Using that $c_{t}=0$ at the boundary we infer
\begin{align*}
\frac{d}{dt} L_{E}(c(t)) 
&=\int_{0}^{1} \frac{c_{u}}{|c_{u}|} \cdot c_{tu} du\\&= \left[\frac{c_{u}}{|c_{u}|} \cdot c_{t}\right]_{0}^{1}- \int_{0}^{1}\vec{k} \cdot c_{t} |c_{u}| du =- \int_{0}^{1}| c_{t}|^{2} |c_{u}| du +\lambda(t)\int_{0}^{1}  \vec{n} \cdot c_{t} |c_{u}|du.
\end{align*}
Note that for a smooth solution  of the flow we have
\[
\lambda(t)=- \frac{\int_{0}^{1} k(t,u) |c_{u}| du }{L_{E}(c(t))} =- \frac{\int_{0}^{1}\vec{k}(t,u) \cdot \vec{n}(t,u) |c_{u}| du }{L_{E}(c(t))} = -\frac{\int_{0}^{1}c_{t} \cdot \vec{n} |c_{u}| du }{L_{E}(c(t))} +\lambda (t),
\]
so that we infer (compare with Remark~\ref{rem1.2} for a geometric interpretation)
\[
0=\int_{0}^{1}c_{t} \cdot \vec{n} |c_{u}| du = \frac{d}{dt} \left( \int_{0}^{1} \frac{1}{2} \begin{pmatrix}
\gamma^{1}(u)\\ \gamma^{2}(u) 
\end{pmatrix} \cdot\begin{pmatrix}
-\gamma^{2}_{u}(u)\\\gamma^{1}_{u}(u) 
\end{pmatrix} du  \right)+\frac{1}{2} [c_{t} \cdot c^{\perp}]_{0}^{1} 
\]
and 
\[
\frac{d}{dt} L_{E}(c(t)) =- \int_{0}^{1}| c_{t}|^{2} |c_{u}| du +\lambda(t)\int_{0}^{1}  \vec{n} \cdot c_{t} |c_{u}|du=- \int_{0}^{1}| c_{t}|^{2} |c_{u}| du  \leq 0. 
\]
\end{rem}

\begin{rem}
The Euclidean evolution \eqref{ACCSF1}, \eqref{ACCSF2} has been studied analytically in different settings subject to different boundary  conditions: for instance for smooth closed (i.e.\ periodic) curves (where $\int_{0}^{1}k |c_{u}|du =2\pi m$, for $m \in \mathbb{Z}$) see \cite{Gage}; for the study of the  area preserving curve shortening flow with Neumann free boundary conditions outside of a convex domain in the plane see \cite{Elena}.
Evolution by mean curvature flow for Legendre submanifolds in a Sazaki manifold is studied in \cite{Smoczyk}. The evolution of closed curves in $\mathcal{H}$ by a flow related to Problem~\ref{HCSF1} is investigated in \cite{Drugan} and \cite{PanSun}: there, due to objectives that differs from ours, the projected curves evolve by the standard curve shortening flow (and their initial closed curves $c_{0}$ cannot be embedded due to the horizontality condition (recall Remark ~\ref{rem1.2})).
\end{rem}

\end{document}